\documentclass[final,onefignum,onetabnum]{siamart190516RT}

\usepackage{lipsum}
\usepackage{amsmath}
\usepackage{amssymb}
\usepackage{amsfonts}
\usepackage{graphicx}
\usepackage{epstopdf}
\usepackage{algorithmic}

\usepackage[percent]{overpic}
\usepackage{subcaption}

\ifpdf
  \DeclareGraphicsExtensions{.eps,.pdf,.png,.jpg}
\else
  \DeclareGraphicsExtensions{.eps}
\fi

\newsiamremark{remark}{Remark}
\newsiamremark{hypothesis}{Hypothesis}
\crefname{hypothesis}{Hypothesis}{Hypotheses}
\newsiamthm{claim}{Claim}

\headers{Oscillatory networks: Insights from piecewise-linear modeling}{S. Coombes, M. Sayli, R. Thul, R. Nicks, M. A. Porter, and Y-M. Lai}

\title{Oscillatory networks: Insights from piecewise-linear modeling\thanks{Submitted to the editors \today{}.
\funding{This work was supported by the Engineering and Physical Sciences Research Council [grant numbers EP/P007031/1 and  EP/V04866X/1].}}}

\author{Stephen Coombes\thanks{School of Mathematical Sciences, University of Nottingham, Nottingham, NG7 2RD, UK 
  (\email{stephen.coombes@nottingham.ac.uk}, \url{https://www.maths.nottingham.ac.uk/plp/pmzsc/}).}
\and Mustafa \c{S}ayli \thanks{School of Mathematical Sciences, University of Nottingham, Nottingham, NG7 2RD, UK 
  (\email{mustafa.sayli1@nottingham.ac.uk})}
\and R\"udiger Thul \thanks{School of Mathematical Sciences, University of Nottingham, Nottingham, NG7 2RD, UK 
  (\email{ruediger.thul@nottingham.ac.uk})}
\and Rachel Nicks \thanks{School of Mathematical Sciences, University of Nottingham, Nottingham, NG7 2RD, UK 
  (\email{rachel.nicks@nottingham.ac.uk})}
\and Mason A. Porter \thanks{Department of Mathematics, University of California, Los Angeles, 90095, USA; Santa Fe Institute, Santa Fe, NM, 87501, USA  
(\email{mason@math.ucla.edu})}
\and Yi Ming Lai\thanks{Faculty of Medicine \& Health Sciences, University of Nottingham, Nottingham, NG7 2HA, UK 
  (\email{Yi.Lai1@nottingham.ac.uk})}
}

\usepackage{amsopn}

\ifpdf
\hypersetup{
  pdftitle={Oscillatory networks: Insights from piecewise-linear modeling},
  pdfauthor={S. Coombes, M. Sayli, R. Thul, R. Nicks, M. A. Porter, and Y-M. Lai}
}
\fi

\usepackage{mathrsfs}
\usepackage{xspace}
\usepackage{cite}
\usepackage[normalem]{ulem}

\renewcommand{\d}{{\rm d}}
\newcommand{\e}{{\rm e}}
\newcommand{\D}{{\rm D}}

\newcommand{\FD}[2]{\frac{\d #1}{\d #2}}
\renewcommand{\vec}[1]{\mathbf{#1}}

\newcommand{\ca}{$\mathrm{Ca}^{2+}$\xspace}

\newcommand{\RSet}{{\mathbb{R}}}
\newcommand{\CSet}{{\mathbb{C}}}
\newcommand{\TSet}{{\mathbb{T}}}
\DeclareMathSymbol{\SSet}{\mathalpha}{AMSb}{"53}
\DeclareMathSymbol{\ZSet}{\mathalpha}{AMSb}{"5A}
\DeclareMathOperator{\sgn}{sgn}

\newcommand{\etal}{\emph{et al.}\xspace}

\newcommand{\map}[1]{\textcolor{blue}{#1}}

\begin{document}

\maketitle

\begin{abstract}
There is enormous interest --- both mathematically and in diverse applications --- in understanding the dynamics of coupled oscillator networks. The real-world motivation of such networks arises from studies of the brain, the heart, ecology, and more.
 It is common to describe the rich emergent behavior in these systems in terms of complex patterns of network activity that reflect both the connectivity and the nonlinear dynamics of the network components. Such behavior is often organized around phase-locked periodic states and their instabilities.  However, the explicit calculation of periodic orbits in nonlinear systems (even in low dimensions) is notoriously hard, so network-level insights often require the numerical construction of some underlying periodic component. In this paper, we review powerful techniques for studying coupled oscillator networks. We discuss phase reductions, phase--amplitude reductions, and the master stability function for smooth dynamical systems. We then focus in particular on the augmentation of these methods to analyze piecewise-linear systems, for which one can readily construct periodic orbits. This yields useful insights into network behavior, but the cost is that one needs to study nonsmooth dynamical systems. The study of nonsmooth systems is well-developed when focusing on the interacting units (i.e., at the node level) of a system, and we give a detailed presentation of how to use
  \textit{saltation operators}, which can treat the propagation of perturbations through switching manifolds, to understand dynamics and bifurcations at the network level. We illustrate this merger of tools and techniques from network science and nonsmooth dynamical systems with applications to neural systems, cardiac systems, networks of electro-mechanical oscillators, and cooperation in cattle herds.
\end{abstract}

\begin{keywords}
Coupled oscillators, networks, phase reduction, phase--amplitude reduction, master stability function, network symmetries, piecewise-linear oscillator models, nonsmooth dynamics, saltation operators
\end{keywords}

\begin{AMS}
34C15, 49J52, 90B10, 92C42, 91D30, 49J52.
\end{AMS}

\tableofcontents

\section*{Dedication}

We dedicate this paper to the memory of our dear friend and colleague Yi Ming Lai.  Although he began with us on the journey to write this paper, which in part reviews some of his research activity in recent years, sadly he did not end that journey with us. RIP Yi Ming Lai 1988--2022.


\section{Introduction\label{sec:intro}}

Real-world networks --- such as those in the brain, the heart, and ecological systems --- exhibit rich emergent behavior.
The observed complex patterns of network activity reflect both the connectivity and the nonlinear dynamics of the network components~\cite{Porter2014}.
 The science of networks \cite{Newman2018} has proven especially fruitful in probing the role of connectivity, as exemplified by \cite{Traud2011}. However, overly focusing on network connectivity can downplay the crucial role of dynamics, and even the investigation of dynamical processes on networks has often focused 
 on a few types of situations \cite{porter2016}, such the spread of infectious diseases \cite{RMP2015} and synchronization in coupled oscillators \cite{Arenas2006}.
 This is perhaps not too surprising, given the significant challenges of understanding even low-dimensional dynamical systems. However, for some time, there has been an appreciation in the applied sciences of the benefits of studying 
 complex systems in the form of networks of
 piecewise-linear (PWL) and possibly discontinuous dynamical systems. 
 
 There is a long history of PWL modeling throughout engineering --- particularly in electrical engineering \cite{Acary2011} and mechanical engineering \cite{di2008bifurcations} --- that has now begun to pervade other disciplines, including the social sciences, finance, and biology \cite{bernardo2008piecewise,Champneys2008}. 
 In neuroscience, the McKean model is a classical example \cite{McKean70} of a PWL system. In the McKean mode, one replaces the cubic nullcline
 of the FitzHugh--Nagumo model \cite{FN-scholar} for action-potential (i.e., nerve-impulse) generation with a PWL function that preserves the original shape, allowing explicit calculations that one cannot perform with the original smooth system.
 At its heart, PWL modeling allows one to obtain analytical insight into a nonlinear model by (1) breaking down its phase space into regions
 in which trajectories obey linear dynamical systems and (2) patching these together across the boundaries between the regions.
 The approach can also handle discontinuous dynamical systems, such as those that arise naturally when modeling impacting mechanical oscillators, integrate-and-fire (IF) models of spiking neurons, and 
 cardiac oscillators with both state-dependent and time-dependent
 switching \cite{Thul2010}.  
 Although PWL modeling is a beautifully simplistic modeling perspective, the
  loss of smoothness precludes the use of many results from the standard toolkit of smooth dynamical systems \cite{bernardo2008piecewise},
  and one must be careful to correctly determine conditions for the existence, uniqueness, and stability of solutions.  

An important perspective in the applied dynamical-systems community is that the piecewise nature of models 
is a much more generally applicable feature for many modern applications in science than the smooth dynamical-systems approach that has dominated to date \cite{Glendinning2015}.  
We refer to the switches and discontinuities in such models as \textit{threshold elements}.
The explicit analysis of PWL models at the network level builds on results at the level of individual nodes (e.g., individual oscillatory units),
in disciplines ranging from engineering to biology, to reap 
benefits for understanding network states. This approach
opens up a new frontier in network science to address the role of node dynamics in the interrelationships between the structure and function of real-world networks \cite{Havlin2012}.  

Throughout the present review, we illustrate the above modeling approach
with applications to biological networks in neuroscience and cardiology. We also illustrate these ideas with explorations of other systems, including Franklin Bells and coordinated behavior in cow herds.

We consider networks of $N$ identical oscillators of
the general form
\begin{equation}
	\dot{x} \equiv \FD{}{t} {x}_i = {f}({x}_i) + {g}_i({x}_1, {x}_2, \ldots, {x}_N)\,, \quad i \in \{1,\ldots, N\} \,, \quad {x}_i \in \RSet^m
\label{Network1}
\end{equation}
and show how to gain insight into emergent network dynamics when the vector field ${f}$ (i.e., the local dynamics) is PWL and the interactions are pairwise. Each oscillator is associated with a node of a structural network (which, most traditionally, takes the form of a graph \cite{Newman2018}), and each interaction is associated with an edge of that network.
With only pairwise interactions, the coupling function is
\begin{equation}
	{g}_i({x}_1, {x}_2, \ldots, {x}_N) =  \sigma \sum_{j=1}^N  w_{ij}  {G} ({x}_i,{x}_j) \, ,
\label{Networkinterraction}
\end{equation}
where $G({x}_i,{x}_j)$ is the dynamics that expresses the coupling between nodes $i$ and $j$, the relative strength of this interaction is $w_{ij}$,
and
$\sigma$ sets the overall network coupling strength.
One achieves insight into network behavior with a merger of techniques that have been developed for nonsmooth systems (see, e.g., \cite{Makarenkov2012}), as exemplified in the books of di Bernardo \textit{et al}.~\cite{bernardo2008piecewise}, Acary \textit{et al.} \cite{Acary2011}, and Jeffrey \cite{Jeffrey2016} for low-dimensional systems with discontinuous behavior and by network-science tools, especially weakly-coupled oscillator theory \cite{Hoppensteadt97} and the master stability function \cite{Pecora2013}
--- that have been developed to describe phase-locked states 
(i.e., states in which all pairs of oscillators are frequency-locked with a constant phase lag between each pair) and their bifurcations.

Our paper proceeds as follows.  In \cref{sec:pwl}, we present the types of PWL models --- including PWL continuous, Filippov, and impacting systems --- that we use as nodes of a network.
We partition the phase space of these PWL models using \textit{switching manifolds}. We give a method to construct periodic orbits, and we describe and employ an extension of Floquet theory to nonsmooth systems to determine a criterion for the stability of a periodic orbit. We use
\textit{saltation operators} to describe the propagation of perturbations through the switching manifolds.
In \cref{sec:iso}, we present a reduction technique that allows one to describe a limit-cycle oscillator in terms of a scalar phase variable and additional variables that encode
directed distances. By again exploiting saltation operators, we show how to calculate the infinitesimal phase and amplitude responses for PWL models. We illustrate this approach for some PWL neuron models. 
We first examine
weakly coupled systems.
 In \cref{sec:phaseoscnetworks}, we consider
 phase-only network descriptions (i.e., dropping the amplitude coordinates) and we also describe the relevant phase-interaction function. 
We highlight the usefulness of a phase-oscillator network description using a combination of theory (specifically, about the stability of phase-locked network states) and numerical simulations, with a focus on neural networks.\footnote{When we write ``neural networks'', we are referring to networks in neuroscience, as opposed to the use of the term ``neural networks'' in contexts such as machine learning.}
In \cref{PWLPhaseAmplitudenetworkfromalism}, we examine phase--amplitude networks, for which one needs
more functions to fully specify all of the interactions between units.
 We use a simple two-node network to highlight the dangers of an overreliance on only phase information and emphasize the benefits of using phase--amplitude coordinates to correctly predict phase-locked behavior for moderate values of the network coupling strength $\sigma$.  
 We then consider strongly coupled systems, for which we do not expect to obtain good predictions of system behavior from approximations of the network dynamics through either phase-only reductions or phase--amplitude network reductions.
 In \cref{sec:strongly}, we develop a theory of phase-locked states in networks of identical PWL oscillators without recourse to any approximation.  In essence, this theory is based on an extension of the master stability function to nonsmooth systems. We use saltation operators to develop this extension. We apply this theory to a  variety of distinct systems, with a focus on synchronous network states and solutions that can arise when a synchronous state loses stability.
 Finally, in \cref{sec:discussion}, we summarize our paper and then briefly discuss extensions and further applications of the methodology in it for analyzing the dynamics of coupled-oscillator networks.

\section{Piecewise-linear oscillators\label{sec:pwl}}

Planar PWL systems \cite{dumortier2006qualitative,shaw1983periodically,freire2012canonical} have dynamics
on two regions (i.e., ``zones''), with a line of discontinuity between those regions. The dynamics of planar PWL systems can be complicated, but they are tractable to study.
Therefore, we start by considering them. We describe the dynamics in the two zones by the variable $x = (v,w)^{\top}\in \nolinebreak \mathbb R^2$, which satisfies
\begin{equation}\label{pwl1}
	\FD{x}{t}=\left\{
	\begin{array}{cc}
		f_{1} \equiv A_{1}x + b_{1} & \mathrm{if} \; x \in R_{1}\,\,\,  \\
		f_{2} \equiv A_{2}x + b_{2} & \mathrm{if} \; x \in R_{2}\,,
	\end{array}\right.
\end{equation}
where $A_{1,2}\in \mathbb R^{2\times 2}$ are constant matrices and $ \  b_{1,2}\in \mathbb R^{2} $ are constant vectors. The regions $R_{1}$ and $R_{2}$ are
\begin{equation}\label{pwlregions}
	R_{1}=\{ x \in \mathbb R^{2} |\ h(x)>0\}  \quad  \mathrm{and} \quad  R_{2} = \{x \in \mathbb R^{2} |\ h(x)<0\} \,,	
\end{equation}
where the \textit{indicator function} $h:\mathbb R^{2} \to \mathbb R^{}$ is
\begin{equation}\label{pwlindicator}
	h(x) = v - a \,. 	
\end{equation}
Switching events occur when $h(x) = 0$, which holds on the \emph{switching manifold} $\Sigma=\{ x \in \mathbb{R}^2 | \ v=a\}$.
The condition $h(x(t_i)) = 0$ implicitly yields the event times $t = t_i$, with $i \in \ZSet$.
If an equilibrium point
exists in the region $R_{\mu}$, one determines its stability by the eigenvalues of $A_{\mu}$, with $\mu \in \{1,2\}$.
When relevant, it is simple to partition phase space into more regions and to thereby incorporate further switching manifolds, so we describe only the simplest situation of two regions of phase space. However, in \cref{pwlMorrisLecar}, we give an example of a system with three
switching manifolds.

Planar PWL systems of the form \cref{pwl1} have been studied for many years and can have rich dynamics.
 For example, Freire \textit{et. al} \cite{freire1998bifurcation} considered continuous systems with two zones and proposed a canonical form that captures many interesting oscillatory behaviors, and Llibre \textit{et. al} \cite{llibre2019limit} studied the existence and maximum number of limit cycles in systems with a discontinuity.
 Planar PWL systems can have almost all types of dynamics that occur in smooth nonlinear dynamical systems, and they can also
 support bifurcations that are not possible in smooth systems~\cite{di2008bifurcations,bernardo2008piecewise}. However, in comparison to smooth systems, the knowledge of bifurcations in PWL systems is largely limited to specific examples \cite{colombo2012bifurcations}.
 Nevertheless, we can
 start to develop a picture of the theory of bifurcations in PWL systems by gathering results from the differential inclusions of Filippov \cite{filippov2013differential}, the ``C bifurcations'' of Feigin \cite{feigin1994forced,di1999local}, and the nonsmooth equilibrium bifurcations of Andronov \emph{et al.}~\cite{andronov2013theory}. Examples of well-known bifurcations
 that arise from discontinuities include grazing bifurcations, sliding bifurcations, and discontinuous saddle--node bifurcations~\cite{bernardo2008piecewise,Harris2015}.

One of the key advantages of PWL modeling is that it allows one to derive closed-form expressions
for periodic orbits\footnote{Every periodic orbit that we consider in this paper is also a limit cycle, so we use the terms ``periodic orbit'' and ``limit cycle''
 interchangeably.} \cite{ponce2014bifurcations}. However, the analysis of such dynamics is not trivial because one needs to match the
solution pieces from
separate linear regimes. Deriving conditions for matching dynamics from different regions
typically necessitates the explicit knowledge of the \textit{times-of-flight} (i.e., the time that is spent by the flow in a zone of phase space before reaching the switching manifold) in each region.
 Essentially, we solve the system \cref{pwl1} in each of its linear zones using matrix exponentials and demand continuity of solutions to construct orbits of the full nonlinear flow. To clarify how to implement this procedure, we denote a trajectory in zone $R_{\mu}$ by $x^{\mu}$ and solve \cref{pwl1} to obtain $x^{\mu}(t,t_0)=x(t,t_0; A_{\mu},b_{\mu})$ using the solution form
\begin{equation}\label{pwl10}
	x(t,t_0;A,b) = G(t - t_0;A)x(t_0) + K(t - t_0;A)b \,,
\end{equation}
where $t_0$ is the initial time, $t > t_0$, and
\begin{equation}\label{GandK}
	G(t;A)=\e^{At}\,, \quad  K(t;A) = \int_{0}^{t}G(s;A)\mathrm{d} s = A^{-1}[G(t;A) - I_2] \,,
\end{equation}
where $I_m$ is the $m \times m$ identity matrix.
One can construct a closed orbit (i.e., a periodic orbit) by connecting two trajectories. One starts from initial data $x(0) = (a,w(0))^{\top}$, which lies on the switching manifold,
in each zone.
One then writes
\begin{equation}\label{per1}
	x(t)=
	\begin{cases}
		x^1(t,0) & \mathrm{if} \; t\in [0,T_1] \\
		x^2(t,T_1) & \mathrm{if} \; t\in (T_1,T] \,,
	\end{cases}
\end{equation}
for some $T > T_1 > 0$. We obtain
a periodic orbit
by requiring that $x$ have period $T$ (i.e., be $T$-periodic). The times $T_{i}$, with $i \in \{1,2\}$ and $T_2 = T - T_1$, gives the times-of-flight between switching events.
To complete the procedure, we must determine the unknowns $(T_1,T_2,w^1(0))$ by simultaneously solving a system of three equations:
$a=v^1(T_1)$, $a=v^2(T_2)$, and $w^2(T_2)=w^1(0)$.
This is easy to do using
 a numerical method for root finding, such as \texttt{fsolve} in {\sc Matlab}, along with a method to compute matrix exponentials (e.g., \texttt{exmp} in {\sc Matlab}). Alternatively, one can readily perform explicit calculations of $G(t;A)$ and $K(t;A)$~\cite{Coombes2008}.

One can classify PWL systems into three different types, depending on their degree of discontinuity \cite{bernardo2008piecewise,leine2013dynamics}. These three types of PWL systems are as follows.
\medskip
\begin{description}
\item[Continuous PWL systems.] These systems have continuous states and continuous vector fields (i.e., $f_1(x) = f_2(x)$) but discontinuities in the first derivative or higher derivatives of the right-hand side functions (i.e.,  $
\partial^n f_1/ \partial x^n \neq \partial^n f_2/ \partial x^n$ for an integer $n \geq 1$),
across the switching manifold. These systems have a degree of smoothness of $2$ or more, but their Jacobian matrices are different on different sides of a switching manifold (i.e., $\D f_1(x) \neq \D f_2(x)$).
\medskip
\item[Filippov systems \cite{filippov2013differential}.]
These systems have continuous states but vector fields that are different on different sides of a switching manifold (i.e., $f_1(x) \neq f_2(x)$).
These systems have a degree of smoothness of $1$.
 The vector field of the system \cref{pwl1} is not defined on the switching manifold $\Sigma=\{ x \in \mathbb{R}^2 | \ h(x)=0\}$. One completes the description of the dynamics on the switching manifold
with a set-valued extension $f(x)$. The extended dynamical system is
\begin{equation}\label{NonsmoothODEinclusion}
	\dfrac{\mathrm{d}x_{}}{\mathrm{d} t} \in f(x)= \left\lbrace \begin{aligned}  & f_1(x) \quad \mathrm{if} \quad x \in R_1 \\
		& \mathrm{\overline{co}} \left\lbrace   f_1(x),f_2(x) \right\rbrace  \quad  \mathrm{if} \quad  x \in \Sigma_{} \\
		& f_2(x) \quad \mathrm{if}  \quad  x \in R_2 \,,
	\end{aligned} \right.
\end{equation}
where $\mathrm{\overline{co}} (\mathcal{A})$ denotes the smallest closed convex set that contains $\mathcal{A}$. In \eqref{NonsmoothODEinclusion}, we have
\begin{equation}\label{ConvexInc}
	\mathrm{\overline{co}} \left\lbrace   f_1(x),f_2(x) \right\rbrace=  \left\lbrace \varsigma f_1(x)+(1-\varsigma)f_2(x) \,\,\,  \text{for all} \,\,\,  \varsigma \in [0,1] \right\rbrace \,,
\end{equation}
where $\varsigma$ (which has no physical meaning) is a parameter that defines the convex combination.
 The extension (i.e., convexification) of the discontinuous system \cref{pwl1} into a convex differential inclusion \cref{NonsmoothODEinclusion} is known as the \textit{Filippov convex method} \cite{filippov2013differential}. If $\langle  \nabla h_{} , f_1 \rangle \langle  \nabla h_{} , f_2 \rangle <0$, a Filippov system can have
 \textit{sliding} motion \cite{jeffrey2011geometry,jeffrey2011sliding}, with $\dot{h}_{}= \nabla h_{}\cdot f=0$, along a switching manifold.\footnote{We use $\langle \cdot, \cdot \rangle$ and $\cdot$ interchangeably to denote the standard vector inner product.}
We then have
\begin{equation}\label{Slidingdynamicconstant}
	\varsigma= \frac{\nabla h_{}\cdot f_2}{\nabla h_{}\cdot (f_2-f_1)} \,.
\end{equation}
\medskip
\item[Impacting systems (i.e., impulsive systems).] These systems have instantaneous discontinuities (i.e., ``jumps'') in a solution at the switching boundary $\Sigma_{}$ that are governed by a smooth jump operator (i.e., a ``switch rule'') $x(t^+)=\mathcal{J}_{}(x(t^-))$,
where $t^-$ denotes the time immediately before the impact and $t^+$ denotes the time immediately after the impact.
These systems have a degree of smoothness of $0$.
 The jump operator $\mathcal{J}_{}$ is often called an \textit{impact rule} (or an \textit{impact law}), and the discontinuity boundary $\Sigma_{}$ is often called an \textit{impact surface}. Depending on the properties of
 $\mathcal{J}_{}$, many different types of dynamics can occur.
To further understand the behavior of impacting systems, see \cite{brogliato1999nonsmooth,brogliato2000impacts,di2008bifurcations,bernardo2008piecewise}.
\end{description}
\medskip

\noindent
To illustrate this classification, we now briefly introduce five different models, each of which has oscillatory behavior and can be written in the form \cref{pwl1}. We defer the detailed form of these models to \cref{Pwlmodels}. In the present section, we emphasize the qualitative aspects of each model with plots of their nullclines and typical periodic orbits (which we construct using the method that we described in the present section).

\medskip
\begin{description}

\item[Absolute model (continuous) [see \cref{Fig:Pwl2pieces}(a)]].
The vector field is continuous across the switching boundary, although its Jacobian is not. The equilibrium point in zone $R_1$ is an unstable focus, and the equilibrium point in zone $R_2$ is a stable focus.
A nonsmooth Andronov--Hopf bifurcation \cite{simpson2008unfolding,jeffrey2018hidden,simpson2007andronov} occurs when an equilibrium crosses from $R_1$ to $ R_2$ and the eigenvalues of the Jacobian jump across the imaginary axis.

\medskip
\item[PWL homoclinic model (continuous) [see \cref{Fig:Pwl2pieces}(b)]].
There is a saddle point for $x \in R_{2}$ and an unstable focus for $x \in R_{1}$, with a vector field that crosses the switching boundary in a  continuous manner. There is a homoclinic orbit that tangentially touches the unstable and stable eigendirections of the saddle point in $R_{2}$. This orbit encloses the unstable focus in $R_{1}$.
See \cite{xu2013homoclinic} for a detailed discussion of the conditions that ensure existence of a limit cycle or a homoclinic orbit.

\medskip
\item[PWL Morris--Lecar model (continuous) [see \cref{Fig:TheMorrisLecarperiodic}]].
\label{pwlMorrisLecar}
The Morris--Lecar model is a planar conductance-based single-neuron model that captures many important features (such as low firing rates) of neuronal firing \cite{morris1981voltage}.  
One can then simplify it to obtain
a PWL system with four zones and three switching manifolds \cite{Coombes2008}.
(By contrast, our other examples have two zones and one switching manifold.) For full details, see \cref{Pwlmodels}.

\medskip
\item[McKean model (Filippov) [see \cref{Fig:Pwl2pieces}(c)]].
The McKean model is a well-known planar PWL model for action-potential generation~\cite{McKean70}. There are two varieties of McKean model. One of them has a
PWL approximation
of a cubic nonlinearity (to capture the behavior of the FitzHugh--Nagumo model), with the associated nullcline
broken into three pieces.
In the other variety, the PWL approximation of the cubic nonlinearity has two pieces
\cite{tonnelier2003mckean}.
We discuss the latter, which requires a set-valued extension on the switching manifold [see \cref{NonsmoothODEinclusion}--\cref{ConvexInc}].
For some parameter values, a stable periodic orbit coexists with a stable
equilibrium point (i.e., an attracting focus). In all of those situations, they are separated by an unstable sliding periodic orbit.

\medskip
\item[Planar IF model (impacting) [see \cref{Fig:Pwl2pieces}(d)]].
In the planar IF model, which is a single-neuron model, whenever the voltage variable $v$ reaches a firing threshold $v_{\text{th}}$, the system resets according to $x(t^+) = \mathcal{J}(x(t^-))\equiv (v_\text{r},w(t^-) + \kappa/\tau )$.
Namely, the voltage $v$ \textit{resets} to $v_\text{r}$ and the recovery variable $w$ is \textit{kicked} by the amount $\kappa/\tau$, where $\kappa$ is the kick strength and $\tau$ is the time scale of the recovery variable.

\begin{figure}[!h]
	\centering
	\includegraphics[width=5.5cm,angle=0]{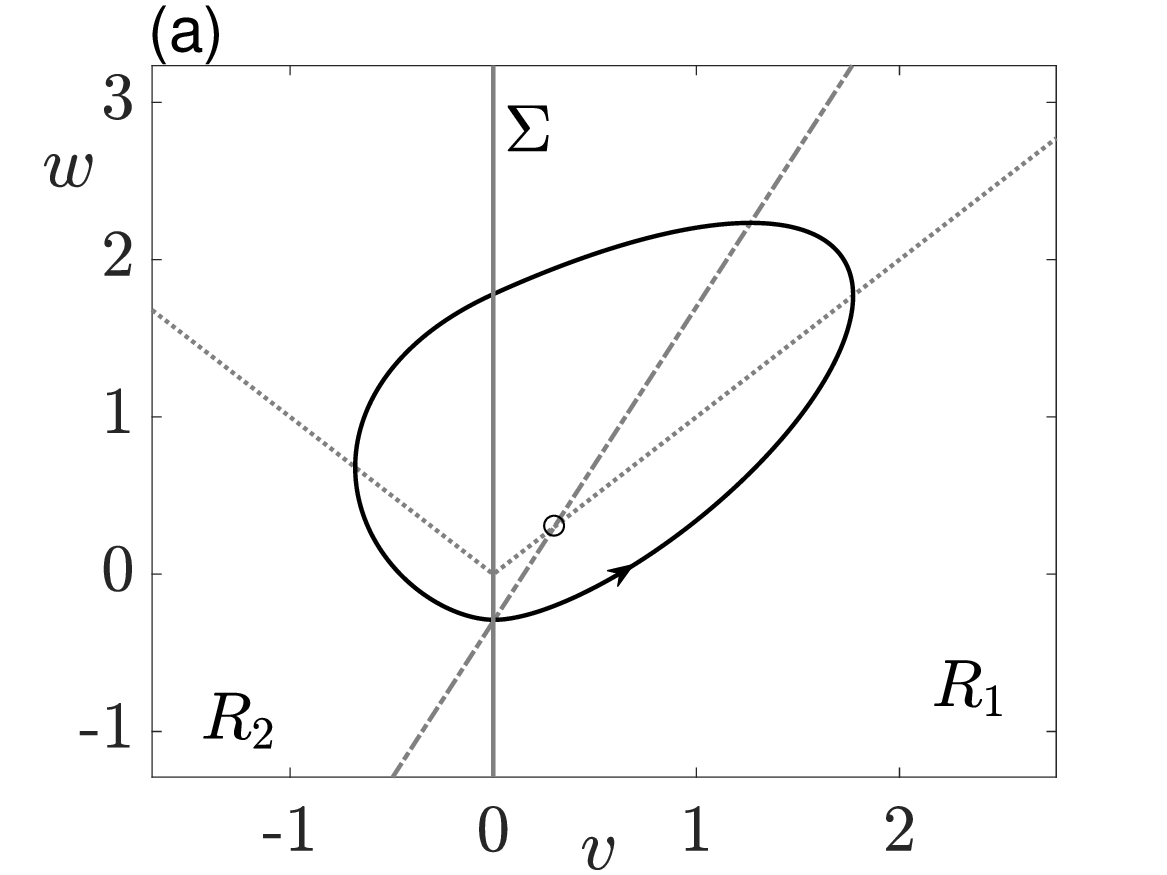}
	\includegraphics[width=5.5cm,angle=0]{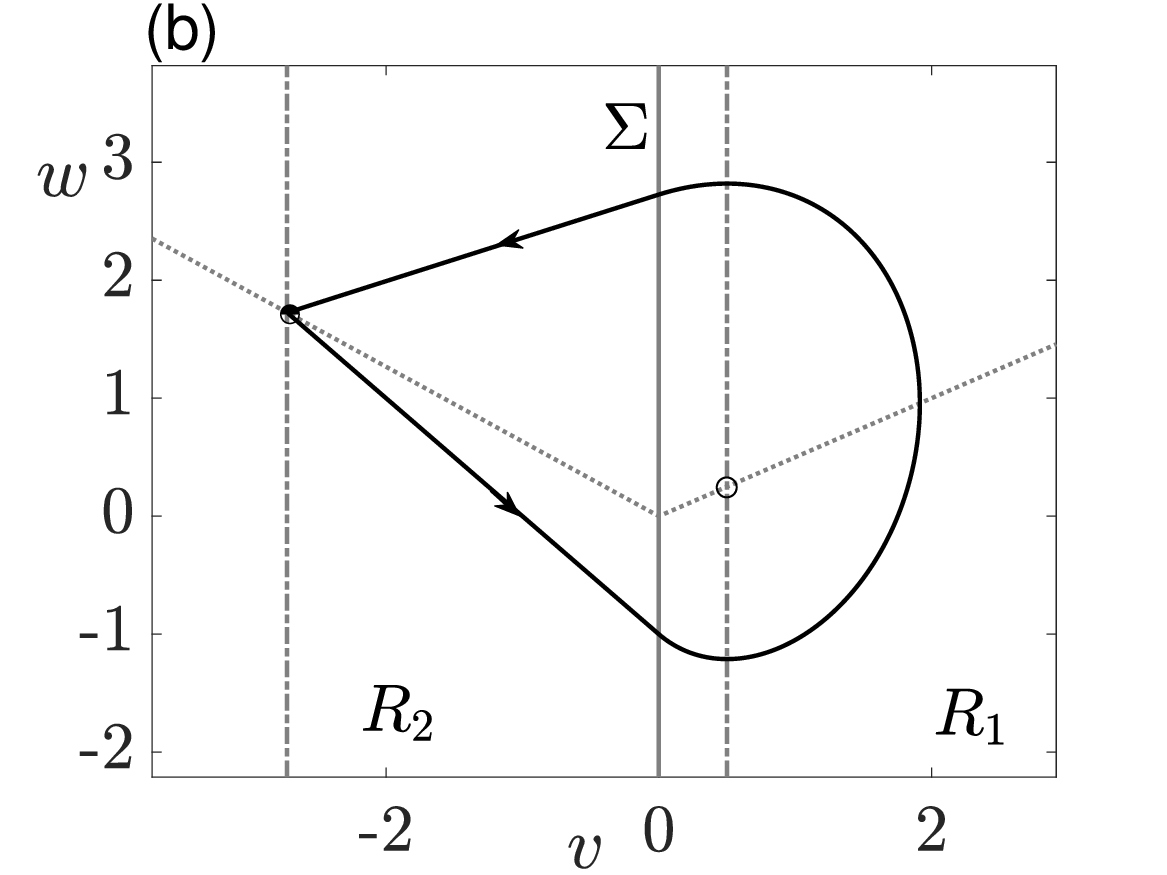}
	\includegraphics[width=5.5cm,angle=0]{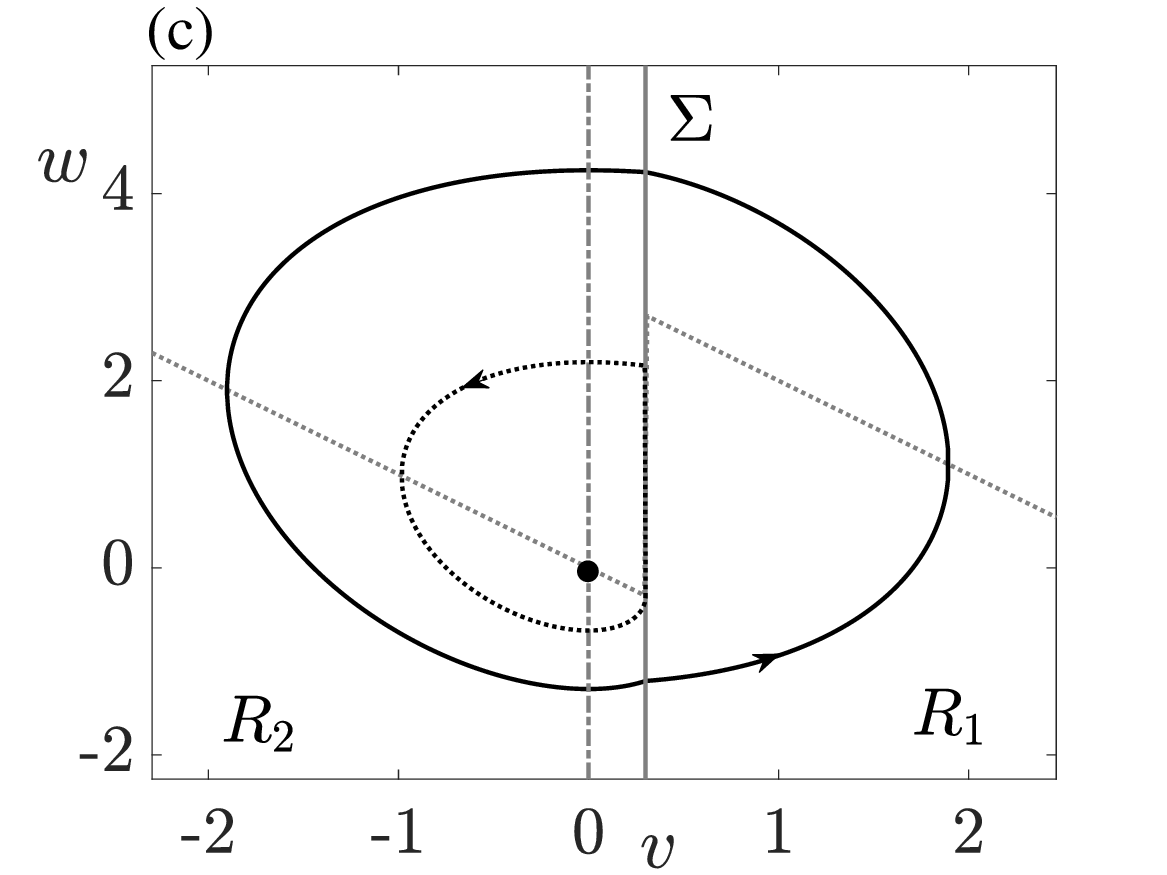}
	\includegraphics[width=5.5cm,angle=0]{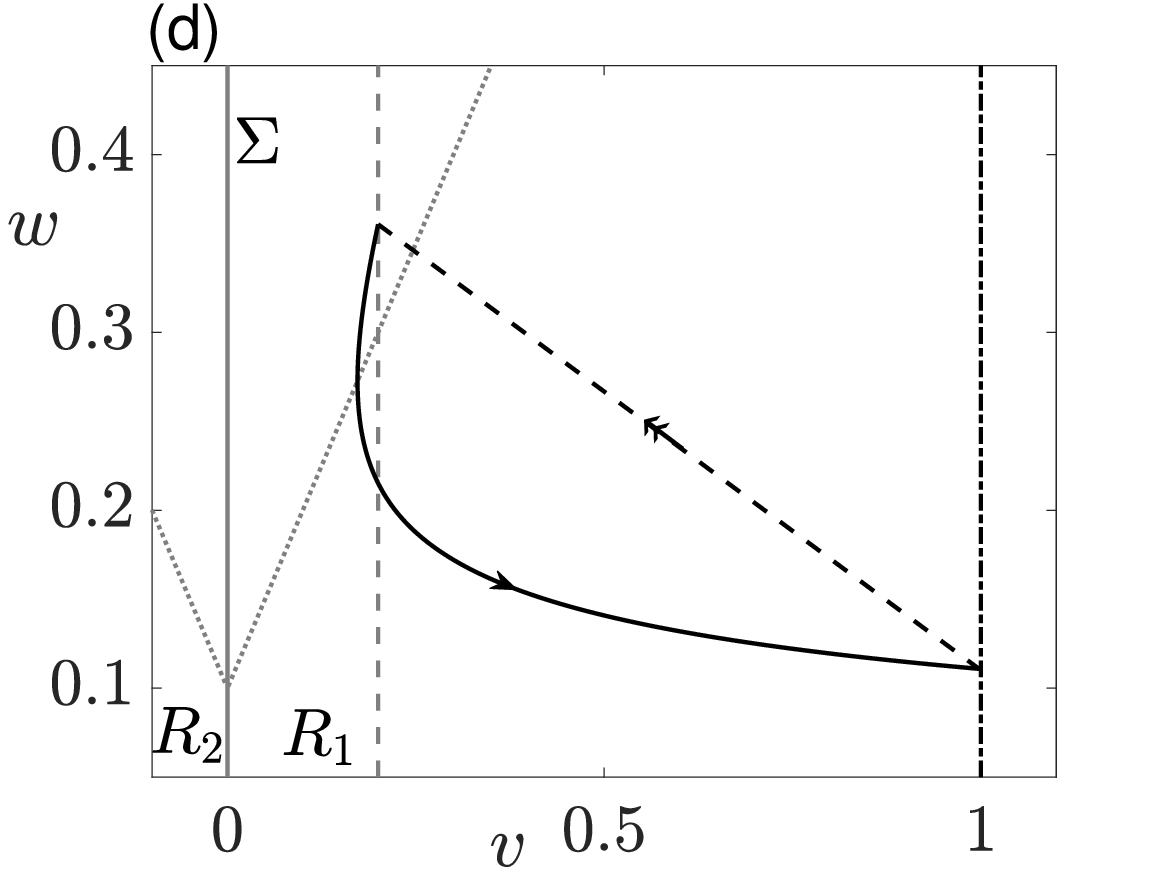}
	\caption{Nullclines and periodic orbits in a variety of planar PWL models. The region $R_1$ (respectively, $R_2$) is the zone with $v>a$ (respectively, $v<a$). We show the stable (respectively, unstable) periodic orbits with solid (respectively, dotted) black curves.  We show the $v$-nullcline (i.e., the curve $\dot{v} = 0$) with a dotted gray curve and the $w$-nullcline (i.e., the curve $\dot{w} = 0$) with a dashed--dotted gray curve. We indicate the switching manifold ($v = a$) with a solid gray line.
	(a) Absolute model. The unstable equilibrium point, which we indicate with an unfilled circle, is in the zone $R_1$. The parameter values are $a=0$, $\overline{w}=-0.1$, $\overline{v}=0.1$, and $d=0.5$.
	(b) PWL homoclinic model.  The repelling focus, which we indicate with an unfilled circle, is in zone $R_1$. The saddle point, which we indicate with a half-filled circle, is in zone $R_2$.
		 The parameter values are $a=0$, $\delta_{1}=2$, $\delta_{2} = -0.3667$, $\tau_{1}=0.5$, and $\tau_{2}=-0.6333$.
	(c) McKean model.  The unstable periodic orbit is of sliding
	type. The stable equilibrium point, which we indicate by a filled black circle and is a focus, is in the zone $R_2$.
	 The parameter values are $a=0.3$, $b=2$, $\gamma=1$, and $I=3$.
	(d) Planar IF model.  We indicate the firing threshold with a dashed--dotted black line and indicate the reset value with a dashed gray line.
	The parameter values are $v_{\text{th}}=1$, $v_{r}=0.2$, $a_{w}=0$, $b_{w}=-1$, $a_{1}=1$, $a_{2}=-1$, and $I=0.1$.
	For further details about these models, see \cref{sec:pwl} and \cref{Pwlmodels}.
	}
	\label{Fig:Pwl2pieces}
\end{figure}

\begin{figure}[hbt!]
	\centering
	\includegraphics[height=7cm]{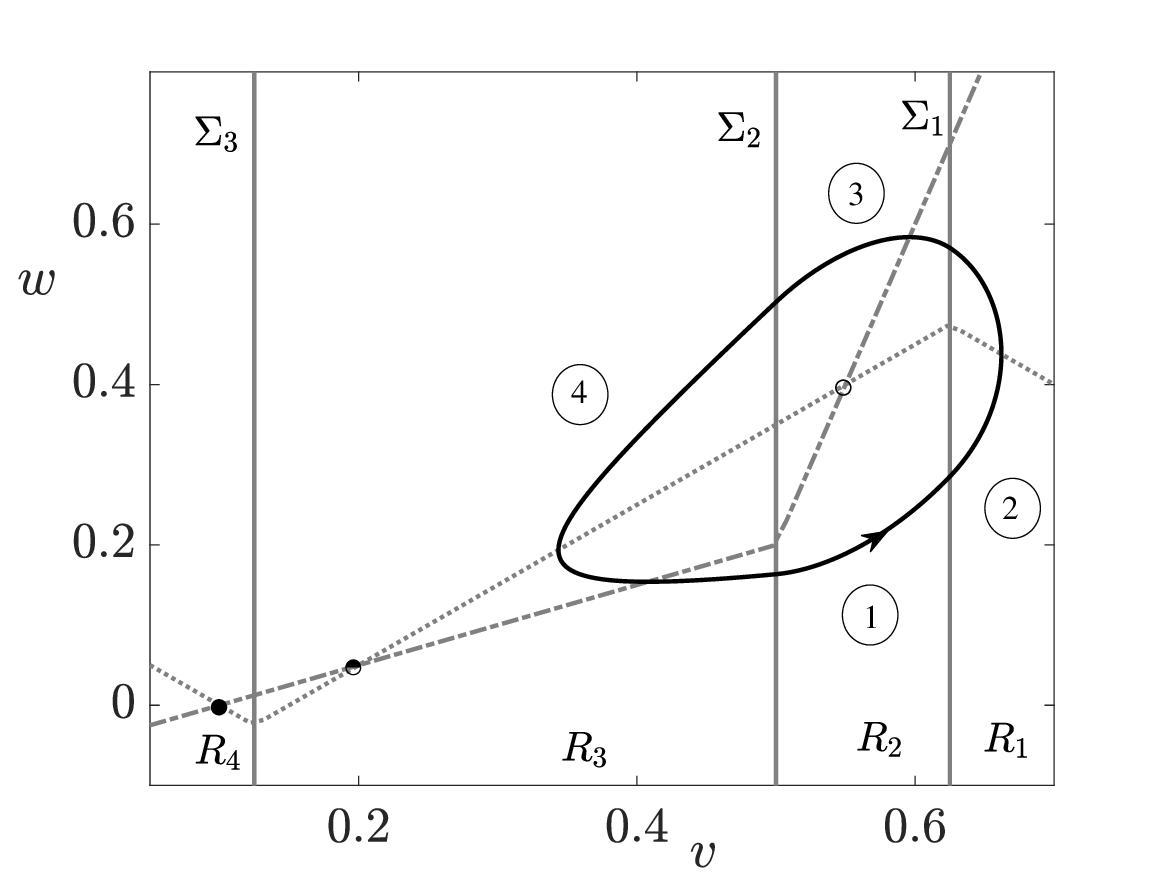}
	\caption{Phase plane of the piecewise-linear Morris--Lecar model, with a stable periodic orbit in black. The periodic orbit has four pieces, with the first and third pieces in $R_2$, the second piece in $R_1$, and the fourth piece in $R_3$. We show the $v$-nullcline with a dotted gray line, the $w$-nullcline with a dashed--dotted gray line, and the switch manifolds $\Sigma_1$, $\Sigma_2$, and $\Sigma_3$ with solid gray lines.
	The nullclines are piecewise-linear approximations of those of the original smooth Morris--Lecar model. The open black circle indicates an unstable equilibrium point, the half-filled black circle indicates a saddle point, and the filled black circle indicates a stable equilibrium point (which is in zone $R_4=\{x\in \mathbb R^{2} |\ v <a/2\}$).
	The parameter values are $C = 0.825$, $I = 0.1$, $a = 0.25$ $b = 0.5$, $b^* = 0.2$,  $\gamma_1 = 2$, and $\gamma_2 = 0.25$.
	For further details this model, see \cref{sec:pwl} and \cref{Pwlmodels}.
	}
	\label{Fig:TheMorrisLecarperiodic}
\end{figure}
\end{description}


\subsection{Floquet theory for nonsmooth systems}
\label{FlqTheory}

Floquet theory \cite{perko2013differential} is a popular and well-developed technique to study the stability and bifurcations of periodic orbits of smooth dynamical systems ${\d}x/{\d}t = f(x)$, where $x \in\mathbb R^{m}$ and $f(x)$ is a continuously differentiable function. If we write a $T$-periodic solution in the form $x^\gamma(t)$, the
variational equation for
this solution is
\begin{equation}
	\dot{\Phi} = {\D} f(x^\gamma (t)) \Phi\,, \quad \Phi(0) = I_m \,.
\label{variational}
\end{equation}
Equation \eqref{variational} has an associated \textit{monodromy} matrix  $\Phi(T)$.  The eigenvalues of $\Phi(T)$, which are $\lambda_k = \e^{\kappa_k T}$ for all $k \in \{0,\ldots,m-1\}$, are the so-called ``Floquet multipliers'' of the limit cycle, and the values $\kappa_k$ are their associated ``Floquet exponents".
For a planar system, for which $x \in \mathbb R^{2}$, one of the Floquet multipliers is equal to $1$ (corresponding to perturbations that are tangent to the periodic orbit) and the other is $\lambda_{\text{smooth}} = \exp(\kappa_{\text{smooth}} T)$, where
\begin{equation}\label{flq2}
	\kappa_{\text{smooth}}=\frac{1}{T}\int_0^{T} \operatorname{Tr} \, (\D f(x^\gamma (t))) \, \d t \, .
\end{equation}
One determines the stability of periodic orbits
from the sign of $\kappa_{\text{smooth}}$. An orbit is linearly stable if $\kappa_{\text{smooth}} < 0$ and unstable if $\kappa_{\text{smooth}} > 0$.

For dynamical systems with nonsmooth or even discontinuous vector fields, one cannot directly use standard Floquet theory \cite{jordan2007nonlinear,klausmeier2008floquet}. It is also necessary to carefully evolve a perturbation across the switching boundaries. We revisit the adaptation of standard Floquet theory (for non-sliding periodic orbits) to PWL systems \cite{bernardo2008piecewise,CoombesThul:2016} of the form  $\dot{x}=A_{\mu}x(t)+b_{\mu}$, where $A_{\mu}\in \mathbb R^{m\times m}$, $b_{\mu} \in \mathbb R^{m}$, and the phase space has $P$ distinct regions $R_\mu$ (with $\mu \in \{1, \dots, P\}$).
Switching events have associated
indicator functions $h_{\mu}(x)$. They occur when $h_{\mu}(x(t_{i}))=0$ and have switching times $t_{i}$, with $i \in \ZSet$.
The state of the system immediately after the switch event is $x(t_{i}^{+})=\mathcal{J}_{\mu}(x(t_{i}^{-}))$, where $\mathcal{J}_{\mu}:\mathbb{R}^m \rightarrow \mathbb{R}^m$ is the switch rule, $t_{i}^{\pm}=\lim_{\epsilon \rightarrow 0^+} (t_{i} \pm \epsilon)$, and $x(t_{i}^{-})$ denotes the state immediately before the switch event.
We construct a periodic orbit $x^{\gamma}(t)$ is by patching
solutions (built from matrix exponentials) across the boundaries of the regions $R_\mu$.

Away from switching events, the variational equation for
a periodic orbit is
\begin{equation}\label{PwlFloquetVar}
	\FD{}{t} \delta x=A_{\mu} \, \delta x \quad \mathrm{for}\  \  x^{\gamma}(t) + \delta x(t) \in R_{\mu} \,,
\end{equation}
where $\delta x(t)$ is a perturbation of the periodic orbit.
The evolution of perturbations in each region is governed by the matrix exponential form $\delta x(t) = G(t-t_0;A_\mu) \delta x (t-t_0)$, where $t>t_0$ and $t_0$ denotes the time at which the trajectory crosses into region $R_\mu$.  To map perturbations across a switching manifold, we use a \textit{saltation operator} \cite{Muller:1995,fredriksson2000normal}.
This allows us to evaluate perturbations during the boundary crossing in which either the solution or the vector field (or both) has a discontinuity. M\"{u}ller \cite{Muller:1995} used saltation operators to calculate Lyapunov exponents of discontinuous systems and Fredriksson and Nordmark \cite{fredriksson2000normal} used them in a normal-form derivation for impact oscillators.
See \cite{Kong2023} for a recent review of saltation operators and their use in engineering.
In our context, saltation operators admit an explicit matrix construction of the form
\begin{equation}
\label{Delatxplusforminchapter}
	S_\mu(t_{i}) = \mathrm{D}\mathcal{J}_{\mu}(x^{\gamma}(t^-_{i})) +
	\frac{[\dot{x}^{\gamma}(t_{i}^{+})-\mathrm{D}\mathcal{J}_{\mu}(x^{\gamma}(t^-_{i}))\dot{x}^{\gamma}(t_{i}^{-})][\nabla_{x} h_{\mu}(x^{\gamma}(t^-_{i}))]^{\top}}{\nabla_{x} h_{\mu}(x^{\gamma}(t^-_{i})) \cdot \dot{x}^{\gamma}(t_{i}^{-})} \,.
\end{equation}
We derive \cref{Delatxplusforminchapter} in \cref{Saltationproof}.

Equation \cref{Delatxplusforminchapter} allows us to write
\begin{equation}
	\delta x(t_i^+) = S_\mu(t_{i}) \delta x(t_i^-)\,, \quad x^{\gamma}(t_i^-)+\delta x(t_i^-) \in R_{\mu}
\label{saltate}
\end{equation}
to describe how perturbations are mapped across a switching manifold at the boundary of region $R_\mu$.  Combining \cref{PwlFloquetVar} and \cref{saltate} allows us to evaluate $\delta x(t)$ over one oscillation period $T$ using $M$ separate times-of-flight. We thus write $T = \sum_{i=1}^M T_i$, with $\delta x(T)=\Psi \delta x(0)$, where $\Psi$ is the product
\begin{equation}\label{MonodromyExplicitinchapter}
	\Psi= S(t_M)G(T_{M}) S(t_{M-1})G(T_{M-1}) \times \cdots \times S(t_2)G(T_2)S(t_1)G(T_1) \,,
\end{equation}
where $G(T_i) = G(T_i;A_{\mu(i)})$ and $S(t_i) = S_{\mu(i)}(t_i)$.
The index $\mu(i) \in \{1,\ldots,P\}$ indicates the region that the periodic orbit is in at time $t_i^-$.
The periodic orbit is linearly stable if all of its nontrivial eigenvalues (i.e., Floquet multipliers) of the matrix $\Psi$ have moduli less than $1$ and equivalently if the corresponding Floquet exponents ($\kappa_k = \ln(\lambda_k)/T$) all have negative real parts.
One (trivial) eigenvalue of $\Psi$ is equal to $1$, corresponding to perturbations that are tangential to the periodic orbit.
For planar systems, one calculates the lone nontrivial Floquet exponent using the formula
\begin{equation}\label{Fexponentformula}
	\kappa=\frac{1}{T} \sum^{M}_{i =1}\left[T_{i} \operatorname{Tr} A_{\mu(i)} + \ln | \det S(t_i) | \right] \,.
\end{equation}
The logarithmic term in \cref{Fexponentformula} reflects the contribution of discontinuous switching to the stability of an orbit. If $S = I_2$ (i.e., there is no saltation), the logarithmic term vanishes and we recover the formula \cref{flq2} for a smooth system.  In \cref{PWLPlanarFloquet}, we derive the Floquet-exponent formula \cref{Fexponentformula} for planar PWL systems. We use this formula
to compute the stability of periodic orbits in all numerical studies of single-oscillator
PWL models.

\section{Isochrons and isostables} \label{sec:iso}

We now examine networks of interacting PWL oscillators. We start by generalizing results from the theory of weakly coupled systems of smooth oscillators.

The theory of weakly coupled oscillators allows us to obtain insights into the phase relationships between the nodes of a network~\cite{Hoppensteadt97}.  Historically, the theory of weakly coupled oscillators
has focused on phase-reduction techniques using the notion of \textit{isochrons}, which extend the phase variable for a limit-cycle attractor to its basin of attraction \cite{Winfree1967,Guckenheimer1975}. More recent research has emphasized the importance of distance from a limit cycle using \textit{isostable} coordinates (which we call ``isostables'' as shorthand terminology) \cite{Guillamon2009,Mauroy2013,Wilson2016,Mauroy2018}.
Employing isochrons and isostables yield
reductions to phase networks and phase--amplitude networks, respectively, although the theory for the latter is far less developed than the theory for the former.

To introduce the concepts of an isochron and an isostable, it is sufficient to consider the dynamical system $\dot{x} = f(x) + g(t)$, with $x \in \RSet^m$.

\subsection{Phase response and amplitude response}
\label{sec:response}
\begin{figure}[!h]
	\centering
	\includegraphics[height=6cm]{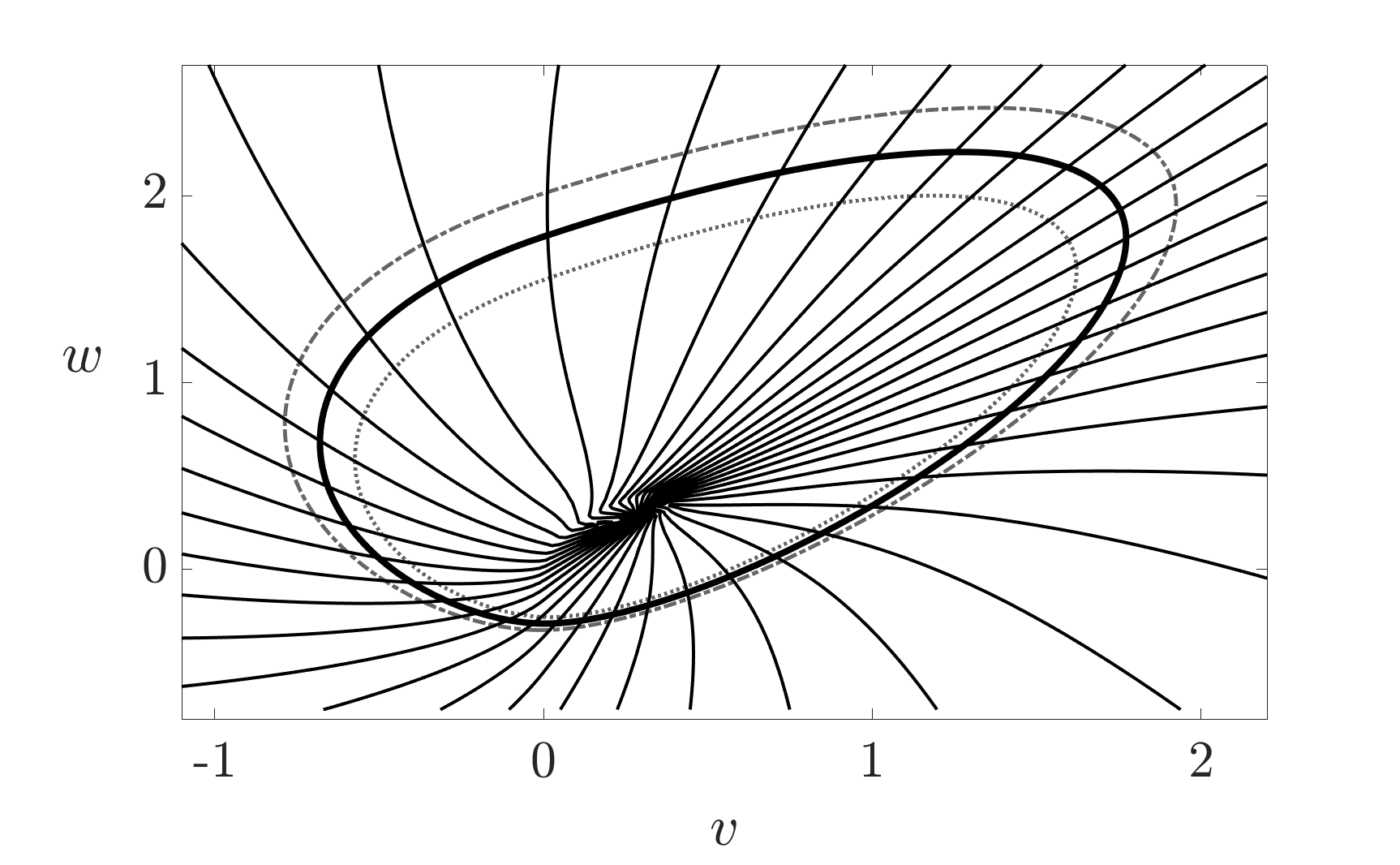}
	\caption{Isochrons and isostables for a stable periodic orbit $x^{\gamma}(t)$ (thick black curve), with Floquet exponent $\kappa \approx -0.1534$, of the absolute model.
	 We calculate the isochrons, which we show as thin black curves, using the numerical technique in \cite{Mauroy2012}.
	  We compute the isostables, $\psi = 0.04$ (gray dotted curve) and $\psi=-0.04$ (gray dashed--dotted curves), in the neighborhood of the limit cycle using the method that we describe in \cref{sec:Phaseamplitudedynamics}. This yields $x(t) = x^{\gamma}(t)+\psi p(t)$, where $p(t)$ is the Floquet mode.
The parameter values are the same as those in \cref{Fig:Pwl2pieces}(a).
}
	\label{Fig:PWLp(t)sFourmodel}
\end{figure}

Consider a $T$-periodic hyperbolic
limit cycle for the case $g(t) = 0$.
Following P\'erez-Cervera \emph{et al.} \cite{Cervera2018}, we parametrize the limit cycle and its $(m-1)$-dimensional stable invariant manifold by writing
\begin{equation}
	\FD{}{t} \theta = \omega \,, \quad \FD{}{t} \psi_k = \kappa_k \psi_k \,, \quad k \in \{1,\ldots,m-1\} \,,
\end{equation}
where $\omega = 2 \pi /T$ and $\kappa_k$ is the $k$th Floquet exponent of the limit cycle.
The dynamics for $\theta$ is uniform rotation, and the dynamics for $\psi_k$ is contraction at a rate of $\kappa_k$.
There exists an analytic map $K: \TSet \times \RSet^{m-1} \rightarrow \RSet^m$ such that $x= K(\theta, \psi_1, \ldots, \psi_{m-1})$ \cite{Cabre2005}.
From the map $K$, we define a scalar function $\Theta(x)$ that assigns a phase to any point in a neighborhood $\Omega$ of the limit cycle. The function $\Theta(x) = \theta$ if there exists $\psi_k \in \RSet$ such that $x = K(\theta, \psi_1, \ldots, \psi_{m-1})$. This function satisfies $\Theta(x(t)) = \Theta(x(t_0)) + \omega (t-t_0)$, and the isochrons are the level curves of $\Theta(x)$.  An isochron extends the notion of a phase (which occurs on a cycle) to the neighborhood $\Omega$.
Similarly, we define a set of functions $\Sigma_k(x)$ that assign a value of the amplitude variable to a point $x \in \Omega$ by setting $\Sigma_k(x) = \psi_k$ if there exists $\theta \in \TSet$ such that $x = K(\theta, \psi_1, \ldots, \psi_{m-1})$.  This function satisfies $\Sigma_k(x(t)) = \Sigma_k(x(t_0)) \e^{\kappa_k(t-t_0)}$, and the isostables are the level curves of $\Sigma_k(x)$.
Intuitively, one can consider each $\psi_{k}$ coordinate to be a signed distance from the limit cycle in a direction that is specified by $v_k$, which is the right eigenvector of $\Phi(T)$ with corresponding eigenvalue $\lambda_k$. See \cite{kvalheim2019existence,wilson2019optimal} for more details. As an illustration, we show a limit cycle
of the absolute model
in \cref{Fig:PWLp(t)sFourmodel} along with some isochrons and isostables in its neighborhood.

Knowledge of isochrons and isostables
allows us to compute
corresponding changes in phase and amplitude under a small perturbation of
$x$ to $x + \Delta x$. The change in phase is $\Delta \Theta(x) = \Theta(x+\Delta x) - \Theta(x) \approx \nabla_x \Theta(x) \cdot \Delta x$, and the change in amplitude is $\Delta \Sigma_k(x) = \Sigma_k(x+\Delta x) - \Sigma_k(x) \approx \nabla_x \Sigma_k(x) \cdot \Delta x$.
It is challenging to determine the map $K$, although it is not necessary to know it to compute
the ($m$-dimensional) \textit{infinitesimal} phase response $\mathcal{Z}$ and amplitude response $\mathcal{I}_k$, which are
\begin{equation}
	\mathcal{Z} \equiv  \nabla_{x^\gamma} \Theta(x)\,, \quad \mathcal{I}_k \equiv \nabla_{x^\gamma} \Sigma_k(x) \,.
\end{equation}
We obtain the infinitesimal phase response (iPRC\footnotemark{}) $\mathcal{Z}$ as the $T$-periodic solution of
the adjoint equation
\begin{equation} \label{Zadjoint}
	\FD{}{t} \mathcal{Z} = - {\rm D} f (x^\gamma (t))^\top \mathcal{Z}\,,
\end{equation}
with the normalization condition $\mathcal{Z}(0) \cdot f (x^\gamma (0)) = \omega$ \cite{ermentrout2010mathematical,ermentrout1991multiple,Hoppensteadt97}.
Similarly, the infinitesimal isostable responses (iIRC\footnotemark[\value{footnote}]) $\mathcal{I}_k$ satisfy
the adjoint equation
\begin{equation} 	\label{Iadjoint}
	\FD{}{t} \mathcal{I}_k =\left(\kappa_{k} I_m-\D f\left(x^{\gamma}(t)\right)^{\top}\right) \mathcal{I}_{k} \,,
\end{equation}
with the normalization condition $\mathcal{I}_k(0) \cdot v_k = 1$, where $v_k$ is the right eigenvector that is associated with the $k$th Floquet exponent of the monodromy matrix \cite{Wilson2016,wilson2018greater,monga2019phase}.

\footnotetext{The ``C'' in iPRC (and iIRC) is a historical hangover from the
phrase ``infinitesimal phase response \textit{curve}'', even though
the phase response and amplitude response are vector-valued functions.
}

For a nonsmooth system, one needs
to augment the above adjoint equations for $\mathcal{Z}$ (see \cref{Zadjoint}) and $\mathcal{I}_k$ (see \cref{Iadjoint}) to examine the behavior at any event time.
For example, Coombes \textit{et al.}~\cite{Coombes2012} determined the discontinuous iPRC for the planar PWL integrate-and-fire (IF) model by enforcing normalization conditions on both sides of a switching manifold.
Additionally, for piecewise-smooth systems, Park \textit{et al.} \cite{park2018infinitesimal} and Wilson \cite{wilson2019isostable} developed a jump operator to map the iPRC through an event by using the above normalization condition and certain linear matching conditions.
This jump operator is equal to
the inverse transpose of the saltation matrix, and related studies \cite[Chapter 5]{Coombes2023} have also made this observation.
 Using a similar approach, Chartrand \textit{et al}. \cite{chartrand2019synchronization} constructed a discontinuous iPRC for the resonate-and-fire model and Shirasaka \textit{et al.} \cite{shirasaka2017phaseAddLATER} showed how to analyze ``hybrid dynamical systems'', which include both continuous and discrete state variables \cite{aihara2010theory}.
  Ermentrout \textit{et al.} \cite{ermentrout2019recent} computed the iPRC of the Izhikevich neuron using a mixture of a jump operator and numerical computations. Wang \textit{et al.}~\cite{wang2019shape} determined the iPRC for several
 planar nonsmooth systems for a limit cycle with sliding dynamics. To do this, they used a modified saltation matrix and then related it to the
 the jump operator at the point where a sliding motion begins
 and terminates. Wang \textit{et al.} subsequently applied their approach to
 neuromechanical control problems~\cite{Wang2022}.

\begin{figure}[!h]
	\flushbottom
	\centering
	\includegraphics[width=4.5cm]{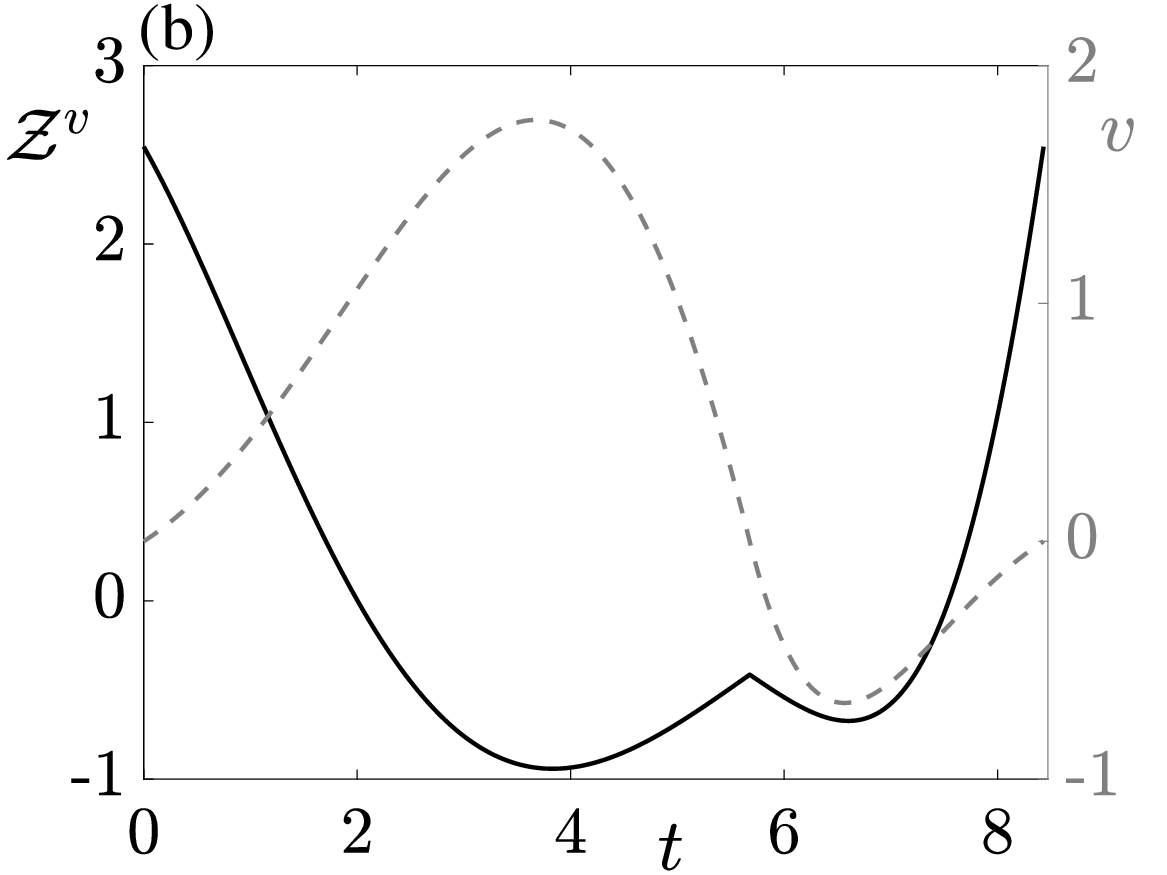}
	\includegraphics[width=4.5cm]{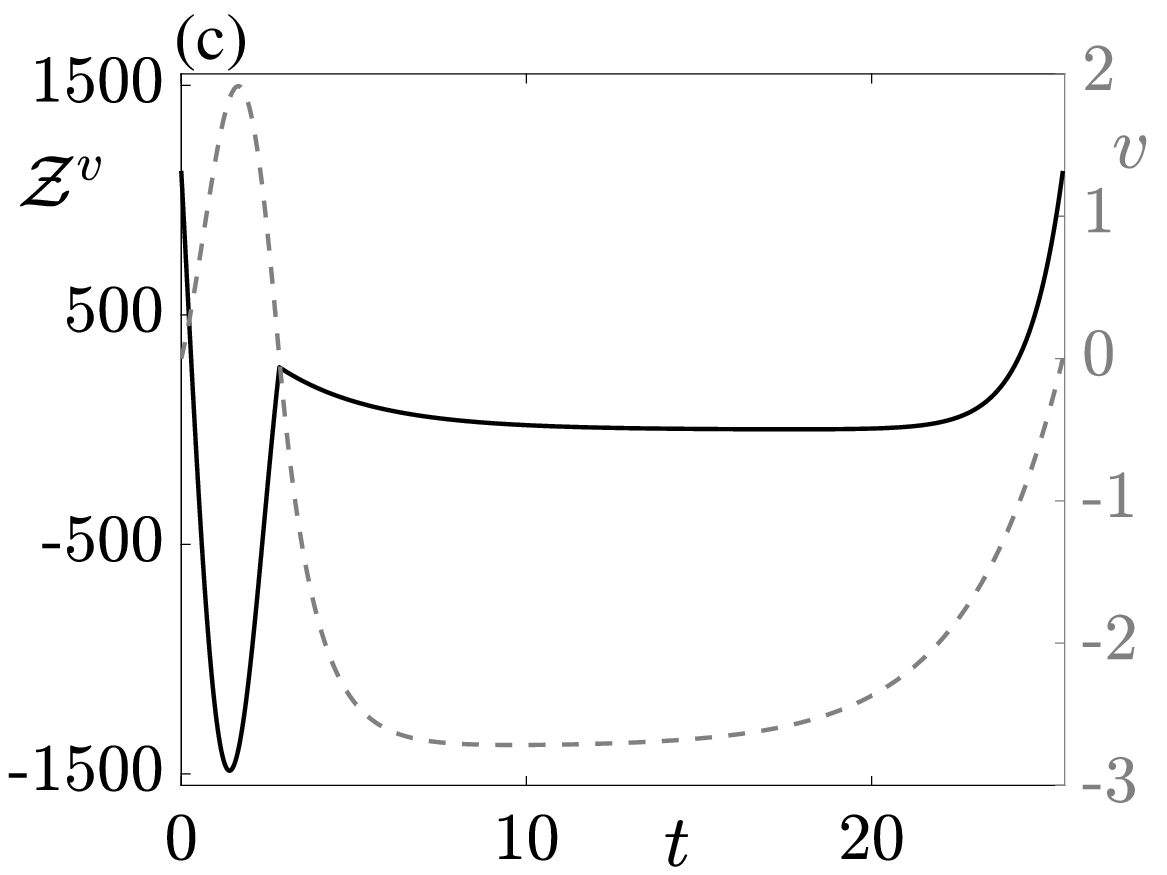}
	\includegraphics[width=4.5cm]{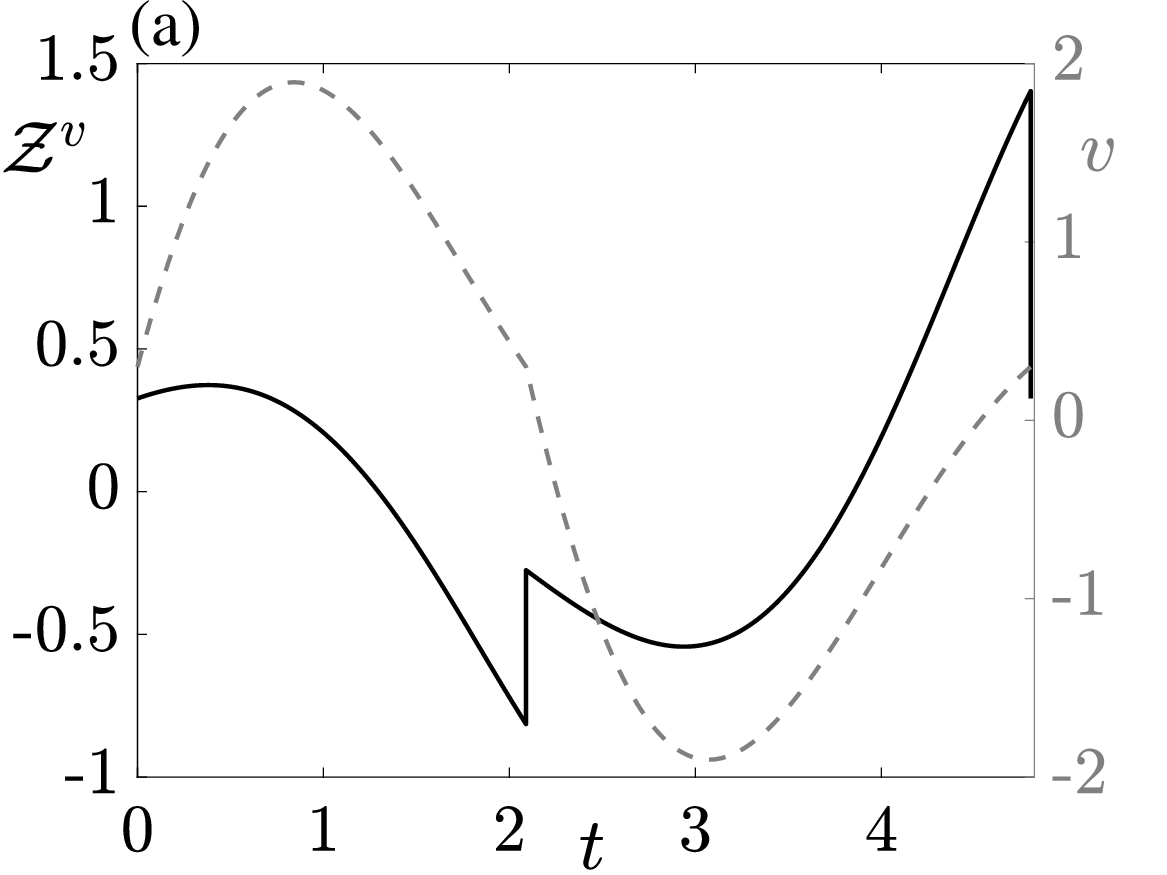}
	\includegraphics[width=4.5cm]{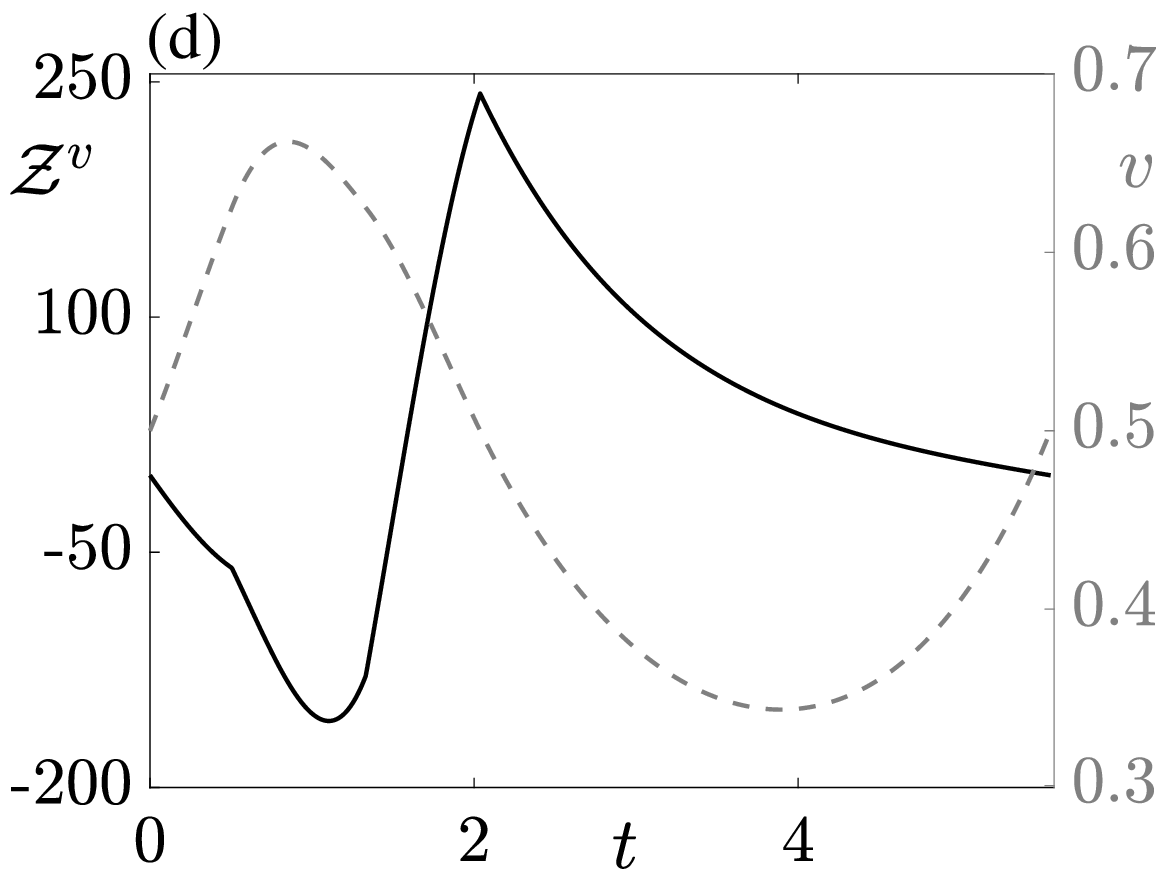}
	\caption{The $v$-component of the iPRC (solid back curve) and underlying shape of the periodic $v$-component (gray dashed curve) for
	(a) the absolute model with the same parameters as in \cref{Fig:Pwl2pieces}(a),
	(b) the PWL homoclinic model with the same parameters as in \cref{Fig:Pwl2pieces}(b),
	(c) the McKean model with the same parameters as in \cref{Fig:Pwl2pieces}(c), and
	(d) the PWL Morris--Lecar model with the same parameters as in \cref{Fig:TheMorrisLecarperiodic}.
	}
	\label{Fig:PWLPRCsFourmodel}
\end{figure}

\begin{figure}[!h]
	\centering
	\includegraphics[width=4.5cm]{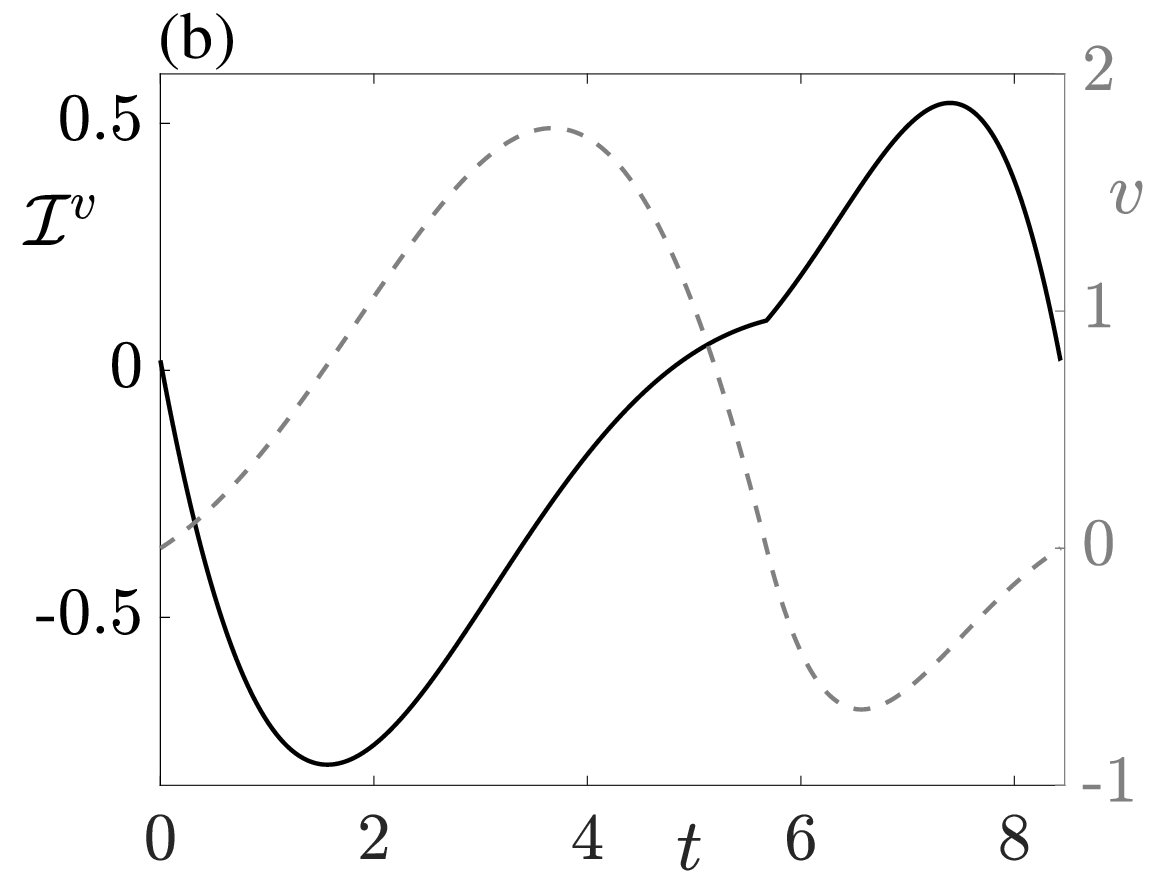}
	\includegraphics[width=4.5cm]{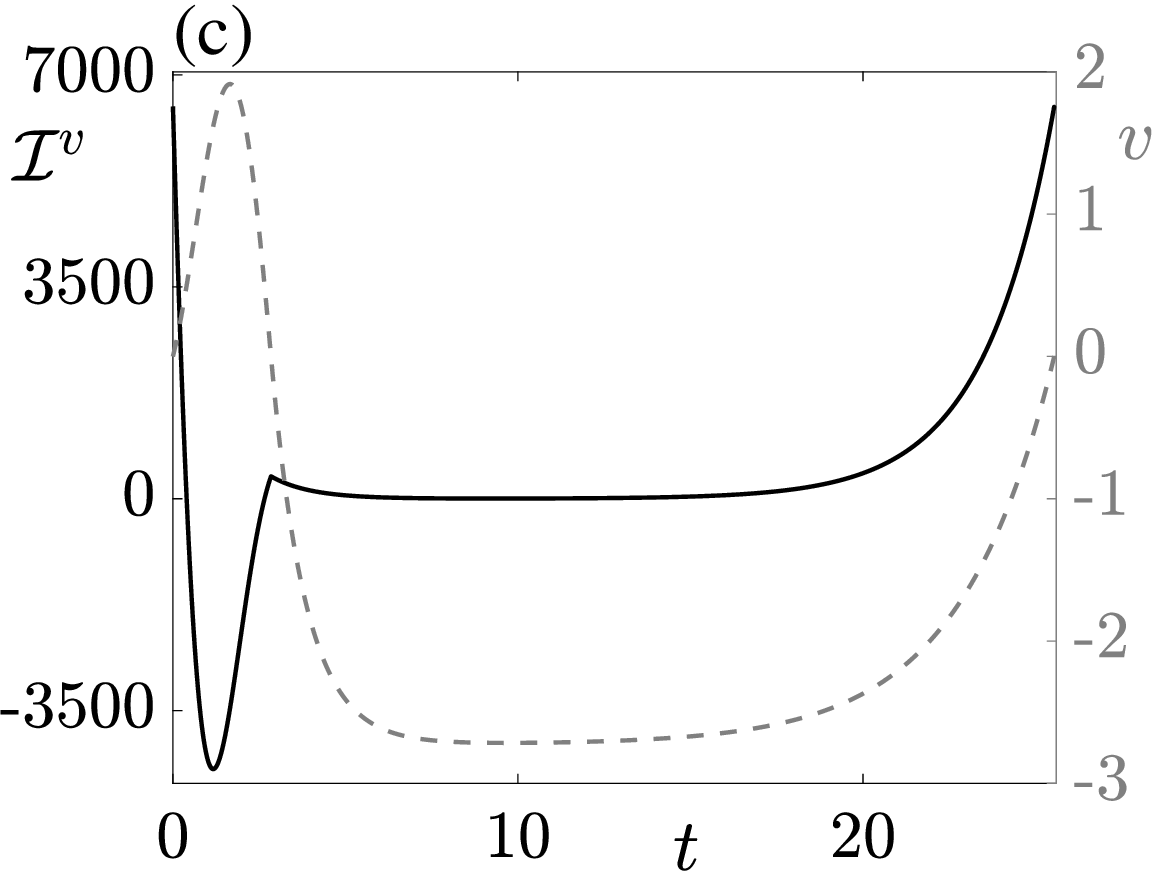}
	\includegraphics[width=4.5cm]{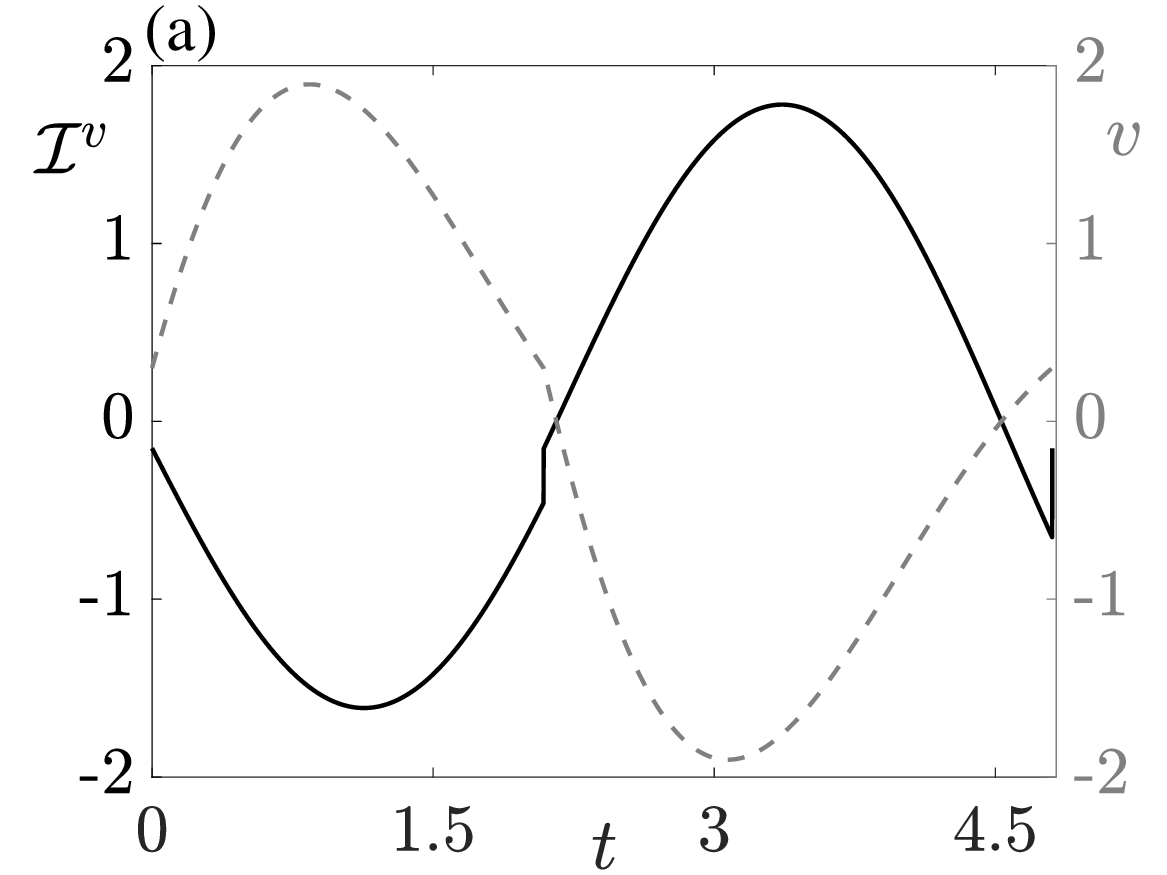}
	\includegraphics[width=4.5cm]{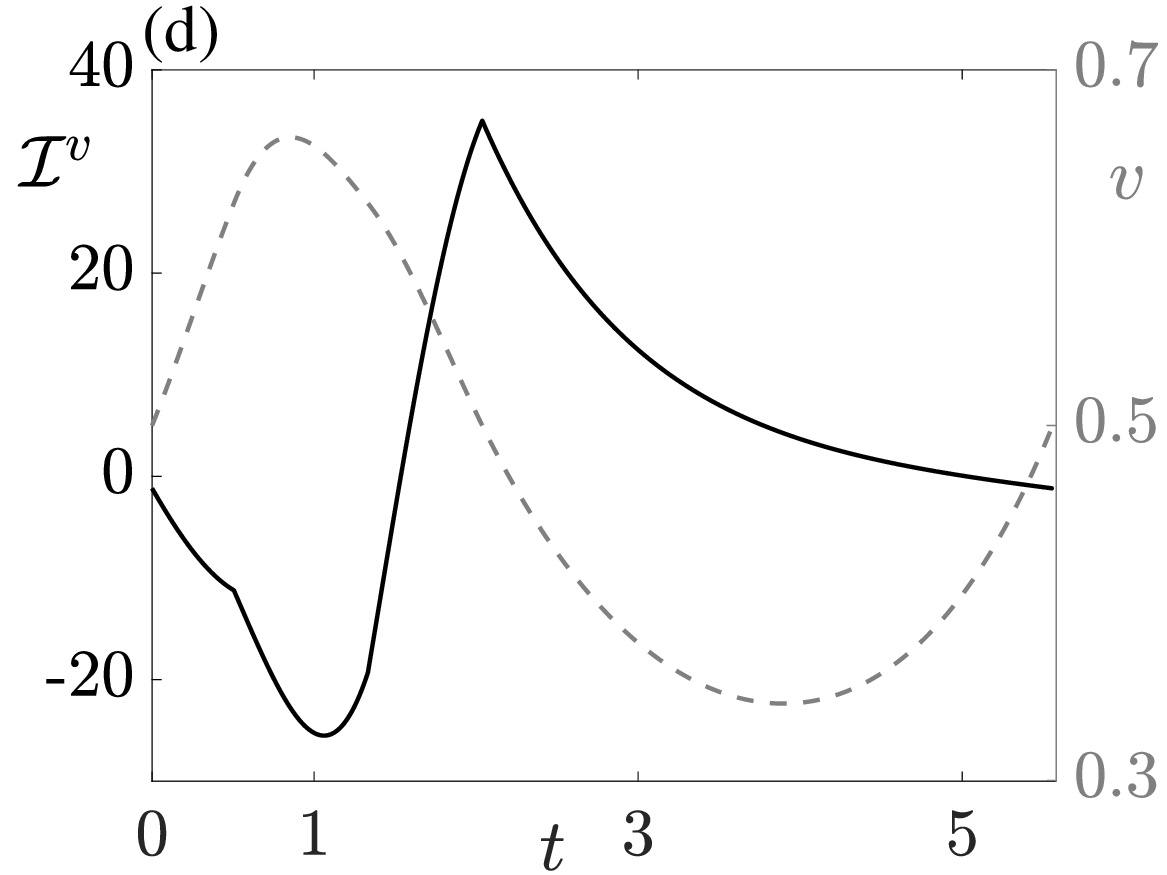}
	\caption{The $v$-component of the iIRC (solid back curve) and underlying shape of the periodic $v$-component (dashed gray curve) for
	(a) the absolute model with the same parameters as in \cref{Fig:Pwl2pieces}(a),
	(b) the PWL homoclinic model with the same parameters as in \cref{Fig:Pwl2pieces}(b),
	(c) the McKean model with the same parameters as in \cref{Fig:Pwl2pieces}(c), and
	(d) the PWL Morris--Lecar model with the same parameters as in \cref{Fig:TheMorrisLecarperiodic}.
	}
	\label{Fig:PWLIRCsFourmodel}
\end{figure}

Suppose that one has a matrix representation of
the iPRC's jump operator of the form $\mathcal{R}^{\top}\mathcal{Z}^+=\mathcal{Z}^-$, where $\mathcal{Z}^{-}$ denotes the iPRC immediately before an event and $\mathcal{Z}^{+}$ denotes the iPRC immediately after it. It is then perhaps simplest to construct the jump operator
by enforcing normalization across the switching manifold.
This balancing of normalization conditions at an event time $t_i$ requires $\langle \mathcal{Z}^+, \dot{x}^\gamma(t_i^+) \rangle = \langle \mathcal{Z}^-, \dot{x}^\gamma(t_i^-) \rangle$, so
\begin{equation}\label{connstantdeltatheatashowNorm}
	\langle\ \mathcal{Z}^+,\dot{x}^{\gamma}(t_{i}^{+})\rangle=\langle\ \mathcal{R}^{\top} \mathcal{Z}^+, \dot{x}^{\gamma}(t_{i}^{-})\rangle = \langle \mathcal{Z}^+, \mathcal{R}^{} \dot{x}^{\gamma}(t_{i}^{-}) \rangle \,,
\end{equation}
which yields
\begin{equation}\label{connstantdeltatheatafinalNorm}
	\langle\ \mathcal{Z}^+, \dot{x}^{\gamma}(t_{i}^{+}) -\mathcal{R}^{} \dot{x}^{\gamma}(t_{i}^{-}) \rangle=0 \,.
\end{equation}
\Cref{connstantdeltatheatafinalNorm} holds for any
$\mathcal{Z}^+$. Therefore, $\dot{x}^{\gamma}(t_{i}^{+}) =\mathcal{R}^{} \dot{x}^{\gamma}(t_{i}^{-})$. Additionally, the action of the saltation matrix on $\dot{x}^{\gamma}(t_{i}^{-})$ satisfies $\dot{x}^{\gamma}(t_{i}^{+})=S(t_i) \dot{x}^{\gamma}(t_{i}^{-})$. To see this, we multiply equation \cref{Delatxplusforminchapter} on the right by $\dot{x}^{\gamma}(t_{i}^{-})$ to obtain
\begin{align}
\label{FbeforeandAfterrelation}
	S(t_{i})\dot{x}^{\gamma}(t_{i}^{-}) &= \mathrm{D}\mathcal{J}_{\mu(i)}(x^{\gamma}(t^-_{i}))\dot{x}^{\gamma}(t_{i}^{-}) \nonumber \\
	&\quad + \frac{[\dot{x}^{\gamma}(t_{i}^{+})-\mathrm{D}\mathcal{J}_{\mu(i)}(x^{\gamma}(t^-_{i}))\dot{x}^{\gamma}(t_{i}^{-})][\nabla_{x} h_{\mu(i)}(x^{\gamma}(t^-_{i}))]^{\top}\dot{x}^{\gamma}(t_{i}^{-})}{\nabla_{x} h_{\mu(i)}(x^{\gamma}(t^-_{i})) \cdot \dot{x}^{\gamma}(t_{i}^{-})}\nonumber \\&= \mathrm{D}\mathcal{J}_{\mu(i)}(x^{\gamma}(t^-_{i}))\dot{x}^{\gamma}(t_{i}^{-})+\dot{x}^{\gamma}(t_{i}^{+})-\mathrm{D}\mathcal{J}_{\mu(i)}(x^{\gamma}(t^-_{i}))\dot{x}^{\gamma}(t_{i}^{-})\nonumber \\&=
	\dot{x}^{\gamma}(t_{i}^{+}) \,.
\end{align}
This implies that $\mathcal{R}^{} = S$, which in turn yields
\begin{equation}
	\mathcal{Z}^+=(S^{\top}(t_i))^{-1}\mathcal{Z}^- \,.
\end{equation}
An analogous
argument for the iIRC gives
\begin{equation}
\label{saltinI}
	\mathcal{I}_{k}^+=(S^{\top} (t_i))^{-1}\mathcal{I}_{k}^-.
\end{equation}
All that remains is to determine $\mathcal{Z}$ and $\mathcal{I}_{k}$ between events.  As usual, the PWL nature of \cref{Zadjoint} and \cref{Iadjoint} implies that
one can use matrix exponentials to obtain closed-form solutions.  For example,
the iPRC $\mathcal{Z}$ and iIRC $\mathcal{I}_{k}$ of the McKean, absolute, and homoclinic models are
\begin{equation}
\label{AdjointsolutionZ}
	\mathcal{Z}(t) = \begin{cases}
		G(t;-A_{1}^{\top})\mathcal{Z}(0) \,,  & 0 \leq t < T_1 \\
		G(t-T_1;-A_{2}^{\top})(S_{1}^{\top})^{-1} G(T_1;-A_{1}^{\top})\mathcal{Z}(0) \,, & T_1 \leq t < T
	\end{cases}
\end{equation}
and
\begin{equation}\label{AdjointsolutionI}
	\mathcal{I}(t) =\begin{cases}
		G(t;Q_{1})\mathcal{I}(0) \,, & 0 \leq t < T_1 \\
		G(t-T_1;Q_{2})(S_{1}^{\top})^{-1} G(T_1;Q_{1};)\mathcal{I}(0) \,, & T_1 \leq t < T \,,
	\end{cases}
\end{equation}
where $Q_{\mu}=(\kappa I_2-A_{\mu}^{\top})$.
One still needs to determine the initial data
$\mathcal{Z}(0)$ and $\mathcal{I}(0)$. To do this, one satisfies the normalization condition and the requirement that responses are periodic.  For example, for \cref{AdjointsolutionZ}, one needs\\
$\mathcal{Z}(0) = (S_{2}^{\top})^{-1}G(T_2;-A_{2}^{\top})(S_{1}^{\top})^{-1}G(T_1;-A_{1}^{\top})\mathcal{Z}(0)$ and $\mathcal{Z}_1(0) \dot{v}^\gamma(0)+\mathcal{Z}_2(0) \dot{w}^\gamma(0) = \omega$.  One then solves this pair of simultaneous linear equations (e.g., using Cramer's rule, as was done in \cite{Coombes2008}) to determine the initial data $\mathcal{Z}(0) = (\mathcal{Z}_1(0),\mathcal{Z}_2(0))$. One analogously determines $\mathcal{I}(0)$ using $\mathcal{I}(0)=(S_{2}^{\top})^{-1}G(T_2;Q_{2})(S_{1}^{\top})^{-1}G(T_1;Q_{1})\mathcal{I}(0)$ and $\mathcal{I}_k(0) \cdot v_k = 1$. One can follow the same procedure for models with as many regions as desired (e.g., for the
PWL Morris--Lecar model, which has four regions).
In \Cref{Fig:PWLPRCsFourmodel} and \cref{Fig:PWLIRCsFourmodel}, we show plots of iPRCs and iIRCs, respectively, that we construct using this method for several PWL models.
{Sayli \textit{et al.}~\cite{Sayli2021PWL} used direct numerical computations to confirm the shapes of these responses.
The similarity between the shapes of some iPRCs and iIRCs, such as that between \cref{Fig:PWLPRCsFourmodel}(d) and \cref{Fig:PWLIRCsFourmodel}(d) for the PWL Morris--Lecar model, was seen previously in studies
of certain smooth models \cite{Gonzales2019}. Indeed, comparing the responses that we have constructed
with those for smooth models~\cite{Gonzales2019,monga2021augmented} illustrates that a PWL approach can successfully capture the qualitative response features of their smooth counterparts.}

\subsection{Phase--amplitude dynamics}
\label{sec:Phaseamplitudedynamics}

With the results from \cref{sec:response}, we are in a position to construct the phase dynamics and amplitude dynamics for weak forcing with $g \neq 0$.
In the neighborhood of a stable limit cycle, we expand the gradients of $\Theta(x)$ and $\Sigma(x)$ and write
\begin{align}
		\nabla_{(x^{\gamma}+\Delta x)} \Theta(x)
	&= \mathcal{Z}(\theta)+H_{\Theta, x^{\gamma}} \Delta x+\mathcal{O}\left(|| \Delta x||^{2}\right) \,, \label{Sec:bacgroundIsostableSOrderTaylorphase} \\
	\nabla_{(x^{\gamma}+\Delta x)} \Sigma_k(x) &= \mathcal{I}_{k}(\theta)+H_{\Sigma_{k}, x^{\gamma}} \Delta x+\mathcal{O}\left(\parallel\Delta x\parallel^{2}\right) \,,  \label{Sec:bacgroundIsostableSOrderTaylorAplitude}
\end{align}
where $H_{\Theta, x^{\gamma}}$ and $H_{\Sigma_{k}, x^{\gamma}}$ are
the Hessian matrices of second derivatives of $\Theta$ and $\Sigma_{k}$, respectively, evaluated at the limit cycle $x^{\gamma}$.
Close to a periodic orbit, we use Floquet theory \cite{perko2013differential} to write
\begin{equation}\label{Sec:bacgroundDelatx}
	\Delta x\left(\theta, \psi_{1}, \ldots, \psi_{m-1}\right)=\sum_{k=1}^{m-1}\psi_{k} p_{k}(\theta/\omega) \,,
\end{equation}
where $p_k(t)= \e^{-\kappa_k t} \Phi(t) v_k$.

Using the chain rule, we see
that $\dot{\theta} =\nabla_{(x^{\gamma}+\Delta x)} \Theta(x) \cdot \dot{x}$ and $\dot{\psi}_k =\nabla_{(x^{\gamma}+\Delta x)} \Sigma_k(x) \cdot \dot{x}$ in the neighborhood of the limit cycle. Therefore,
equations
\cref{Sec:bacgroundIsostableSOrderTaylorphase}, \cref{Sec:bacgroundIsostableSOrderTaylorAplitude}, and \cref{Sec:bacgroundDelatx} yield a phase--amplitude approximation of the full dynamics that is accurate to second order. This approximation is
\begin{align}
	\FD{\theta}{t} &=\omega+\left( \mathcal{Z}(t)+\sum_{k=1}^{m-1}\left[\mathcal{B}^{k}(t) \psi_{k}\right]  \right) \cdot g(t)\label{Sec:bacgroundReducedfinalPhase} \,,  \\
	\FD{\psi_{k}}{t} &=\kappa_{k} \psi_{k}+\left( \mathcal{I}_{k}(t)+\sum_{l=1}^{m-1}\left[\mathcal{C}_{k}^{l}(t)\psi_{l}\right] \right) \cdot g(t) \,, \label{Sec:bacgroundReducedfinalAmplitude}
\end{align}
where we define the notation $\mathcal{B}^{k}(t) \equiv H_{\Theta, x^{\gamma}} p_{k}(t)$ and $\mathcal{C}_{k}^{l}(t) \equiv H_{\Sigma_{k}, x^{\gamma}}p_{l}(t)$ and we enforce the conditions
\begin{align}
\label{Sec:backBModifiedNormal}
	-\mathcal{Z}(\theta(t))^{\top} \D f\left(x^{\gamma}(t)\right) p_{k}(t)&=f\left(x^{\gamma}(t)\right)^{\top} \mathcal{B}^{k}(t) \,, \\
	\label{Sec:backCmodifiedNormal}
	\mathcal{I}_{k}(\theta(t))^{\top}\left(\kappa_{k} I_m-\D f\left(x^{\gamma}(t)\right)\right) p_{l}(t)&=f\left(x^{\gamma}(t)\right)^{\top} \mathcal{C}_{k}^{l}(t) \,.
\end{align}
Following Wilson and Ermentrout~\cite{wilson2019isostable}, one can show that $\mathcal{B}^{k}$ and $\mathcal{C}_{k}^{l}$ satisfy
\begin{align}
	\label{Sec:backBmanupilated}
		\FD{}{t} \mathcal{B}^{k}&=-\left(\D f^{\top}\left(x^{\gamma}(t)\right)+\kappa_{k} I_m \right) \mathcal{B}^{k} \,, \\
		\FD{}{t} \mathcal{C}_{k}^{l}&=-\left(\D f^{\top}\left(x^{\gamma}(t)\right)+\left(\kappa_{l}-\kappa_{k}\right) I_m\right) \mathcal{C}_{k}^{l} \,,
\label{Sec:backCmodifiedEvolve}
\end{align}
and we have used the fact that the Hessian of a vector field vanishes for PWL dynamical systems.
Importantly, because the system is PWL, we again use matrix exponentials to construct explicit formulas for $p_k(t)$ (by first solving the variational equation \cref{variational} for $\Phi(t)$), $\mathcal{B}^{k}(t)$, and $\mathcal{C}_{k}^{l}(t)$ (which are all $T$-periodic), being mindful to incorporate appropriate jump conditions.

As we show in \cref{JumpinB}, the jump condition on $\mathcal{B}$ for the transition across a switching manifold is
\begin{align}\label{BSOrderTaylorphasePlanar8}
	\mathcal{B}^+=(S^{\top}(t_i))^{-1}\mathcal{B}^-+C^{-1} (t_i) \eta (t_i) \,,
\end{align}	
where we have suppressed the $k$ indices, $C(t_i)$ and $\eta(t_i)$ are
\begin{equation}
\label{sigmavectorB}
	C(t_i)=\begin{bmatrix}
		\dot{v}^{\gamma}(t_{i}^{+}) & \dot{w}^{\gamma}(t_{i}^{+}) \\
		0 & 1
	\end{bmatrix}\,, \quad
	\eta(t_i)=\begin{bmatrix}
		\mathcal{Z}^- \cdot (A_{\mu(i)}p(t^-_{i}))-\mathcal{Z}^+ \cdot (A_{\mu(i+1)}p(t^+_{i})) \\
		\frac{p^{v}(t^{-}_{i})}{\dot{v}^{\gamma}(t_{i}^{-})}(A_{\mu(i)}^{\top}\mathcal{Z}^--A_{\mu(i+1)}^{\top}\mathcal{Z}^+)\cdot (0,1)
	\end{bmatrix}
\end{equation}
for a planar system, and $p^v$ denotes the $v$-component
of $p$.
Similarly, the jump condition for $\mathcal{C}$ for the transition across a switching boundary is
\begin{align}
\label{CSOrderTaylorphasePlanar8}
	\mathcal{C}^+=(S^{\top}(t_i))^{-1}\mathcal{C}^-+C^{-1} (t_i)\zeta(t_i) \,,
\end{align}	
where we have again suppressed the $k$ indices and
\begin{equation}\label{rhovectorC}
	\zeta(t_i)=\begin{bmatrix}
		\mathcal{I}^+ \cdot [(\kappa_{} I_2-A_{\mu(i+1)})p(t^+_{i})]-\mathcal{I}^- \cdot [(\kappa_{} I_2-A_{\mu(i)})p(t^-_{i})] \\
		\frac{p^{v}(t^{-}_{i})}{\dot{v}^{\gamma}(t_{i}^{-})}\left[(A_{\mu(i)}^{\top}-\kappa I_2)\mathcal{I}^--(A_{\mu(i+1)}^{\top}-\kappa I_2)\mathcal{I}^+\right]\cdot (0,1)
	\end{bmatrix} \,.
\end{equation}

\begin{figure}[!h]
	\centering
	\includegraphics[width=5.5cm]{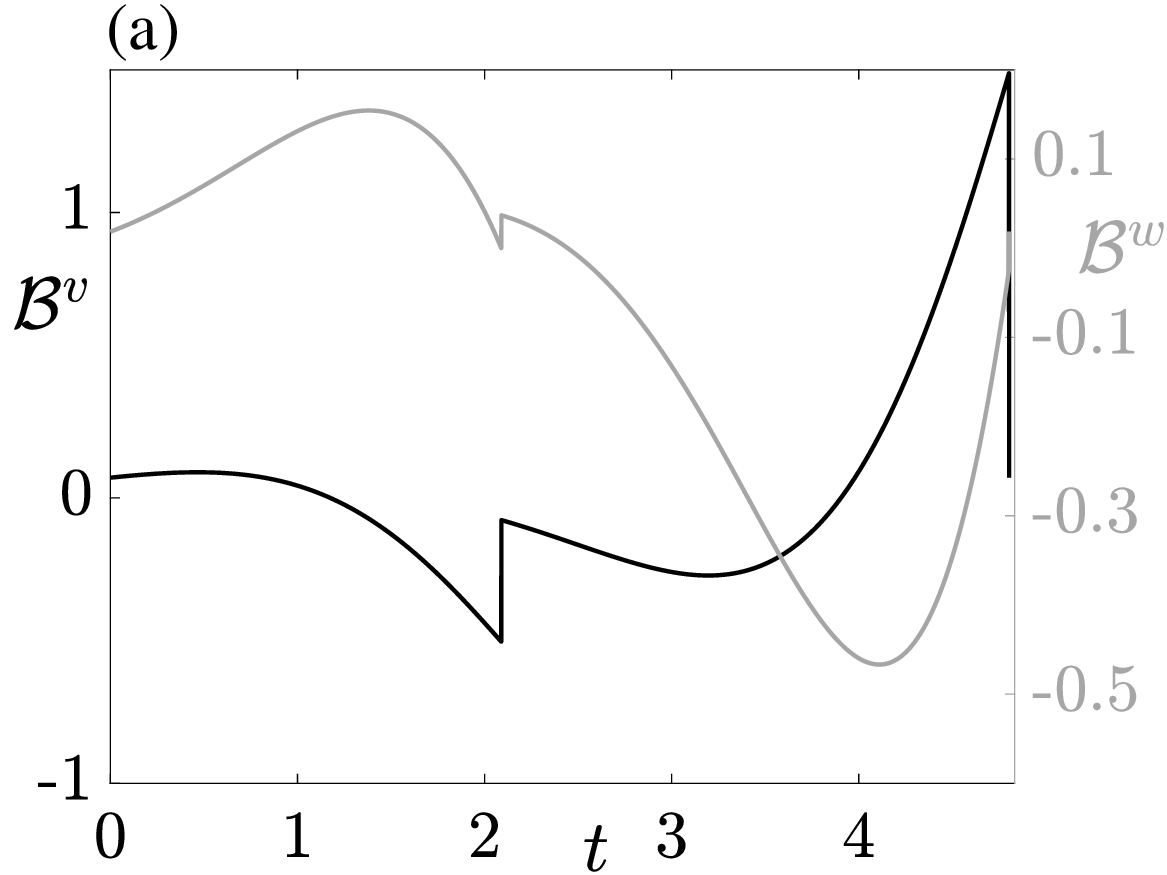}
	\includegraphics[width=5.5cm]{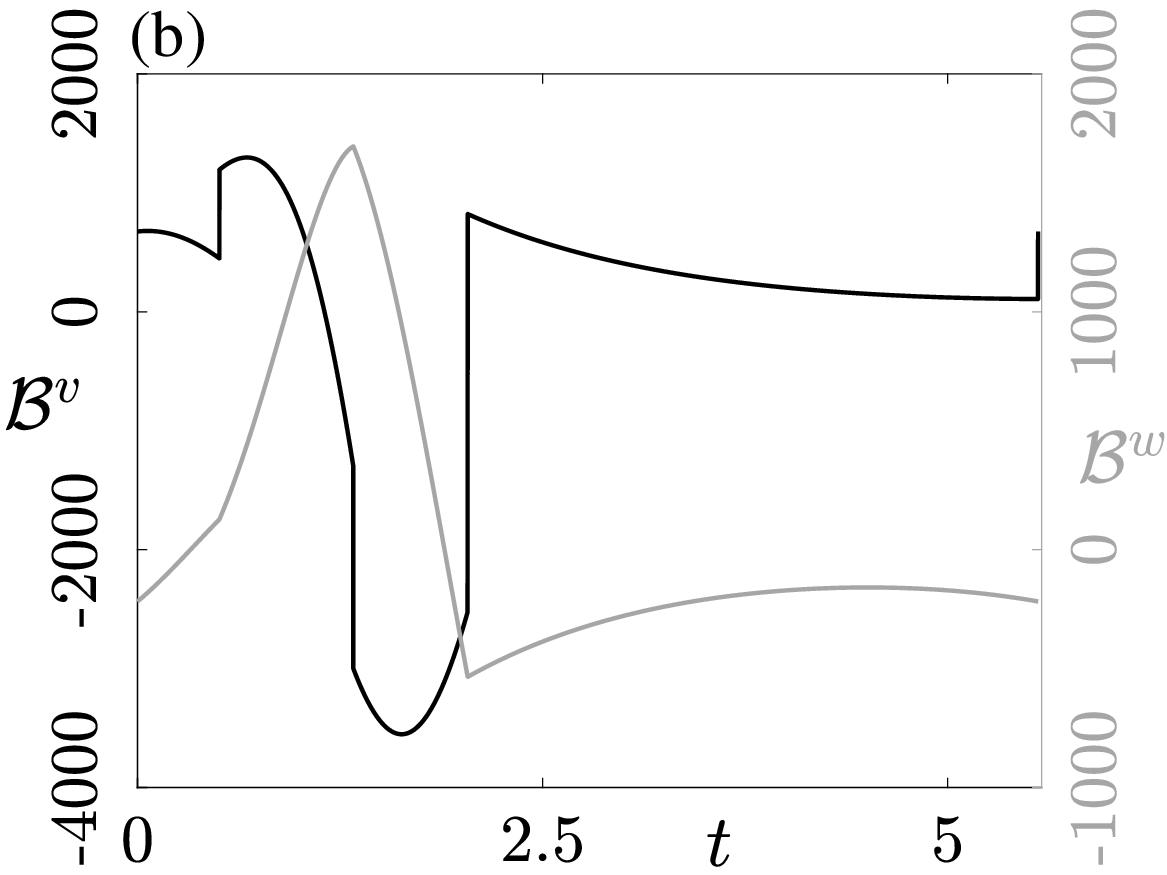}
	\caption{A plot of 
	$\mathcal{B}(t)$, with the $v$ and $w$ components on the left and right vertical axes, respectively. (a) The McKean model with the same parameters as in \cref{Fig:Pwl2pieces}(c). (b) The PWL Morris--Lecar model with the same parameters as in \cref{Fig:TheMorrisLecarperiodic}.
	}
	\label{Fig:PWLB(t)s}
\end{figure}

\begin{figure}[!h]
	\centering
	\includegraphics[width=5.5cm]{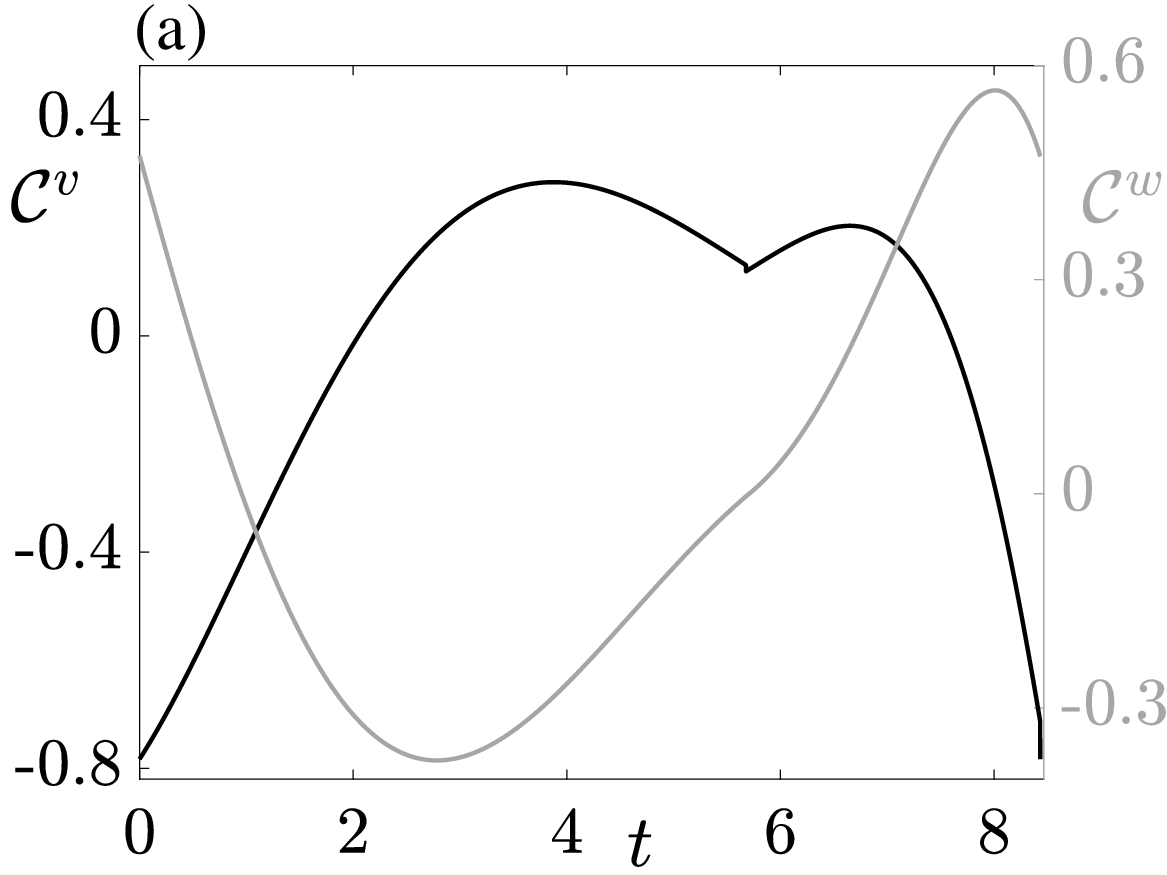}
	\includegraphics[width=5.5cm]{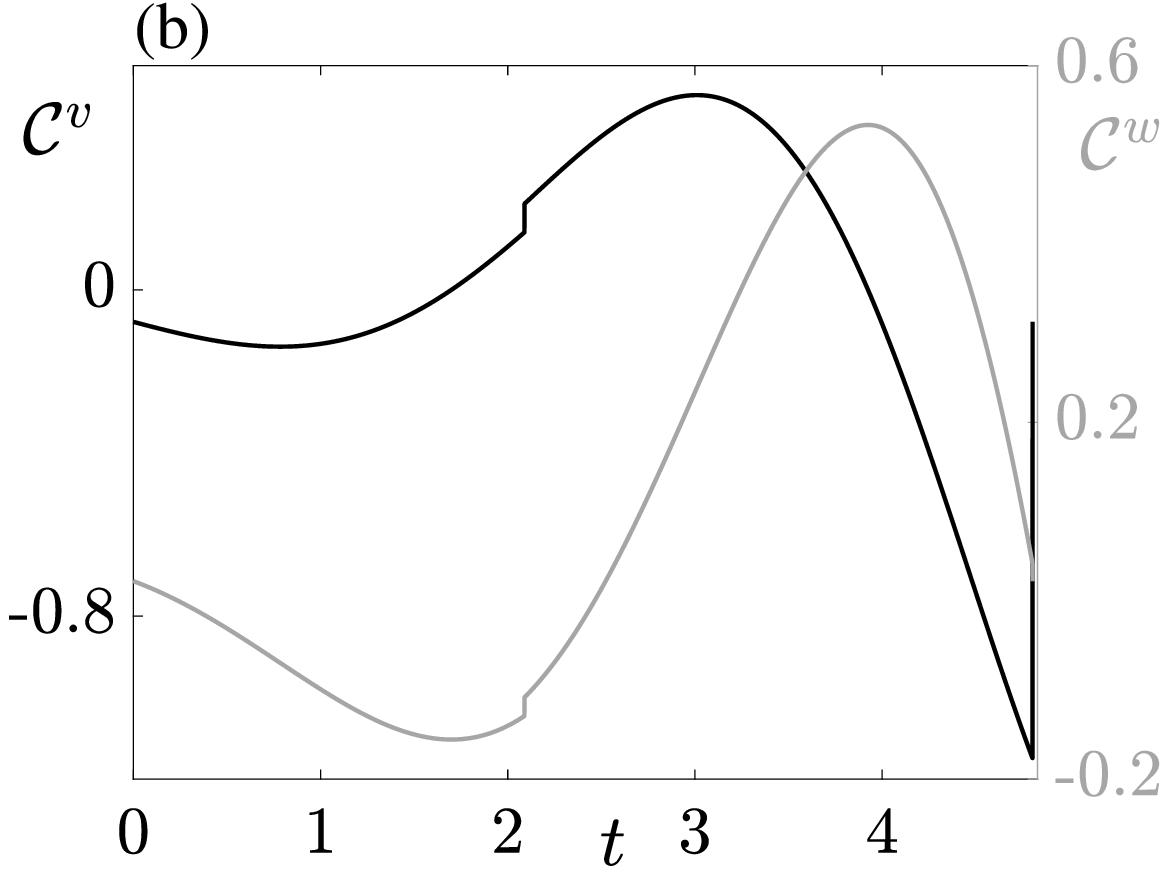}
	\caption{A plot of 
	$\mathcal{C}(t)$, with the $v$ and $w$ components on the left and right vertical axes, respectively. (a) The absolute model with the same parameters as in \cref{Fig:Pwl2pieces}(a). (b) The McKean model with the same parameters as in \cref{Fig:Pwl2pieces}(c).
	}
	\label{Fig:PWLC(t)s}
\end{figure}

For example, for PWL models with two zones (such as the McKean model), the above method
yields the following explicit formulas:
\begin{equation}\label{AdjointsolutionP}
	p(t) = \begin{cases}
	\e^{-\kappa t} G(t;A_1) \widetilde{v} \,, & 0 \leq t < T_1 \\
	\e^{-\kappa t} G(t-T_1;A_2) S(t_1) G(T_1;A_1)  \widetilde{v} \,, & T_1 \leq t < T \,,
	\end{cases}
\end{equation}
where $\widetilde{v}$ is the right eigenvector that is associated with the nontrivial Floquet exponent and
\begin{equation}\label{AdjointsolutionB}
	\mathcal{B}(t) = \begin{cases}
		G(t;K_1)\mathcal{B}(0) \,, &  0 < t \leq T_1 \\
		G(t-T_1;K_{2})[(S^{\top}(t_1))^{-1}G(T_1;K_{1})\mathcal{B}(0) + C_{1}^{-1}(t_1)\eta(t_1)]  \,, & T_1 \leq t <T \,,
		\end{cases}
\end{equation}
where $K_\mu= -(A_\mu^\top + \kappa I_2)$ and
\begin{equation}\label{AdjointsolutionC}
	\mathcal{C}(t) = \begin{cases}
		G(t;-A_{1}^{\top})\mathcal{C}(0) \,, & 0 < t \leq T_1 \\
		G(t-T_1;-A_{2}^{\top})[(S^{\top}(t_1))^{-1}G(T_1;-A_{1}^{\top})\mathcal{C}(0)+C_{1}^{-1}(t_1)\zeta(t_1)] \,, & T_1 \leq t <T \,.
	\end{cases}
\end{equation}
For $\mathcal{B}(t)$ and $\mathcal{C}(t)$, one can determine initial data in an analogous fashion as for equation \cref{AdjointsolutionZ} by simultaneously enforcing the periodicity constraints and conditions \cref{Sec:backBModifiedNormal} and \cref{Sec:backCmodifiedEvolve}. For further details, see \cite{Sayli2021PWL}. In \cref{Fig:PWLB(t)s} and \cref{Fig:PWLC(t)s}, we show example plots of $\mathcal{B}(t)$ and $\mathcal{C}(t)$ that we obtain with the above approach.
.

We are now ready
to examine how to use the phase and amplitude
to describe the dynamics of networks of the form \cref{Network1}.

\section{Phase-oscillator networks\label{sec:phaseoscnetworks}}

We first consider the case of strong attraction to a
limit cycle. Therefore, to leading order, we do not need to
consider amplitude coordinates.
Using \cref{Sec:bacgroundReducedfinalPhase}, we take a
leading-order approximation of \cref{Network1} with \cref{Networkinterraction} and $|\sigma| \ll 1$ to obtain
\begin{equation}
	\FD{}{t} \theta_i = \omega + \sigma \mathcal{Z}(\theta_i/\omega) \cdot \sum_{j=1}^N w_{ij} G(x^\gamma(\theta_i/\omega),x^\gamma(\theta_j/\omega)) \, , \quad i \in \{1,\ldots, N\} ,
\end{equation}
with $\theta_i \in [0, 2 \pi)$.
This reduced dynamical system
evolves on $\TSet^N$, whereas
the original dynamical system
 evolves on $\RSet^{N \times m}$.
  We obtain a further (and pragmatic) reduction to a model in terms of phase
 {differences} (rather than products of phases)
after
{averaging} over one oscillation period.
See, e.g., \cite{ermentrout1991multiple} and the review \cite{Ashwin2016}. We obtain the Kuramoto-like model \cite{Kuramoto1984}
\begin{equation}
	\FD{}{t} \theta_i = \omega  + \sigma \sum_{j=1}^N w_{ij }H(\theta_j-\theta_i)\,, \quad H(\theta) = \frac{1}{T} \int_0^{T} \mathcal{Z}(t) \cdot G(x^\gamma(t), x^\gamma(t+ \theta/\omega)) \, \d t  \,,
		\label{eq:phasedynamics}
\end{equation}
where the \textit{phase-interaction function} $H$ is $2 \pi$-periodic. We write it
as a Fourier series $H(\theta) = \sum_{n \in \ZSet} H_n \e^{{\rm i} n \theta}$, where the complex Fourier coefficients $H_n$ take the form $H_n = \mathcal{Z}_n \cdot G_{-n}$ and $\mathcal{Z}_n$ and $G_{n}$ are the corresponding vector Fourier coefficients of $\mathcal{Z}$ and $G$, respectively.
{For computationally useful representations of the coefficients,
see \cite{Coombes2008}.}

Using \cref{eq:phasedynamics}, it is straightforward to construct relative equilibria (which correspond to oscillatory network states) and determine their stability in terms of both local dynamics and structural connectivity {\cite{Ermentrout1992}}.
The structural connectivity is encoded in a graph (i.e., a structural network) with weighted adjacency matrix $A$ (i.e., coupling matrix) with entries $w_{ij}$.
For a graph of $N$ nodes, one specifies the connectivity pattern by an \emph{adjacency matrix} $w \in \RSet^{N \times N}$ (which is sometimes also called a ``coupling matrix'' or a ``connectivity matrix'') with entries $w_{ij}$. The spectrum of the graph is the set of eigenvalues of $w$. This spectrum also determines the eigenvalues of the associated combinatorial graph Laplacian $\mathcal{L}$. We denote the eigenvalues of $w$ by $\lambda_l$, with $l \in \{0,\ldots,N-1\}$; we denote the corresponding right eigenvectors by $u_l$. 

For a
phase-locked state $\theta_i(t) = \Omega t + \phi_i$ (where $\phi_i$ is the constant phase of each oscillator), one determines stability in terms of the eigenvalues of the Jacobian matrix $\widehat{H}(\Phi)$ of \cref{eq:phasedynamics}, where
$\Phi = (\phi_1, \ldots, \phi_N)$ and its components are
\begin{equation}
	[\widehat{H}(\Phi)]_{ij}=\sigma [H'(\phi_j-\phi_i) w_{ij} - \delta_{ij} \sum_{k=1}^{N} H'(\phi_k-\phi_i) w_{ik}] \,.
	\label{eq:H_Jacobian}
\end{equation}
The globally synchronous steady state, $\phi_i = \phi$ for all $i$, exists in a network with a phase-interaction function that vanishes at the origin (i.e., $H(0) = 0$) or for a network with a constant row sum (i.e., $\sum_j w_{ij} = \text{constant}$ for all $i$).
Using the Jacobian \eqref{eq:H_Jacobian}, synchrony is linearly stable if $\sigma H'(0) > 0$ and all of the eigenvalues of the structural network's combinatorial \textit{graph Laplacian} \cite{Newman2018}
\begin{equation}
	[\mathcal{L}]_{ij} \equiv -w_{ij} + \delta_{ij} \sum_{k=1}^N w_{ik}
	\label{Laplacian}
\end{equation}
lie in the right-hand side of the complex plane.
Because the eigenvalues of a graph Laplacian all have the same sign (apart from a single 0 value), stability is determined entirely by the sign of $\sigma H'(0)$.

\begin{figure}[!h]
	\centering
	\includegraphics[width=12cm]{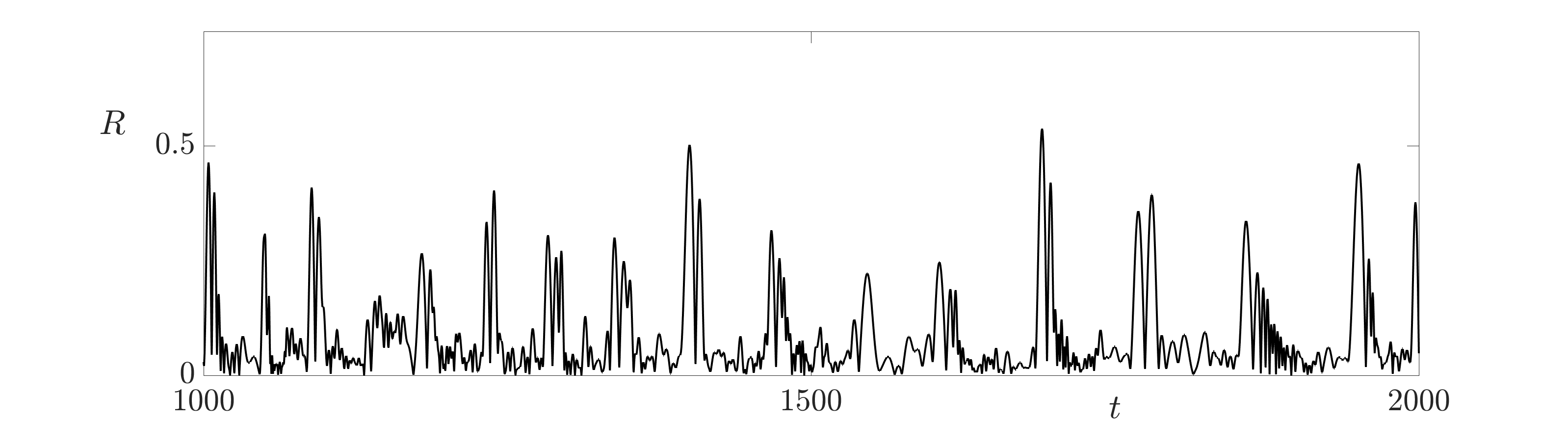}
	\caption{Evolution of the Kuramoto order parameter for the dynamics of a weakly coupled phase-oscillator network of
	$N = 1000$ PWL Morris--Lecar neurons with linear voltage coupling.
	The evolution of the Kuramoto order parameter $R = |N^{-1} \sum_{j=1}^N \e^{{\rm i} \theta_j}|$ illustrates that the system fluctuates
	between unstable states of synchrony ($R = 1$)  and asynchrony ($R=0$).
	In the upper-left panel of \cref{Fig:Hshapes}, we show the phase-interaction function $H(\theta)$ is $H_1(\theta)$.
		}
	\label{Fig:PWLMLweak}
\end{figure}

In a globally coupled network with $w_{ij}=1/N$, the graph Laplacian $\mathcal{L}$ has one 0 eigenvalue, and $(N-1)$
degenerate eigenvalues at $-1$, so synchrony is stable if $\sigma H'(0) > 0$.
In a globally connected network, one also expects the splay state $\phi_i = 2 \pi i/N$ to exist generically \cite{ashwin1992dynamics}. Additionally, in the limit $N \rightarrow \infty$, the eigenvalues
 to determine stability
 are related to the Fourier coefficients of $H$ by the equation $\lambda_n = - 2 \pi {\rm i} n \sigma H_{-n}$ \cite{Kuramoto1984}. To illustrate these results in a concrete setting, it is informative to consider a globally coupled network of PWL Morris--Lecar neurons.
  In this case, $w_{ij} = 1/N$ and $G(x_i,x_j) = v_j - v_i = v((\theta_j-\theta_i)/\omega)- v(0)$ for some common orbit $v(t)$.
    This yields $H(\theta) = \sum_{n \in \ZSet} \mathcal{Z}_n^v v_{-n} [\e^{-2 \pi {\rm i} n \theta}-1]$, where the superscript $v$ denotes voltage component and we can readily calculate the Fourier coefficients of the phase response $\mathcal{Z}$ and orbit $v$ for a PWL system~\cite{Coombes2008}.
 In the upper-left panel of \cref{Fig:Hshapes}, we show a plot of the phase-interaction function. By visually inspecting the plot, we see that $H'(0) < 0$.
 Therefore, for $\sigma > 0$, the synchronous state is unstable.
{See \cite{Kopell2002} for a geometric argument for why synchrony is unstable for gap-junction coupling when
the uncoupled oscillators are near a homoclinic bifurcation (as is the case here).} A numerical calculation of this splay state's eigenvalues also illustrate that the synchronous state is unstable.
Direct numerical simulations with large networks of oscillators
 illustrate an interesting large time-scale rhythm for which the Kuramoto synchrony order parameter $R = |N^{-1} \sum_{j=1}^N \e^{{\rm i} \theta_j}|$ fluctuates (possibly chaotically) between the value $R = 1$ for complete synchrony and the value $R = 0$ \cite{Han1995,Coombes2008}. In \cref{Fig:PWLMLweak}, we illustrate these dynamics.

\subsection{An application to the structure--function relationship in large-scale brain dynamics}

The weakly-coupled-oscillator theory that we described in \cref{sec:phaseoscnetworks} is natural for
 exploring relationships between the brain’s structural connectivity (SC) and the associated supported neural activity (i.e., its function). 
  There are studies of the SC of the human brain,
 and graph-theoretic approaches have revealed a variety of features, including
  a small-world architecture \cite{Bassett2006}, hub regions and cores \cite{Oldham2019}, rich-club organization \cite{Betzel2016}, a hierarchical-like 
  modular structure \cite{Sporns2016}, and cost-efficient
  wiring \cite{Betzel2017}.
One can evaluate the emergent brain activity that SC supports using functional-connectivity (FC) network analyses \cite{bassett2017}, which describe patterns of temporal 
coherence 
between brain regions.
Researchers have associated disruptions in SC network and FC networks with many psychiatric and neurological diseases \cite{Bassett2009}.  A measure of FC that is especially appropriate for network models of the form \cref{eq:phasedynamics} is the pairwise phase coherence
\begin{equation}
	R_{ij} = \left | \lim_{t \rightarrow \infty} \frac{1}{t } \int_0^t \cos(\theta_i(s)-\theta_j(s)) \, \d s \right | \,.
\label{FC}
\end{equation}

\begin{figure}[htbp]
\centering
\begin{subfigure}[t]{0.49\textwidth}
\centering
\begin{overpic}[clip,trim=-0.2in 0.0in -0.1in -0.3in,width=\linewidth]{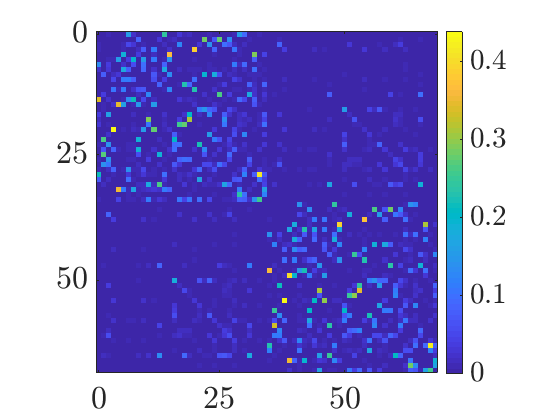}
\put(1,67){(a)}
\end{overpic}
\end{subfigure}
\begin{subfigure}[t]{0.49\textwidth}
\centering
\begin{overpic}[clip,trim=-0.2in 0.0in -0.1in -0.3in,width=\linewidth]{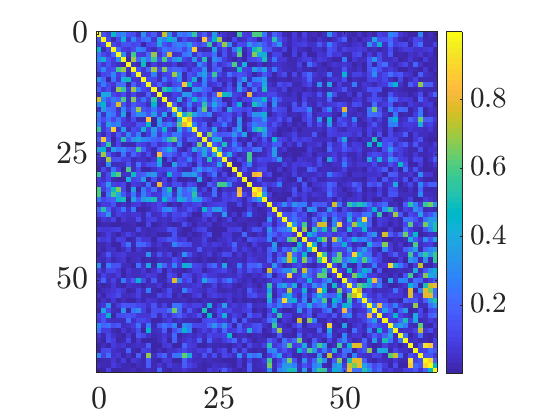}
\put(1,67){(b)}
\end{overpic}
\end{subfigure}
\caption{(a) A structural-connectivity (SC) matrix from diffusion magnetic-resonance-imaging (MRI) data that was made available through the Human Connectome Project \cite{Essen2013}. The data was processed using a probabilistic-tractography approach \cite{Glasser2013}. In this data set, the pairwise connectivity strength in a dense 60,000-network is equal to the fraction of streams that propagate from each voxel $i$
on the white/gray matter boundary and terminate at voxel $j$. This network was
then parcellated to create a 68-node network, which we employ. See~\cite{Tewarie2019} for further details. 
(b)  A functional-connectivity (FC) matrix that uses
the phase-coherence measure \cref{FC} from a phase-oscillator network with the SC 
pattern in (a) and the biharmonic phase-interaction function \cref{biharmonic}. The parameter values are $\sigma=1$, $\omega=1$, $a = 0.1$, and $r = r_c - 1/2 = \cos(a)/2 - 1/2$.
\label{Fig:SC-FC}
}
\end{figure}

Models of interacting 
{neural masses} yield natural choices of the phase-interaction function~\cite{Hlinka2012,Forrester2019}. 
For simplicity, we use a biharmonic phase-interaction function \cite{Hansel1993}
\begin{equation}
	H(\theta) = -\sin (\theta - 2\pi a) + r \sin(2 \theta)
\label{biharmonic}
\end{equation}
to illustrate how SC can influence FC.  
Using the results of the present section, we find that the stability boundary for the synchronous state is $H'(0) = 0$, which yields $r = r_c = \cos(a)/2$.  
Direct simulations of the phase oscillator network \cref{eq:phasedynamics} using human connectome data (parcellated into 68 brain regions) beyond the point of instability of the synchronous state reveal 
very rich patterns of pairwise coherence \cref{FC}.
These complicated FC dynamics reflect the fact that all eigenmodes of the graph Laplacian $\mathcal{L}$ are unstable,
leading to network dynamics that mixes all of these states. In \cref{Fig:SC-FC}, we show a plots of the emergent FC matrix from the network model \cref{eq:phasedynamics} and 
the interactions that are prescribed by the SC matrix.


\subsection{Dead zones in networks with refractoriness}

Consider once again a phase-oscillator network \eqref{eq:phasedynamics} that one obtains from a phase reduction of a weakly coupled system \cref{Network1} with coupling function \cref{Networkinterraction}.
The phase-interaction function $H(\theta)$ has a
\textit{dead zone} $U$ if $U \subset [0, 2\pi)$ is an open interval on which $H(U) = 0$ \cite{Ashwin2019}. Let $\mathrm{DZ}(H)$ denote the union of all dead zones of $H$.
When the phase difference $\theta_j - \theta_i \in \mathrm{DZ}(H)$, oscillator $i$ does not respond to changes in oscillator $j$ because the connection between them is temporarily absent. Therefore,
dead zones of the interaction function $H$ lead to an effective decoupling of network nodes for certain network states $\theta = (\theta_1, \ldots, \theta_N)$. For a
state $\theta$, the \textit{effective interaction graph} is a subgraph of the underlying structural network (which is defined by the connection strengths $w_{ij}$) that include only the edges $j \to i$ for which $\theta_j - \theta_i \notin \mathrm{DZ}(H)$.
Along solution trajectories, the effective interaction graph evolves with time. The set of subgraphs of the underlying structural network (i.e., graph) that are realizeable
 by trajectories of the system depends on the dead zones of the coupling function $H$. Ashwin \emph{et al.}~\cite{Ashwin2019} explored the interplay between dead zones of coupling functions and the realization of particular effective interaction graphs, and they began to explore how the dynamics of the associated coupled system
corresponds to changing effective interactions along a trajectory.

Because one derives the phase-oscillator network \cref{eq:phasedynamics} from the original nonlinear-oscillator network described by \cref{Network1}, \cref{Networkinterraction}, it is natural to examine 
conditions on the nonlinear oscillator dynamics that yield
a dead zone of the phase-interaction function $H$.
Both the coupling function $G$ and the iPRC $\mathcal{Z}$ influence whether or not there is an open interval $A \in [0, 2\pi)$ with
$H(\theta) = 0$ for all $\theta \in A$.
See Ashwin \emph{et al.}~\cite{Ashwin2021} for conditions for dead zones for relaxation oscillators with a separable coupling that acts only through one component of $x_i\in \mathbb{R}^m$. Here a coupling function $G$ is separable if it can be written as
\begin{equation}
G(x_i, x_j)= G^{in}(x_j) \odot G^{res}(x_i),
\end{equation} where $G^{in}$ is an input function and $G^{res}$ is the response function and $v \odot w$ denotes the Hadamard (element-wise) product of the vectors $v$ and $w$.
They also showed that pulsatile coupling can yield $\xi$-approximate dead zones $U_\xi$, on which $\sup\{ |H(\theta)| : \theta \in U_\xi\} \leq \xi$.

\begin{figure}[!h]
	\centering
	\includegraphics[height=5.5cm]{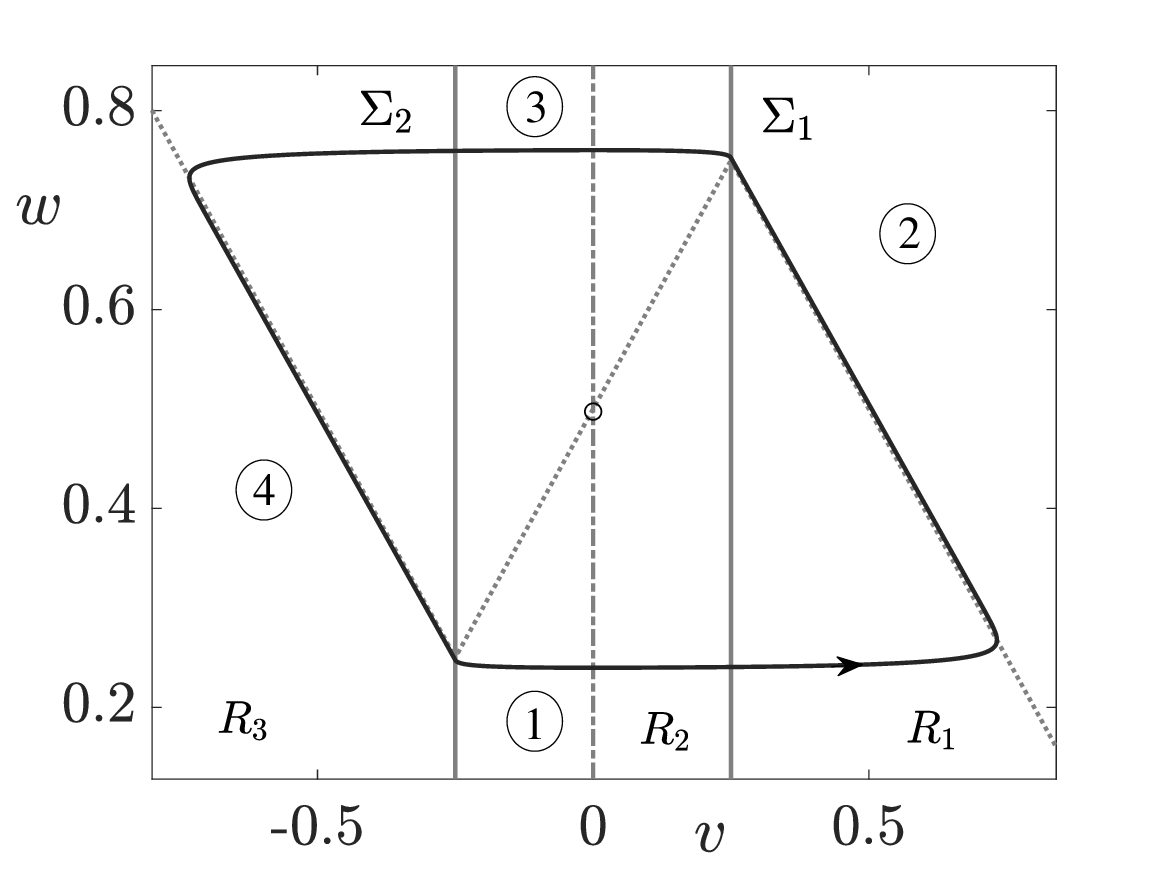}
	\caption{The phase plane for the McKean model has a $v$-nullcline
	that is a piecewise-linear approximation of a cubic (dotted line) and a linear $w$-nullcline
	(dashed--dotted line). The parameters are $C = 0.01$, $I = 0$, $\gamma = 0$, and $a = -0.5$. The solid black curve indicates a stable periodic orbit.
	}
	\label{Fig:continuousPWLML}
\end{figure}

In the present section, we show how to obtain $\xi$-approximate dead zones in the phase-interaction function for networks of synaptically coupled PWL neuron models with refractoriness.
Many
neural oscillators have a refractory period after emitting
an action potential (i.e., a nerve impulse).
During this time, the neuron does not respond to input. For a phase-oscillator, during the refractory period, input does not cause the oscillator phase to advance thereby preventing further
firing events. Therefore, the iPRC, $\mathcal{Z}(t)$ is approximately $0$ for one or more intervals $(t_1, t_2) \subset [0, T)$, where $T$ is the oscillator period. An example of a planar PWL model with such a refractory period is the continuous McKean model with ``three pieces"~\cite{McKean70, Coombes2008}. This is a PWL caricature of the FHN model, with the $v$-nullcline broken into three pieces, 
which partitions the phase space into three zones, with two switching manifolds. The dynamics of the system satisfy
\begin{equation}\label{this}
	C\dot{v}=\rho(v)-w+I \,, \quad \dot{w} = v-\gamma w\,,
\end{equation}
where $C > 0$, $\gamma \ge 0$, and (to approximate the cubic $v$-nullcline) $\rho(v)$ is given by equation~\cref{fMcKeannadMorrisLecar}. The dynamical system \cref{this} has
a stable periodic orbit when there is a single unstable equilibrium on the center
branch of the cubic $v$-nullcline. In Figure~\ref{Fig:continuousPWLML}, we show the phase portrait for parameter values with a stable periodic orbit. See \cref{Pwlmodels} for further details about the model.

\begin{figure}[!h]
	\centering
	\includegraphics[height=5cm]{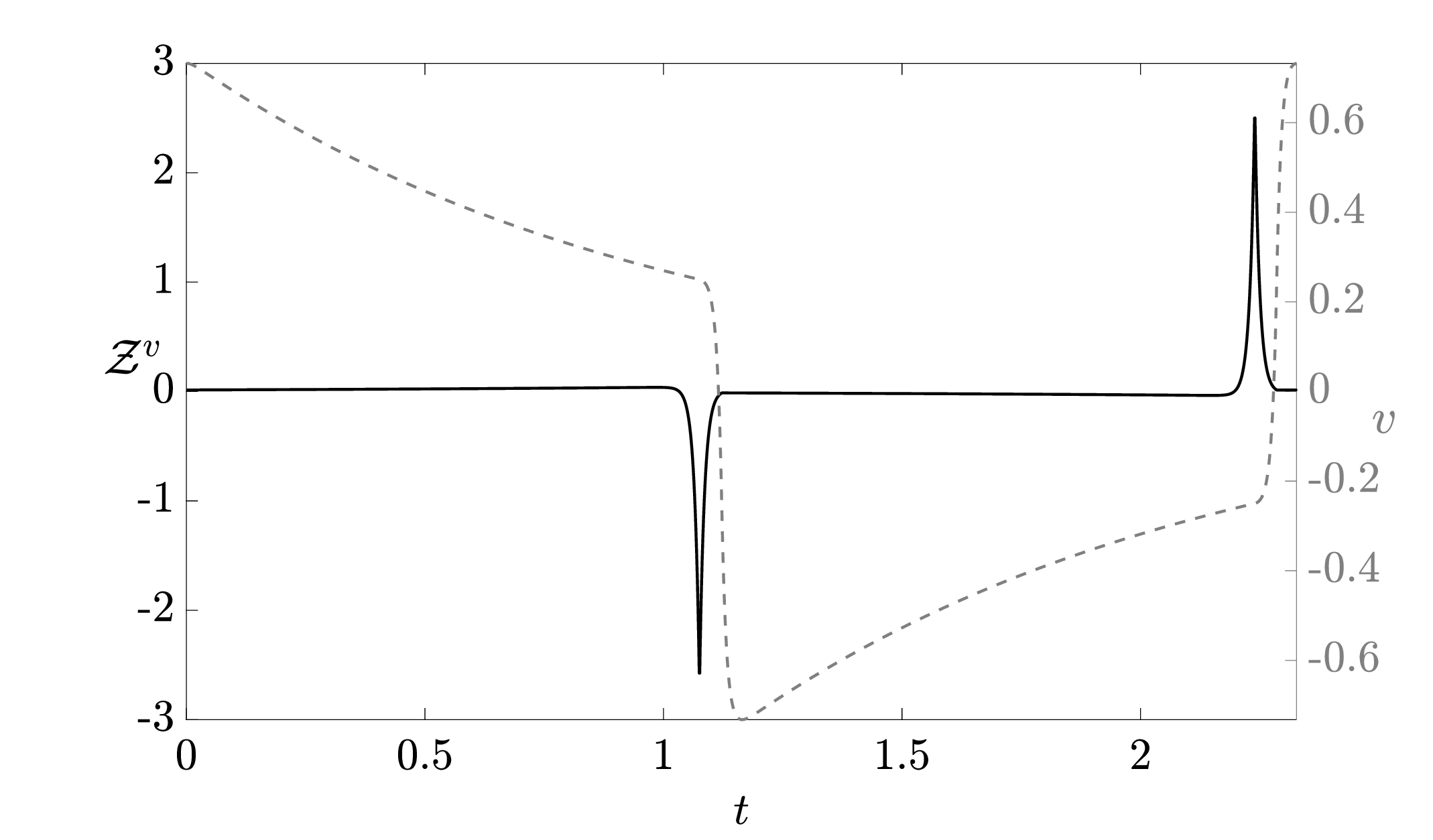}
	\caption{The $v$ component $\mathcal{Z}^v$ of the iPRC of the continuous McKean model (solid black curve) and the corresponding value of $v$ on the limit cycle (dashed gray curve). The parameter values are $C=0.01$, $I = 0$, $\gamma = 0$, and $a = -0.5$.  Observe that $\mathcal{Z}^v$ is approximately $0$ for much
	of the limit cycle, with peaks corresponding to locations near the switching planes. 
	}
	\label{fig:FHNiPRC}
\end{figure}

For $C \ll 1$, the dynamics of the voltage
$v$ are fast and the dynamics of the recovery variable $w$ are slow.
Therefore, as $C \rightarrow 0$, the system spends most of its time
on the left and right branches of the $v$-nullcline,
with fast switching between the two branches.
Consequently, $\mathcal{Z}^v$ (i.e., the $v$ component of the iPRC) is approximately $0$ for much
 of the limit cycle, with peaks corresponding to locations near the switching planes. In Figure~\ref{fig:FHNiPRC}, we show $\mathcal{Z}^v$ (and $v$ on the limit cycle) for $C=0.01$.
{In the singular $C \rightarrow 0$ limit, the iPRC for this model is discontinuous \cite{Izhikevich2000,Coombes2001}.}

We now compute the phase-interaction function for a network of $N$ synaptically coupled continuous McKean neurons with time-dependent forcing:
\begin{equation}
	C\dot{v}_i=\rho(v_i)-w_i+ I + \sigma_{j=1}^N \sum_j w_{ij} s_j(t) \,, \quad \dot{w}_i = v_i-\gamma w_i\,, \quad  i \in \{1, \ldots, N\} \,.
\end{equation}

Suppose that the synaptic input from neuron $j$ takes the standard ``event-driven''
form
\begin{equation}\label{eventsynapse}
	s_j(t) =  \sum_{p \in \ZSet} \eta(t-t_j^p) \,,
\end{equation}
where $t_j^p$ denotes the $p$th firing time of neuron $j$ and the causal synaptic filter $\eta$ describes the shape of the post-synaptic response.

\begin{figure}[!h]
	\centering
	\includegraphics[height=5cm]{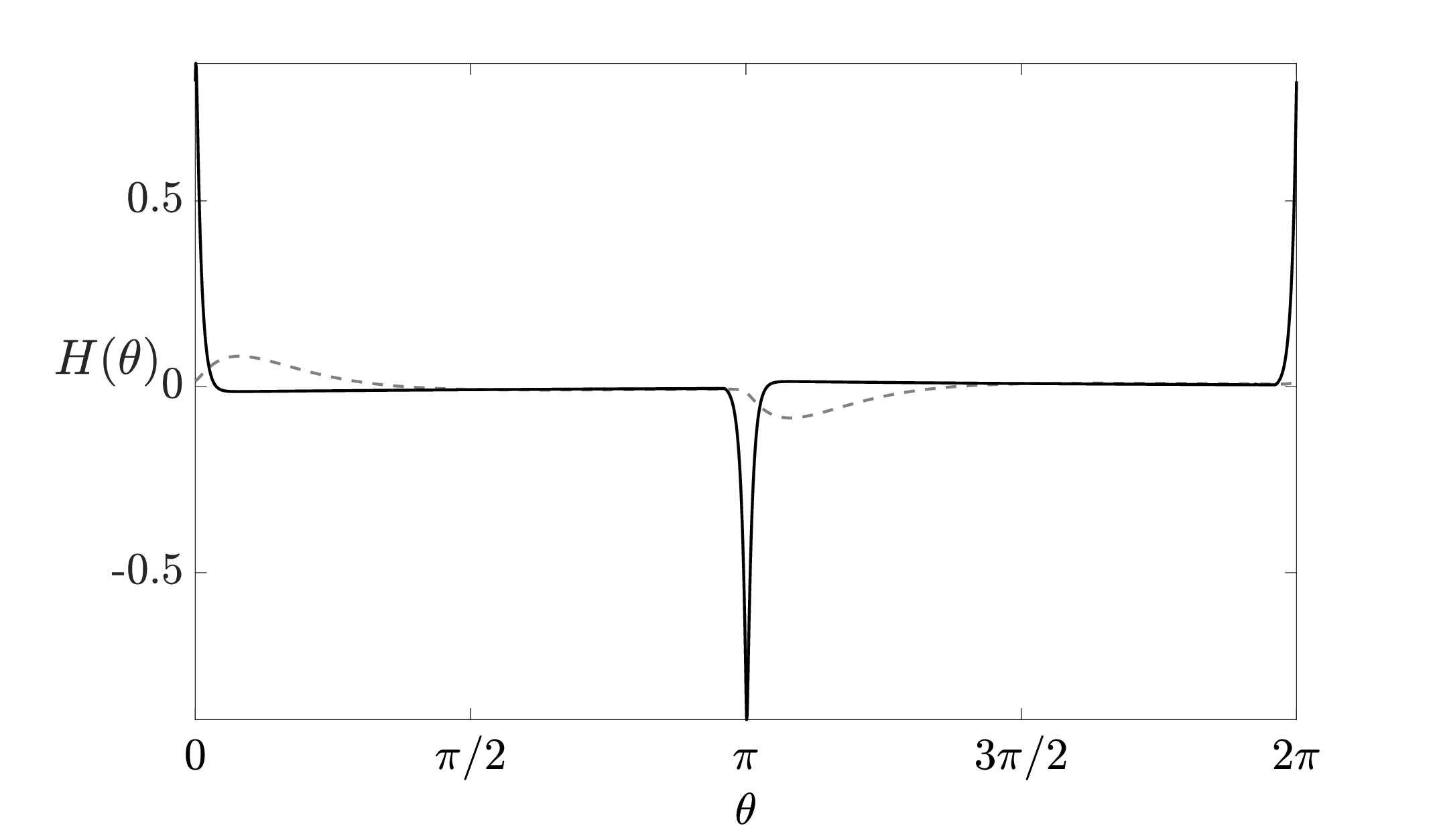}
	\caption{The phase-interaction function $H(\theta)$ for the continuous McKean model with synaptic coupling for $\alpha=1000$ (black curve, fast synapses) and $\alpha=10$ (gray curve, slower synapses). Synapse-firing events occur when $v = 0.6$ and $\dot{v} > 0$.
	The other parameter values are $I = 0$, $\gamma = 0$, and $a = -0.5$. The larger value of $\alpha$ results in a larger dead zone of $H(\theta)$.
	 }
	\label{fig:FHNinteraction}
\end{figure}

For a phase-locked system, one writes the firing times $t_j^p$
as $t_j^p = pT - \phi_j/\omega$ for a phase offset $\phi_j \in [0, 2\pi)$. Therefore, the phase-interaction function is
\begin{equation}
	H(\theta) = \frac{1}{T} \int_0^{T} \mathcal{Z}^v(\omega t-\theta)  P(\omega t) \, \d t = \frac{1}{2\pi}\int_0^\infty \mathcal{Z}^v(u-\theta) \eta (u/\omega)\, \d u \,,
\end{equation}
where $P(\phi) = \sum_{p \in \ZSet} \eta(\phi/\omega - pT)$. Because $\mathcal{Z}^v$ is $2\pi$-periodic, we can write $\mathcal{Z}^v(u) = \sum_{n \in \ZSet} \mathcal{Z}^v_n \e^{{\rm i} nu}$, where $\mathcal{Z}^v_n = (2\pi)^{-1} \int_0^{2\pi} \mathcal{Z}^v(u) e^{-{\rm i} nu} \, \d u$. Consequently,
\begin{equation}
	H(\theta) = \frac{1}{T} \sum_{n \in \ZSet} \mathcal{Z}^v_{-n} \e^{{\rm i} n \theta} \hat{\eta}(n/T) \,,
\end{equation}
where $\hat{\eta}(k) = \int_0^\infty \e^{-2\pi {\rm i} k t} \eta(t) \, \d t$ is the Fourier transform of the causal synaptic filter $\eta$.
That is, $H(\theta) = \sum_{n\in \ZSet} H_n \e^{{\rm i} n \theta}$, where $H_n = \mathcal{Z}^v_{-n} \hat{\eta}(n/T)/T$.
If we adopt the common choice
\begin{equation}\label{alphafunction}
	\eta(t) = \alpha^2 t \e^{-\alpha t} \Theta(t)  \,,
\end{equation}
where $\alpha>0$ and $\Theta$ is a Heaviside step function, then
\begin{equation}
	H_n = \frac{\alpha^2\mathcal{Z}^v_{-n}}{T(\alpha + 2\pi i n/T)^2} \,.
\end{equation}
In Figure~\ref{fig:FHNinteraction}, we show the phase-interaction function $H(\theta)$ for two values of $\alpha$.
This interaction function has two large dead zones. The larger the value of $\alpha$, the larger the dead zones of $H(\theta)$. For the chosen parameter values, the dead zones of $H$ are
symmetric for large values of $\alpha$ \cite{Ashwin2019}. (That is, $-\theta \in \mathrm{DZ}(H)$ if $\theta \in \mathrm{DZ}(H)$, then $-\theta \in DH(H)$. Such a coupling function is ``dead-zone symmetric''.) This symmetry places restrictions on
the effective interaction graphs which can be realized by the trajectories of the model. For example, if $H$ is dead-zone symmetric, then all of the effective interaction graphs for $H$ are undirected \cite[Proposition 3.7]{Ashwin2019}.
In the limit $\alpha \to \infty$, we also observe that $H(\theta)$ is a scaled version of the $v$ component of the iPRC.
This follows from
 $H(\theta) =\mathcal{Z}^v(-\theta)/T$ and the fact that $\lim_{\alpha \rightarrow \infty} \eta(t) = \delta(t)$, giving pulsatile coupling.

\section{Phase--amplitude networks\label{PWLPhaseAmplitudenetworkfromalism}}

We now consider the second-order
approximation \cref{Sec:bacgroundReducedfinalPhase}--\cref{Sec:bacgroundReducedfinalAmplitude}
that allows us to use
both phase and amplitude coordinates to treat oscillatory network dynamics. In contrast to the phase-only approach
in \cref{sec:phaseoscnetworks}, there has been much less work on the theory and applications of phase--amplitude networks although this is now 
growing, as exemplified by the work 
in \cite{wilson2019sirev}. 
Because of both this and to facilitate our exposition, we focus
on a small network of two identical planar oscillators with
linear coupling through the $v$ component. Pairs and larger networks of linearly coupled smooth Morris--Lecar neurons are considered in Nicks \emph{et al.} \cite{Nicks2023} where conditions for linear stability of various phase-locked states in globally coupled phase-amplitude networks are also derived. Our discussion parallels the one in Ermentrout \emph{et al.} \cite{ermentrout2019recent} for a smooth model of synaptically coupled thalamic neurons \cite{rubin2004high}.

Specifically, we consider equation \cref{Network1} with $N = 2$ oscillators, a coupling strength of $|\sigma| \ll 1$, and
\begin{equation}\label{phaseisostablenetworkPWLcoupling1}
	\begin{aligned}
		g_1(x_1,x_2) &= \sigma[v_2-v_1,0]^{\top} \,, \\
		g_2(x_1,x_2) &= \sigma[v_1-v_2,0]^{\top} \,.
	\end{aligned}
\end{equation}
To determine the form of $g(t)$ in the phase--amplitude equations \cref{Sec:bacgroundReducedfinalPhase}--\cref{Sec:bacgroundReducedfinalAmplitude} to obtain the corresponding phase-amplitude reduction of the network equations \cref{Network1}, we write
$v_{i}(t)=v^{\gamma}\left(\theta_{i}(t)\right)+\psi_{i}(t) p^v\left(\theta_{i}(t)\right)$ and assume that the amplitudes $\psi_i$ are $\mathcal{O}(\sigma)$.  Substituting this expression into \cref{Sec:bacgroundReducedfinalPhase} and \cref{Sec:bacgroundReducedfinalAmplitude} and
keeping terms up to order $\mathcal{O}(\sigma^2)$
yields the phase-amplitude reduced network equations
\begin{align}\label{phaseisostablenetworkPWLcoupling4}
		\dot{\theta}_{1} &=\omega+\sigma\left[h_{1}\left(\theta_{1}, \theta_{2}\right)+\psi_{1} h_{2}\left(\theta_{1}, \theta_{2}\right)+\psi_{2} h_{3}\left(\theta_{1}, \theta_{2}\right)\right]  \,,  \\
		\dot{\psi}_{1} &=\kappa \psi_{1}+\sigma\left[h_{4}\left(\theta_{1}, \theta_{2}\right)+\psi_{1} h_{5}\left(\theta_{1}, \theta_{2}\right)+\psi_{2} h_{6}\left(\theta_{1}, \theta_{2}\right)\right] \,, \notag \\
		\dot{\theta}_{2} &=\omega+\sigma\left[h_{1}\left(\theta_{2}, \theta_{1}\right)+\psi_{2} h_{2}\left(\theta_{2}, \theta_{1}\right)+\psi_{1} h_{3}\left(\theta_{2}, \theta_{1}\right)\right] \,,  \notag \\
		\dot{\psi}_{2} &=\kappa \psi_{2}+\sigma\left[h_{4}\left(\theta_{2}, \theta_{1}\right)+\psi_{2} h_{5}\left(\theta_{2}, \theta_{1}\right)+\psi_{1} h_{6}\left(\theta_{2}, \theta_{1}\right)\right] \,, \notag
\end{align}
where we give the detailed forms of the doubly $2\pi$-periodic functions $h_1, \ldots, h_6$ in \cref{sec:PhasAmpInteraction}.

\begin{figure}[!h]
	\centering
		\includegraphics[height=3.2cm]{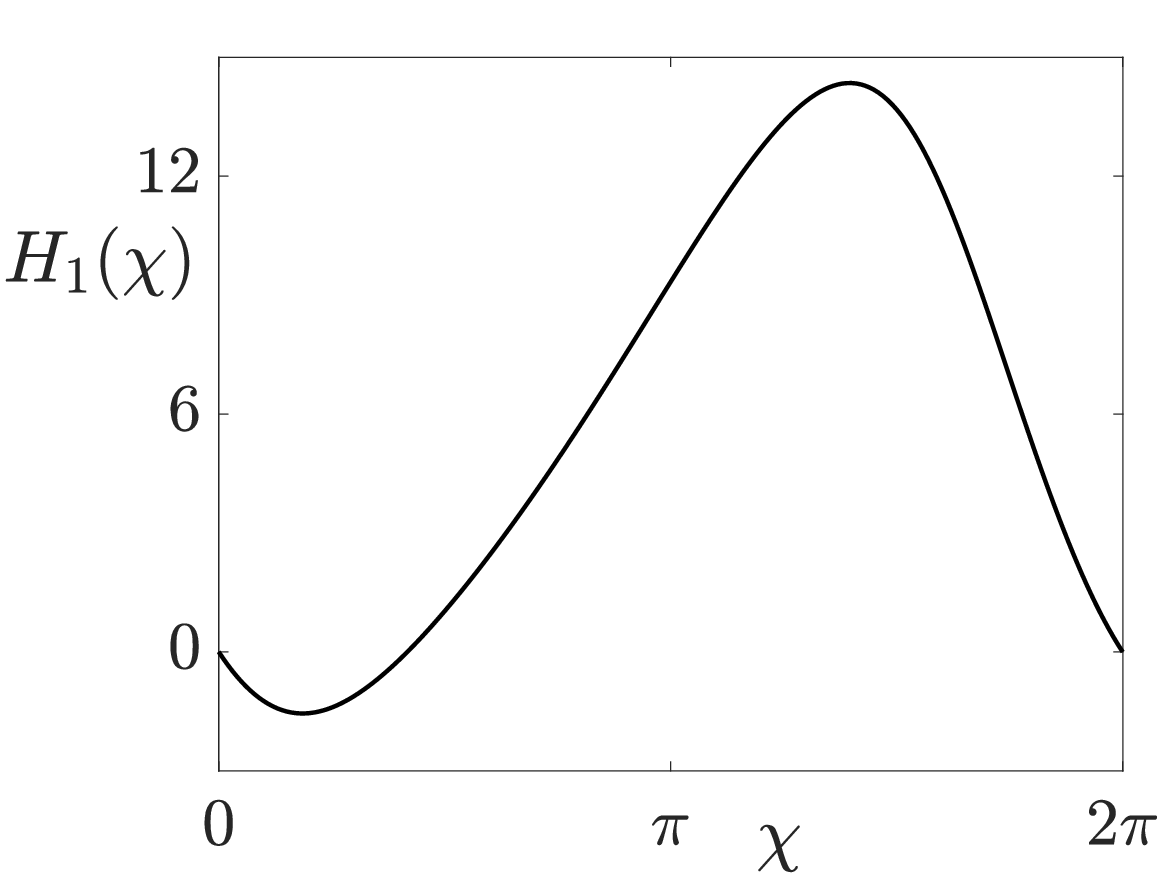}
	    \includegraphics[height=3.2cm]{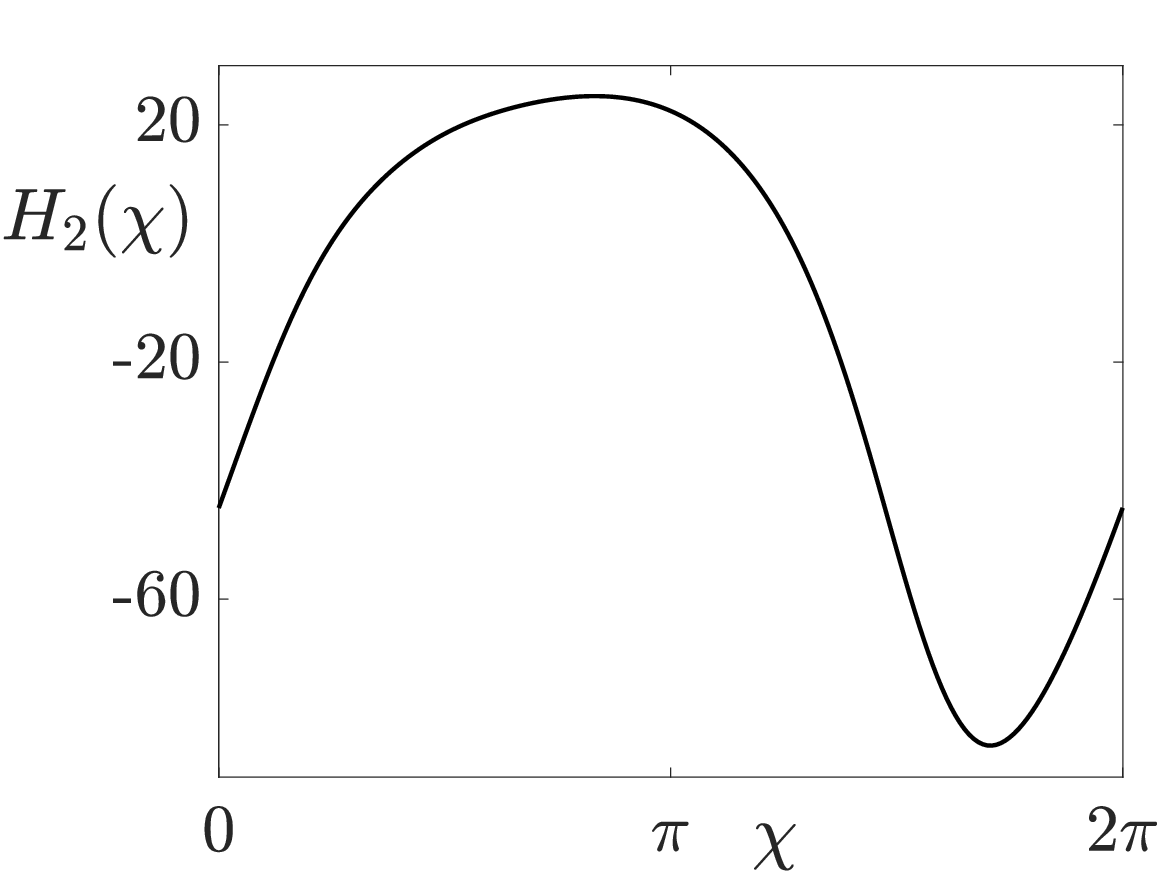}
		\includegraphics[height=3.2cm]{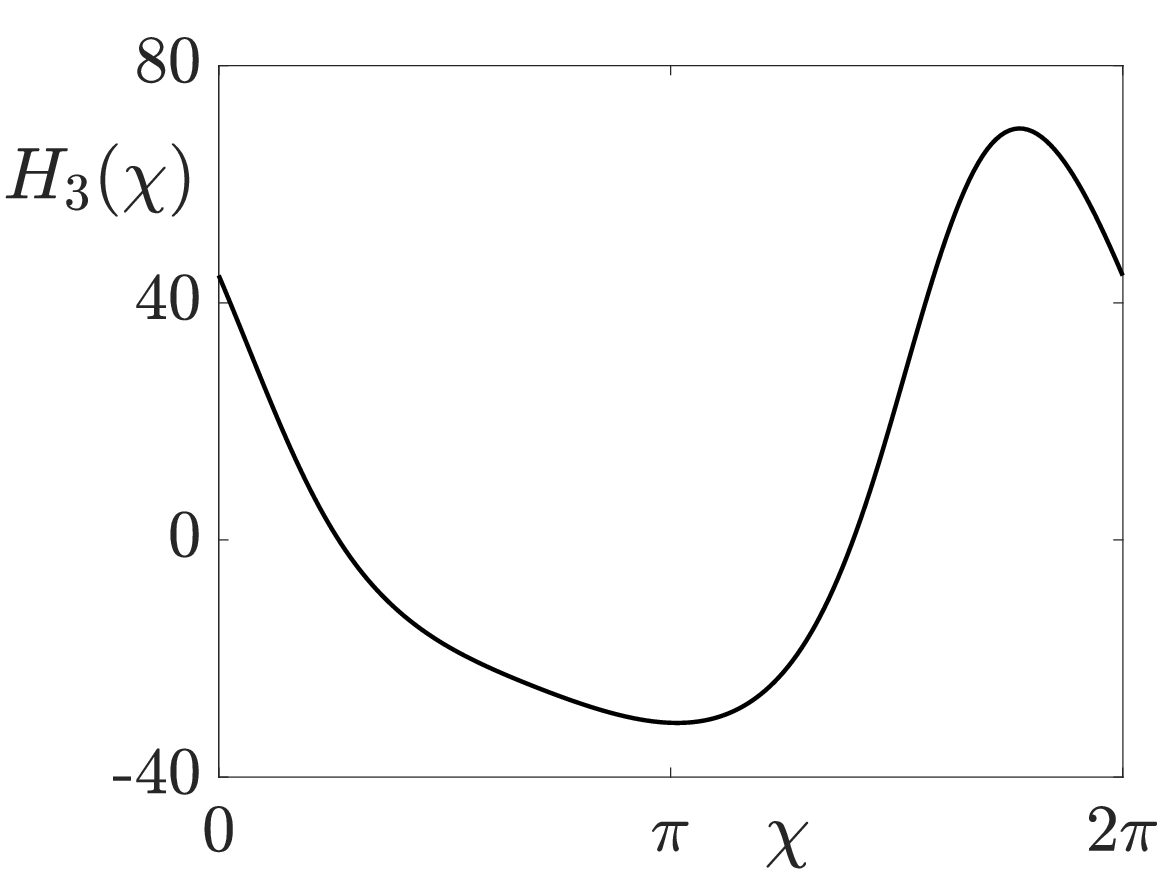}
		\includegraphics[height=3.2cm]{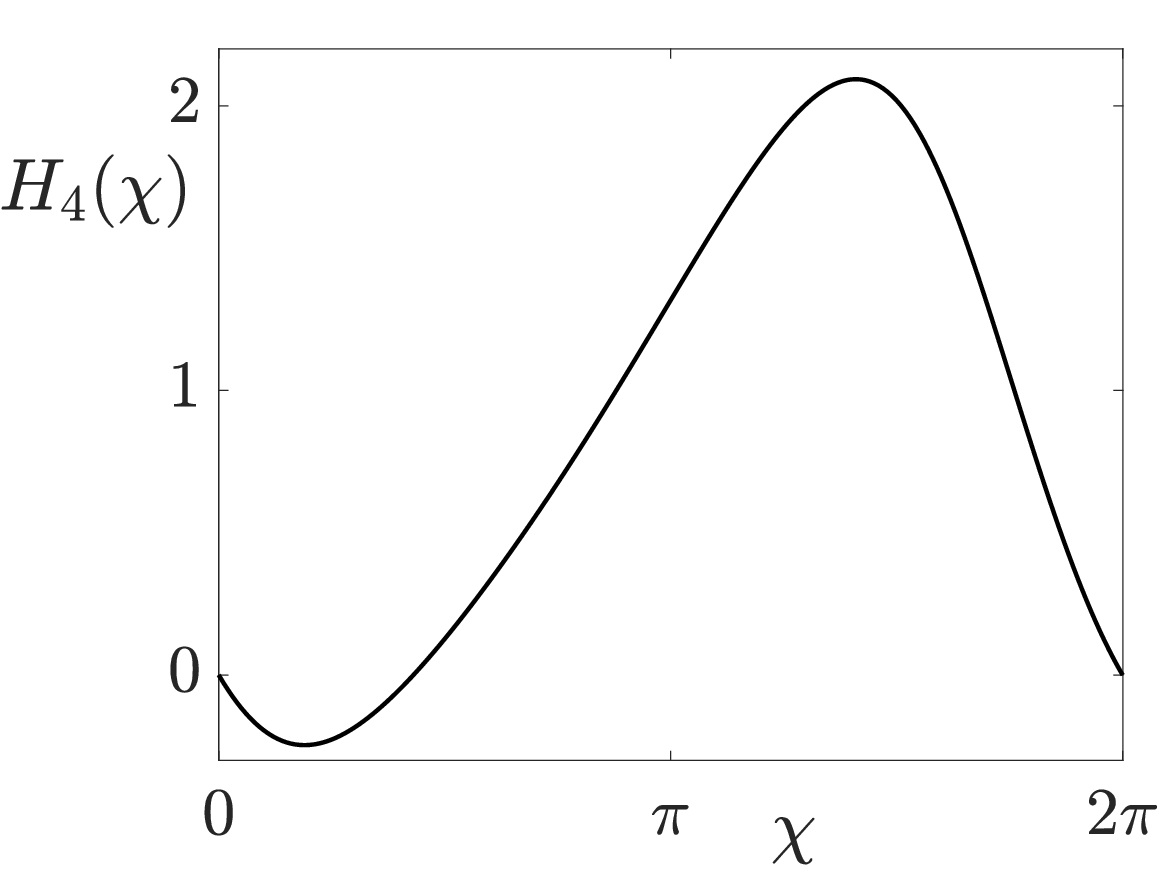}
		\includegraphics[height=3.2cm]{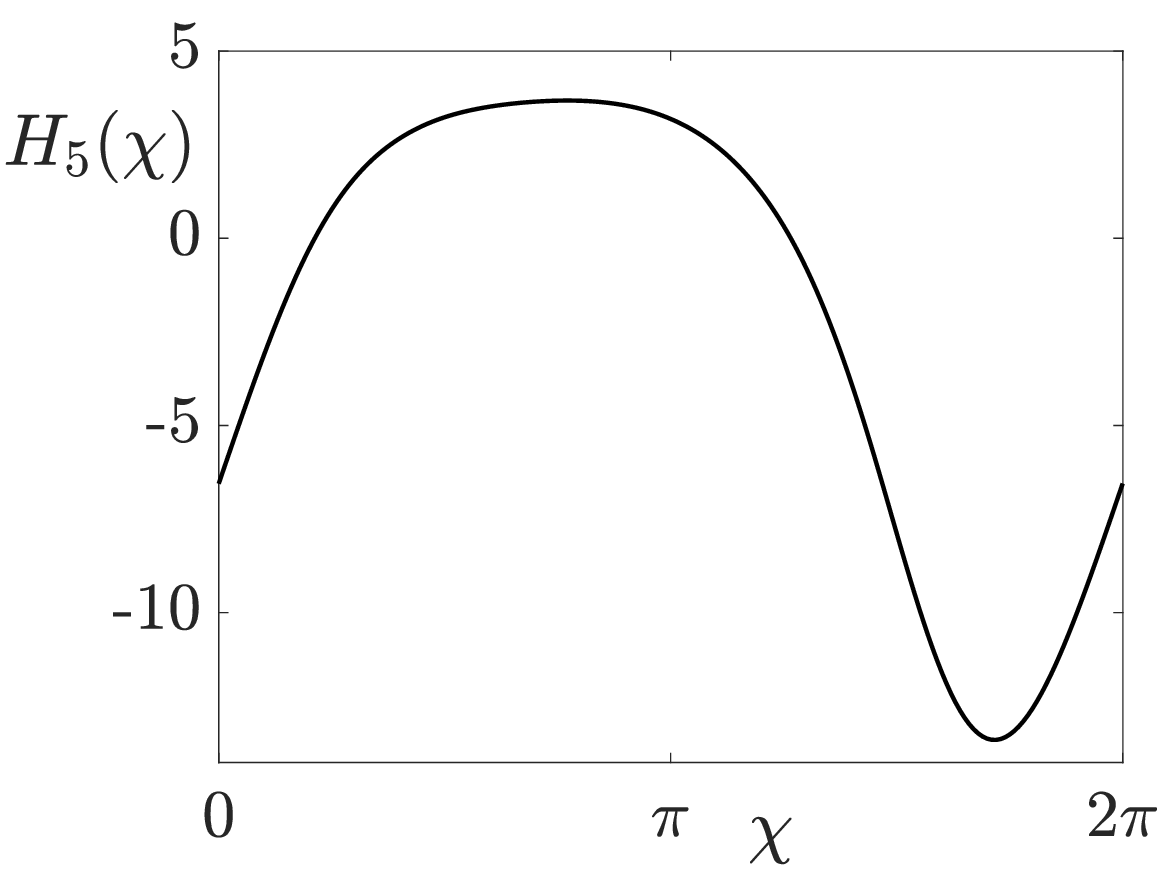}
		\includegraphics[height=3.2cm]{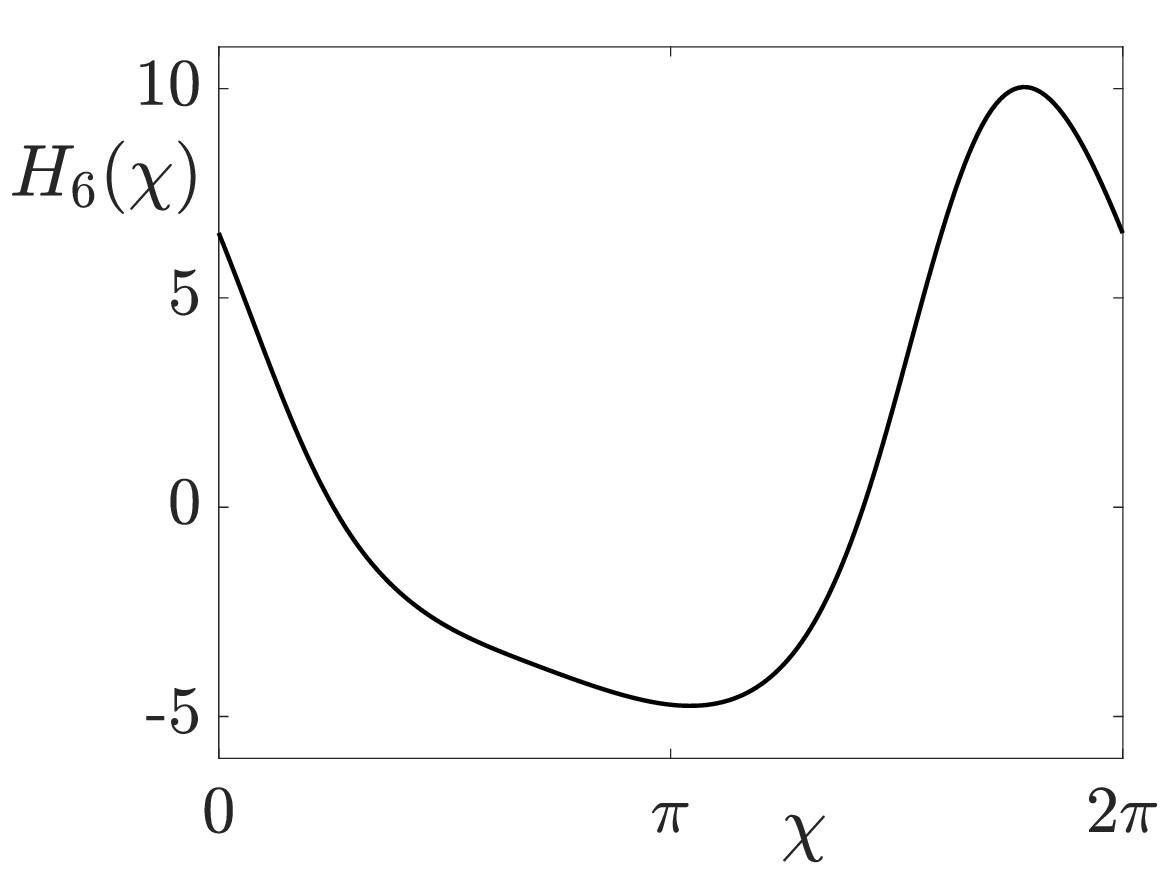}	
\caption{The interaction functions $H_1,\ldots,H_6$ for the PWL Morris--Lecar model.  The parameters are the same as those in \cref{Fig:TheMorrisLecarperiodic}.}
	\label{Fig:Hshapes}
\end{figure}

To further reduce the system \cref{phaseisostablenetworkPWLcoupling4} to a phase-difference form, we use averaging (see \cref{sec:phaseoscnetworks}) and
write $H_{i}(y)=(2\pi)^{-1} \int_{0}^{2\pi} h_{i}(s,y+s)\mathrm{d}s$ and $\chi \equiv \theta_{2} - \theta_{1}$. This yields
\begin{align}\label{phaseisostablenetworkPWLcoupling7}
		\displaystyle \dot{\chi} &= \sigma\left[H_1(-\chi)-H_1(\chi)+\psi_1\left(H_3(-\chi)-H_2(\chi) \right)+\psi_2\left( H_2(-\chi)-H_3(\chi)\right)\right] \,, \\
		\displaystyle \dot{\psi}_1 &= \kappa \psi_1+\sigma\left[H_4(\chi)+\psi_1 H_5(\chi)+\psi_2 H_6(\chi) \right] \,,  \notag \\
		\displaystyle \dot{\psi}_2 &= \kappa \psi_2+\sigma\left[H_4(-\chi)+\psi_2 H_5(-\chi)+\psi_1 H_6(-\chi) \right] \,. \notag
\end{align}
In \cref{Fig:Hshapes}, we show the
six interaction functions $H_1,\ldots, H_6$ for the PWL Morris--Lecar model. We compute these functions using the Fourier representation that we described in \cref{sec:phaseoscnetworks}. Note that these six functions are all that is needed to describe the phase-amplitude reduced dynamics of networks of any finite size $N$ \cite{Park2021, Nicks2023}.

For the synchronous $0$-amplitude solution $[\chi,\psi_1,\psi_2]^{\top}=[0,0,0]^{\top}$, the Jacobian of \cref{phaseisostablenetworkPWLcoupling7} has the form
\begin{equation}\label{Jacobsecondorder}
J = \begin{bmatrix} -2\sigma H_1'(0) & 2\sigma H_3(0) & -2\sigma H_3(0)  \\
		\sigma H_4'(0) & \kappa- \sigma H_6(0) & \sigma H_6(0)  \\
		- \sigma H_4'(0) &\sigma H_6(0) & \kappa- \sigma H_6(0)   \end{bmatrix} \,
\end{equation}
where we have used the fact that linear coupling gives $H_2(0) = -H_3(0)$ and $H_5(0) = -H_6(0)$. All eigenvalues have negative real part, so
 the synchronous solution is linearly stable when $\kappa < 0$ (which we assume to obtain
 a stable periodic orbit), $\kappa < 2\sigma(H_1'(0) + H_6(0))$, and $H_1'(0)(\kappa - 2\sigma H_6(0)) + 2\sigma H_3(0)H_4'(0) < 0$. Reducing 
 to the phase-only description
 by taking $H_2, \ldots, H_6 \equiv 0$ recovers
  the result that the synchronous solution is linearly stable when $\sigma H_1'(0) > 0$. One can similarly determine stability
  conditions for the antisynchronous state (for which
  the phase difference between the two oscillators is $\chi = \pi$). 
  For the antisynchronous state, the shared orbit satisfies
   $\psi_1= \psi_2 = \psi$, where $\psi = -\sigma H_4(\pi)/(\kappa + \sigma(H_5(\pi) + H_6(\pi)))$, which is constant (so that the orbit coincides with an isostable of the node dynamics) \cite{Nicks2023}.

\begin{figure}[!h]
	\centering
\includegraphics[width=10cm]{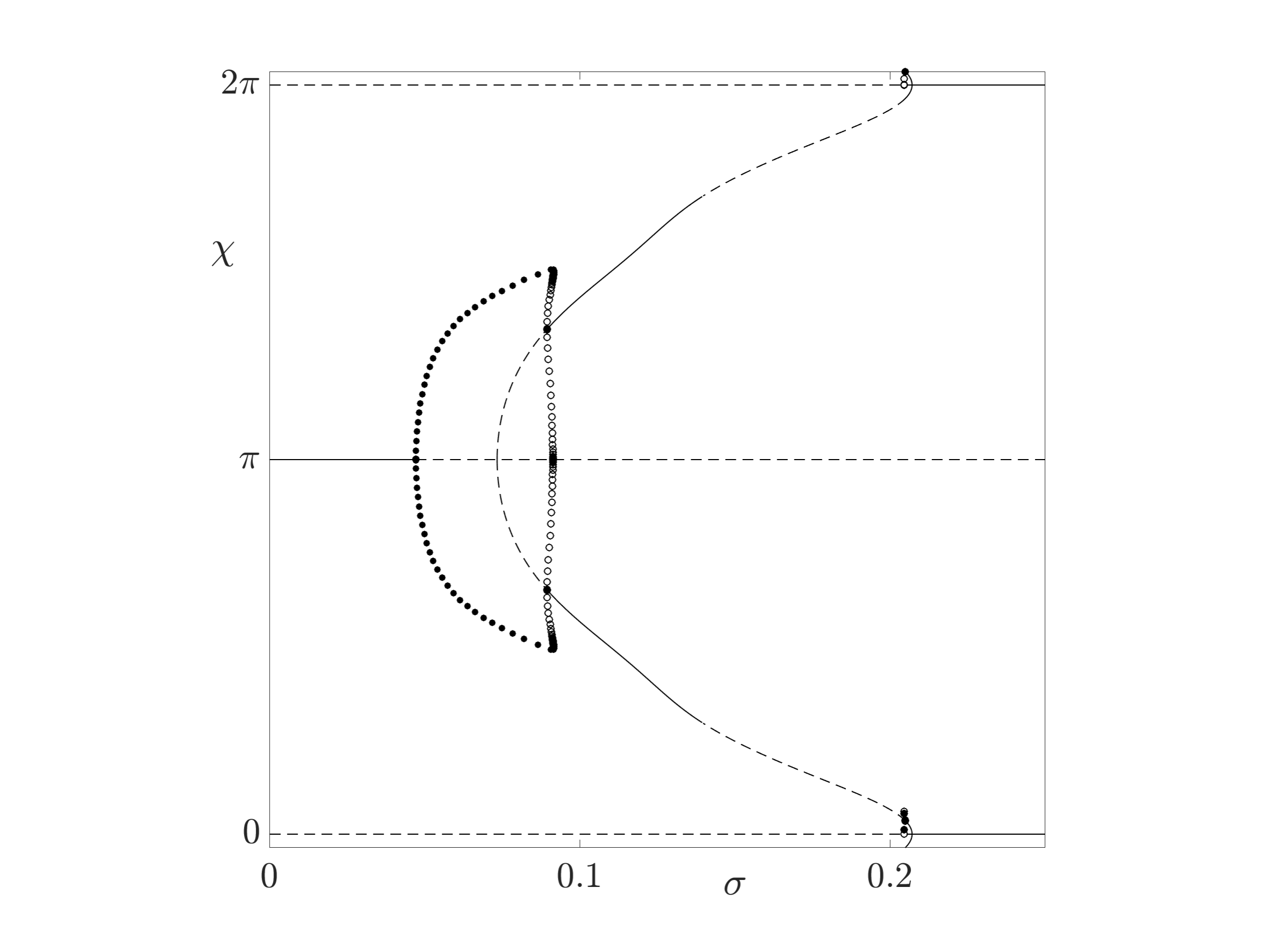}
\caption{Bifurcation diagram for two linearly coupled PWL Morris--Lecar models (see \cref{ML}) showing the phase difference $\chi$ under variation in the overall coupling strength $\sigma$. The solid (respectively, dashed)  curves indicate stable (respectively, unstable) branches of steady-state solutions. The filled (respectively, empty) circles indicate the amplitude of stable (respectively, unstable) periodic orbits.
Two branches with phase difference $\chi \neq 0,\pi$ meet at $\sigma \lessapprox 0.15$. These branches both terminate in a limit point, so the
apparent change of stability is just an artifact of this coincidence.
For $0.15 \lessapprox \sigma \lessapprox 0.2$, we do not obtain any stable solution branches.
The parameters are the same as those in \cref{Fig:TheMorrisLecarperiodic}.}
	\label{Fig:PWLMLBif}
\end{figure}

Importantly, for both solutions, the phase-only reduction does not predict any bifurcations from changing
$\sigma>0$, whereas the phase--amplitude approach does allow this possibility. This is the case because both the 
eigenvalues of \cref{Jacobsecondorder} and of the Jacobian for the antisynchronous state have a richer dependence on the coupling strength $\sigma$.
See \cref{Fig:PWLMLBif} for an
 interesting bifurcation diagram for the PWL Morris--Lecar model that we obtain by varying $\sigma$.
  We see that we can restabilize the synchronous state by increasing $\sigma$ when $\sigma \gtrapprox 0.2$.
    Moreover,
    at smaller values of $\sigma$, stable periodic orbits arise from a Andronov--Hopf bifurcation of the antisynchronous state. 
     In one region, for which $0.15 \lessapprox \sigma \lessapprox 0.2$, our analysis predicts that there are no stable solution branches.
     Direct numerical simulations (see \cref{Fig:PMLNetwoksimulations}) of the full model \cref{ML} confirm this prediction.


Although the qualitative predictions of the phase--amplitude formalism are better than those of the phase-only formalism, it remains to be seen if these predictions can also give successful quantitative insights.
We explore this issue in \cref{sec:strongly}.

\section{Strongly coupled oscillator networks}
\label{sec:strongly}

In previous sections, we explored how collective behavior (such as phase-locked states) arises in
weakly coupled networks.
We considered the dynamics of the system \cref{Network1} on
a reduced phase space that is given by the Cartesian product of each oscillator's phase and possibly a subset of the oscillator amplitudes.
However, the assumption of weak coupling is not valid in many real-world situations.
There are far fewer results for strongly coupled oscillator networks than for weakly coupled oscillator networks, and the former are often restricted to special states such as synchrony \cite[Chapter 7]{Coombes2023}.

One popular approach to obtain insights into the behavior of strongly coupled oscillators in the context of
smooth dynamical systems is the master-stability-function (MSF) approach. The MSF approach\footnote{At least on occasion, MSF approaches were used before they were invented officially in the 1990s. See Segel and Levin \cite{segel1976} (a conference-proceeding paper from 1976).} of
 Pecora and Carroll \cite{Pecora1998} to assess the stability of synchronous states of a network in terms of the spectral properties of the network's adjacency matrix is exact. It does not rely on any approximations, aside from those in numerical implementations (to construct periodic orbits and compute
 Floquet exponents).
 In the present section, we describe how to augment this MSF approach for PWL systems using the saltation operators that we described in \cref{FlqTheory}. For PWL systems, one can use semi-analytical approaches (with numerical computations only for times of flight between switching manifolds) instead of the numerical computations (i.e., simulations of differential equations) that one uses for
 smooth nonlinear systems.\footnote{Recently, Corragio \emph{et al.}~\cite{Corragio2021} used an alternative approach for systems with a so-called ``$\sigma$-QUAD property'' (which includes many discontinuous neural, genetic, and impact networks) to prove global asymptotic convergence to synchronization in networks of piecewise-smooth dynamical systems.}

\subsection{The master stability function for nonsmooth systems}
\label{subsec:strongly}

To introduce the MSF formalism, we start with an arbitrary connected network of $N$ coupled identical oscillators
\cref{Network1},
\cref{Networkinterraction} with $G(x_i,x_j) = \mathcal{H}(x_i) - \mathcal{H}(x_j)$.
The output for each oscillator
is determined by a vector function $\mathcal{H}:\ \RSet^{m} \rightarrow \RSet^{m} $ (which can be either linear or nonlinear). The network dynamics are
\begin{equation}
\label{MSFnetworkEqn}
	\FD{}{t} x_i =  f\left(x_{i}\right)-\sigma \sum_{j=1}^{N} \mathcal{L}_{i j} \mathcal{H}\left(x_{j}\right) \,,
\end{equation}
where, the matrix $\mathcal{L} \in \RSet^{N \times N}$, with entries $\mathcal{L}_{i j}$, is the graph Laplacian
\cref{Laplacian}.
By construction, the matrix $\mathcal{L}$ has $0$ row sums. The $N-1$ constraints $x_{1}(t)=x_{2}(t)=\cdots=x_{N}(t)=s(t)$
define the invariant \textit{synchronization manifold}, where $s(t)$ is a solution in $ \RSet^{m}$ of the associated uncoupled system.
That is, $\mathrm{d}s(t)/\mathrm{d}t=f(s(t))$.
Any motion that begins
 on the synchronization manifold remains there, so the associated
 synchronized state is \textit{flow-invariant}.

When all oscillators are initially on the synchronization manifold with identical initial conditions, they always remain synchronized. To assess the stability of a synchronized state, we perform a linear stability analysis by inserting a perturbed solution $x_{i}(t)=s(t)+\delta x_{i}(t)$ into \cref{MSFnetworkEqn} to obtain the variational equation
\begin{equation}\label{MSFnetworkVar}
	\frac{\mathrm{d}\delta x_{i}}{\mathrm{d} t}=\D f(s) \delta x_{i}-\sigma \D \mathcal{H}(s) \sum_{j=1}^{N} \mathcal{L}_{i j} \delta x_{j} \,,
\end{equation}
where $\D f(s)\in \RSet^{m \times m}$ and $\D \mathcal{H}(s)\in \RSet^{m \times m}$, respectively, denote the Jacobians of $f(s)$ and $\mathcal{H}(s)$, which one evaluates at the synchronous solution $s(t)$. We introduce $U = (\delta x_1,\delta x_2, \ldots, \delta x_N)\in \RSet^{mN}$ and use the tensor product (i.e., Kronecker product) $\otimes$ for matrices to write the variational equation as
\begin{equation}\label{MSFnetworkVarTensor}
	\frac{\mathrm{d U}}{\mathrm{d} t}= \left[I_N\otimes \D f(s)-\sigma( \mathcal{L}\otimes \D \mathcal{H}(s))\right]U \,.
\end{equation}

We organize the normalized right eigenvectors of $\mathcal{L}$ into a matrix $P$ such that $P^{-1}\mathcal{L}=\Lambda P^{-1}$, with $\Lambda=\mathrm{diag}(\lambda_1,\lambda_2,\ldots,\lambda_N)$, where $\lambda_{\eta}$ (with $\eta \in \{1,\ldots, N\}$) are the corresponding eigenvalues of $\mathcal{L}$. We introduce a new variable $Y$ using the linear transformation $Y = (P\otimes I_m)^{-1}U$ to obtain a block-diagonal system
\begin{equation}\label{MSFnetworkVarTensorEig}
	\frac{\mathrm{d Y}}{\mathrm{d} t}= \left[ I_N\otimes \D f(s)-\sigma( \Lambda \otimes \D \mathcal{H}(s))\right] Y \,,
\end{equation}
where $I_N$ is the $N\times N$ identity matrix. This yields
a set of $N$ decoupled $m$-dimensional equations,
\begin{equation}\label{MSFnetworkVarEigfirst}
	\frac{\mathrm{d} \xi_{l}}{\mathrm{d} t}=\left[\D f(s)-\sigma \lambda_{l} \D \mathcal{H}(s)\right] \xi_{l}\,, \quad l \in \{1, \ldots, N\} \,,
\end{equation}
that are parametrized by the eigenvalues of the graph Laplacian $\mathcal{L}$. The Jacobians
$\D f(s)$ and $\D \mathcal{H}(s)$ are independent of the block label $l$.
Because the row sums of $\mathcal{L}$ are $0$, there is always a $0$ eigenvalue $\lambda_{1} = 0$, with a corresponding eigenvector $[1,1, \ldots, 1]^\top$ that characterizes a perturbation that is tangential to the synchronization manifold. The remaining $N-1$ transversal perturbations (which are associated with the other $N-1$ solutions of equation \cref{MSFnetworkVarEigfirst}) must damp out for the synchronous state to be linearly stable.
In general, some eigenvalues $\lambda_{l}$ of $\mathcal{L}$ may be complex. (For example, this can occur when the adjacency
matrix is not symmetric.)
This leads us to consider the system
\begin{equation}\label{MSFnetworkVarEigfinal}
	\frac{\mathrm{d}\xi}{\mathrm{d} t}=[\D f(s)-\mu \D \mathcal{H}(s)] \xi \,, \quad \mu \in \CSet \,, \quad \xi \in \CSet^{m} \,.
\end{equation}
All of the individual variational equations in the system \cref{MSFnetworkVarEigfirst} have the same structure as that of the system \cref{MSFnetworkVarEigfinal}. The only difference is that
$\mu=\sigma \lambda_l$.
Equation \cref{MSFnetworkVarEigfinal} is the so-called \textit{master variational equation}. To determine its stability, we calculate its largest Floquet exponent \cite{guckenheimer2013nonlinear} as a function of $\mu$. The resulting function is the so-called \textit{master stability function} (MSF). More explicitly, for a given $s(t),$ the MSF is
the function that maps the complex number $\mu$ to the largest Floquet exponent of the dynamical system \cref{MSFnetworkVarEigfinal}. The synchronized state of a network of coupled oscillators is linearly stable if the MSF is negative at $\mu=\sigma \lambda_{l}$, where $\lambda_{l}$ ranges over all eigenvalues of the matrix $\mathcal{L}$ except for $\lambda_{1} = 0$.

The Laplacian form of the coupling in equation \cref{MSFnetworkEqn} guarantees that there exists a
synchronous state. However, other forms of coupling are also natural. For example, consider
\begin{equation}\label{MSFnetworkMoreGeneral}
	\dfrac{\mathrm{d}x_{i}}{\mathrm{d} t}= f\left(x_{i}\right)+\sigma \sum_{j=1}^{N} w_{i j} \mathcal{H}\left(x_{j}\right) \,.
\end{equation}
Substituting $x_{i}(t)=s(t)$, with $i \in \{1, 2, \ldots, N$\}, into equation \cref{MSFnetworkMoreGeneral} yields
\begin{equation}\label{MSFnetworkMoreGeneralSync}
	\dfrac{\mathrm{d}s}{\mathrm{d} t}=   f\left(s\right)+\sigma \mathcal{H}\left(s\right) \sum_{j=1}^{N} w_{i j} \,.
\end{equation}
To guarantee that all oscillators have the same behavior, we assume that
$\sum_{j=1}^{N} w_{i j} = \mathrm{constant}$ for all $i$. If the constant is $0$, then we say that the system is \textit{balanced} \cite{deneve2016efficient,van1996chaos,rubin2017balanced,roudi2007balanced}.
In a balanced network, the existence of a synchronous network state is independent of the interaction parameters,
 so varying these parameters cannot induce any nonsmooth bifurcations (arising from a change of the orbit shape and its possible tangential intersection with a switching manifold).

One can apply the MSF framework to chaotic systems, for which one calculates Liapunov exponents instead of Floquet exponents \cite{PecoraCarroll:1997,Pecora1998,dabrowski2012largest}. One can also generalize the MSF formalism to network settings in which the coupling between oscillators includes a time delay \cite{Kyrychko2014,dahms2012cluster}. A synchronous solution is a very special network state, and more elaborate types of behavior can occur. An example is a
``chimera state''(see \cite{CoombesThul:2016, lodato2007synchronization}), in which some oscillators are synchronized but others behave asynchronously~\cite{panaggio2015}.
The original MSF approach allows one to investigate the stability of
networks of identical oscillators, but it has been extended to study stability in networks of almost identical oscillators \cite{sun2009master}. For other discussions of the MSF formalism and its applications, see \cite{Ashwin2016,ArenasGuileraMoreno:2008,Pecora2013,Porter2014}.

One cannot directly apply the MSF methodology to networks of nonsmooth oscillators, and it is desirable to extend it to such systems. We first review a technique that adapts the MSF to PWL systems \cite{CoombesThul:2016}, and we then apply this approach to the models in \cref{sec:pwl}. We seek to show how the linear stability of the synchronous solution changes under variations of the coupling strength in networks of coupled oscillators.

For networks of the form \cref{MSFnetworkEqn} with \textit{linear} vector functions $\mathcal{H}$ (including the 
``linear diffusive case''
$\mathcal{H}(x^{1}, x^{2}, x^{3}\map{,} \ldots, x^{m})=(x^{1},0, 0, \ldots, 0)$) that one builds
from PWL systems of the form \cref{pwl1}, both $\D f(s)$ and $\D \mathcal{H}(s)$ are piecewise-constant matrices.
Therefore, in each
region $R_{\mu}$, equation \cref{MSFnetworkVarEigfinal} takes the form
\begin{equation}\label{MSFnetworkVarEigfinalPWL2}
	\frac{\mathrm{d}\xi_{\mu}}{\mathrm{d} t}=[A_{\mu}-\beta J] \xi_{\mu}\,, \quad \beta \in \CSet \,,
\end{equation}
where $J=\D \mathcal{H}(s)$ and $\xi_{\mu} \in \CSet^{m}$. We solve \cref{MSFnetworkVarEigfinalPWL2} using matrix exponentials.
This yields $\xi_{\mu} (t)=G\left(t;A_{\mu}-\beta J \right)\xi_{\mu} (0)$, where $G(t;A)$ is given by equation \cref{GandK}, although we need to be careful when evolving perturbations through the switching manifolds. Using the
notation $U=(\delta x_1, \delta x_2, \ldots, \delta x_N) \in\mathbb R^{Nm}$, at each event time $t_{i}$, we write $U^{+}=(I_N \otimes S(t_{i}))U^{-}$. We then use
the transformation $Y=(P\otimes I_m)^{-1}U$  and obtain $Y^{+}=(I_N\otimes S(t_i))Y^{-}$, which has $m\times m$ dimensions and an $N-$block structure.
The action of the saltation operator on each block is $\xi_{}(t_{i}^+)= S(t_{i}) \xi_{}(t_{i}^-)$. We use the technique in \cref{Saltationproof} to treat perturbations across a switching boundary. After one period of motion (with $M$ switching events), this yields $\xi(T) = \Psi \xi(0)$, where
\begin{equation}
\label{Gamma}
	\Psi=S(t_M)  G (T_{M} ) S(t_{M-1})  G (T_{M-1} )  \times \cdots \times S(t_2) G(T_{2}) S(t_1) G( T_{1})
\end{equation}
and $G(T_i) = G(T_i; A_{\mu(i)}-\beta J )$ (see \cref{MonodromyExplicitinchapter}).
For PWL systems, all of the individual variational equations, which take the form \cref{MSFnetworkVarEigfirst}, have the same structure as that of the system \cref{MSFnetworkVarEigfinalPWL2}. The only difference is that now there is an additional term $\beta=\sigma \lambda_l$. Therefore, by choosing a reasonable value of $\beta$ in the complex plane, we can determine the stability of \cref{MSFnetworkVarEigfirst} by checking
that the MSF of \cref{MSFnetworkVarEigfinalPWL2} is negative for each $\beta=\sigma \lambda_l$.
Alternatively, we can calculate $\Psi$ for each $l$; we use the notation $\Psi(l)$ to emphasize this. We then obtain that the synchronous state is linearly stable if the periodic solution of a single oscillator is linearly stable and the eigenvalues of $\Psi(l)$, for each $l \in \{2, \dots, n\}$ lie within the unit disc.

The power of the MSF approach is that it allows one to treat the stability of synchronous states for
{all possible networks}. One first computes the MSF and then uses the spectrum of the chosen network to determine stability. Unlike in weakly-coupled-oscillator theory, one can perform the MSF stability analysis without making any approximations.

\subsection{MSF versus weakly-coupled-oscillator theory for systems of $N=2$ oscillators}
\label{PWLthreemethodcomarisons}

Before we present applications of the augmented MSF to a few example nonsmooth systems, we compare and contrast this exact approach to results from weakly-coupled-oscillator theory
without focusing too much on network structure. Consider a simple reciprocal network (i.e., all coupling is bidirectional) of two nodes with an interaction that is described by \cref{phaseisostablenetworkPWLcoupling1}. 
The nonzero eigenvalue of the graph Laplacian $\mathcal{L}$ of this network is $+2$.
For the phase-only description, the synchronous state is linearly stable if $\sigma H'(0) > 0$ (see \cref{sec:phaseoscnetworks}). For the phase--amplitude description, the synchronous state is linearly stable if the
three eigenvalues of \cref{Jacobsecondorder} are all in the left-hand side of the complex plane (see \cref{PWLPhaseAmplitudenetworkfromalism}). For the exact approach of the present section, the synchronous state is linearly stable if
the MSF is negative at $2 \sigma$ (see \cref{subsec:strongly}).
Using the same oscillator parameters as those in \cref{Fig:Pwl2pieces} and \cref{Fig:TheMorrisLecarperiodic}, we find that weakly-coupled-oscillator theory sometimes fails to capture the behavior that is predicted by the exact (MSF) approach. For the McKean and absolute models, all three approaches give the same qualitative prediction that the synchronous
state is linearly stable for small positive $\sigma$ (i.e., for weak coupling) and this stability persists for larger $\sigma$ (i.e., strong coupling).
\begin{figure}[h]
	\centering
		\includegraphics[width=6cm]{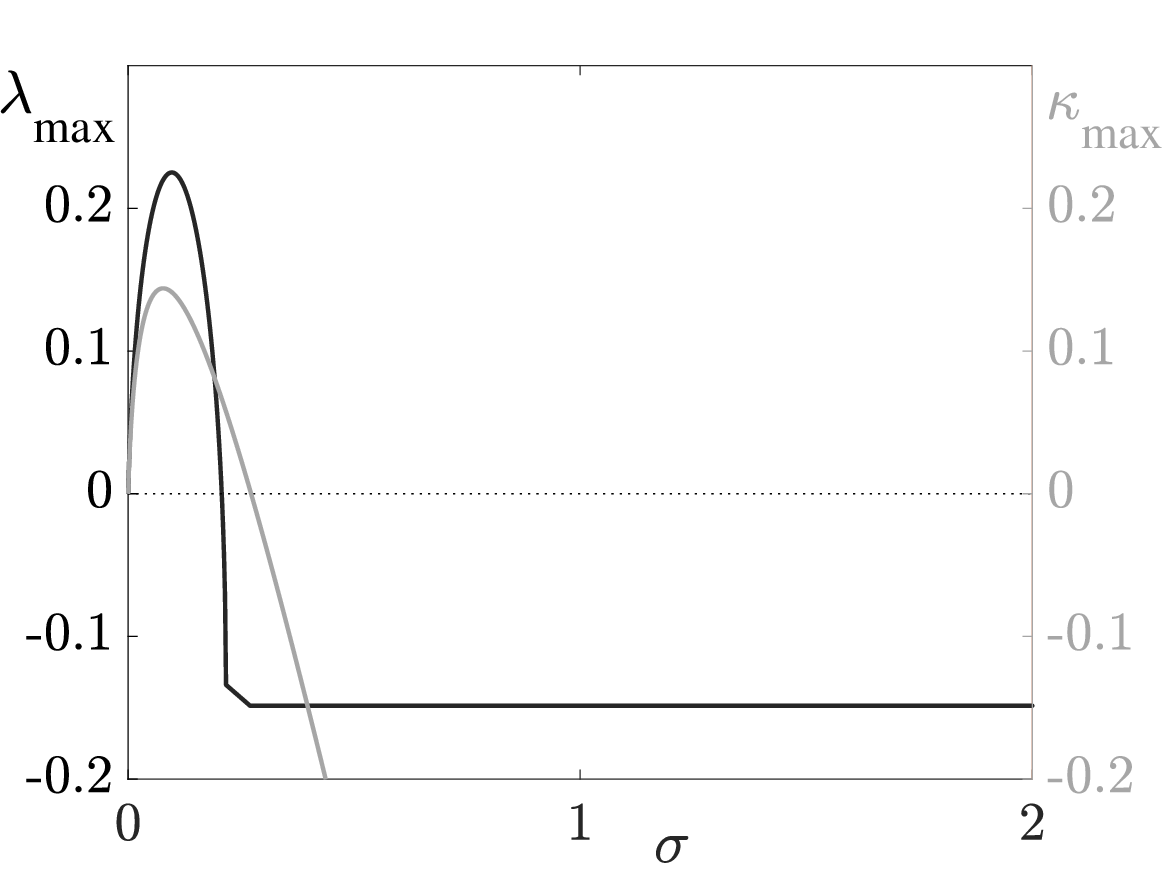}
	\caption{Predictions of linear stability of the synchronous state for the PWL Morris--Lecar model (see \cref{ML}) on a reciprocal two-oscillator network using the phase--amplitude approximation (from weakly-coupled-oscillator theory) and MSF approaches.
	We plot the real part of the largest eigenvalue $\lambda_{\text{max}}=\lambda_{\text{max}}(\sigma)$ of the Jacobian from the phase--amplitude reduction, which predicts that the synchronous state restabilizes at $\sigma_c  \approx 0.2071$ as one increases the coupling strength $\sigma$. The largest Floquet exponent from the MSF approach is $\kappa_{\text{max}}=\kappa_{\text{max}}(\sigma)$, which gives a more accurate prediction of restabilization at $\sigma_c \approx 0.272$.
	}
	\label{Fig:PMLNetwork}
\end{figure}

\begin{figure}[h]
	\centering
	\includegraphics[height=2.3cm]{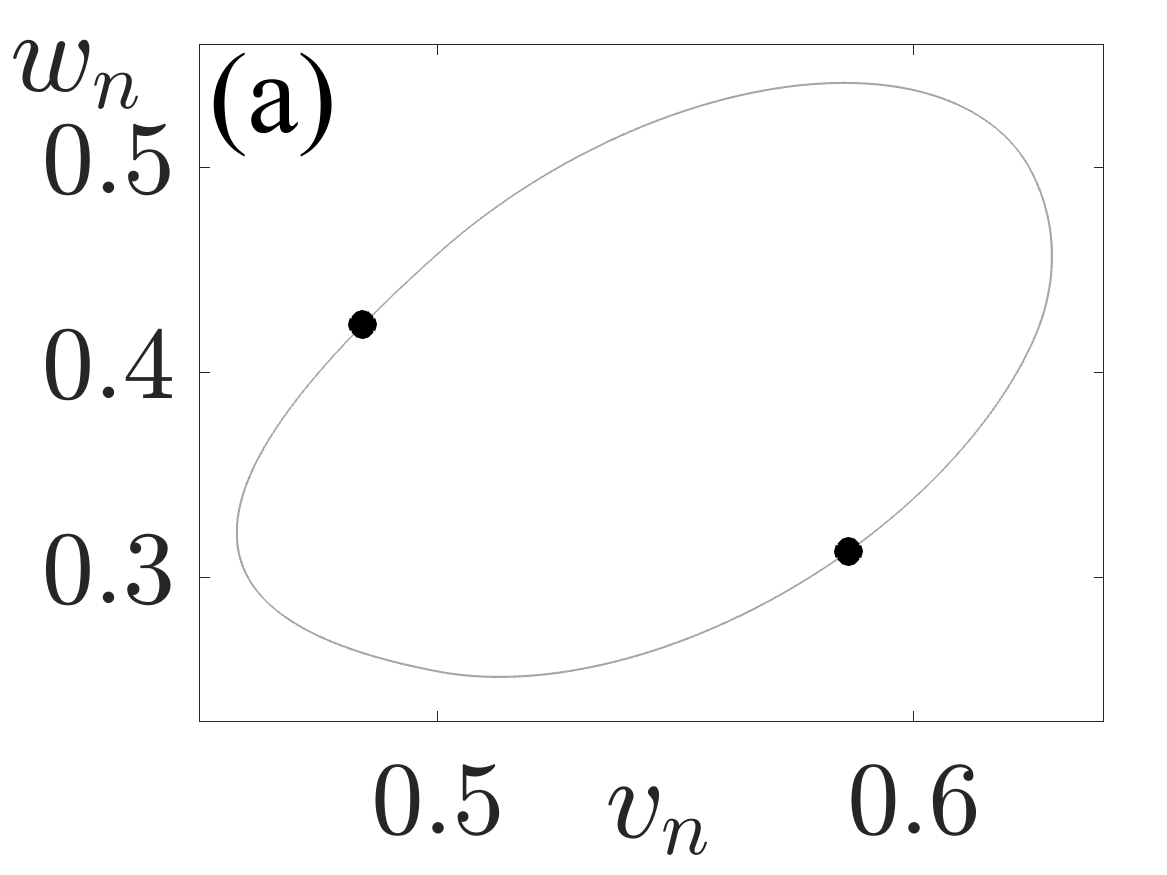}
	\includegraphics[height=2.3cm]{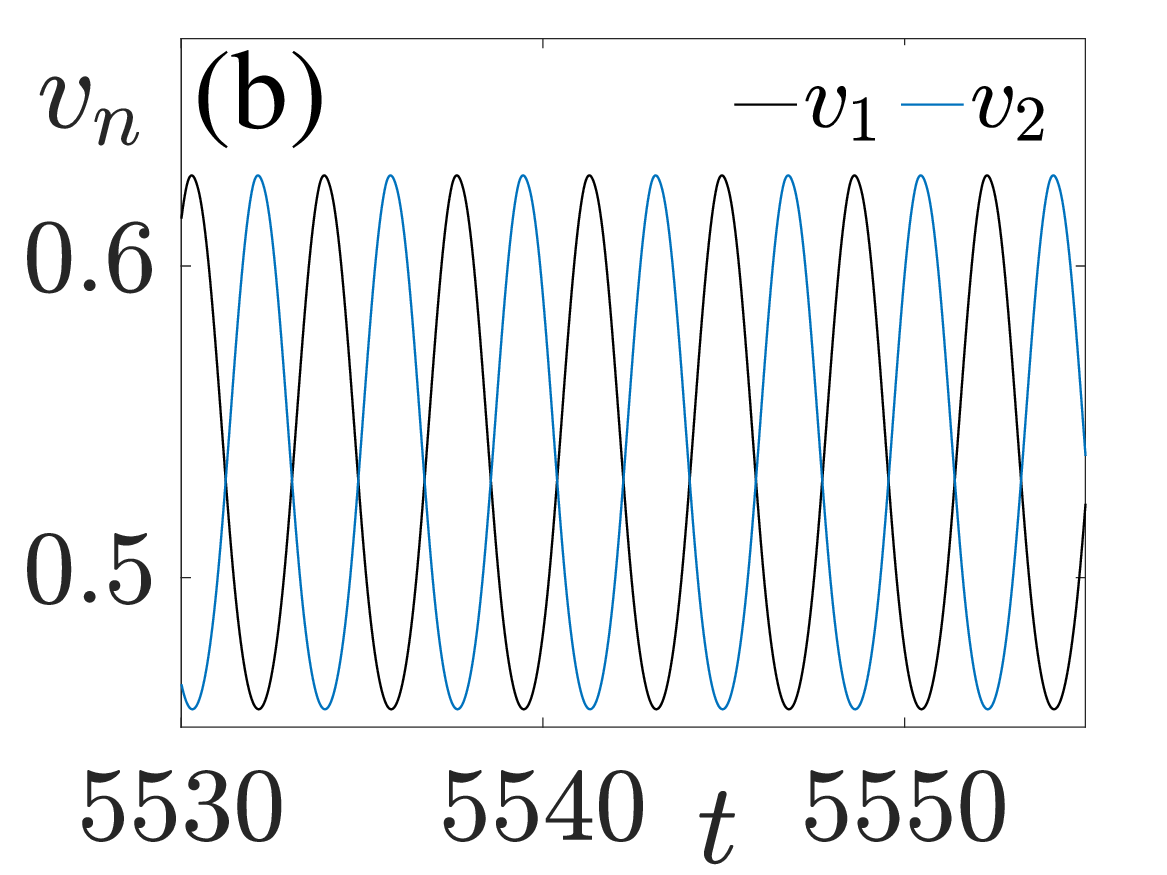}
	\includegraphics[height=2.3cm]{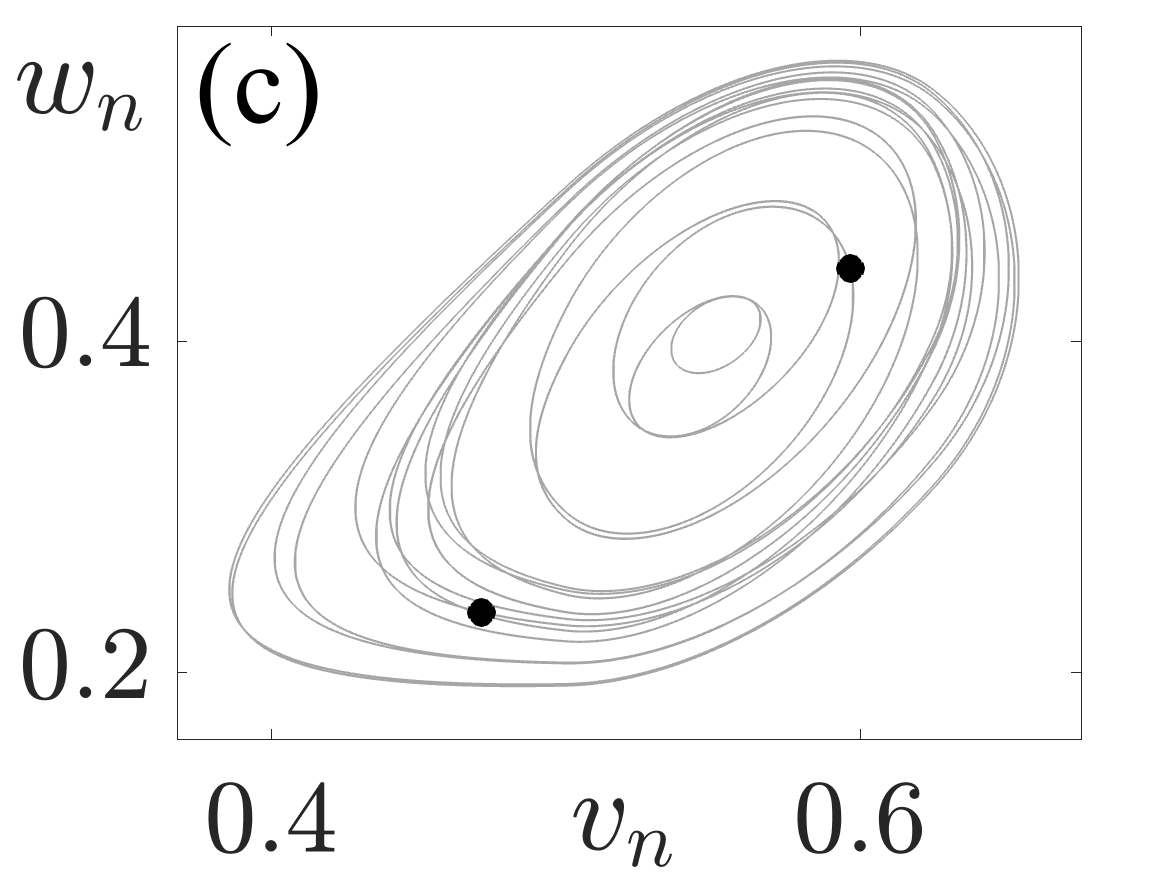}
	\includegraphics[height=2.3cm]{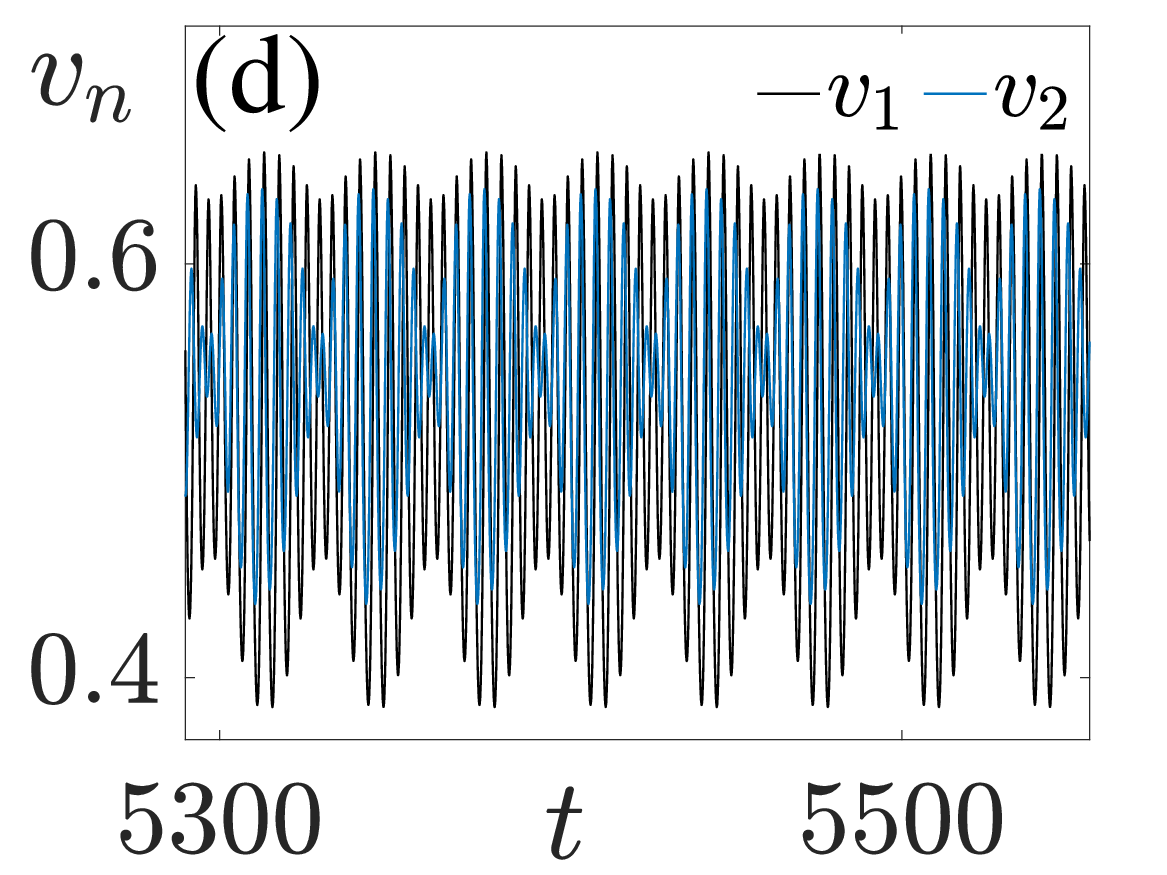}
	\includegraphics[height=2.3cm]{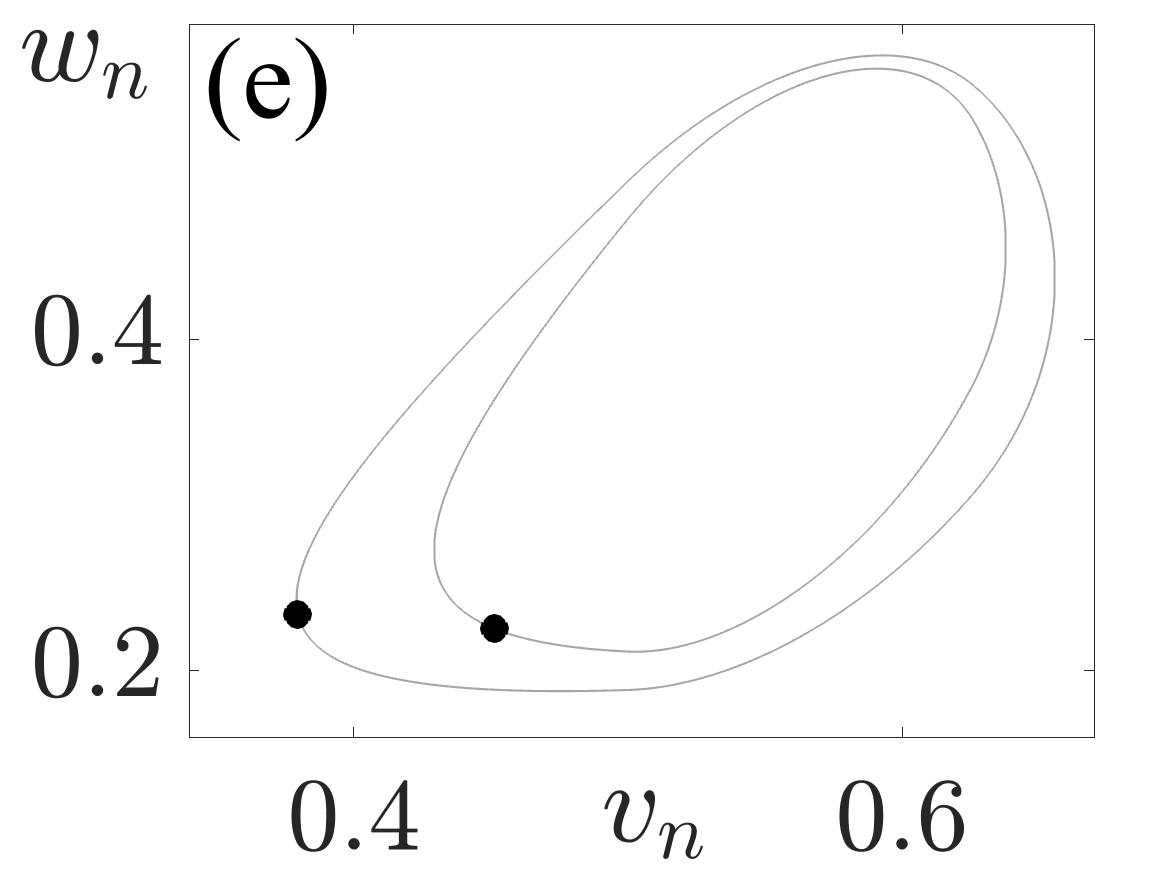}
	\includegraphics[height=2.3cm]{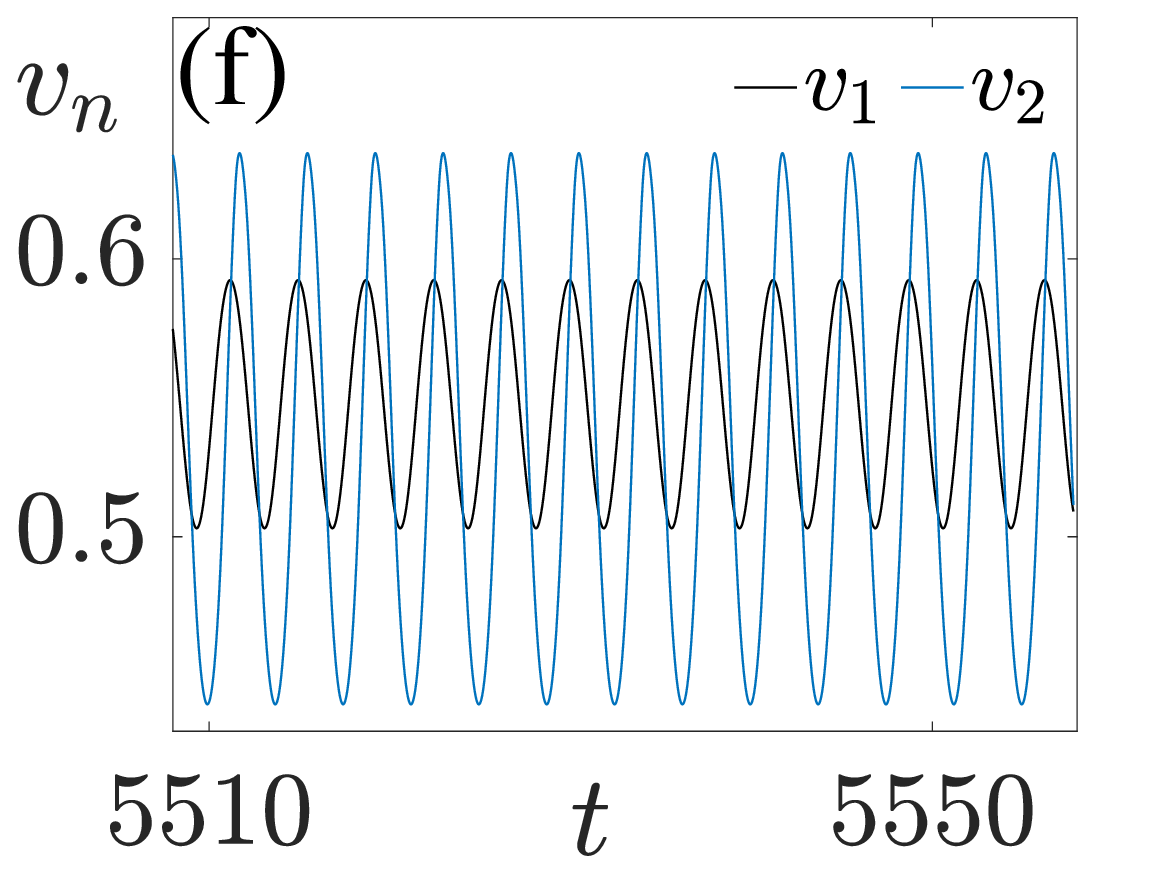}
	\includegraphics[height=2.3cm]{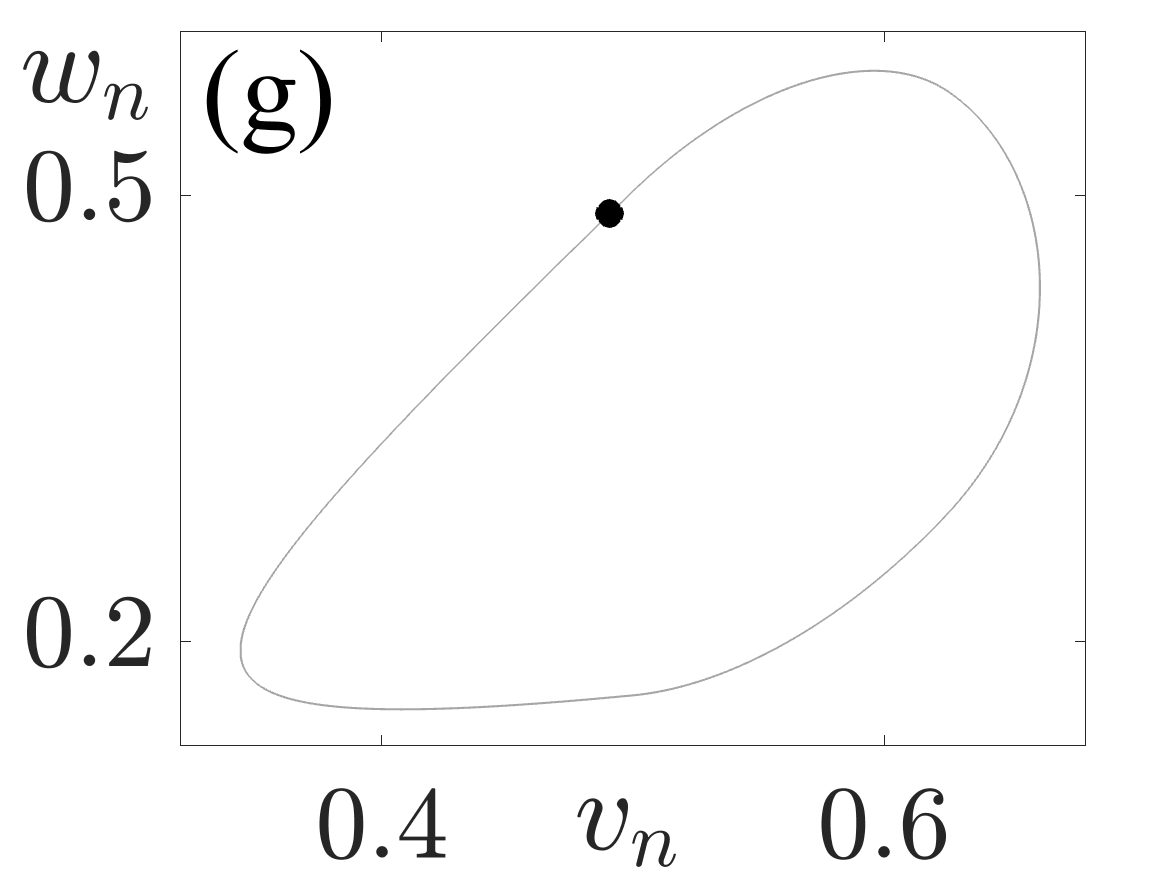}
	\includegraphics[height=2.3cm]{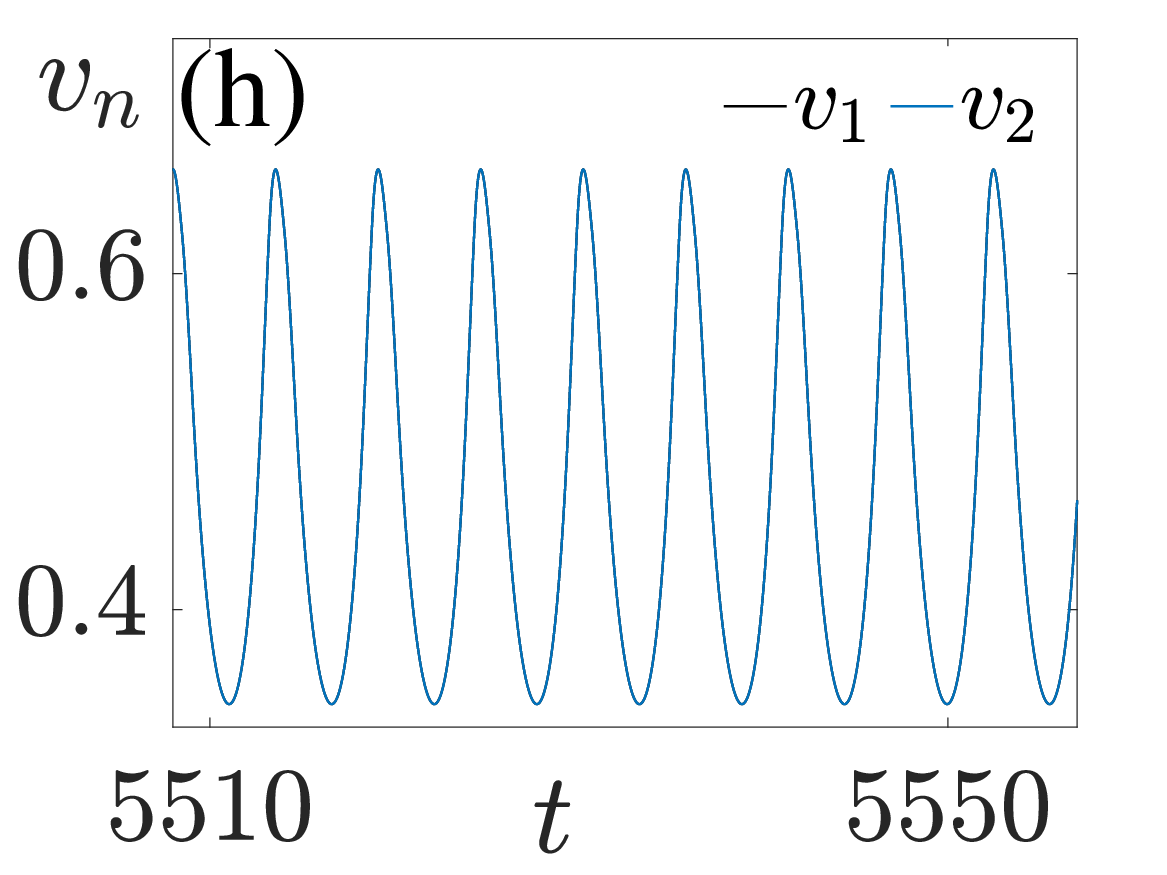}
	\caption{Direct numerical simulations of two reciprocally coupled
	PWL Morris--Lecar oscillators (see \cref{ML}) for coupling strengths of
	(a,b) $\sigma = 0.1$, (c,d) $\sigma = 0.18$, (e,f) $\sigma = 0.25$, and (g,h) $\sigma = 0.28$. In panels (a, c, e, g), we show network activity in the $(v_n,w_n)$ plane.
	In panels (b, d, f, h), we show the corresponding time series for $v_1$ and $v_2$.
	For all $\sigma \gtrapprox 0.272$, the synchronous
	state is always stable. For $\sigma \lessapprox 0.272$, we observe different frequency-locked patterns. The oscillator
	parameters are the same as those in \cref{Fig:TheMorrisLecarperiodic}.
		\label{Fig:PMLNetwoksimulations}}
\end{figure}

For the PWL Morris--Lecar model (see \cref{ML}),
the prediction from the phase-only approximation is that synchrony is always unstable for weak positive coupling. By contrast, the phase--amplitude approximation and MSF approach predict that synchrony can restabilize with increasing coupling strength $\sigma$, although they predict somewhat different values for the critical coupling strength $\sigma=\sigma_c$ at which the network restabilizes.
In \Cref{Fig:PMLNetwork}, we plot the real part $\kappa_{\text{max}} = \text{Re} (\ln(\text{MSF}(\beta)))/T$ (where $\beta=2 \sigma$) of the largest Floquet exponent from the MSF
as a function of $\sigma$. In the same figure, we plot the the real part $\lambda_{\text{max}} = \max \{ \text{Re}(\lambda_{1}), \mathrm{Re}(\lambda_{2}), \mathrm{Re}(\lambda_{3}) \}$ of the largest eigenvalue from the phase--amplitude approximation.
 The phase--amplitude prediction is that $\sigma_c \approx 0.2071$, whereas the (exact) MSF prediction is that $\sigma_c \approx 0.272$.
The phase-only theory is incorrect qualitatively, the phase--amplitude theory is correct qualitatively, and the MSF approach (which agrees with direct numerical simulations) is correct both qualitatively and quantitatively.
In \cref{Fig:PWLMLBif}, we explored
the behavior of the PWL Morris--Lecar model in the phase--amplitude reduction for coupling strengths $\sigma \in (0,\sigma_c)$ (where synchrony is unstable) using bifurcation analysis, which
predicts the existence of a stable antisynchronous state (for which there is a relative phase of $\pi$ between the two oscillators) and of frequency-locked states (i.e., states in which oscillators are synchronized at the same frequency)
of different amplitudes.
 Direct numerical simulations (see \cref{Fig:PMLNetwoksimulations}) confirm these predictions.

\begin{figure}[h]
	\centering
	\includegraphics[height=4cm]{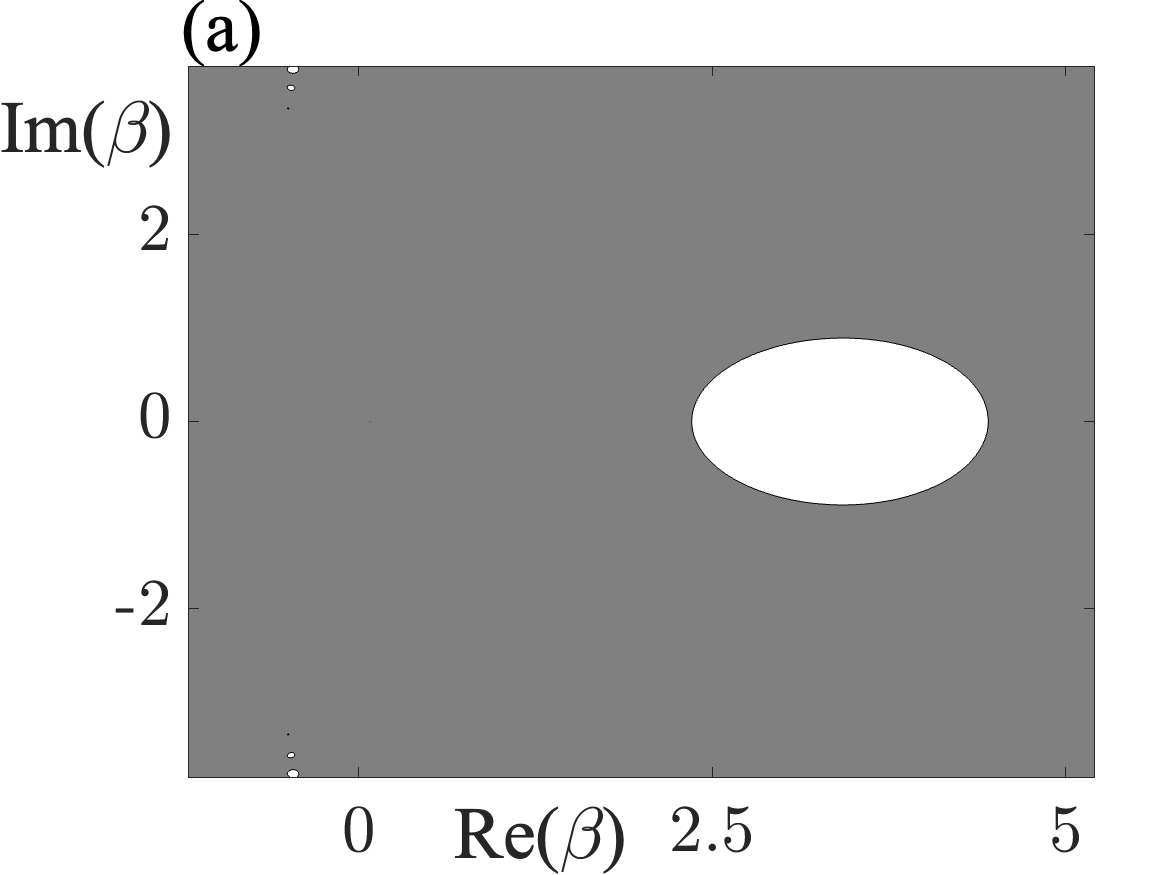} \hspace*{1cm}
	\includegraphics[height=4cm]{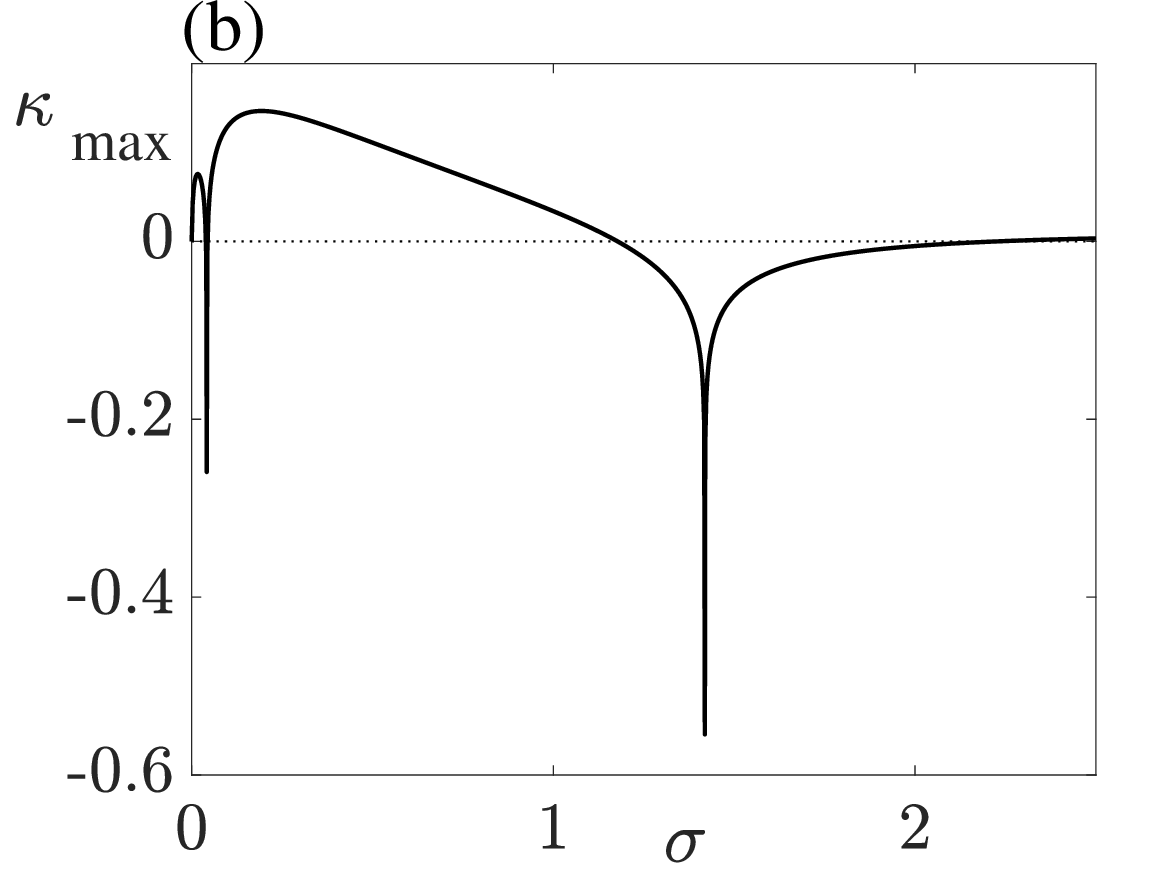}
	\caption{Predictions from the MSF approach of the stability of the synchronous state in the PWL homoclinic model for a reciprocal two-node network.
	 (a) The MSF has a nontrivial structure, with two ellipsoidal regions where it is negative (which we color in white). The region near the origin is very small and not visible with the employed scales.  It is easier to see this small region of stability in (b), in which we plot the real part of the largest network Floquet exponent for $\sigma = \beta/2$ with $\beta \in \RSet$.
One region of stability is $0.0395 \lessapprox \sigma \lessapprox 0.0439$ and the other is $1.178 \lessapprox \sigma \lessapprox 2.226$. The
oscillator parameters are the same as those in \cref{Fig:Pwl2pieces}.
}
	\label{Fig:HomoclinicNetwork}
\end{figure}

For the PWL homoclinic model, both the phase-only approximation and the phase--amplitude approximation predict that synchrony is always unstable for weak positive coupling $\sigma$ in a two-oscillator reciprocal network.
These predictions are both inconsistent with the MSF prediction (which agrees with direct numerical simulations) of two windows of positive coupling with stable synchronous states.
In \Cref{Fig:HomoclinicNetwork}, we show a plot of the MSF that reveals a nontrivial structure, with two ellipsoidal regions where it is negative. (There is a very small ellipsoidal region near the origin that is not visible with the employed scales.)
We also show a slice through $\beta$ along the real axis that illustrates where the real part of the largest network Floquet exponent is negative, generating two regions in which the synchronous state is stable.

\subsection{A brief note about graph spectra}
\label{sec:spectra}

As we have seen in our discussions, the spectrum of a graph is important for determining the stability of the synchronous state in both the weakly-coupled-oscillator and MSF approaches. We thus briefly discuss the spectra of a few simple but notable types of graphs. See \cite{vanMieghem2012} for a thorough exploration of graph spectra.

For a network (i.e., a graph) of $N$ nodes, one specifies the connectivity pattern by a coupling matrix $w \in \RSet^{N \times N}$ (which is often called an ``adjacency matrix'') with entries $w_{ij}$. The spectrum of the graph is the set of eigenvalues of the matrix $w$. This spectrum also determines the eigenvalues of the associated combinatorial graph Laplacian $\mathcal{L}$. In our discussion, we denote the eigenvalues of $w$ by $\lambda_l$, with $l \in \{0,\ldots,N-1\}$, and we denote the corresponding right eigenvectors by $u_l$.

\medskip

\begin{description}
\item[Global.]
The simplest type of network with global coupling has adjacency-matrix entries
$w_{ij} = N^{-1}$. The associated network is fully connected
with homogeneous coupling.
The matrix $w$ has an eigenvector $(1,1,\ldots,1)$ with eigenvalue $\lambda_0=1$ and $N-1$ degenerate eigenvalues $\lambda_l=0$, for $l \in \{1,\ldots,N-1\}$, with corresponding eigenvectors $u_l$ that satisfy the constraint $\sum_{l=0}^{N-1} u_l = 0$.

\medskip

\item[Star.]
A star network has a hub-and-spoke structure, with a central oscillator that is adjacent to $N-1$ leaf nodes (which are not adjacent to each other). Star networks arise in computer-network topologies in which one
central computer acts as a conduit to transmit messages (providing a common connection point for all nodes through a hub). This star-graph architecture has the adjacency matrix
\begin{equation}\label{star}
	w =
\begin{bmatrix}
0 & 1/K & 1/K  & \cdots & 1/K \\
1 & 0 & 0   & \cdots & 0 \\
1 & 0 & 0   & \cdots & 0 \\
\vdots & \vdots & \vdots  & \ddots & \vdots \\
1 & 0 & 0   & \cdots & 0 \\
\end{bmatrix}
\end{equation}
for some constant $K$. If $K = N - 1$, the matrix $w$ has an eigenvalue $\lambda_0=1$ with corresponding eigenvector $[1,1,\ldots,1]^\top$, an eigenvalue $\lambda_1 = -1$ with corresponding eigenvector $[-1,1,\ldots,1]^\top$, and
$N-2$ degenerate eigenvalues $\lambda_l=0$, for $l \in \{3,\ldots,N-1\}$, with corresponding eigenvectors $u_l$ of the form $[0,u_1,\ldots, u_{N-1}]^\top$ that satisfy the constraint $\sum_{l=1}^{N-1} u_l= 0$.

\medskip

\item[Circulant.]
A circulant network's adjacency matrix has entries $w_{ij} = w_{|i-j|}$. Its rows are
shifted versions of the column vector $[w_0, \ldots, w_{N-1}]^top$.
Its eigenvalues are $\lambda_l = \sum_{j=0}^{N-1} w(|j|) \omega_l^j$, where $\omega_l=\exp(2 \pi {\rm i} l /N)$ is an $N$th root of unity. The eigenvectors are $u_l = [1,\omega_l, \omega_l^2, \ldots, w_l^{N-1}]^\top$.
\end{description}

\subsection{Network symmetries and cluster states}
\label{ssec:symmetries}

Perfect global synchronization is just one of many states that can emerge in networks of oscillators.
Indeed, one expects instabilities of the synchronous state to generically yield ``cluster states'', in which subpopulations synchronize, but not necessarily with each other.
Such cluster synchronization has been relatively well-explored in phase-oscillator networks \cite{Ashwin07, Brown03}, although less is known about it in networks of limit-cycle oscillators. For this more general scenario, researchers have made progress in networks with symmetry or when the coupling has a linear diffusive (i.e., Laplacian) form \cite{golubitsky2023,Golubitsky1986, Pogromsky2002, Pogromsky2008, Belykh2000}.

Pecora \textit{et al}. \cite{Pecora2014} and Sorrentino \textit{et al}. \cite{Sorrentino2016} extended the MSF approach (see \cref{subsec:strongly}) to analyze the stability of cluster states that stem either from network symmetries or from
Laplacian coupling (see equations \cref{MSFnetworkEqn} and \cref{MSFnetworkMoreGeneral}). Cluster states arise naturally in networks with symmetry, and cluster synchronization can also occur in networks without symmetry when some of the nodes
have synchronous input patterns \cite{Golubitsky2016}. For networks of identical oscillators that satisfy equation \cref{MSFnetworkEqn}, a symmetry of the network is a permutation $\gamma$ of the nodes that does not change the governing equations.
These permutations are precisely the ones that satisfy $M_\gamma \mathcal{L} = \mathcal{L} M_\gamma$, where $\mathcal{L} \in \RSet^{N \times N}$ is the graph Laplacian \cref{Laplacian} and $M_\gamma$ is the $N\times N$ permutation matrix for the permutation $\gamma \in S_N$.
The network symmetries form a group $\Gamma \subseteq S_N$ that is
isomorphic to the group of automorphisms of the graph that underlies
the network.

For a given adjacency
matrix,
one can identify the automorphism group $\Gamma$ using computational-algebra routines (such as those that are implemented in SageMath \cite{sagemath}).
One can then apply the algorithms in \cite{Sorrentino2016} to enumerate all possible cluster states for the associated network structure. Some of these correspond to isotropy subgroups\footnote{A subgroup of a Lie group $\Gamma$ is an isotropy subgroup for the action of $\Gamma$ on a vector space $V$ if it is the largest subgroup that leaves invariant some vector in $V$ \cite{golubitsky2023, Golubitsky1988}.}
$\Sigma \subseteq \Gamma$ and thus arise from network symmetries. 
The orbit under $\Sigma$ of node $i$ is the set $\{ \gamma(i) \ : \ \gamma \in \Sigma\}$. The orbits permute subsets of nodes among each other and thereby partition the nodes into clusters. Nodes that are part of the same orbit (i.e., in the same cluster) have synchronized dynamics $x_{\gamma(i)} \equiv x_i$ for any $\gamma \in \Sigma$ (see \cite[Thm III.2]{Golubitsky2016}). Isotropy subgroups that are conjugate in $\Gamma$ 
lead to cluster states with identical existence and stability criteria \cite{Golubitsky1988, Ashwin2016}. The remaining possible cluster states arise from the specific
choice of Laplacian coupling. One can determine them using an algorithm that considers whether or not merging two clusters in a state that is determined by symmetry yields a dynamically valid state (i.e., whether or not it yields consistent equations of motion when $x_i$ is the same
for all nodes in the merged cluster). Sorrentino \textit{et al.} \cite{Sorrentino2016} referred to such cluster states as ``Laplacian clusters".
See \cite{Sorrentino2016} for a detailed explanation of the algorithm to determine these clusters, and see \cite{Nicks2018} for an illustration of this algorithm. One can automate
this algorithm using computer-algebra tools \cite{sagemath}.

The above steps yield a list of possible cluster states. 
The existence and stability of these states depends on the
 node dynamics $f$, the output function $\mathcal{H}$, and the coupling strength $\sigma$. The presence of symmetry in a system imposes constraints on the form of the Jacobian matrix, which one can use to greatly simplify stability calculations. For periodic cluster states that one predicts from symmetry, there are well-established methods for stability calculations in symmetric systems
 to block-diagonalize the Jacobian and
generalize the MSF formalism \cite{Golubitsky1988, golubitsky2023}. Sorrentino \textit{et al.} \cite{Sorrentino2016} extended these techniques to Laplacian cluster states. We follow \cite{Sorrentino2016, Nicks2018} and summarize this analysis.

Consider a periodic cluster state that arises from a network symmetry with the corresponding isotropy subgroup $\Sigma \subseteq \Gamma$. The fixed-point subspace of $\Sigma$
is $\Upsilon = \mbox{Fix}(\Sigma)$, which is the synchrony subspace of the cluster state. The cluster state consists of $M$ clusters $\mathcal{C}_k$, with $k \in \{1, \ldots, M\}$, where $M = \operatorname{dim}( \operatorname{Fix} (\Sigma))= \dim(\Upsilon)$. Let $s_k(t)$ denote the synchronized state of nodes in cluster $\mathcal{C}_k$ and
recall the notation of \cref{subsec:strongly}. The variational equation of \cref{MSFnetworkEqn} about the cluster state is
\begin{align}\label{eq:ClusterVariational}
	\FD{U}{t}  =  \left[ \sum_{k=1}^M E^{(k)} \otimes \D f(s_k(t)) - \sigma \sum_{k=1}^M \left( \mathcal{L} E^{(k)}\right) \otimes \D\mathcal{H}(s_k(t)) \right ] U \,,
\end{align}
where $E^{(k)}$ is the diagonal $N \times N$ matrix with entries $E_{ii}^{(k)} = 1$ if $i \in \mathcal{C}_k$ and $E_{ii}^{(k)} = 0$ otherwise. To determine the stability of the periodic cluster state, we need to compute the Floquet exponents of equation \cref{eq:ClusterVariational}. We block-diagonalize the variational equation \cref{eq:ClusterVariational} using the system's symmetries to
simplify this task. 
One can decompose the action of $\Sigma$ on the phase space $\mathbb{R}^{Nm}$ 
into a collection
of irreducible representations of $\Sigma$ (i.e., the most trivial invariant subspaces under the action of $\Sigma$). Some of these subspaces are isomorphic to each other; we combine these subspaces to obtain ``isotypic components" \cite{Golubitsky1988, golubitsky2023}.
Each isotypic component is invariant under the variational equation \eqref{eq:ClusterVariational}, so one can determine the Floquet exponents by considering the restriction of this equation
to each isotypic component.
Therefore, the decomposition puts the variational equations into block-diagonal form. We then compute Floquet exponents for each block to determine the stability of the cluster state. See~\cite{Golubitsky1988} for a detailed discussion of the process of isotypic decomposition and its use in stability computations. Pecora \textit{et al.} \cite{Pecora2014} presented an explicit algorithm to (1) determine the isotypic decomposition for a given cluster state from symmetry and (2) compute a transformation matrix $Q$ so that $\mathcal{L}^\prime = Q \mathcal{L} Q^{-1}$ is block diagonal. Applying this transformation to the variational equation \eqref{eq:ClusterVariational} yields a block-diagonal system of equations:
\begin{align}\label{eq:blockVariational}
	\FD{V}{t} =  \left[\sum_{k=1}^M J^{(k)} \otimes \D f(s_k(t))  - \sigma \sum_{k=1}^M \left( \mathcal{L}^\prime J^{(k)}\right) \otimes \D\mathcal{H}(s_k(t)) \right] V \,,
\end{align}
where $V(t) = (Q \otimes I_m) U(t)$ and $J^{(k)} = Q E^{(k)} Q^{-1}$. The isotypic component of the trivial representation is $\mbox{Fix}(\Sigma) = \Upsilon$, which is the synchronization manifold. This gives an $M \times M$ block in $\mathcal{L}^\prime$ that corresponds to perturbations within the synchronization manifold; one of the Floquet exponents will be $0$ and the remaining $M-1$ correspond to intercluster perturbations. 
The remaining blocks correspond to the isotypic components of other irreducible representations of $\Sigma$. When the node-space representation has $l\geq 1$ isomorphic copies of a particular irreducible representation, we obtain
a block of size $l \times l$.
Such a block corresponds to a perturbation that is
transverse to the synchronization manifold (intracluster perturbations); the associated Floquet multipliers determine the stability
under a synchrony-breaking perturbation. For a cluster state to be linearly stable, all Floquet exponents (except the one that is always $0$) must have a negative real part.

For a periodic Laplacian cluster state, the synchronization manifold is an invariant subspace, but it is not the fixed-point subspace of any subgroup of $\Gamma$. However, we can still block-diagonalize the Laplacian $\mathcal{L}$ so that the top-left block corresponds to perturbations within the synchronization manifold. To do this, we use the algorithm of Sorrentino \textit{et al}. \cite{Sorrentino2016}. Suppose that we start with a cluster state from symmetry with isotropy group $\Sigma$ that has $M$ clusters and a variational equation that is
block-diagonalized by the matrix $Q$.
Suppose that we merge two clusters in this state to obtain a Laplacian cluster state.
Upon this merger, the dimension of the synchronization manifold decreases by $1$ and the dimension of the transverse manifold increases by $1$. We obtain new coordinates on the synchronization manifold by transforming the new synchronization vector in the node-set coordinates (this vector has $1$ entries in the position of each node in the new merged cluster and $0$ entries everywhere else) into the coordinates of the block-diagonalization of the cluster state with isotropy group $\Sigma$.
The orthogonal complement of the new synchronization vector gives the new transverse direction. We normalize the resulting vectors and use
them as rows of an orthogonal matrix $Q^\prime$ whose other rows satisfy $Q^\prime_{ij} = \delta_{ij}$.
The matrix $\chi=Q^\prime Q$ block-diagonalizes $\mathcal{L}$ to a matrix $\mathcal{L}^{\prime \prime}$ that has a top-left block of size
$(M-1) \times (M-1)$. Therefore, the transformation matrix $\chi$ block-diagonalizes the variational equation for the Laplacian cluster state, facilitating the ability to determine
 both the $m(M-1)$ Floquet exponents within the synchronization manifold and the $m(M+1)$ transverse Floquet exponents. This process for computing the required matrix $\chi$ is illustrated with examples in \cite{Sorrentino2016} and \cite{Nicks2018}.

\begin{figure}[t]
\begin{center}
\includegraphics[width=8cm]{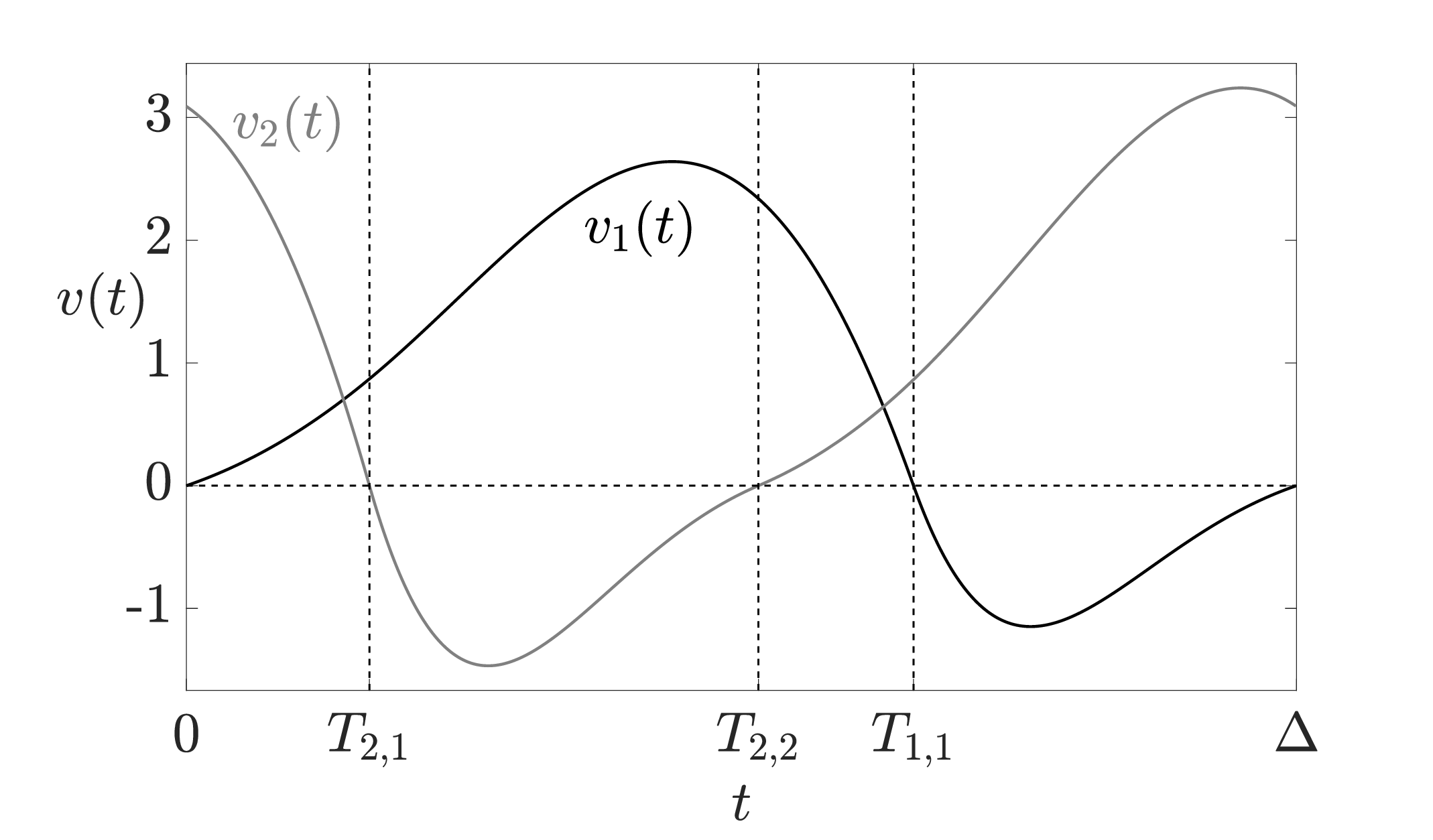}
\caption{The $v$ components of the orbits $s_1$ and $s_2$ over one period. One needs to solve seven nonlinear algebraic equations to determine the unknown initial data $w_1(0)$, $v_2(0)$, and $w_2(0)$ and the switching times $T_{1,1}$, $T_{2,1}$, $T_{2,2}$, and $T_{1,2}= \Delta$.
 }
\label{Fig:switching}
\end{center}
\end{figure}

For PWL systems of the form \cref{MSFnetworkEqn} with linear 
vector function $\mathcal{H}$, it is relatively straightforward to construct the periodic orbits $s_k(t)$ for a
cluster state and to determine its stability by applying the modified Floquet theory (which accounts for the lack of smoothness of the dynamics) of \cref{FlqTheory}
to the block-diagonalized system. For example, suppose that we have a small network of linearly coupled oscillators whose dynamics satisfy the absolute PWL model (see \cref{Fig:Pwl2pieces}(a)). As an illustration, consider the five-node network in
\cite{Sorrentino2016} with graph Laplacian matrix
 \begin{equation}\label{eq:Laplacian5}	
	 \mathcal{L} = \begin{bmatrix} 3 &-1 &0 & -1 & -1\\ -1 & 3 &-1 & 0 & -1  \\ 0 &-1 &3 & -1 & -1\\-1 & 0 &-1 & 3 & -1  \\ -1 & -1 & -1 & -1 & 4\end{bmatrix} \,.
 \end{equation}
 The network supports a Laplacian cluster state with clusters $\mathcal{C}_1 = \{ 1,3,5 \}$ and $\mathcal{C}_2 = \{ 2,4\}$ \cite{Sorrentino2016, Nicks2018}. For this cluster state, ${x}_1 = {x}_3 = {x}_5 = s_1$ and ${x}_2 = x_4= s_2$, where $x_i = [v_i, w_i]^\top$ for $i \in \{1,\ldots, 5\}$ and the invariant-subspace equations have the form $\dot{\vec{s}} = A_{\mu_1, \mu_2} \vec{s} + b_{\mu_1, \mu_2}$, where $\vec{s} = [s_1, s_2]^\top$ and
\begin{equation}
	A_{\mu_1, \mu_2} = \begin{bmatrix}
		A_{\mu_1} - 2 \sigma \D\mathcal{H} & 2 \sigma \D\mathcal{H}\\
		3 \sigma \D\mathcal{H} & A_{\mu_2} -3 \sigma \D\mathcal{H}
				\end{bmatrix} \,, \ \
b_{\mu_1, \mu_2} = \begin{bmatrix}
b_{\mu_1} \\
b_{\mu_2}
\end{bmatrix}\,, \quad \mu_i=  \begin{cases} 1\,, & v_i>0 \\ 2\,, & v_i<0 \end{cases}
\end{equation}
and we define $A_{1}$, $A_2$, $b_{1}$, and $b_2$ in \cref{tab:pwl2piecesmodels}. Also let
$\mathcal{H}(x) = [v, 0]^\top$ so that the coupling acts only through the first component.
This is a $4$-dimensional PWL system with two switching planes, $v_1=0$ and $v_2=0$. One can construct the periodic orbit on the $4$-dimensional synchronous manifold by following the method that we outlined in \cref{sec:pwl}. Starting from the initial data $\vec{s}(0) = [0, w_1(0), v_2(0), w_2(0)]^\top$, we now have to solve a system of seven nonlinear algebraic equations for $w_1(0)$, $v_2(0)$, and $w_2(0)$ and the four switching times $T_{1,1}$, $T_{2,1}$, $T_{2,2}$, and $T_{1,2}= \Delta$ (see \cref{Fig:switching}).


 With the block-diagonalization of the variational equation \eqref{eq:blockVariational}, one uses
the initial data and switching times to explicitly compute the Floquet multipliers of the periodic orbit.  One can compute Floquet multipliers that correspond to perturbations within the synchronization manifold without using the block-diagonalization. We have
\begin{equation}
	\FD{}{t}\delta {\vec{s}} = A_{\mu_1, \mu_2}\delta\vec{s} \,,
\end{equation}
which one can solve using matrix exponentials, being careful to use saltation matrices to evolve perturbations through switching manifolds. After one period, $\delta\vec{s}(\Delta) = \Psi_s\delta\vec{s}(0)$, where $\Psi_s$ is the monodromy matrix on the synchronization manifold.
Considering all evolutions and transitions through switching manifolds, we obtain
\begin{align}
	\Psi_s = S_{12} \mathcal{E}_{2,1}(\Delta- T_{1,1}) S_{11}  \mathcal{E}_{1,1}(T_{1,1}-T_{2,2}) S_{22}  \mathcal{E}_{1,2}(T_{2,2}-T_{2,1})S_{21}  \mathcal{E}_{1,1}(T_{2,1}) \,,
\end{align}
with saltation matrices
\begin{align}
	S_{ij} &= P_i \otimes S_i(T_{i,j}) \,, \notag  \\
	P_1 &= \begin{bmatrix}
 1& 0 \\
0 & 0
\end{bmatrix}\,,
	\,\, P_2 = \begin{bmatrix}
0& 0 \\
0 & 1
\end{bmatrix} \,, \notag \\
\label{eq:saltation}
	 S_i(t) &= \begin{bmatrix}
\dot{v}_i(t^+)/\dot{v}_i(t^-) & 0 \\
(\dot{w}_i(t^+)-\dot{w}_i(t^-))/\dot{v}_i(t^-) & 1
\end{bmatrix}\,, \quad i \in \{1,2\}\,, \\
	\mathcal{E}_{\mu_1, \mu_2}(t) &= \e^{A_{\mu_1, \mu_2}t}\,, \qquad \mu_i \in \{ 1, 2\} \,.
\end{align}
The Floquet multipliers for perturbations within the synchronization manifold are the eigenvalues of the monodromy matrix $\Psi_s$. One of these eigenvalues is always $1$, corresponding to perturbations along the periodic orbit.

The block-diagonalization of $\mathcal{L}$ for the cluster state that we have been discussing is \cite{Sorrentino2016, Nicks2018}
\begin{equation}
	\mathcal{L}'' =\begin{bmatrix} 3 &-\sqrt{6} &0 & 0 & 0\\ -\sqrt{6} & 2 &0 & 0 & 0  \\ 0 &0 &5 & 0 & 0\\0 & 0 &0 & 3 & 0  \\ 0 & 0 & 0 & 0 & 3\end{bmatrix} \,.
\end{equation}
 In the directions that are transverse to the synchronization manifold, this block-diagonalisation yields the following three decoupled Floquet problems:
\begin{equation}
\begin{aligned}
	\dot{V}_3 &= (D f(s_1) -3\sigma D \mathcal{H})V_3\,,  \\
	\dot{V}_4 &= (D f(s_2) -3\sigma D \mathcal{H})V_4\,,  \\
	\dot{V}_5 &= (D f(s_1) -5\sigma D \mathcal{H})V_5\,,
\end{aligned}
\end{equation}
which (as usual) one can solve using
matrix exponentials and saltation matrices. This yields $V_i(\Delta) = \Psi_{t_i} V_i(0)$, where the monodromy matrices $\Psi_{t_i}$ for the transverse directions are
\begin{equation}
\begin{aligned}
	\Psi_{t_3} &= S_1(\Delta) \mathcal{E}_L^3(\Delta-T_{1,1}) S_1(T_{1,1}) \mathcal{E}_R^{3}(T_{1,1}) \, ,   \\
	\Psi_{t_4} &=\mathcal{E}_R^3(\Delta- T_{2,2}) S_2(T_{2,2}) \mathcal{E}_L^3(T_{2,2}-T_{2,1}) S_2(T_{2,1}) \mathcal{E}_R^3(T_{2,1}) \,,  \\
	\Psi_{t_5} &= S_1(\Delta) \mathcal{E}_L^5(\Delta-T_{1,1}) S_1(T_{1,1}) \mathcal{E}_R^5(T_{1,1}) \,,
\end{aligned}
\end{equation}
and
\begin{equation}
	\mathcal{E}_\mu^{\beta}(t) = \e^{(A_\mu-\beta\sigma D\mathcal{H})t} \,.
\end{equation}
The eigenvalues of the mondromy matrices $\Psi_{t_i}$, with $i \in {3,4,5}$ give the Floquet multipliers for directions that are transverse to the synchronization manifold.

The change of basis from $U$ coordinates to $V$ coordinates has no effect on the action of the saltation matrices. (Recall that $V = (Q \otimes I_m) U$.)
To
evolve $U$ through a discontinuity, we write
$U^+ = S U^-$, where
\begin{equation}
	S = \sum_{k=1}^M E^{(k)} \otimes S_k \,.
\end{equation}
Therefore, $V^+ =\widehat{S}  V$, where
\begin{equation}
	\widehat{S} = (Q \otimes I_m) S (Q\otimes I_m)^{-1}  = \sum_{k=1}^M (Q E^{(k)} Q^{-1} \otimes S_k) = \sum_{k=1}^M (J^{(k)} \otimes S_k)\,.
\end{equation}

Because the vector field of the absolute model is continuous, all saltation matrices are the identity matrix.
One then does an algebraic calculation to show that the cluster state is stable for the choice of parameters in \cref{Fig:Pwl2pieces}(a). One finds bifurcations of the periodic orbit by determining when the Floquet multipliers leave the unit disk. As one varies the parameters, the order of the times at which trajectories cross the switching planes can also change. One constructs bifurcation diagrams by similarly treating all types of cluster states from network symmetries and Laplacian clustering.

\begin{figure*}[t]
\begin{center}
\includegraphics[width=12cm]{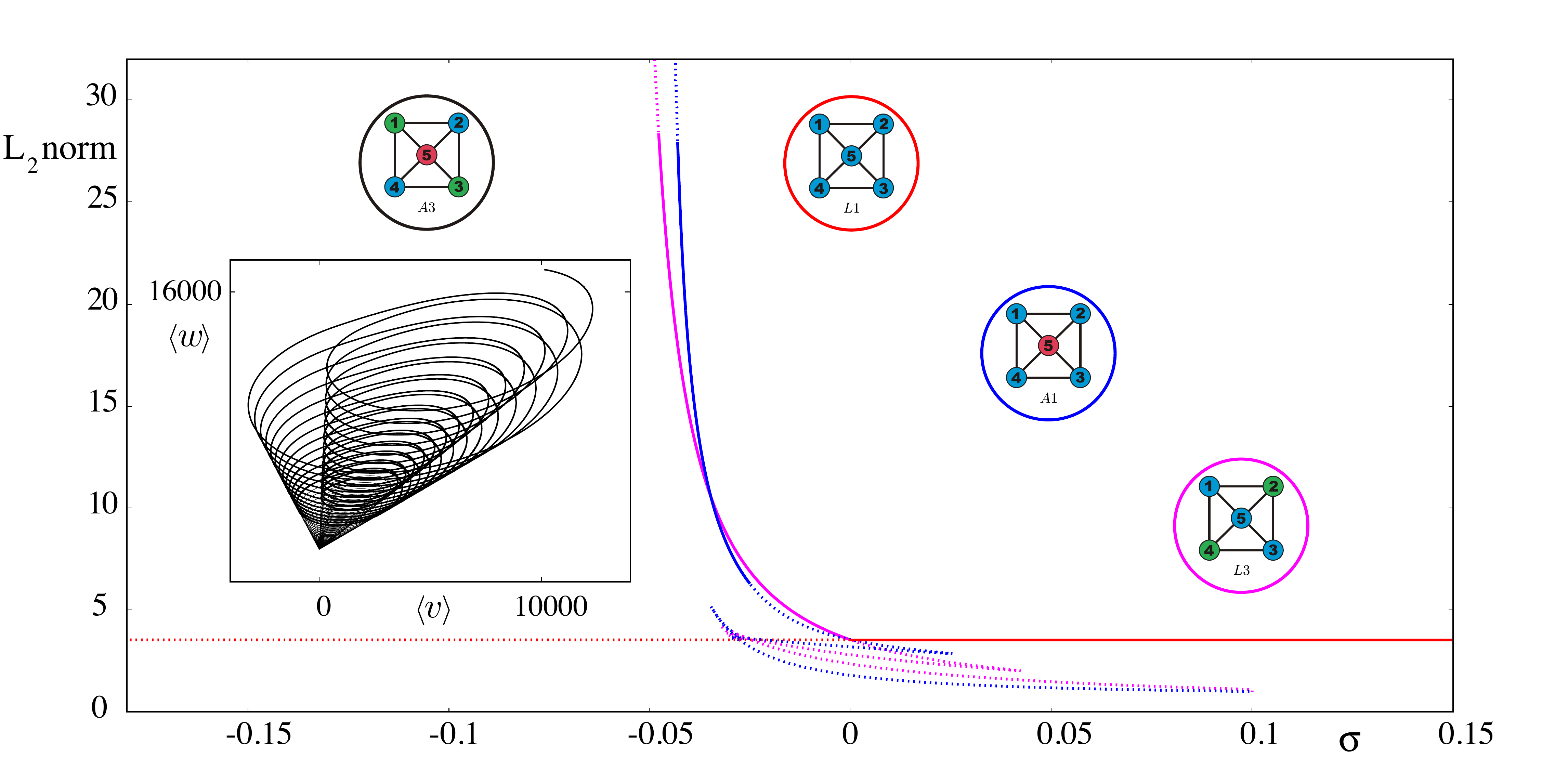}
\caption{Bifurcations between cluster states from varying the coupling strength $\sigma$ in a five-node network of absolute-model oscillators for the same parameters as in \cref{Fig:Pwl2pieces}(a). We indicate stable periodic orbits with solid curves and unstable solutions with dotted curves. We color each branch according to the circle with its associated cluster state.
We show only one branch of the $L3$ solutions, where we expect a branch of conjugate solutions with identical stability properties.
In the inset, we show the mean-field dynamics $[\langle v \rangle, \langle w \rangle ] = \sum_{i=1}^5 [v_i, w_i]/5$ for $\sigma= -0.05$. This shows an $A3$ cluster state with dynamics that blow up in finite time.
 This behavior dominates for $\sigma \lessapprox -0.0477$, where the $L3$ branch loses stability. All of the depicted bifurcations from stable states are tangent bifurcations, in which a Floquet multiplier passes through the value $+1$.
}
\label{Fig:AbsoluteBifDiag}
\end{center}
\end{figure*}

For the absolute model with the choice of parameters in \cref{Fig:Pwl2pieces}(a) and interaction function $\mathcal{H}(x)= [v,0]^\top$, we show the bifurcations from varying the
coupling strength $\sigma$ in \cref{Fig:AbsoluteBifDiag}.
All of the bifurcations from stable states are tangent bifurcations \cite{kuznetsov2013elements},
in which a Floquet multiplier passes through the value $+1$.


One can use the above approach to determine the stability of cluster states in any network of PWL nodes; see \cite{Nicks2018} for more examples. The computational difficulty of applying the MSF approach for a cluster state scales with the number of clusters in the state and with the number of switching planes in the PWL model of the individual oscillators. It does \emph{not} scale with the size (i.e., the number of nodes) of a network. Finally, we note that one can view synchrony as a single-cluster state, for which the above methodology reduces to the standard MSF approach in \cref{subsec:strongly}.

\subsection{An application to synaptically coupled, spiking neural networks}

It is common to model spiking neural networks using integrate-and-fire (IF)
neurons. Coombes \emph{et al.}~\cite{Coombes2012} explored the nonsmooth nature of systems of IF neurons.
The MSF approach has been used to study synaptically coupled networks of nonlinear (specifically, adaptive exponential) IF neurons \cite{Ladenbauer2013}, for which one uses numerical computations to obtain periodic orbits.
 Nicks \textit{et al.}~\cite{Nicks2018} showed how to make analytical progress on
 the dynamics of PWL planar IF neurons.

We follow Nicks \textit{et al.} \cite{Nicks2018} and consider a network of $N$ synaptically coupled planar IF neurons with
the time-dependent forcing
$I \rightarrow I + \sigma \sum_j w_{ij} s_j(t)$.
The synaptic input from neuron $j$ takes the standard event-driven form \cref{eventsynapse}
We adopt the common choice of a continuous $\alpha$-function, so that $\eta(t)$ is \cref{alphafunction}. We can then express $s_i(t)$ as the solution to the impulsively forced linear system
\begin{equation}
	\left ( 1 + \frac{1}{\alpha} \FD{}{t} \right ) s_i = u_i\,, \quad
	\left ( 1 + \frac{1}{\alpha} \FD{}{t} \right ) u_i = \sum_{p \in \ZSet} \delta(t-t_i^p)\,.
\label{history}
\end{equation}

We exploit the linearity of the synaptic dynamics between firing events to write
the network model in the form \cref{MSFnetworkMoreGeneral} with
$\dot{x}_i = f(x_i)$, where $x_i = (v_i,w_i,s_i,u_i)$ and $f$ has the form \cref{pwl1}, with
\begin{equation}
	A_{1,2}=\begin{bmatrix}
a_{1,2} & -1 & 0 & 0 \\
a_w/\tau & b_w/\tau & 0 & 0 \\
0 & 0 & -\alpha & \alpha \\
0 & 0 & 0 & -\alpha
	\end{bmatrix}
\end{equation}
and $b_1=[I,0,0,0]^\top=b_2$, and one applies the jump operator $\mathcal{J} (x_i) =  (v_{\text{r}}, w_i +\kappa/\tau, s_i, u_i+ \alpha)$ whenever $h(x_i)=v_i-v_{\text{th}} = 0$.
The vector function that specifies the interaction is $\mathcal{H}(x_i) = [s_i,0,0,0]^\top$.

\begin{figure}
\begin{center}
\includegraphics[width=8cm]{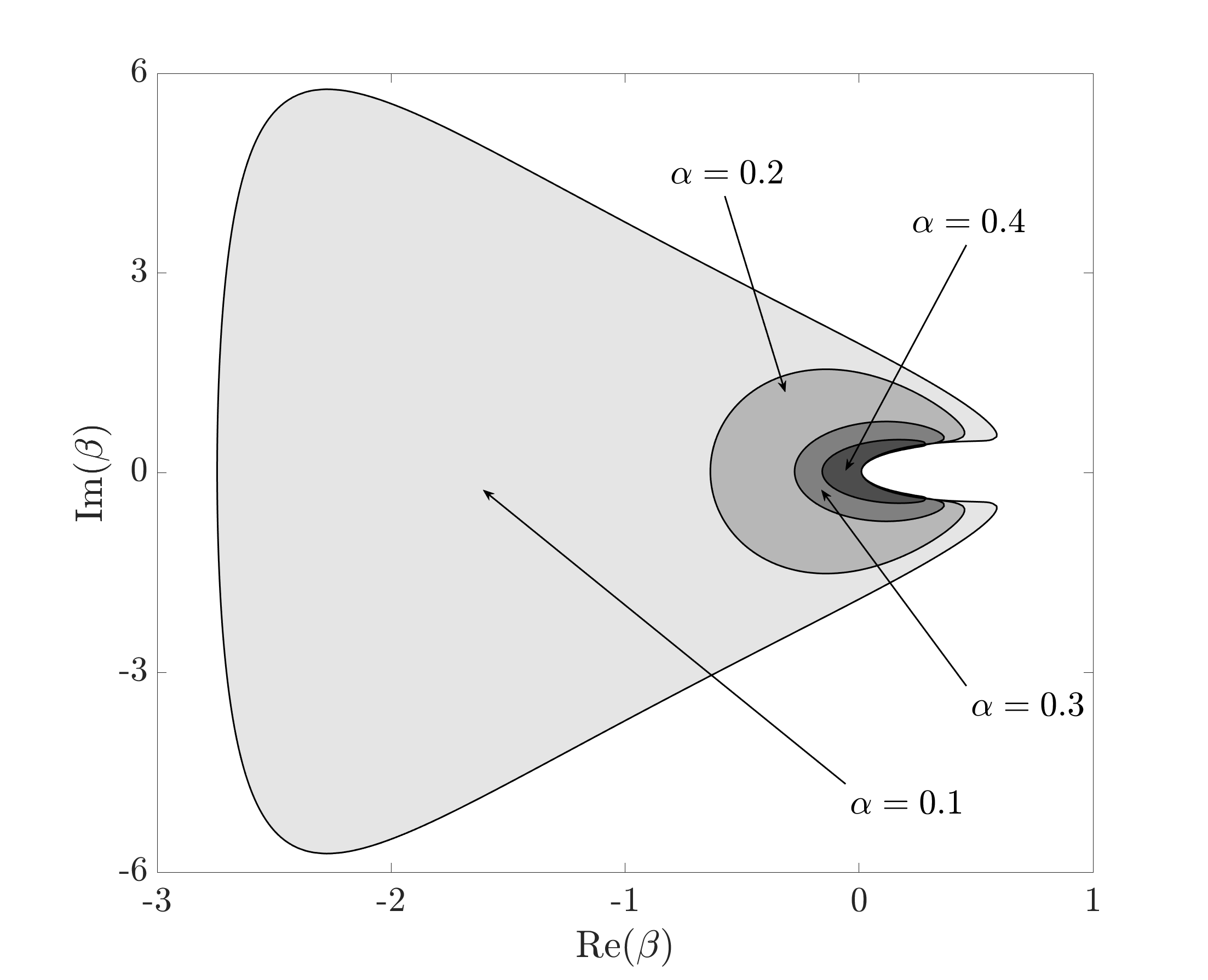}
\caption{The MSF for a network of synaptically coupled, planar integrate-and-fire (IF) neurons for the synchronous tonic orbit in \cref{Fig:Pwl2pieces}(d). The shaded regions indicate where the MSF is negative for various values of the synaptic rate parameter $\alpha$. The largest depicted region is for $\alpha=0.1$, with progressively smaller areas for $\alpha=0.2$, $\alpha = 0.3$, and $\alpha = 0.4$.  The synchronous solution is stable if all of the eigenvalues of $\sigma w$ lie within a shaded area for a given value of $\alpha$.
\label{Fig:MSFIF}}
\end{center}
\end{figure}

For a synchronous orbit of the type
in \cref{Fig:Pwl2pieces}(d) (so that a trajectory only visits the region of phase space that is described by $A_2$ and the $T$-periodic trajectory satisfies the constraints $v(T)=v_\text{th}$, $w(0)=w(T)+\kappa/\tau$, $s(0)=s(T)$, and $u(0)=u(T)+\alpha$), we only need to consider saltation at firing events and the saltation matrix takes the explicit form
\begin{equation}
	S(t) = \begin{bmatrix}
\dot{v}(t^+)/\dot{v}(t^-) & 0 & 0 & 0\\
(\dot{w}(t^+)-\dot{w}(t^-))/\dot{v}(t^-) & 1 & 0 & 0 \\
(\dot{s}(t^+)-\dot{s}(t^-))/\dot{v}(t^-) & 0 & 1 & 0\\
(\dot{u}(t^+)-\dot{u}(t^-))/\dot{v}(t^-) & 0 & 0 & 1\\
	\end{bmatrix} \,.
\end{equation}
See \cref{Saltationproof} for the general formula for the saltation operator of a PWL system.
In this case, the expression for $\Psi$ in \cref{Gamma} reduces to
\begin{equation}
	\Psi = S(T) \exp \{ (A_2 + \beta \D \mathcal{H}) T \} \, ,
\label{Gamma1}
\end{equation}
where $\beta = \sigma \lambda_l$ and $\lambda_l$ is the $l$th eigenvalue of $w$.
The matrix $\D \mathcal{H}$ is a constant matrix with entries $[\D \mathcal{H}]_{ij} = 1$ if $i=1$ and $j=3$ and $[\D \mathcal{H}]_{ij} = 0$ otherwise.
Therefore, using equation \cref{Gamma1} and the prescription in \cref{subsec:strongly}, we are able to construct MSF (see \cref{Fig:MSFIF}).


\begin{figure}[hbtp]
\centering
\begin{subfigure}[t]{0.32\textwidth}
\centering
\begin{overpic}[clip,trim=-0.2in 0.0in -0.1in -0.3in,width=\linewidth]{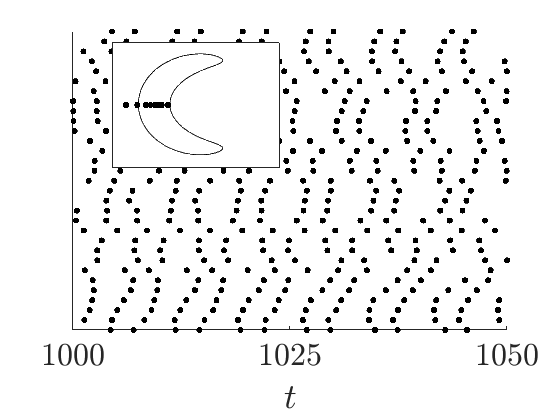}
\put(1,67){(a)}
\end{overpic}
\end{subfigure}
\begin{subfigure}[t]{0.32\textwidth}
\centering
\begin{overpic}[clip,trim=-0.2in 0.0in -0.1in -0.3in,width=\linewidth]{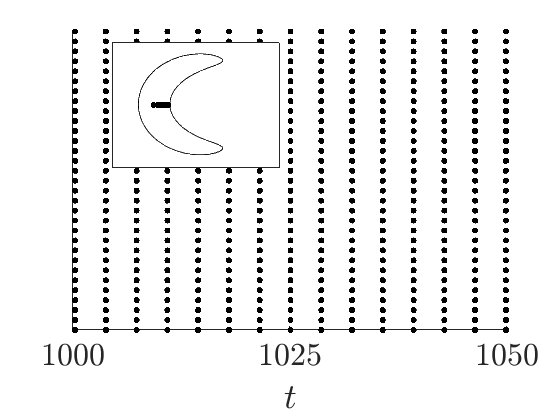}
\put(1,67){(b)}
\end{overpic}
\end{subfigure}
\begin{subfigure}[t]{0.32\textwidth}
\centering
\begin{overpic}[clip,trim=-0.2in 0.0in -0.1in -0.3in,width=\linewidth]{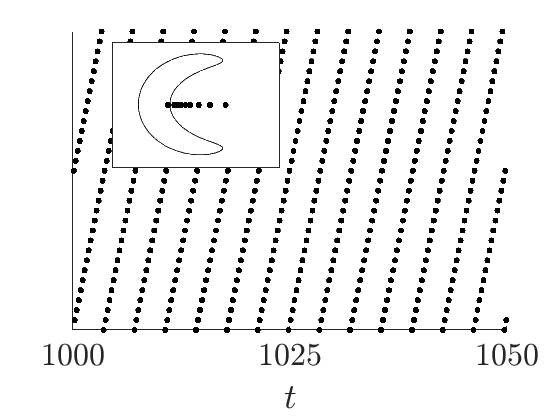}
\put(1,67){(c)}
\end{overpic}
\end{subfigure}
\caption{Raster plot of spike times from direct numerical simulations of a network of synaptically coupled, planar integrate-and-fire (IF) neurons with $N=31$ oscillators, $d=3$, and $\alpha=0.4$. A raster plot allows us to convey neuron-by-neuron variations in spike times.
In the inset, we plot the MSF and superimpose the eigenvalues of $\sigma w$. In (a), $\sigma=-0.1$ and synchrony is unstable. In (b), $\sigma=-0.025$ and synchrony is stable. In (c), $\sigma=0.1$ and synchrony is unstable.
The predicted instability borders (at $\sigma=0$ and $\sigma \approx -0.05$) are in good agreement with the predictions from the nonsmooth MSF analysis.  For $\sigma > 0$, the typical pattern of firing activity beyond an instability of the synchronous state is a periodic traveling wave. For $\sigma <0$, a spatiotemporal pattern emerges via a period-doubling instability of the firing times.
\label{Fig:pwlIFsim}}
\end{figure}

As a particular realization of a network architecture that guarantees synchrony, we use a balanced ring network with odd $N$ and $w_{ij} = w(|i-j|)$.
We calculate the distances $|i - j|$ modulo $(N-1)/2$ and use the decay rate $w(x) = (1-a|x|/d) \e^{-|x|/d}$. We choose the parameter $a$  so that $\sum_{j=1}^N w_{ij} = 0$ (a balance condition) for the network
size $N$ and a scale $d$.
The eigenvalues $\lambda_l$ of the associated (symmetric and circulant) adjacency
matrix are real and given by $\lambda_l = \sum_{j=0}^{N-1} w(|j|) \omega_l^j$.
The balance condition enforces $\lambda_0=0$.  Additionally, $\lambda_{N-l} = \lambda_l$ for $l \in \{1,\ldots, (N-1)/2\}$, so any excited pattern 
(which arises from an instability) 
is given by a combination $e_m+e_{-m} = 2 \, \text{Re} (e_m)$ for some
$m \in \{1, \ldots, (N-1)/2\}$.
Given the shape of the MSF function in \cref{Fig:MSFIF}, one determines the value of $m$ using $\lambda_m= \max_{l} \lambda_l$. In \cref{Fig:pwlIFsim}, we compare direct simulations of a network versus the predictions of the MSF.  When the network's eigenvalues lie within the region where the MSF is negative, small perturbations of synchronous initial data decay
away and the system settles to a synchronous periodic orbit, as expected.  When one of the eigenvalues crosses the $0$ level set of the MSF from negative to positive, two types of instabilities emerge.
 One of them leads to a
 spatiotemporal pattern of spike doublets (i.e., a burst of two spikes), which arise because an eigenvalue of $\Gamma$ leaves the unit disk at $-1$ (through a period-doubling bifurcation), and the other yields a periodic traveling wave (with asynchronous firing) because an eigenvalue of $\Gamma$ leaves the unit disk at $+1$ (through a tangent bifurcation).


\subsection{An application to neural-mass networks}

The human brain has roughly \(10^{11}\) neurons and roughly \(10^{15}\) synapses.
Although there is general consensus that the synaptic interactions between neurons drive brain dynamics, these astronomical numbers prohibit the construction, analysis, and simulation of an entire brain network that is built from single-neuron models such as the absolute model or the Morris--Lecar model (see {\cref{sec:pwl}}). Instead, it is instructive to coarse-grain neural behavior by grouping neurons and studying the interactions between these groups. This idea led to neural-mass models~\cite{Ashwin2016}, which describe the average dynamics of large populations of neurons.

One of the most influential neural-mass models is the Wilson--Cowan model~\cite{Wilson73,Wilson72}
\begin{equation}
	\FD{u}{t}  = -u + F(I_u+ w^{uu}u-w^{vu}v)\,, \quad \tau \FD{v}{t} = - v + F(I_v + w^{uv} u - w^{vv}v) \,,
\label{eq:WC}
\end{equation}
where \(u\) and \(v\), respectively, indicate the activity of excitatory and inhibitory populations of neurons. A firing-rate function \(F(x)\), which researchers often take to have a sigmoidal shape, mediates the interactions between the two populations. The quantities \(I_{u,v} \) represent background inputs, and \(w^{\alpha \beta}\) (with $\alpha, \beta \in \{u,v\}$) denote connection strengths between populations.
 The positive constant \(\tau\) encodes the relative time scale between the dynamics of the two populations.

To make analytical progress, Coombes \etal \cite{Coombes20181} considered a PWL firing-rate function of the form
\begin{equation}
	F(x) = \begin{cases}
0 \,, & x \leq 0 \\
\epsilon^{-1}x \,,& 0 < x < \epsilon \\
1 \,,& x \geq \epsilon \,.
	\end{cases}
\label{pwlrate}
\end{equation}
With this choice, it is straightforward to compute periodic orbits of the dynamical system \eqref{eq:WC} and to determine their linear stability using the techniques that we described in \cref{sec:pwl}. One can think of the system \eqref{eq:WC} as modeling an appropriately chosen brain region, so coupling oscillators that satisfy \eqref{eq:WC} lets one investigate the dynamics of interacting brain regions. By introducing the coupling matrices \(\mathcal{W}^{\alpha \beta} \in \mathbb R^{N\times N}\), with \(\alpha,\beta \in \{u,v\} \), we obtain
a network of \(N\) interacting oscillators
with dynamics
\begin{align}
	\FD{u_i}{t}  &= -u_i + F \left (I_u+ \sum_{j=1}^N \mathcal{W}^{uu}_{ij}u_j-\sum_{j=1}^N \mathcal{W}^{vu}_{ij}v_j \right ) \,, \\
\tau \FD{v_i}{t} &= - v_i + F \left (I_v + \sum_{j=1}^N \mathcal{W}^{uv}_{ij} u_j - \sum_{j=1}^N \mathcal{W}^{vv}_{ij}v_j \right) \,, \qquad i \in \{1,\ldots,N\} \,.
\label{eq:WCnetwork}
\end{align}
Although the dynamical system \cref{eq:WCnetwork} is not exactly in the form that we described in \cref{subsec:strongly}, one can analyze this network using essentially the same MSF techniques. For simplicity and to guarantee the existence of a synchronous network state, we impose the
{row-sum} constraint \(\sum_{j=1}^N \mathcal{W}^{\alpha,\beta}_{ij} = w^{\alpha \beta}\) for \(\alpha,\beta \in \{u,v\} \). These row-sum constraints are natural for networks arranged on a ring, because the coupling matrix is circulant (see \cref{sec:spectra}). The synchronous network state satisfies $[u_i(t), v_i(t)] = [u(t), v(t)]$ for all $i \in \{1,\ldots,N\}$, where $[u(t), v(t)]$ satisfies equations
\cref{eq:WC}.

It is convenient to introduce the vector
 $X=[u_1, v_1, u_2, v_2, \ldots , u_N, v_N]^\top \in \RSet^{2N}$ and change variables by writing $Y = \mathcal{W} X + C$, where
\begin{equation}
	\mathcal{W} =
\mathcal{W}^{uu} \otimes \begin{bmatrix} 1 & 0 \\ 0 & 0	\end{bmatrix}
- \mathcal{W}^{vu} \otimes \begin{bmatrix} 0 & 1 \\ 0 & 0	\end{bmatrix}
+ \mathcal{W}^{uv} \otimes \begin{bmatrix} 0 & 0 \\ 1 & 0	\end{bmatrix}
- \mathcal{W}^{vv} \otimes \begin{bmatrix} 0 & 0 \\ 0 & 1	\end{bmatrix} \,,
\end{equation}
the matrix $C = \vec{1}_N \otimes [I_u, I_v]^\top$, and $\vec{1}_N$ is an $N$-dimensional vector with all entries equal to $1$.
One can then succinctly describe the switching manifolds by the relations $Y_i=0$ and $Y_i=\epsilon$, and the dynamics takes the form
\begin{equation}\label{eq:WCN_dY}
	\FD{}{t} Y = \mathcal{A} (Y- C) + \mathcal{W} \mathcal{J} F(Y)\,,\quad \mathcal{J}=I_N \otimes \begin{bmatrix}
1 & 0 \\
0 & 1/\tau
\end{bmatrix}  \,,
\end{equation}
where \( \mathcal{A} = - \mathcal{W} \mathcal{J} \mathcal{W}^{-1}\). We denote the synchronous solution by $\overline{Y}(t)= (\overline{U}(t), \overline{V}(t)$, $\overline{U}(t), \overline{V}(t),  \ldots, \overline{U}(t), \overline{V}(t))$, with
\begin{equation}
	\begin{bmatrix}
\overline U(t) \\ \overline V(t)
	\end{bmatrix}
=
	\begin{bmatrix}
w^{uu} &-w^{vu}\\
w^{uv} &-w^{vv}
	\end{bmatrix}
	\begin{bmatrix}
\overline u(t) \\ \overline v(t)
	\end{bmatrix}
+
	\begin{bmatrix}
I_u \\ I_v
	\end{bmatrix}\,,
\end{equation}
and consider small perturbations such that $Y = \overline{Y} + \delta Y$. We thereby obtain
\begin{equation}
	\FD{}{t} \delta Y = \mathcal{A}\,  \delta Y + \mathcal{W} \mathcal{J} \, \D F(\overline{Y})\, \delta Y  \,,
\label{eq:perturbN}
\end{equation}
where $\D F(\overline{Y})$ is the Jacobian of $F$ evaluated along the periodic orbit.

As we showed in \cref{subsec:strongly}, we need appropriately diagonalize \eqref{eq:perturbN}. Suppose that we can diagonalize all \(\mathcal{W}^{\alpha \beta}\) with respect to the same basis, and let $P=[e_1 ~e_2 ~\ldots ~e_N]$ be the matrix whose columns consist of the basis vectors. Such simultaneous diagonalization is feasible for circulant matrices, which naturally obey the above row-sum constraint.
Let \(\{\nu_j^{\alpha\beta}\}\), with $j \in \{1,\ldots,N\}$, denote the eigenvalues of \(\mathcal{W}^{\alpha \beta}\). We then write
\begin{equation}
	(P \otimes I_2)^{-1}\mathcal{W} (P \otimes I_2) = \operatorname{diag}(\Lambda_1,\Lambda_2, \ldots, \Lambda_N) \equiv \Lambda \, ,
\end{equation}
where
\begin{equation}
	\Lambda_p = \begin{bmatrix} \nu^{uu}_p & - \nu^{vu}_p \\
\nu^{uv}_p & - \nu^{vv}_p
\end{bmatrix} \,, \quad p \in \{1,2,\ldots, N\} \,.
\end{equation}
Additionally, $(P \otimes I_2)^{-1}\mathcal{A} (P \otimes I_2) = -\Lambda (I_N \otimes J) \Lambda^{-1}$.

Consider perturbations of the form $\delta Z = (P \otimes I_2)^{-1}\delta Y$. Equation \cref{eq:perturbN} then implies that the linearized dynamics satisfies
\begin{equation}
	\FD{}{t} \delta Z = \Lambda (I_N \otimes J)\left [ -\Lambda^{-1}  + (I_N \otimes \D)\right ] \delta Z \,,
\label{deltaZ}
\end{equation}
where $\D \in \RSet^{2 \times 2}$ is the Jacobian of $(F(\overline{U}),F(\overline{V}))$. The matrix $\D$ is a
piecewise-constant matrix that is nonzero only if either $0 < \overline{U}(t) < \epsilon$ or $0  < \overline{V}(t) < \epsilon$.
Analogously to
\cref{MSFnetworkVarEigfirst}, equation \cref{deltaZ} has a block structure in which the dynamics in each of $N$ $2 \times 2$ blocks satisfies
\begin{equation}
	\FD{}{t} \xi = [A_p + \Lambda_p J \D]\xi  \,, \quad p \in \{1,\ldots, N\} \,, \quad \xi \in \RSet^2 \,,
\label{cool}
\end{equation}
with $A_p= - \Lambda_p  J \Lambda^{-1}_p$.

The problem that is defined by \cref{cool} is time-independent between switching manifolds, so
one can construct a solution in a piecewise fashion from matrix exponentials and write $\xi(t) = \exp [(A_p + \Lambda_p J \D)t] \xi(0)$. One can then construct a perturbed trajectory over one period of oscillation in the form $\xi(\Delta)=\Gamma_p \xi(0)$, where $\Psi_p \in \RSet^{2 \times 2}$ is
\begin{equation}
	\Psi_p= \e^{A_p\Delta_8}\e^{A^-_p(\epsilon)\Delta_7}\e^{A_p\Delta_6}\e^{A^-_p(\epsilon)\Delta_5}\e^{A_p\Delta_4}\e^{A^+_p(\epsilon)\Delta_3}\e^{A_p\Delta_2}\e^{A^-_p(\epsilon)\Delta_1} \,,
\label{eq:Gamma}
\end{equation}
with
\begin{equation}
	A^\pm_p(\epsilon)=\left( A_p + \epsilon^{-1} \Lambda(p) J T^\pm \right)\,,\quad
T^+=   \begin{bmatrix} 1 & 0 \\ 0 & 0 \end{bmatrix}\,, \quad
	T^-  =\begin{bmatrix} 0 & 0 \\ 0 & 1 \end{bmatrix}\,.
\end{equation}
Note
(and see equation \cref{Gamma}) that there is no saltation (i.e., $S = I$).

\begin{figure}[htbp]
\begin{center}
\includegraphics[width=5cm]{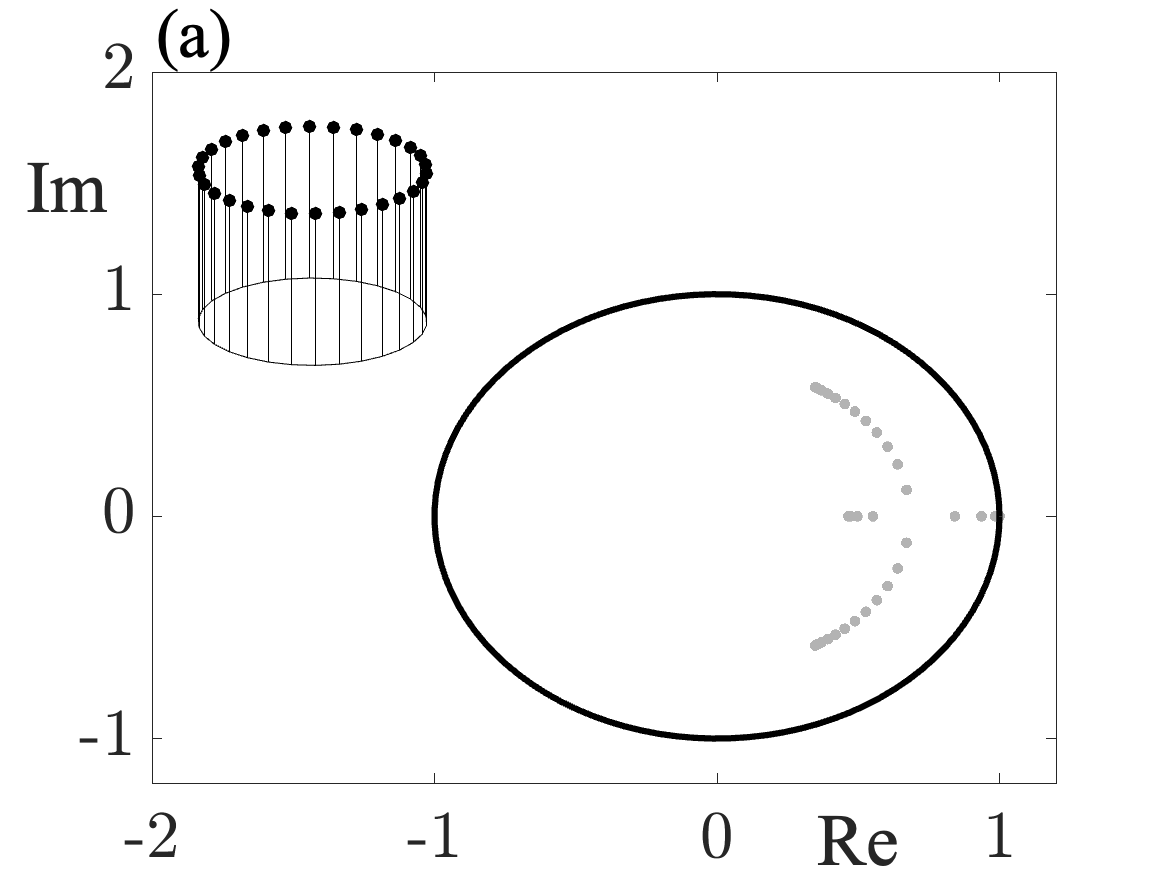} \hspace*{.75cm}
\includegraphics[width=5cm]{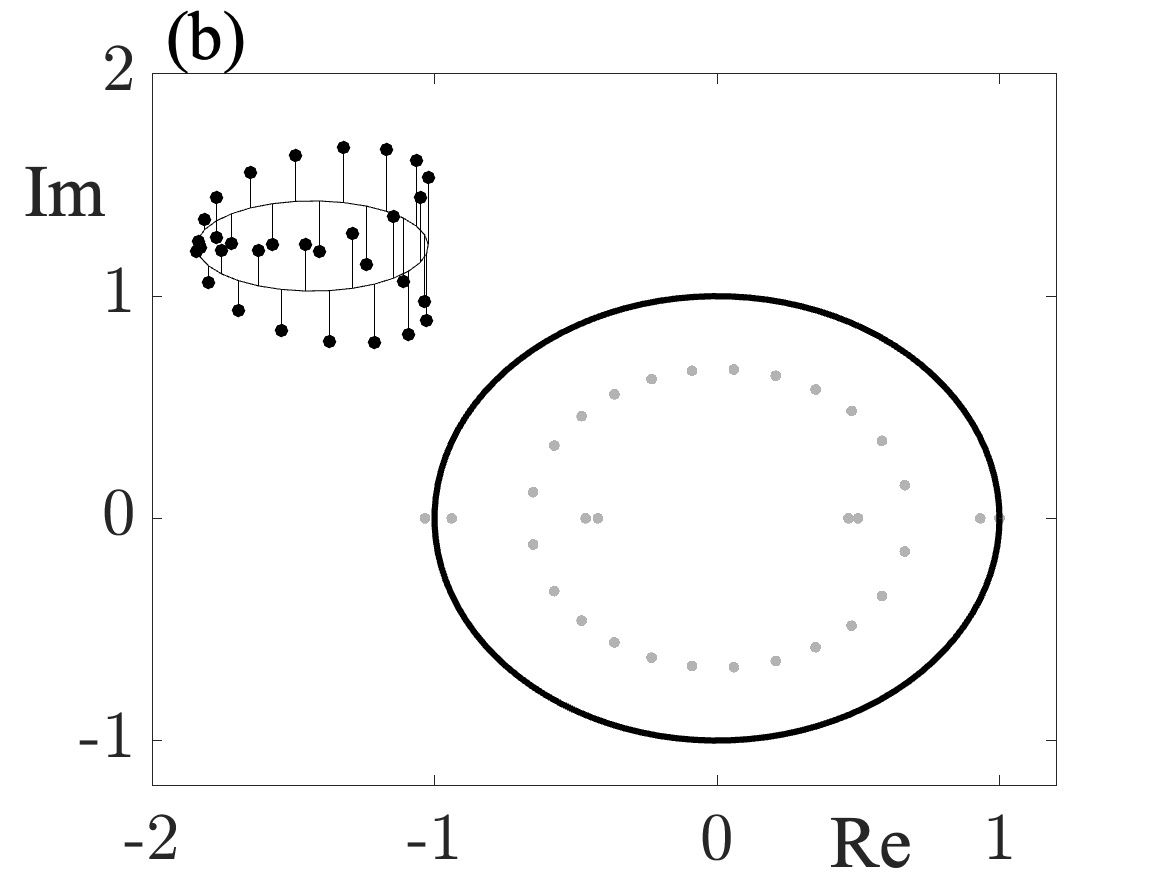}
\caption{Plots of the eigenvalues of \cref{eq:Gamma} for a
ring of $N = 31$ Wilson--Cowan oscillators, with
$\sigma_{\alpha \beta}=\sigma$ for all $\alpha$ and $\beta$.
The parameter values are $\epsilon=0.04$, $\tau=0.6$, $I_u=-0.05$, $I_v=-0.3$, $w^{uu}=1$, $w^{vu}=2$, $w^{uv}=1$, and $w^{vv}=0.25$ for coupling strengths of (a) $\sigma=0.15$ and (b) $\sigma=0.191$. The inset of (a) shows the synchronous network state, and the inset of (b) shows the eigenvector that is associated with the eigenvalue that lies outside the unit disk.
\label{Fig:PWLSpectrumN31}
}
\end{center}
\end{figure}

As an illustration, consider a network of Wilson--Cowan oscillators
on a ring graph with an odd number of nodes.  Let $\operatorname{dist}(i,j)=\min\{|i-j|, N-|i-j|\}$ be the distance between nodes $i$ and $j$. We then define a set of exponentially decaying connectivity matrices
\begin{equation}
	\mathcal{W}^{\alpha \beta}_{ij} = w^{\alpha \beta} \frac{ \e^{-\operatorname{dist}(i,j)/ \sigma_{\alpha \beta}}}{\sum_{j=0}^{N-1} \e^{-\operatorname{dist}(0,j)/\sigma_{\alpha \beta}}} \,.
\end{equation}
In this example, are four circulant matrices; they are parametrized by the four quantities
$\sigma_{\alpha \beta}$, which respect the row-sum constraints $\sum_{j=1}^N \mathcal{W}^{\alpha \beta}_{ij} = w^{\alpha \beta}$.
In \Cref{Fig:PWLSpectrumN31}, we plot the eigenvalues of $\Psi_p$ for $p \in \{1,\ldots, N\}$ for two different parameter choices. In \cref{Fig:PWLSpectrumN31}(a), all of the eigenvalues (excluding the one at \((1,0)\) that arises from time-translation invariance) lie within the unit disk. In \cref{Fig:PWLSpectrumN31}(b), one eigenvalue leaves the unit disk along the negative real axis.
This latter scenario predicts an instability of the synchronous state. In the inset of panel (b), we show the eigenvector that corresponds to the eigenvalue that crosses
to the outside of the unit disk.
 The prediction of this instability is in excellent agreement with direct numerical simulations \cite{Coombes20181}.

\subsection{An application to cardiac alternans}

One can conceptualize a beating heart as a network of muscle cells
in which each heartbeat results from their coordinated contraction and subsequent relaxation.
 Because the dynamics of an organ results
 from the orchestrated behavior of individual cells, a major avenue in cardiac research is
 the investigation of the dynamical repertoire of individual cardiac muscle cells \cite{Qu2014}.
 Molecular changes at the individual-cell level can yield irregular cell behavior, which then feeds forward to pathological heart
 dynamics (such as cardiac arrhythmias). A vital component that determines the behavior of cardiac muscle cells is the intracellular calcium (\ca) concentration. In basic terms, rises and falls of the cytosolic \ca concentration are responsible for muscle contraction, and irregularities and
 abnormalities of the intracellular \ca dynamics have been linked to a plethora of cardiac pathologies \cite{Landstrom2017}.

The intracellular \ca concentration in cardiac muscle cells has rich spatiotemporal patterns that arise from the interplay of diffusively coupled calcium-release units (CRUs). One can decompose each CRU into
compartments with \ca fluxes between them, so
a cardiac muscle cell corresponds to a network of networks. In other words, each node of the cellular network is itself a network (i.e., a CRU).
Using a 5-dimensional PWL representation of a well-established cardiac \ca model \cite{Thul2010}, one can express the dynamics of a network of \(N\) CRUs as
 \begin{equation}
	\FD x t = A x + F(t) + \mathcal{L} \otimes H x\,,
\label{eq:general_network}
\end{equation}
where $x = (x_1, x_2, \ldots, x_N)$ is a $5N$-dimensional vector. Each entry $x_\mu$, with $\mu \in \{1,\ldots, N$\}, is
the 5-dimensional state vector of a single CRU. The matrix $A \in \mathbb R^{5N \times 5N}$ is constant and block diagonal with entries in the set \(\{A_i\}\).  The constant matrices \(A_i \in \mathbb R^{5 \times 5} \) 
are associated with a single CRU, analogously to the matrices $A_1$ and $A_2$
in equation \cref{pwl1}.
As usual, the matrix $\mathcal L \in \mathbb N^{N \times N}$ denotes the combinatorial graph
Laplacian matrix of the network and the matrix $H \in \mathbb R^{5 \times 5}$ encodes which variables are coupled and how strongly they are coupled.
The time-dependence $F(t)=1_N \otimes v(t) \in \mathbb R^{5N}$ distinguishes the present example from the other examples in this section. 
The explicit time-dependent drive \(v(t) \), which is \(\Delta\)-periodic, models an experimental condition that is known as a voltage clamp, which is used often to disentangle the different cellular mechanisms that contribute to the complex spatiotemporal patterns of the intracellular \ca concentration of cardiac cells.

Because of the explicit time-dependence in equation \cref{eq:general_network}, the switching manifolds are not only state-dependent (as in all of the previous examples in this section), but some of them are also time-dependent. This leads to a 
system in which any trajectory is determined by a sequence of state-dependent and time-dependent switches \cite{Thul2010,lai2019master,Veasy2019}. As demonstrated in \cref{sec:pwl}, one can readily compute the synchronous network state \(s(t) \) of \cref{eq:general_network} using matrix exponentials. One can then linearize equation \cref{eq:general_network} around the synchronous network state $s(t)$ by using the ansatz $x(t) = 1_N \otimes s(t) + \delta x$ and following the general approach in \cref{sec:strongly}. 
Analogously to equation \cref{MSFnetworkVarEigfirst}, this yields 
\begin{equation}
	\FD{\xi_l}{t} = \left [A_i - \lambda_l H \right ]\xi_l \,.
\label{eq:MSF_components}
\end{equation}
where \(\xi_l \in \mathbb R^5\) and \(\lambda_i\) are the eigenvalues of \(\mathcal L\). Because we are perturbing from the synchronous network state, we assume that all CRUs have the same associated \(A_i\) to obtain \cref{eq:MSF_components}.

\begin{figure}[h!]
\includegraphics[width=0.5\textwidth]{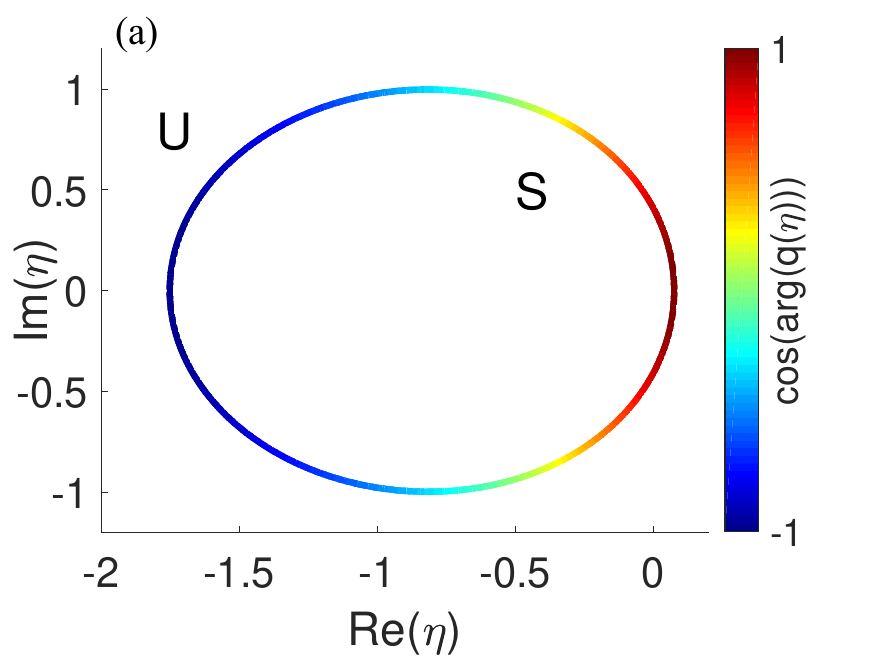}
\includegraphics[width=0.5\textwidth]{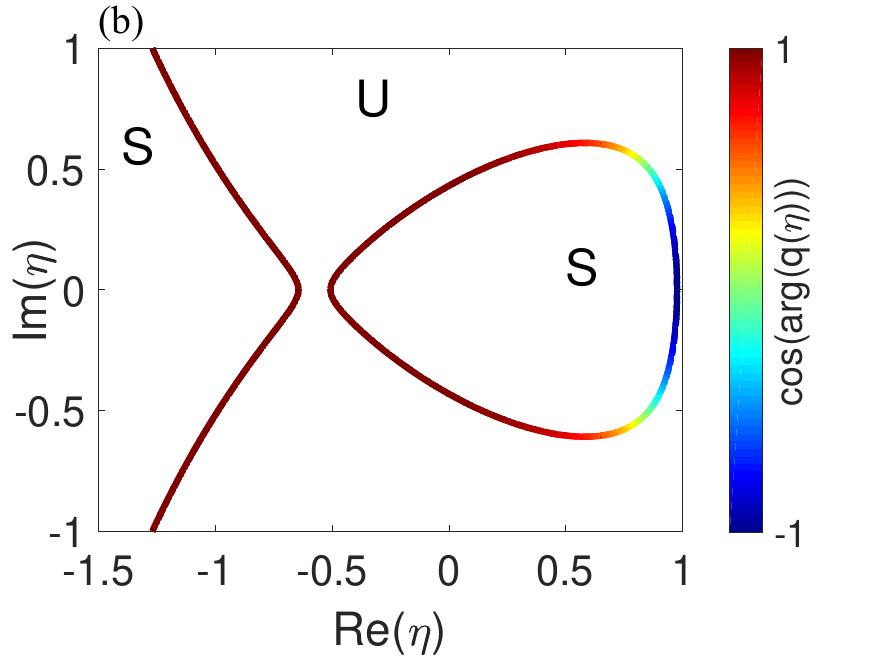}
\caption{Zero contour of the MSF for $\Delta = 0.9$ and two different coupling regimes. 
The MSF is negative in regions that we label by ``S" and positive in regions that we label by ``U". 
The color indicates the value of $\cos \left ( \arg \left ( q \left (\eta \right)\right)\right)$, where $q$ is the largest eigenvalue of $\Psi(\eta)$ and $\Psi$ is given by 
\cref{eq:Gamma-cardiac}.
The synchronous solution is stable if all of the points with $\eta = \lambda_l$ lie in the region S, where $\lambda_l$ are the eigenvalues of the graph Laplacian $\mathcal{L}$ (which also incorporates the coupling strength of the network). For more details, see \cite{lai2019master}.
}
\label{fig:cardiac MSF}
\end{figure}

The dynamical system
 \cref{eq:MSF_components} is continuous at the switching manifolds, so one can obtain its solution
 using matrix exponentials. Let $\Delta_i$ denote the time-of-flight for when the dynamics are associated with  $A_i$, let $\Delta=\sum_i \Delta_i$
 denote the period of the synchronous state, and let $\xi(0)$ denote an initial perturbation.
 According to the relation \cref{Gamma} and noting that there is no saltation (i.e., $S=I$), the perturbation after one period is $\xi_l(\Delta)=\Psi(\lambda_l) \xi_l(0)$, where
\begin{equation}
	\Psi(\lambda_l)=\exp\left[(A_m -  \lambda_l H)\Delta_N \right] \times \cdots \times \exp\left[(A_1 -  \lambda_l H)\Delta_1 \right]\,.
\label{eq:Gamma-cardiac}
\end{equation}
As we showed in \cref{subsec:strongly}, one can use the relation \cref{eq:Gamma-cardiac} to construct the MSF. In \Cref{fig:cardiac MSF}, we illustrate that the topology of the MSF can vary substantially across different coupling regimes. In the left panel of \cref{fig:cardiac MSF}, the zero contour of the MSF forms a closed loop; the MSF is negative inside the loop and positive outside it.
 On the contrary, in the right panel of \cref{fig:cardiac MSF}, there are
 two distinct regions in which the MSF is negative. 
 The colors 
 reveal that if the MSF changes sign along the real axis, then the only instabilities are either a period-doubling bifurcation (i.e., a $-1$ bifurcation) or a tangent bifurcation (i.e., a $+1$ bifurcation).

\begin{figure}[h!]
\centering
\includegraphics[width=0.9\textwidth]{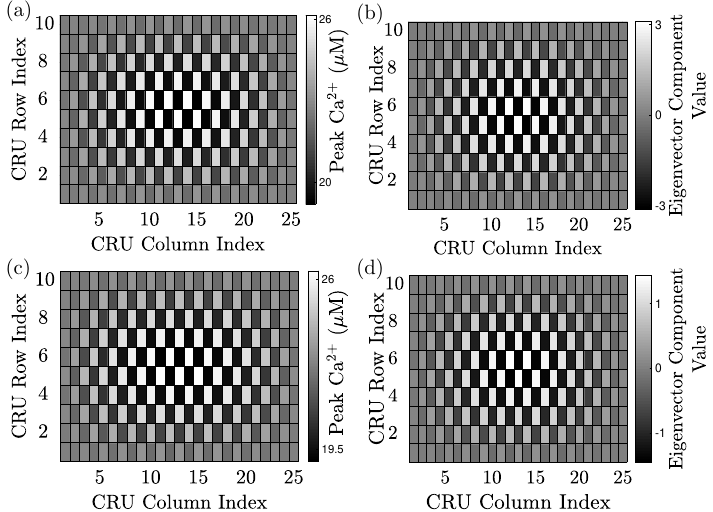}
\caption{Instabilities of the synchronous network state induced by (a, b) a period-doubling bifurcation and (c, d) a tangent bifurcation. Panels (a, c)  show the peak \ca concentration during one period $\Delta$ in one of the CRU compartments, while panels (b, d) depict the eigenvectors that correspond to the only eigenvalue that lies outside the unit disk in for the patterns in panels (a, c).  For more details, see \cite{lai2019master}.}
\label{fig:2d-TP_0_9}
\end{figure}

As was demonstrated by Lai \emph{et al.} \cite{lai2019master}, MSF plots like those in \cref{fig:cardiac MSF} allow one to understand
abrupt changes in spatial \ca patterns from small 
changes of a single parameter. In \Cref{fig:2d-TP_0_9}, we illustrate the
behavior that emerges when the synchronous solution
destabilizes via a period-doubling bifurcation or a tangent
bifurcation. 
In \Cref{fig:2d-TP_0_9}(a), we show the peak \ca concentration during one period of the $\Delta$-periodic drive $v(t)$ for a period-doubling bifurcation. The CRUs are arranged on a regular two-dimensional grid, so one can reference each CRU
by a row and column index. Each small rectangle represents the \ca concentration of a single CRU.
One can clearly see a spatially alternating pattern, which is more pronounced near the center of the figure than it is near the edges.
This spatially alternating pattern also alternates in time: when a CRU has a large peak during one pacing period, it has a small peak during the next pacing period (and vice versa). In other words, each CRU exhibits a period-2 orbit and adjacent CRUs oscillate out of phase with each other. This phenomenon is known as ``subcellular \ca alternans'' and is a precursor to severe cardiac arrhythmia. 
For the behavior
in \cref{fig:2d-TP_0_9}(a), only one eigenvalue lies outside the unit disk. In \cref{fig:2d-TP_0_9}(b), we show the corresponding right eigenvector, which is in excellent agreement with direct numerical simulations. 
When an eigenvalue leaves the unit disk along the positive real axis, one observes a pattern like the one in \cref{fig:2d-TP_0_9}(c). As with the period-doubling bifurcation, the peak \ca concentrations exhibit an alternating spatial pattern. However, in contrast to the period-doubling bifurcation, each CRU follows a period-1 orbit, so the peak amplitude is the same across pacing periods, rather than alternating from one pacing period to the next. In \Cref{fig:2d-TP_0_9}(d), we show the eigenvector that is associated with the only eigenvalue that lies outside the unit disk for the behavior in \cref{fig:2d-TP_0_9}(b). This eigenvector also agrees excellently with direct numerical simulations.

\subsection{An application to Franklin-bell networks}

Benjamin Franklin was one of the leading political figures of his time, and he was also a prolific inventor and scientist. To facilitate his studies into the nature of electricity, he employed lightning as an electrical power source. To be notified when an iron rod outside his house was sufficiently electrified by lighting, Franklin employed what is now known as a ``Franklin bell'' \cite{Labaree:1962}.
A Franklin bell is a metal ball that oscillates between two metal plates, which are driven by electrical charge. A Franklin bell is an example of an impacting system; the ball velocity changes nonsmoothly when it contacts either plate.
A network of $N$ Franklin bells satisfies the dynamical system \cite{ShirokyGendelman:2016,csayli2019synchrony}
\begin{alignat}{2}
\label{originalnetworkModelling}
	\ddot{u}_n+\gamma_1 \dot{u}_n+\gamma_2 u_n+ \sigma \displaystyle \sum_{m=1}^{N} w_{nm}\left( u_{m}-u_{n}\right) &= \sgn(\dot{u}_n) f\,, \quad  t \neq t_{n_{i}}\,, \\
\label{originalnetworkModelling1}
	\dot{u}_n(t_{n_{i}}^{+}) &= -k\dot{u}_n(t_{n_{i}}^{-})\,, \quad \ t = t_{n_{i}}\,,
\end{alignat}
where $u_n$ denotes the position of the ball of the \(n\)th Franklin bell, which is restricted between two impacting manifolds at \(\pm a\). One implicitly determines the time $t_{n_{i}}$ of the $i$th impacting event of the $n$th oscillator using the relation $u_n(t_{n_{i}})= \pm a$. The parameter $\sigma$ is a global coupling strength, and the network structure is encoded by a matrix with elements $w_{nm}$. The constant \(k \in \RSet^+\) is the coefficient of restitution upon impact, $f$ is a constant force (which is determined by a sum of the repelling and attracting electrostatic forces), $\gamma_{1} > 0$ is a damping coefficient, and $\gamma_2 > 0$ sets the natural frequency of the pendulum.

It is convenient to write the dynamical system \cref{originalnetworkModelling}, \cref{originalnetworkModelling1}
as a system of first-order differential equations (i.e., in the standard form of a dynamical system) by introducing the state vector $x_n=[u_n,v_n]^\top$, where $v_n = \dot{u}_n$. This yields
\begin{alignat}{2}
\label{statespaceModellingnetwork}
	\dot{x}_{n} &= F(x_n)+ \sigma \displaystyle \sum_{m=1}^{N} w_{nm}[\mathcal{H}(x_{m})-\mathcal{H}(x_{n})]\,, \quad  \ t \neq t_{n_{i}}\,, \\
	x_n(t_{n_{i}}^{+}) &=  g \left (x_n(t_{n_{i}}^{-})\right) \,, \quad  \ t = t_{n_{i}}\,.
\end{alignat}
The vector field $F: \RSet^2 \rightarrow \RSet^2$ is
$F(x_n)= Ax_n+f_{e_n}$, where
\begin{equation}
\label{Af}
	A=
		\begin{bmatrix}
			0 & 1\\-\gamma_2 & -\gamma_1
		\end{bmatrix}\,, \quad
	f_e =
		\begin{bmatrix}
			0 \\
			f
	\end{bmatrix} \sgn(v) \,.
\end{equation}
The function $\mathcal{H}: \RSet^2 \rightarrow \RSet^2$ is $\mathcal{H}[u,v])=[0,u])^\top$.  The form of the coupling in equation \cref{statespaceModellingnetwork} ensures the existence of the synchronous network state \(s(t)\). To determine its linear stability, we rewrite \cref{statespaceModellingnetwork} using the graph Laplacian $\mathcal L$ (see \cref{subsec:strongly}). The
dynamics between impacts is
\begin{equation}\label{graphlaplacian}
	\dot{x}_{n}(t)=F(x_n(t))-\sigma\displaystyle \sum_{m=1}^{N}\mathcal L_{nm} \mathcal{H}(x_{m})\,,
\end{equation}
which has the form
\cref{MSFnetworkEqn}.
Consequently, after diagonalization, the Floquet problem for the linear stability of the synchronous network state
becomes $\xi_l(\Delta)=\Psi(l)\xi_l(0)$ (with $l \in \{1,\ldots,N\}$), where
\begin{equation}\label{blockstabilitymatrix}
	\Psi(l)=K(t_2)\e^{A_l \Delta_{2}}K(t_{1})e^{A_l \Delta_{1}}\,, \quad A_l=A-\sigma \lambda_l \D \mathcal{H}\,,
\end{equation}
and the saltation operator is
\begin{equation}
	K(t) = \begin{bmatrix}
		-k & 0 \\
		\frac{k \dot{v}(t^-)+\dot{v}(t^+)}{\dot{u}(t^-)} & -k
	\end{bmatrix} \, .
\end{equation}

\begin{figure}[!h]
\centering
\includegraphics[width=4cm]{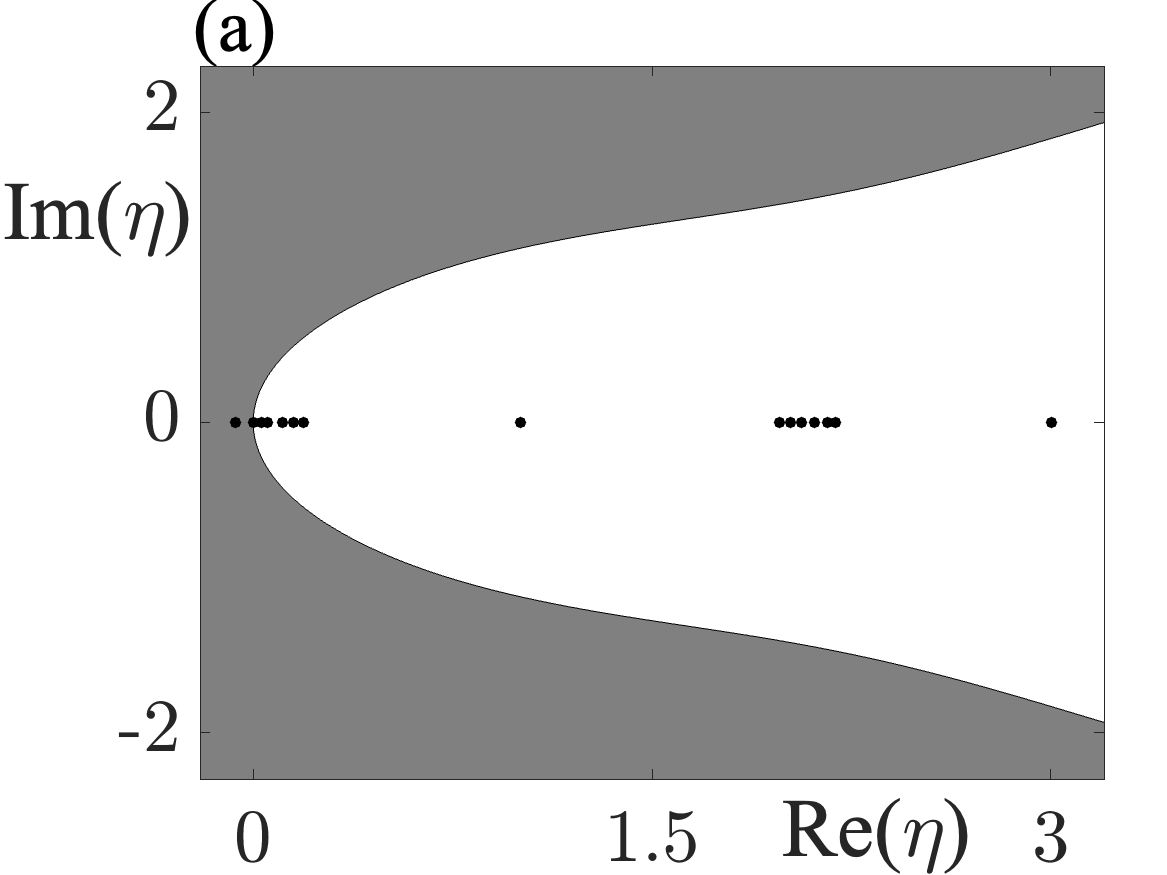}
\includegraphics[width=4cm]{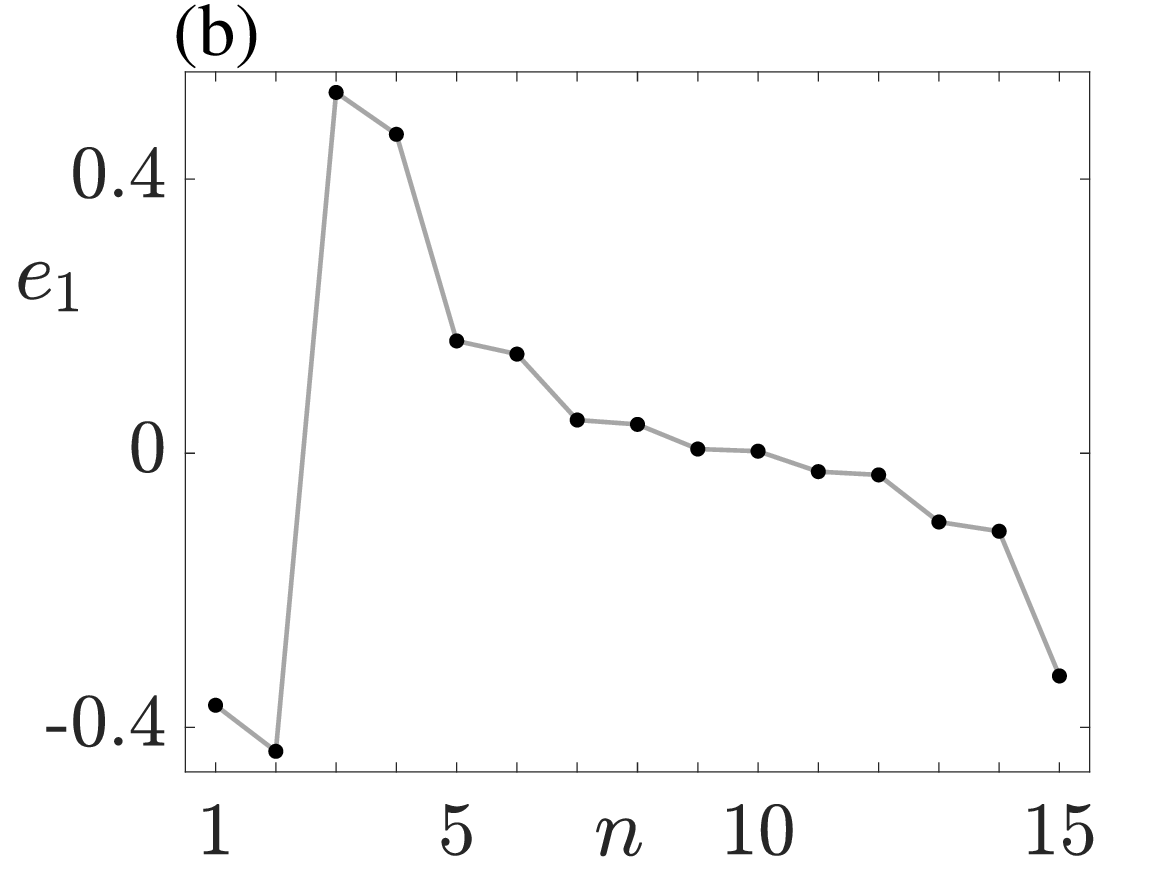}
\includegraphics[width=4cm]{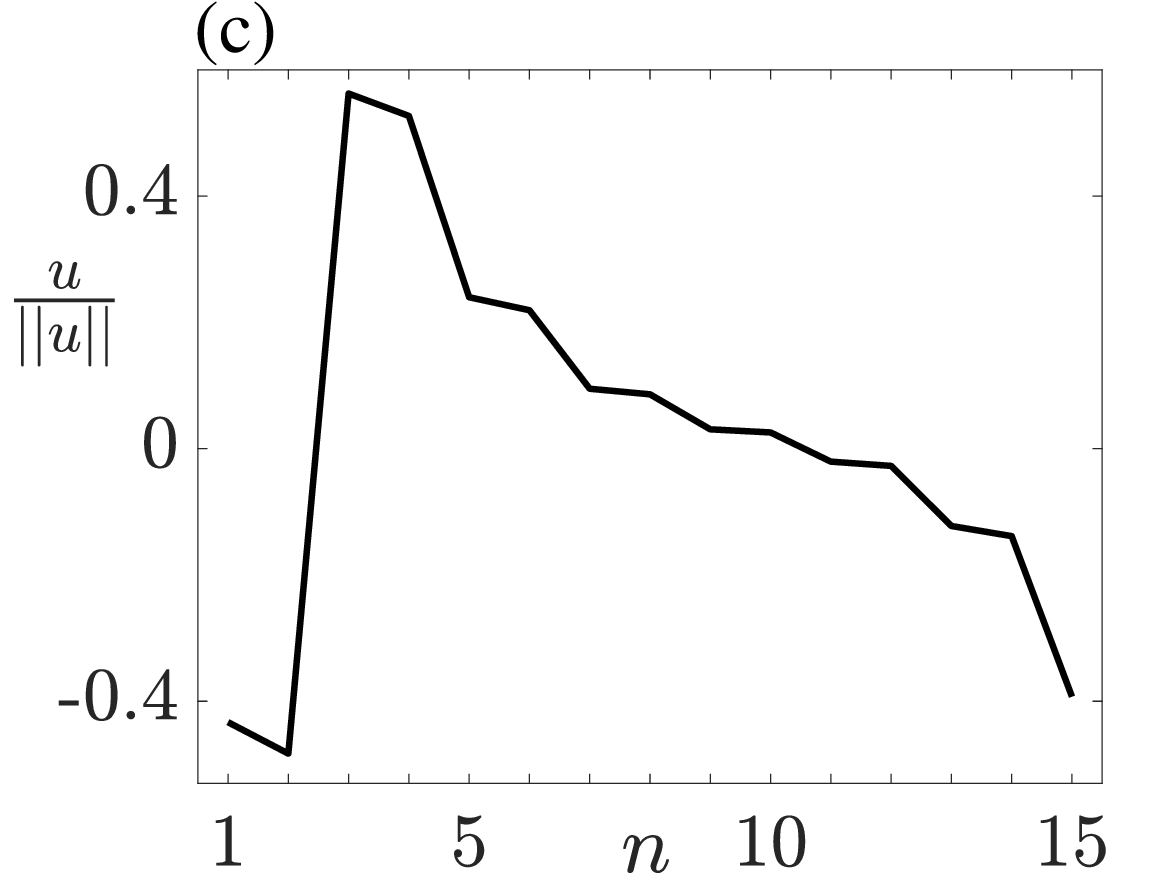}
\caption{(a) MSF and the values of $\eta_l$ (black dots) for a Franklin-bell network with $15$ nodes, where $c_{n}=1$ if $n$ is odd and $c_{n}=0.1$ if $n$ is even, except for $c_2=-0.1$. The MSF is negative in the white region and positive in the gray region. (b) Normalized eigenvector $e_1$ that corresponds to the eigenvalue in panel (a) for which the MSF is negative (the leftmost black dot).
(c) Normalized position $u_n$ at a fixed time for each oscillator $n$ in the network. For the other parameter values, see \cite{csayli2019synchrony}.
}
\label{Fig:Franklin Bell network}
\end{figure}

As we showed in \cref{subsec:strongly}, we obtain the MSF from \cref{blockstabilitymatrix} with the replacement $\sigma \lambda_l \rightarrow \eta \in \CSet$. In \Cref{Fig:Franklin Bell network}, we illustrate the dynamics of a network of 15 Franklin bells when the adjacency
matrix has entries
\begin{equation}
	\label{ringconnectivity}
	w_{nm}=c_{n}(\delta_{n,m-1}+\delta_{N-n+1,1})+c_{n-1}\delta_{n,m+1}+c_{N}\delta_{1,N-m+1}\,, \quad c_n \in \RSet
\end{equation}
for \(n,m \in \{1,2,\ldots,N\}\).
The MSF (see \cref{Fig:Franklin Bell network}(a)) is negative for only one \(\eta_l\), so the synchronous state is unstable.
In \cref{Fig:Franklin Bell network}(b), we show the eigenvector that corresponds to the critical value \(\eta_l\). We see in \cref{Fig:Franklin Bell network}(c) that we obtain excellent agreement with direct numerical simulations. Because the adjacency
matrix
with the entries \cref{ringconnectivity} is symmetric, all of the eigenvalues are real, so \(\eta_l\) is real. See {\c{S}}ayli \emph{et al.} for a discussion of the predictive power of the MSF for an example with complex eigenvalues.

\subsection{An application to coordination in cow herds}

Grazing animals, such as cows,
protect themselves from predators by living in herds \cite{mendel2001}, and synchronizing their behavior (by tending to eat and lie down at the same time) helps them remain together as a herd \cite{rook1991}.
Sun \emph{et al.} \cite{Sun2011} developed a piecewise-linear dynamical system as a simplistic model to study collective in herds of cattle.
One can treat some aspects of their model --- both for a single cow and for a network of cows --- using the focal techniques of the present paper.

Cows are ruminants.
They eat
plant food, swallow, and regurgitate it at some later stage; they then again chew the partly digested plant food.
During the first stage (standing/feeding), they stand up to graze. However, they typically
lie down and ruminate (i.e., chew the cud) in a second stage (lying/ruminating). A cow thus oscillates between two stages. One can construct a simplistic caricature of a cow by separately considering its observable state (eating, lying down, or standing) and its unobservable level of hunger or desire to lie down. Sun \emph{et al.} formulated a model in terms of a
variable $x(t;\theta)$, where $x=(v,w) \in [\delta ,1 ] \times [\delta, 1]$ with a parameter $\delta \in (0,1)$ and an observable state
 $\theta \in \{ \mathcal{E}, \mathcal{R}, \mathcal{S} \}$.  The variables $v$ and $w$, respectively, represent the extent of the desire of a cow to eat and lie down. The
 variable $\theta$ represents the state of a cow, which can be eating ($\mathcal{E}$), lying down ($\mathcal{R}$), or standing ($\mathcal{S}$).  The dynamics in the $(v,w)$ plane is confined to a box, and a cow switches between states whenever a trajectory intersects with the edge of the box.

The dynamics of $x=x(t;\theta)$ takes the simple form $\dot{x} = a(\theta) x$, where
\begin{equation}
	a(\theta) = \begin{bmatrix}
\alpha(\theta) & 0 \\
0 & \beta(\theta)
\end{bmatrix} \,,
\quad
\begin{bmatrix}
\alpha(\mathcal{E}) & \alpha(\mathcal{R}) & \alpha(\mathcal{S})  \\
\beta(\mathcal{E}) & \beta(\mathcal{R})  & \beta(\mathcal{S})
\end{bmatrix}
=
\begin{bmatrix}
-\alpha_2 & \alpha_1 & \alpha_1 \\
\beta_1 & -\beta_2 & \beta_1
\end{bmatrix} \,,
\end{equation}
with hunger parameters $\alpha_{1,2} > 0$ and lying parameters $\beta_{1,2}>0$. The parameter $\alpha_1$ represents the rate of increase of hunger, $\alpha_2$ represents the decay rate of hunger, $\beta_1$ represents the rate of increase of the desire to lie down, and $\beta_2$ represents the decay rate of the desire to lie down. One prescribes switching conditions (at the four edges of the box) using four indicator functions: $h_1(x) = v-1$, $h_2(x)=w-1$, $h_3(x)=v-\delta$, and $h_4(x)=w-\delta$. The model's four state-transition rules take the form $\theta \rightarrow g_i (\theta)$, where $g_1(\theta) = \mathcal{E}$, $g_2 (\theta) = \mathcal{R}$, and $g_3(\theta) = g_4(\theta) = \mathcal{S}$.
If a trajectory intersects the corner of the box, one can apply a state-tiebreaker rule~\cite{Sun2011}, although we will not consider such scenarios.

\begin{figure}[h]
	\centering
\includegraphics[width=8cm]{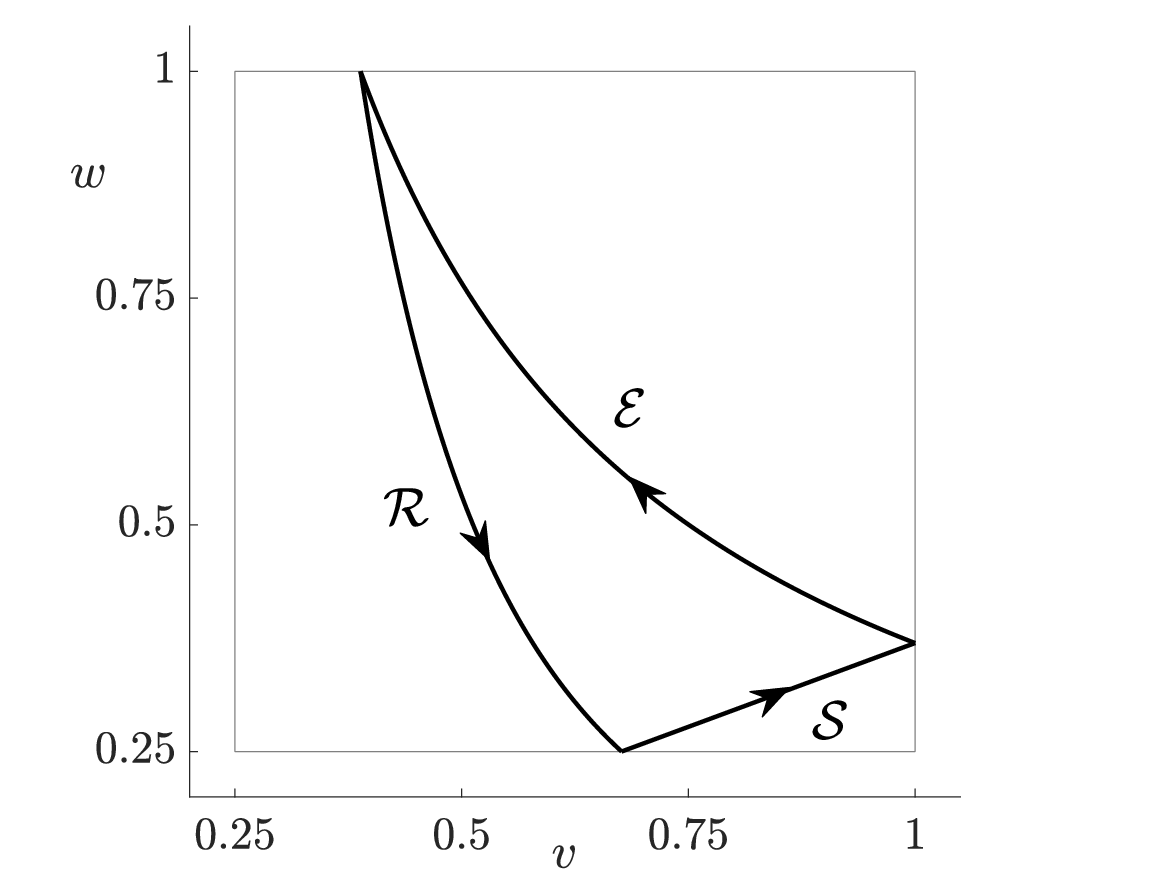}
	\caption{A periodic orbit with a cyclic state-transition sequence $\mathcal{E} \rightarrow \mathcal{R} \rightarrow \mathcal{S}$ in the model for cow dynamics.
	The parameters are $\delta=0.25$, $\alpha_1=0.05$, $\beta_1=0.05$, $\alpha_2=0.95\alpha_1$, and $\beta_2=0.125$.}
	\label{Fig:cow}
\end{figure}

The general prescription in \cref{Delatxplusforminchapter} yields saltation matrices at each of the four possible state transitions. They take the explicit forms
\begin{align}
	S_1(t) = S_3(t) &= \begin{bmatrix}
\dot{v}(t^+)/\dot{v}(t^-) & 0 \\
(\dot{w}(t^+)-\dot{w}(t^-))/\dot{v}(t^-) & 1
		\end{bmatrix}\,, \notag \\
	S_2(t) = S_4(t) &= \begin{bmatrix}
1 & (\dot{v}(t^+)-\dot{v}(t^-))/\dot{w}(t^-) \\
0 & \dot{w}(t^+)/\dot{w}(t^-)
	\end{bmatrix} \,.
\end{align}
For a given state, one readily obtains phase-space trajectories as convex curves
$w = k v^{\beta(\theta)/\alpha(\theta)}$, with $k=w(t_0;\theta)/v(t_0;\theta)^{\beta(\theta)/\alpha(\theta)}$. The associated time evolution is
$(v(t;\theta), w(t;\theta)) = (\e^{\alpha(\theta)t}v(0;\theta), \e^{\beta(\theta)t}w(0;\theta))$ for $t \geq t_0$.
Many periodic orbits are possible, and one can catalog them in terms of state-transition sequences. Sun \emph{et al.} \cite{Sun2011} identified four low-period orbits, which have the following cyclic state-transition sequences: (a) $\mathcal{E} \rightarrow \mathcal{R} \rightarrow \mathcal{E}$, (b) $\mathcal{E} \rightarrow \mathcal{R} \rightarrow \mathcal{S} \rightarrow \mathcal{E}$, (c) $\mathcal{E} \rightarrow \mathcal{S} \rightarrow \mathcal{R} \rightarrow \mathcal{E}$, and (d) $\mathcal{E} \rightarrow \mathcal{S} \rightarrow \mathcal{R} \rightarrow \mathcal{S} \rightarrow \mathcal{E}$.
For example, consider
a periodic orbit of type (b). Starting from $x(0) = (1,w(0;\mathcal{E}))$, we obtain the time-of-flight $T_1 = -\beta_1^{-1}\log w(0;\mathcal{E})$ from the relation $h_2(x(T_1))=0$.  This, in turn, allows one to determine the initial data for the next piece of the trajectory. It is $x(T_1) = (v(T_1;\mathcal{E}),1)$, from which we obtain the time-of-flight $T_2 = -\beta_2^{-1} \log \delta$. The third and final piece of the orbit has initial data $x(T_1+T_2)=(v(T_2;\mathcal{R}),\delta)$ and a time-of-flight of $T_3 = -\alpha_1^{-1} \log v(T_2;\mathcal{R})$. To determine the value of $w(0;\mathcal{E})$, one enforces the periodicity condition $w(T;\mathcal{S})=w(0;\mathcal{E})$, where $T=T_1+T_2+T_3$ is the period of the periodic orbit. One thus obtains
\begin{equation}
	w(0;\mathcal{E}) = \delta^{\frac{1+\frac{\beta_1}{\beta_2}}{1+\frac{\alpha_2}{\alpha_1}}} \,.
\end{equation}
To ensure that $\delta < w(0;\mathcal{E}) < 1$, $\delta < v(T_1;\mathcal{E}) < 1$, and $\delta < v(T_2;\mathcal{R}) < 1$, the trajectory needs to satisfy
 the inequality $(\alpha_2/\alpha_1) \cdot (\beta_2/\beta_1) > 1$ and it also needs to satisfy $\alpha_1^{-1} + \alpha_2^{-1} \geq \beta_1^{-1}+ \beta_2^{-1}$ when $\beta_1 < \alpha_2$.
In \cref{Fig:cow}, we show an example of an orbit that one constructs in this fashion.

To determine the stability of this orbit, we calculate the Floquet exponent
\cref{Fexponentformula} and obtain
\begin{equation}
	\kappa= \kappa_{\text{smooth}}+ \frac{1}{T} \ln \left|   \det S_1(T) \det S_4(T_1+T_2) \det S_2(T_1) \right|  =\frac{1}{T} \ln \left ( \frac{\alpha_2}{\alpha_1} \right ) \,,
\end{equation}
where we use the fact that $\kappa_{\text{smooth}} = \left [T_1 (-\alpha_2+\beta_1) + T_2(\alpha_1-\beta_2) +T_3 (\alpha_1+\beta_1) \right ]/T=0$.
Therefore, the orbit (if it exists) is linearly stable for $q \equiv \alpha_2/\alpha_1 <1$.  The Floquet multiplier is $-q < 0$, so one can lose stability only through a period-doubling instability.

To model a herd (i.e., a network) of $N$ identical cows, we suppose that each cow has an associated variable $x_i(t;\theta_i)$ that evolves according to
\begin{equation}
	\FD{x_i}{t} = a(\theta_i) x_i + \sigma \sum_{j=1}^N w_{ij} \chi (\theta_j) x_j \,, \quad \chi(\theta) = \begin{bmatrix} \chi_{\mathcal{E}}(\theta) & 0 \\ 0 & \chi_{\mathcal{R}}(\theta) \end{bmatrix} \,,
\end{equation}
where $\chi_\psi$ is the indicator function
\begin{equation}
	\chi_\psi(\theta) = \begin{cases}
1\,, & \text{if}~ \theta=\psi \\
	0\,, & \text{otherwise} \,.
	\end{cases}
\end{equation}
For $w_{ij} >0$, this network model describes the case that a cow feels hungrier when it notices other cows eating and feels a greater desire to lie down when it notices other cows lying down.
Because the indicator functions change with time, the network has a time-dependent coupling, so one cannot use our previous MSF for it.
 Nonetheless, the PWL nature of the dynamics allows one to obtain analytical insights into the model's behavior.
  Assuming the row-sum normalization $\sum_{j=1}^N w_{ij}=1$ for all $i$, a synchronous orbit (if it exists) satisfies the equation $x_i(t)=s(t)$ for all $i$, where
$\dot{s} = a(\theta,\sigma)s$ and $a(\theta,\sigma) = a(\theta) +\sigma \chi(\theta)$. Therefore, one can use the approach that we described above for a single cow to
 construct a synchronous orbit under the replacement $a(\theta) \rightarrow a(\theta,\sigma)$.
For perturbations that do not change the order (i.e., the total number of states, including repetitions) in a
  state-transition sequence, linearization around the synchronous state leads to the evolution of the network perturbation $U=(\delta x_1, \delta x_2, \ldots, \delta x_N) \in\mathbb R^{2N}$ over one period of the form $U(T) = \Psi U(0)$, where
\begin{equation}
	\Psi  = K_1(T)\e^{A(\mathcal{S},\sigma) T_3}K_4(T_1+T_2)\e^{A(\mathcal{R},\sigma) T_2}K_2(T_1)\e^{A(\mathcal{E},\sigma) T_1} \,,
\label{bigGamma}
\end{equation}
with $A(\theta,\sigma)=I_N \otimes a(\theta) + \sigma w \otimes \chi(\theta)$ and $K_\mu(t_i) = I_N \otimes S_\mu(t_i)$.
Therefore, the synchronous state is linearly stable if all of the eigenvalues of $\Psi \in \RSet^{2N \times 2N}$ lie within the unit disk.

The above analysis allows one to generate a quantitative answer to the following question: Can herd interactions promote synchrony? Consider
 the choice $q>1$, so that an isolated cow
(i.e., $\sigma=0$) cannot achieve a stable $\mathcal{E} \rightarrow \mathcal{R} \rightarrow \mathcal{S} \rightarrow \mathcal{E}$ cycle. One can numerically calculate the eigenvalues of $\Psi$ \cref{bigGamma} to determine whether or not there is a critical value of $\sigma$ that brings all of the eigenvalues back inside the unit disk.  Numerical calculations for several types of row-normalized networks (e.g., star networks and nearest-neighbor circulant networks)
suggest that this is indeed the case, with a common critical value of $\sigma=\sigma_c$ that is independent of $N$.
For the parameters in \cref{Fig:cow} with $q=1.5$, we find that $\sigma_c \approx 0.025$.

\section{Discussion\label{sec:discussion}}

In this review, we discussed several popular mathematical frameworks for analyzing synchronized
states in coupled networks of identical oscillators.  
We focused on oscillator dynamics that take the form of piecewise-linear (PWL) ordinary-differential-equation models.
This choice allows the semi-analytical construction of periodic orbits without the need to employ
numerical ODE solvers. We 
demonstrated that it is also mathematically tractable to determine the stability of periodic states in networks of such coupled oscillators.
The key augmentation to standard theoretical approaches is the use of saltation operators to treat the nonsmooth nature of the individual oscillator models and the network models. 
We thereby highlighted
the usefulness of combining techniques from smooth dynamical systems --- in particular, weakly-coupled-oscillator theory and the master stability function (MSF) --- with techniques from nonsmooth modeling and analysis to deliver new tools for the analysis of dynamical systems on networks.

Given the prevalence of nonsmooth models in mechanics and biology (as well as in other areas), it is very appealing to further apply and extend these approaches. 
For example, one can apply such methodology 
to networks of scalar-valued
nodes with threshold-linear nonlinearities (of ReLu type, which is now ubiquitous in machine learning \cite{szand2021}),
which have become very popular for developing ideas about 
so-called {``sequential attractors''} \cite{Curto2016,Curto2013,Bel2021,Parmelee2022,Parmelee2022a}. 
Additionally, Cho \textit{et al}. \cite{Cho2017} have connected synchronized cluster states and {``chimera states''} \cite{panaggio2015}
(in which a subpopulation of oscillators synchronizes in an otherwise incoherent sea). Their research was formulated in a smooth setting, and it would be fascinating to explore it using a PWL perspective.

The extension of the methodology to treat various complexities --- including nonidentical oscillators, oscillators with high-dimensional (non-planar) dynamics, excitable systems, coupling delays, adaptive networks (in which a dynamical process on a network is coupled to the dynamics of the network's structure), temporal networks (in which a network’s entities and/or their interactions change with time), and multilayer networks (which can incorporate multiple types of interactions, multiple subsystems, and other complexities), and oscillator networks with polyadic (i.e., beyond pairwise) interactions --- is mathematically interesting and can build naturally on existing inroads on these challenges that have been made for smooth networks \cite{Porter2019}. Relevant studies to extend to a PWL framework include investigations of networks of Kuramoto oscillators with heterogenous frequencies \cite{Ott2008} and modular structures \cite{Skardal2019}, an extension of the MSF for coupled nearly-identical dynamical systems \cite{sun2009master} and dynamical systems with delays \cite{Kyrychko2014,Otto2018}, and extension of coupled-oscillator theory to networks with polyadic interactions (which are sometimes called ``higher-order'' interactions)
\cite{Bick2016,Grilli2017,Lambiotte2019,bick2022}.  
It is also worth extending the analysis of models of noisy PWL oscillators
to networks of such systems \cite{Simpson2011}.  
The further development of techniques for analyzing nonsmooth network dynamics is extremely relevant for systems with switches or thresholds, which arise in models of social-influence-driven opinion changes \cite{Singh2013} and contagions \cite{Lehmann2018}. 
There are numerous outstanding 
challenges in the study of dynamics on networks that may benefit from the perspective of nonsmooth modeling and analysis.


\section*{Acknowledgements}

We thank Thilo Gross for helpful discussions.


\newpage

\appendix

\section{Piecewise-linear models}
\label{Pwlmodels}
In \cref{tab:pwl2piecesmodels}, we summarize the components $A_{\mu}$ and $b_{\mu}$ (with $\mu \in \{1,2\}$) of the PWL models in \cref{sec:pwl} when they are written in the form \cref{pwl1}.

\begin{table} [H]
	\centering
	\begin{tabular}{ |p{3cm}|p{2.5cm}|p{2.5cm} |p{1.7cm} |p{1.7cm}| }
		\hline
		Model & $A_{1}$ & $A_{2}$ & $b_{1}$ & $b_{2}$\\
		\hline
		&&&&\\
		McKean model $\varsigma(w)=(\gamma a+w)/ I$, $w \in[-\gamma a,-\gamma a+I],$ on $v=a$. & $\left[\begin{array}{cc} -\gamma & -1 \\	b & 0	\end{array}\right] $ & $\left[\begin{array}{cc} -\gamma & -1 \\ b & 0 \end{array}\right] $    & $\left[\begin{array}{c}	I \\ 0 \end{array}\right]$ & $\left[\begin{array}{c} 0 \\ 0 \end{array}\right]$   \\ [1em]
		\hline
		&&&&\\
		The absolute model & $\left[\begin{array}{cc}  1 &-1  \\	1 & -d	\end{array}\right] $ & $\left[\begin{array}{cc} -1 & -1 \\  1 & -d  \end{array}\right] $    & $\left[\begin{array}{c}	-a  \\ d \overline{w}-\overline{v}  \end{array}\right]$ & $\left[\begin{array}{c} a  \\ d \overline{w}-\overline{v}  \end{array}\right]$  \\ [1em]
		\hline
		&&&&\\
		PWL homoclinic model & $\left[\begin{array}{cc} \tau_{1}&-1 \\ \delta_{1} & 0	\end{array}\right] $ & $\left[\begin{array}{cc} \tau_{2}&-1 \\ \delta_{2} & 0  \end{array}\right] $    & $\left[\begin{array}{c} 0 \\ -1 \end{array}\right]$ & $\left[\begin{array}{c} 0 \\ -1 \end{array}\right]$   \\
		\hline
		&&&&\\
		Planar integrate-and-fire (IF) model & $\left[\begin{array}{cc}  a_{1} & -1 \\ a_{w}/{\tau}  & b_{w}/{\tau} \end{array}\right] $ & $\left[\begin{array}{cc} a_{2} & -1 \\ a_{w}/{\tau}  & b_{w}/{\tau}  \end{array}\right] $    & $\left[\begin{array}{c}	I \\ 0  \end{array}\right]$ & $\left[\begin{array}{c} I \\ 0 \end{array}\right]$   \\
		\hline
	\end{tabular}
	\caption{Components of the examined models in the form \cref{pwl1}. We complete the definition of the McKean model by using the Filippov convention.
	}
	\label{tab:pwl2piecesmodels}
\end{table}

\noindent
The dynamics of the PWL Morris--Lecar model with three zones is
\begin{equation} \label{ML}
	C\dot{v}=\rho(v)-w+I \,, \quad \dot{w}=g(v,w) \,,
\end{equation}
with a continuous $\rho(v)$ (to approximate a cubic $v$-nullcline) of the form
\begin{equation}\label{fMcKeannadMorrisLecar}
	\rho(v) = \left\{
	\begin{array}{lll}
		-v  & \mathrm{if} \ \ \  v<a/2  \\
		v-a  & \mathrm{if} \ \ \ a/2 \leq v\leq (1+a)/2 \\
		1-v  & \mathrm{if} \ \ \ v>(1+a)/2\,,
	\end{array} \right.
\end{equation}
and continuous $g$ function
\begin{equation}\label{fMcKeannadMorrisLecar2}
	g(v,w) = \left\{
	\begin{array}{ll}
		(v-\gamma_1 w+b^*\gamma_1-b)/\gamma_1 & \mathrm{if} \ \ \  v<b \\
		(v-\gamma_2 w+b^*\gamma_2-b)/\gamma_2  & \mathrm{if} \ \ \ v\geq b\,,
	\end{array} \right.
\end{equation}
with the constraints $-a/2<b^*<(1-a)/2$, $a/2<b<(1+a)/2$, $\gamma_2>0$, and $\gamma_1\in \mathbb{R}$.
To construct periodic solutions, such as the one in \cref{Fig:TheMorrisLecarperiodic}, we use the formalism in \cref{sec:pwl}.
We break the periodic orbit into pieces such that each piece is governed by a linear dynamical system.
This is similar to the system \cref{pwl1}, but now the orbit has four distinct pieces that evolve according to $\text{d}x/\text{d}t=A_{\mu}x+b_{\mu}$, with $\mu \in \{1,2,3,4\}$, in three linear regimes $R_1=\{x\in \mathbb R^{2} |\ v >(1+a)/2\}$, $R_2=\{x\in \mathbb R^{2} |\ b<v <(1+a)/2\}$, and $R_3=\{x\in \mathbb R^{2} |\ a/2<v < b\}$. Therefore,
\begin{equation}
	A_{1}= \begin{bmatrix}
		1/C & -1/C \\ 1/\gamma_2  & -1
	\end{bmatrix},\ A_{2}= \begin{bmatrix}
		-1/C & -1/C \\ 1/\gamma_2 & -1
	\end{bmatrix}, \ A_{4}= \begin{bmatrix}
		1/C & -1/C \\ 1/\gamma_1 & -1
	\end{bmatrix}
\end{equation}
and
\begin{equation}
	b_{1}= \begin{bmatrix}
		(I-a)/C  \\ b^*-b/\gamma_2
	\end{bmatrix}, \ b_{2}= \begin{bmatrix}
		(1+I)/C  \\ b^*-b/\gamma_2
	\end{bmatrix},  \ b_{4}= \begin{bmatrix}
		(I-a)/C  \\ b^*-b/\gamma_1
	\end{bmatrix}\,,	
\end{equation}
with $A_{3}=A_{1}$ and $b_{3}=b_{1}$. Let $T_{\mu}$ denote the time-of-flight for each piece, and let $T=\Sigma_{\mu=1}^{4}T_{\mu}$ denote the corresponding period of the orbit. To build a closed orbit, we use the boundary-crossing values of the voltage variable (i.e., $v=b$ and $v=(1+a)/2$) and equation \cref{pwl10}, and we enforce periodicity of the solution. Choosing initial data $x^1(0)=(b,w^1(0))^{\top}$ and enforcing continuity of solutions by using the matching conditions $x^{\mu+1}(0)=x^{\mu}(T_{\mu})$ for $\mu \in \{1,2,3\}$ determines $T_\mu$ and $w(0)$ through the simultaneous solution of the equations
$v^{1}(T_1)=(1+a)/2$, $v^{2}(T_2)=(1+a)/2$, $v^{3}(T_3)=b$, $v^{4}(T_4)=b$, and $w^1(0)=w^4(T_4)$.


\section{Saltation operator\label{Saltationproof}}
Using the notation of \cref{FlqTheory}, we denote a periodic orbit by $x^\gamma$, a perturbed orbit by $\widetilde{x}$, an event time by $t_0$, and a perturbed event time by $\widetilde{t}_0$. We obtain the last two from the equations
$h_\mu(x^\gamma(t_0)) = 0$
and $h_\mu(\widetilde{x}(\widetilde{t}_0)) = 0$, respectively.  The difference between the perturbed and unperturbed events times is $\delta t_0 = \widetilde{t}_0 - t_0$.
The periodic and perturbed states after the switching event are $x^{\gamma}(t_{0}^{+})=\mathcal{J}_{\mu}(x^{\gamma}(t_{0}^{-}))$ and $\tilde{x}(\tilde{t}^{+}_{0})=\mathcal{J}_{\mu}(\tilde{x}(\tilde{t}^{-}_{0}))$, where $\mathcal{J}_{\mu}$ is the switch rule.
Without loss of generality, we consider $\delta t_{0}>0$ so $x^{\gamma}(t)$ and $\tilde{x}(t)$ are on opposite sides of the switching manifold (because $x^{\gamma}(t)$ has already crossed the switching boundary).
We then have that $\tilde{x}(\tilde{t}^{-}_{0})=\tilde{x}(t^{-}_{0}+ \delta t_{0}) \approx x^{\gamma}(t^{-}_{0})+\delta x(t^{-}_{0})+\dot{x}^{\gamma}(t_{0}^{-}) \delta t_{0}$.

We do a first-order Taylor expansion of $\mathcal{J}_{\mu}$ and obtain
\begin{align}\label{Taylorexpofstate}
	\tilde{x}(\tilde{t}^{+}_{0})&=\mathcal{J}_{\mu}(\tilde{x}(\tilde{t}^{-}_{0})) \approx \mathcal{J}_{\mu}(x^{\gamma}(t^{-}_{0})+\delta x(t^{-}_{0})+\dot{x}^{\gamma}(t_{0}^{-}) \delta t_{0}) \nonumber\\
	&\approx \mathcal{J}_{\mu}(x^{\gamma}(t^{-}_{0})) + \mathrm{D}\mathcal{J}_{\mu}(x^{\gamma}(t^-_{0}))[ \delta x(t^{-}_{0})+\dot{x}^{\gamma}(t_{0}^{-}) \delta t_{0}]\nonumber\\
	&\approx x^{\gamma}(t^{+}_{0})+\mathrm{D}\mathcal{J}_{\mu}(x^{\gamma}(t^-_{0}))[ \delta x(t^{-}_{0})+\dot{x}^{\gamma}(t_{0}^{-}) \delta t_{0}] \,,
\end{align}
where $\mathrm{D}\mathcal{J}_{\mu}$ is the Jacobian matrix of $\mathcal{J}_{\mu}$. The first-order Taylor expansion of $h_{\mu}(\tilde{x}(\tilde{t}^{-}_{0}))$ is
\begin{align}\label{TaylorexpofH}
	h_{\mu}(\tilde{x}(\tilde{t}^{-}_{0})) &=h_{\mu}(\tilde{x}(t^{-}_{0}+\delta t_{0}))=h_{\mu}(x^{\gamma}(t^{-}_{0}+\delta t_{0})+\delta x(t^{-}_{0}+\delta t_{0})) \nonumber\\
	& \approx h_{\mu}(x^{\gamma}(t^-_{0})+\dot{x}^{\gamma}(t_{0}^{-}) \delta t_{0})+\nabla_{x} h_{\mu}(x^{\gamma}(t^{-}_{0}+\delta t_{0})) \cdot \delta x(t^{-}_{0}+\delta t_{0})\nonumber \\
	& \approx h_{\mu}(x^{\gamma}(t^-_{0}))+\nabla_{x} h_{\mu}(x^{\gamma}(t^-_{0})) \cdot \dot{x}^{\gamma}(t_{0}^{-}) \delta t_{0}+ \nabla_{x} h_{\mu}(x^{\gamma}(t^-_{0})) \cdot \delta x(t^-_{0}) \,.
\end{align}
Using \cref{TaylorexpofH} and the fact that $h_{\mu}(x^{\gamma}(t_{0}))=0=h_{\mu}(\tilde{x}(\tilde{t}_{0}))$, we obtain
\begin{equation}\label{deltatdiff}
	\delta t_{0}=-\frac{\nabla_{x} h_{\mu}(x^{\gamma}(t^-_{0})) \cdot\delta x(t^-_{0})}{\nabla_{x} h_{\mu}(x^{\gamma}(t^-_{0})) \cdot \dot{x}^{\gamma}(t_{0}^{-})}.
\end{equation}	
We approximate $\tilde{x}(t^{+}_{0})$ as
\begin{align}\label{Delatxplusappxsimate}
	\tilde{x}(t^{+}_{0}) \approx \tilde{x}(\tilde{t}^{+}_{0})- \dot{\tilde{x}}(\tilde{t}^{+}_{0})\delta t_{0} \approx \tilde{x}(\tilde{t}^{+}_{0})-\dot{x}^{\gamma}(t_{0}^{+}+\delta t_{0}) \delta t_{0}
	\approx \tilde{x}(\tilde{t}^{+}_{0})-\dot{x}^{\gamma}(t_{0}^{+}) \delta t_{0} \,.
\end{align}

Using (\ref{Taylorexpofstate}) and (\ref{Delatxplusappxsimate}) yields
\begin{align}\label{DelatxplusAppendix}
	\delta x(t^+_{0}) &=\tilde{x}(t^+_{0})-x^{\gamma}(t^+_{0}) \approx \tilde{x}(\tilde{t}^{+}_{0})-\dot{x}^{\gamma}(t_{0}^{+}) \delta t_{0}-x^{\gamma}(t^+_{0}) \nonumber \\
	& \approx  x^{\gamma}(t^{+}_{0})+\mathrm{D}\mathcal{J}_{\mu}(x^{\gamma}(t^-_{0}))[\delta x(t^{-}_{0})+\dot{x}^{\gamma}(t_{0}^{-}) \delta t_{0}]
	-[x^{\gamma}(t^+_{0})+\dot{x}^{\gamma}(t_{0}^{+}) \delta t_{0}] \nonumber \\
	&=\mathrm{D}\mathcal{J}_{\mu}(x^{\gamma}(t^-_{0}))\delta x(t^{-}_{0})+ [\mathrm{D}\mathcal{J}_{\mu}(x^{\gamma}(t^-_{0}))\dot{x}^{\gamma}(t_{0}^{-})- \dot{x}^{\gamma}(t_{0}^{+}) ]\delta t_{0^{}} \,.
\end{align}
Therefore, using (\ref{deltatdiff}) and (\ref{DelatxplusAppendix}), we write $\delta x(t^+)$ in the form \cref{saltate},	
where $S_\mu(t)$ is the saltation matrix
\begin{align}\label{Delatxplusform}
	S_\mu(t) = \mathrm{D}\mathcal{J}_{\mu}(x^{\gamma}(t^-)) +
	\frac{[\dot{x}^{\gamma}(t^{+})-\mathrm{D}\mathcal{J}_{\mu}(x^{\gamma}(t^-))\dot{x}^{\gamma}(t^{-})][\nabla_{x} h_{\mu}(x^{\gamma}(t^-))]^{\top}}{\nabla_{x} h_{\mu}(x^{\gamma}(t^-)) \cdot \dot{x}^{\gamma}(t^{-})}\,.
\end{align}


\section{The nontrivial Floquet exponent for planar PWL systems\label{PWLPlanarFloquet}}

For planar systems, one eigenvalue of the monodromy matrix $\Psi$ is
$1$. Let $\mathrm{e}^{\kappa T}$ denote the nontrivial Floquet multiplier. Using the relation $\mathrm{det} \, \Psi= \mathrm{e}^{\kappa T} \times 1$ and equation \cref{MonodromyExplicitinchapter}, we obtain
\begin{align}
		\mathrm{e}^{\kappa T} &=\operatorname{det} \left[ S(t_M)G(T_{M}) S(t_{M-1})G(T_{M-1}) \times \cdots \times S(t_2)G(T_2)S(t_1)G(T_1) \right] \nonumber \\
		&=\operatorname{det} S(t_{M}) \times \cdots \times \operatorname{det}S(t_{1})  \operatorname{det} G(T_{M}) \times \cdots \times \operatorname{det}G(T_{1}) \nonumber \\
		&=\operatorname{det} S(t_{M}) \times \cdots \times \operatorname{det}S(t_{1})  \operatorname{det} \mathrm{e}^{A_{M} T_{M}} \times \cdots \times \operatorname{det} \mathrm{e}^{A_{1}T_{1}} \,.
\end{align}
Using the well-known fact $\operatorname{det}\mathrm{e}^{A}=\mathrm{e}^{\mathrm{Tr} A}$,
 we obtain the useful formula
\begin{equation}\label{AppxFexponentformula}
	\kappa=\frac{1}{T} \sum^{M}_{i =1}\left[T_{i} \operatorname{Tr} A_{\mu(i)}+\ln | \det S(t_i) | \right] \,.
\end{equation}

\section{Derivation of the jump condition in $\mathcal{B}(t)$\label{JumpinB}}

Using equation \cref{Sec:bacgroundDelatx}, we consider a perturbed solution of the form $x(t)=x^{\gamma}+\psi p(t)$, where $\psi=\mathcal{O}(\sigma)$ has a small value, that
 crosses the switching manifolds at the perturbed switching times $\tilde{t}^{}_{i}=t_{i}+g_{i}(\psi)$, which we obtain
 by solving
 $h_{\mu}(x(t_{i}+g_{i}(\psi)))=0$. In general, $g_{i}(\psi)$ depends on the geometry of a switching manifold
 and on the displacement $\psi p(t)$. For the PWL models that we consider, one can calculate $g_{i}(\psi)$ explicitly.
   To first order, the Taylor expansion of $h_{\mu}(\tilde{x}(\tilde{t}^{}_{i}))$ is
\begin{equation}
\label{TaylorexpofHB}
		h_{\mu}(\tilde{x}(\tilde{t}^{}_{i}))  \approx
		  h_{\mu}(x^{\gamma}(t^-_{i}))+\nabla_{x} h_{\mu}(x^{\gamma}(t^-_{i})) \cdot \dot{x}^{\gamma}(t_{i}^{-}) g_{i}(\psi)
		+ \nabla_{x} h_{\mu}(x^{\gamma}(t^-_{i})) \cdot \psi p(t^{-}_{i}) \,.
\end{equation}
Because $h_{\mu}(x^{\gamma}(t_{i}))=0=h_{\mu}(\tilde{x}(\tilde{t}_{i}))$, we obtain
\begin{equation}
\label{deltatdiffB}
	g_{i}(\psi)=-\frac{\nabla_{x} h_{\mu}(x^{\gamma}(t^-_{i})) \cdot\psi p(t^{-}_{i})}{\nabla_{x} h_{\mu}(x^{\gamma}(t^-_{i})) \cdot \dot{x}^{\gamma}(t_{i}^{-})}
	=- \frac{\psi p^{v}(t^{-}_{i})}{\dot{v}^{\gamma}(t_{i}^{-})} \,,
\end{equation}	
where $p^v(t)$ is the $v$ component of $p(t)$. \Cref{Sec:bacgroundIsostableSOrderTaylorphase} implies that
\begin{equation}
\nabla_{(x^{\gamma}(\tilde{t}^-_{i})+\psi p(\tilde{t}^-_{i}))} \Theta(x)
	\approx \mathcal{Z}(t^-_{i}+g_{i}(\psi))+\psi \mathcal{B}(t^-_{i}+g_{i}(\psi))
\end{equation}
immediately before the switching event.
We obtain a similar equation
by evaluating equation \cref{Sec:bacgroundIsostableSOrderTaylorphase} at $\tilde{t}^+_{i}=t^+_{i}+g_{i}(\psi)$.

We follow the technique that was proposed by Wilson \cite{wilson2019isostable} to derive a jump condition in the iIRC,
$\mathcal{B}(t)$, and
$\mathcal{C}(t)$ for an $m$-dimensional piecewise-smooth systems with an $(m-1)$-dimensional switching
manifold $\Sigma_{\mu}$ that is transverse to $x^{\gamma}(t)$.
 We make four assumptions. First, for all $k$, the phase function $\psi_{k}(x)$ is continuous
 in an open neighborhood of $x^{\gamma}(t)$. Second, for all $k$, the function $\psi_{k}(x)$
 is at least twice differentiable in the interior of each region $R_{\mu}$. Third, each boundary $\Sigma_{\mu}$ is at least $C^1$ (i.e., continuously differentiable) in an open ball $B(p_{\mu},R)$ that is centered at $p_{\mu}$ (the intersection point of $\Sigma_{\mu}$ and $x^{\gamma}(t)$) with radius $R$.
 It then follows that at each intersection point $p_{\mu}$ that there exists a tangent hyperplane $\Pi$ rhat is spanned by an orthonormal set of $(m-1)$-dimensionals vectors $w_k$ for $k \in \{ 1, \ldots , m-1\}$.
 Fourth, for all $k$, the directional derivatives of $\psi_{k}$
 exist on $\Pi$ in all tangential directions $w_k$ and are identical from both sides. (Otherwise, the associated coordinate of $\psi_{k}$ is not continuous \cite{wilson2019isostable}.)
 For the planar PWL models that we consider, $h_{\mu}(x)=v-a_{\mu}$, where $a_{\mu}$ is a constant. Therefore, $w_{1}=[0,1]^{\top}$.

Using the continuity of $\Theta(x)$ and the fourth assumption
about the phase function $\psi_{k}(x)$, 
approaching from either side of the switching manifold yields
\begin{equation}
	\left (\nabla_{(x^{\gamma}(\tilde{t}^-_{i})+\psi p(\tilde{t}^-_{i}))} \Theta(x) \right ) \cdot w_1 =
	\left (\nabla_{(x^{\gamma}(\tilde{t}^+_{i})+\psi p(\tilde{t}^+_{i}))} \Theta(x) \right ) \cdot w_1 \map{\,,}
\end{equation}
where we drop the subscript on $\psi$ for convenience.
Equivalently,
\begin{equation} \label{BSOrderTaylorphasePlanar2}
	\left [\mathcal{Z}(t^-_{i}+g_{i}(\psi))+\psi \mathcal{B}(t^-_{i}+g_{i}(\psi)) \right ] \cdot w_1 =
	\left [\mathcal{Z}(t^+_{i}+g_{i}(\psi))+\psi \mathcal{B}(t^+_{i}+g_{i}(\psi)) \right ] \cdot w_1 \,.
\end{equation}
We Taylor expand equation \cref{BSOrderTaylorphasePlanar2} in $\psi_{}$ to obtain
\begin{align}\label{BSOrderTaylorphasePlanar3}
	&\left[\mathcal{Z}(t^-_{i})+\left(\left.\frac{\d\mathcal{Z}}{\d t}\right|_{t=t^-_{i}}\right) g_{i}(\psi)+\psi \mathcal{B}\left(t^-_{i}\right)\right] \cdot w_{1} \nonumber \\
	&\qquad=\left[\mathcal{Z}(t^+_{i})+\left(\left.\frac{\d\mathcal{Z}}{\d t}\right|_{t=t^+_{i}}\right) g_{i}(\psi)+\psi \mathcal{B}\left(t^+_{i}\right)\right] \cdot w_{1} +\mathcal{O}\left(\psi^{2}\right) \,.
\end{align}
We set the $\mathcal{O}(\psi^{0})$ terms equal the two sides of equation \cref{BSOrderTaylorphasePlanar3} and use the normalization conditions to obtain
the jump operator for $\mathcal{Z}$ at $t_{i}$. This jump operator is the same one that we obtained
 in \cref{sec:response}. Collecting the $\mathcal{O}(\psi^{})$ terms in \cref{BSOrderTaylorphasePlanar3} yields
\begin{align}\label{BSOrderTaylorphasePlanar67}
	\left[\left(\left.\frac{\d\mathcal{Z}}{\d t}\right|_{t=t^-_{i}}\right) g_{i}(\psi)+\psi \mathcal{B}^-\right] \cdot w_{1}
	=\left[\left(\left.\frac{\d\mathcal{Z}}{\d t}\right|_{t=t^+_{i}}\right) g_{i}(\psi)+\psi \mathcal{B}^+\right] \cdot w_{1} \,.
\end{align}

We use equations \cref{Zadjoint} and (\ref{deltatdiffB}) to rewrite \cref{BSOrderTaylorphasePlanar67} to obtain
\begin{align}\label{BSOrderTaylorphasePlanar5}
	\psi \left[\frac{p^{v}(t^{-}_{i})}{\dot{v}^{\gamma}(t_{i}^{-})}A_{\mu(i)}^{\top}\mathcal{Z}^-+ \mathcal{B}^-\right] \cdot w_{1}
	=\psi \left[\frac{p^{v}(t^{-}_{i})}{\dot{v}^{\gamma}(t_{i}^{-})}A_{\mu(i+1)}^{\top}\mathcal{Z}^++\mathcal{B}^+\right] \cdot w_{1} \,.
\end{align}	
The condition \cref{Sec:backBModifiedNormal} holds on both sides of a switching manifold, so
\begin{align}\label{BSOrderTaylorphasePlanar6}
	\mathcal{Z}^-\cdot (A_{\mu(i)}p(t^-_{i}))+f^-_{\mu(i)}\cdot\mathcal{B}^-=0=\mathcal{Z} \cdot (A_{\mu(i+1)}p(t^+_{i}))+f^+_{\mu(i)}\cdot\mathcal{B}^+ \,,
\end{align}	
where $f^-_{\mu(i)}$ is the vector field evaluated on the limit cycle immediately before a switching event and $f^+_{\mu(i)}$ is the vector field evaluated on the limit cycle immediately after it.
Combining (\ref{BSOrderTaylorphasePlanar6}) and (\ref{BSOrderTaylorphasePlanar5}) yields
\begin{align}\label{BSOrderTaylorphasePlanar7}
	\mathcal{B}^\cdot f^+_{\mu(i)} &= \mathcal{B}^- \cdot f^-_{\mu(i)} +\mathcal{Z}^- \cdot (A_{\mu(i)}p(t^-_{i}))-\mathcal{Z}^+ \cdot (A_{\mu(i+1)}p(t^+_{i}))  \\
	\mathcal{B}^+\cdot w_{1} &= \mathcal{B}^-\cdot w_{1}+\frac{p^{v}(t^{-}_{i})}{\dot{v}^{\gamma}(t_{i}^{-})}\left[A_{\mu(i)}^{\top}\mathcal{Z}^--A_{\mu(i+1)}^{\top}\mathcal{Z}^+\right]\cdot w_{1} \,.
\end{align}	
Therefore, the jump condition on $\mathcal{B}$ during a transition across a switching manifold is
\begin{align}\label{APPBSOrderTaylorphasePlanar8}
	\mathcal{B}^+=(S^{\top}(t_i))^{-1}\mathcal{B}^-+C^{-1}(t_i)\eta(t_i) \,,
\end{align}
where $C(t_i)$ and  $\eta(t_i)$ are given by \cref{sigmavectorB}.


\section{Interaction functions}
\label{sec:PhasAmpInteraction}
The interaction functions in the dynamical system \cref{phaseisostablenetworkPWLcoupling4}
are
\begin{equation}
	\begin{array}{lll}
		h_{1}\left(\omega \theta_{1}, \omega \theta_{2}\right)= \mathcal{Z}^{v}\left(\theta_{1}\right) \left(v^{\gamma}\left(\theta_{2}\right)-v^{\gamma}\left(\theta_{1}\right)\right) \,, \\
		h_{2}\left(\omega\theta_{1}, \omega\theta_{2}\right)=\mathcal{B}^{v}\left(\theta_{1}\right) \left(v^{\gamma}\left(\theta_{2}\right)-v^{\gamma}\left(\theta_{1}\right)\right)- \mathcal{Z}^{v}\left(\theta_{1}\right)p^v(\theta_{1}) \,,\\
		h_{3}\left(\omega\theta_{1}, \omega\theta_{2}\right)=\mathcal{Z}^{v}\left(\theta_{1}\right)p^v(\theta_{2}) \,, \\
		h_{4}\left(\omega\theta_{1}, \omega\theta_{2}\right)=\mathcal{I}^{v}\left(\theta_{1}\right) \left(v^{\gamma}\left(\theta_{2}\right)-v^{\gamma}\left(\theta_{1}\right)\right) \,,\\
		h_{5}\left(\omega\theta_{1}, \omega\theta_{2}\right)=\mathcal{C}^{v}\left(\theta_{1}\right) \left(v^{\gamma}\left(\theta_{2}\right)-v^{\gamma}\left(\theta_{1}\right)\right)- \mathcal{I}^{v}\left(\theta_{1}\right)p^v(\theta_{1}),\\
		h_{6}\left(\omega\theta_{1}, \omega\theta_{2}\right)=\mathcal{I}^{v}\left(\theta_{1}\right)p^v(\theta_{2}) \,,
	\end{array}
\end{equation}
where $\mathcal{Z}^{v}$, $\mathcal{I}^{v}$, $\mathcal{B}^{v}$, and $\mathcal{C}^{v}$ are the $v$ components of the corresponding vectors.


\tableofcontents


\newpage



\begin{thebibliography}{100}

\bibitem{Acary2011}
{\sc V.~Acary, O.~Bonnefon, and B.~Brogliato}, {\em Nonsmooth Modeling and
  Simulation for Switched Circuits}, vol.~69 of Lecture Notes in Electrical
  Engineering, Springer-Verlag, Heidelberg, Germany, 2011.

\bibitem{aihara2010theory}
{\sc K.~Aihara and H.~Suzuki}, {\em Theory of hybrid dynamical systems and its
  applications to biological and medical systems}, Philosophical Transactions
  of the Royal Society A, 368 (2010), pp.~4893--4914.

\bibitem{andronov2013theory}
{\sc A.~A. Andronov, A.~A. Vitt, and S.~E. Khaikin}, {\em Theory of
  Oscillators: Adiwes International Series in Physics}, vol.~4, Elsevier,
  Amsterdam, The Netherlands, 2013.

\bibitem{ArenasGuileraMoreno:2008}
{\sc A.~Arenas, A.~D\'{\i}az-Guilera, J.~Kurths, Y.~Moreno, and C.~Zhou}, {\em
  Synchronization in complex networks}, Physics Reports, 469 (2008),
  pp.~93--153.

\bibitem{Arenas2006}
{\sc A.~Arenas, A.~D\'iaz-Guilera, and C.~J. P\'erez-Vicente}, {\em
  {Synchronization processes in complex networks}}, Physica D, 224 (2006),
  pp.~27--34.

\bibitem{Ashwin2019}
{\sc P.~Ashwin, C.~Bick, and C.~Poignard}, {\em State-dependent effective
  interactions in oscillator networks through coupling functions with dead
  zones}, Philosophical Transactions of the Royal Society A, 377 (2019),
  20190042.

\bibitem{Ashwin2021}
{\sc P.~Ashwin, C.~Bick, and C.~Poignard}, {\em Dead zones and phase reduction
  of coupled oscillators}, Chaos: An Interdisciplinary Journal of Nonlinear
  Science, 31 (2021), p.~093132.

\bibitem{Ashwin2016}
{\sc P.~Ashwin, S.~Coombes, and R.~Nicks}, {\em Mathematical frameworks for
  oscillatory network dynamics in neuroscience}, The Journal of Mathematical
  Neuroscience, 6 (2016), 2.

\bibitem{Ashwin07}
{\sc P.~Ashwin, G.~Orosz, J.~Wordsworth, and S.~Townley}, {\em Dynamics on
  networks of cluster states for globally coupled phase oscillators}, SIAM
  Journal on Applied Dynamical Systems, 6 (2007), pp.~728--758.

\bibitem{ashwin1992dynamics}
{\sc P.~Ashwin and J.~W. Swift}, {\em The dynamics of $n$ weakly coupled
  identical oscillators}, Journal of Nonlinear Science, 2 (1992), pp.~69--108.

\bibitem{Bassett2006}
{\sc D.~S. Bassett and E.~Bullmore}, {\em Small-world brain networks}, The
  Neuroscientist, 12 (2006), pp.~512--523.

\bibitem{Bassett2009}
{\sc D.~S. Bassett and E.~T. Bullmore}, {\em Human brain networks in health and
  disease}, Current Opinion in Neurobiology, 22 (2009), pp.~340--347.

\bibitem{bassett2017}
{\sc D.~S. Bassett and O.~Sporns}, {\em Network neuroscience}, Nature
  Neuroscience, 20 (2017), pp.~353--364.

\bibitem{Bel2021}
{\sc A.~Bel, R.~Cobiaga, W.~Reartes, and H.~G. Rotstein}, {\em Periodic
  solutions in threshold-linear networks and their entrainment}, SIAM Journal
  on Applied Dynamical Systems, 20 (2021), pp.~1177--1208.

\bibitem{Belykh2000}
{\sc V.~N. Belykh, I.~V. Belykh, and M.~Hasler}, {\em Hierarchy and stability
  of partially synchronous oscillations of diffusively coupled dynamical
  systems}, Physical Review E, 62 (2000), pp.~6332--6345.

\bibitem{Betzel2016}
{\sc R.~F. Betzel, S.~Gu, J.~D. Medaglia, F.~Pasqualetti, and D.~S. Bassett},
  {\em Optimally controlling the human connectome: the role of network
  topology.}, Scientific Reports, 6 (2016), p.~30770.

\bibitem{Betzel2017}
{\sc R.~F. Betzel, J.~D. Medaglia, L.~Papadopoulos, G.~L. Baum, R.~Gur, R.~Gur,
  D.~Roalf, T.~D. Satterthwaite, and D.~S. Bassett}, {\em {The modular
  organization of human anatomical brain networks: Accounting for the cost of
  wiring}}, Network Neuroscience, 1 (2017), pp.~42--68.

\bibitem{Bick2016}
{\sc C.~Bick, P.~Ashwin, and A.~Rodrigues}, {\em Chaos in generically coupled
  phase oscillator networks with nonpairwise interactions}, Chaos: An
  Interdisciplinary Journal of Nonlinear Science, 26 (2016), p.~094814.

\bibitem{bick2022}
{\sc C.~Bick, E.~Gross, H.~A. Harrington, and M.~T. Schaub}, {\em What are
  higher-order networks?}, SIAM Review (in press; arXiv:2104.11329),  (2022).

\bibitem{brogliato1999nonsmooth}
{\sc B.~Brogliato}, {\em Nonsmooth Mechanics}, Springer-Verlag, Heidelberg,
  Germany, 1999.

\bibitem{brogliato2000impacts}
{\sc B.~Brogliato}, {\em Impacts in Mechanical Systems: Analysis and
  Modelling}, Springer-Verlag, Heidelberg, Germany, 2000.

\bibitem{Brown03}
{\sc E.~Brown, P.~Holmes, and J.~Moehlis}, {\em Globally coupled oscillator
  networks}, in {Perspectives and Problems in Nonlinear Science: A Celebratory
  Volume in Honor of Larry Sirovich}, E.~Kaplan, J.~Marsden, and
  K.~Sreenivasan, eds., Springer, New York, 2003, pp.~183--215.

\bibitem{Cabre2005}
{\sc X.~Cabr\'e, E.~Fontich, and R.~{De La Llave}}, {\em The parameterization
  method for invariant manifolds {III}: {O}verview and applications}, Journal
  of Differential Equations, 218 (2005), pp.~444--515.

\bibitem{Champneys2008}
{\sc A.~R. Champneys and M.~{di Bernardo}}, {\em Piecewise smooth dynamical
  systems}, Scholarpedia, 3 (2008), 4041.

\bibitem{chartrand2019synchronization}
{\sc T.~Chartrand, M.~S. Goldman, and T.~J. Lewis}, {\em Synchronization of
  electrically coupled resonate-and-fire neurons}, SIAM Journal on Applied
  Dynamical Systems, 18 (2019), pp.~1643--1693.

\bibitem{Cho2017}
{\sc Y.~S. Cho, T.~Nishikawa, and A.~E. Motter}, {\em Stable chimeras and
  independently synchronizable clusters}, Physical Review Letters, 119 (2017),
  084101.

\bibitem{colombo2012bifurcations}
{\sc A.~Colombo, M.~{di Bernardo}, S.~J. Hogan, and M.~R. Jeffrey}, {\em
  Bifurcations of piecewise smooth flows: Perspectives, methodologies and open
  problems}, Physica D, 241 (2012), pp.~1845--1860.

\bibitem{Coombes2001}
{\sc S.~Coombes}, {\em Phase locking in networks of synaptically coupled
  {McKean} relaxation oscillators}, Physica D, 160 (2001), pp.~173--188.

\bibitem{Coombes2008}
{\sc S.~Coombes}, {\em Neuronal networks with gap junctions: A study of
  piece-wise linear planar neuron models}, SIAM Journal on Applied Dynamical
  Systems, 7 (2008), pp.~1101--1129.

\bibitem{Coombes20181}
{\sc S.~Coombes, Y.~M. Lai, M.~{\c{S}}ayli, and R.~Thul}, {\em Networks of
  piecewise linear neural mass models}, European Journal of Applied
  Mathematics, 29 (2018), pp.~869--890.

\bibitem{CoombesThul:2016}
{\sc S.~Coombes and R.~Thul}, {\em Synchrony in networks of coupled non-smooth
  dynamical systems: Extending the master stability function}, European Journal
  of Applied Mathematics, 27 (2016), pp.~904--922.

\bibitem{Coombes2012}
{\sc S.~Coombes, R.~Thul, and K.~C.~A. Wedgwood}, {\em Nonsmooth dynamics in
  spiking neuron models}, Physica D, 241 (2012), pp.~2042--2057.

\bibitem{Coombes2023}
{\sc S.~Coombes and K.~C.~A. Wedgwood}, {\em Neurodynamics: An Applied
  Mathematics Perspective}, vol.~75 of Texts in Applied Mathematics, Springer,
  2023.

\bibitem{Corragio2021}
{\sc M.~Coraggio, P.~{De Lellis}, and M.~{di Bernardo}}, {\em {Convergence and
  synchronization in networks of piecewise-smooth systems via distributed
  discontinuous coupling}}, Automatica, 129 (2021), 109596.

\bibitem{Curto2013}
{\sc C.~Curto, A.~Degeratu, and V.~Itskov}, {\em Encoding binary neural codes
  in networks of threshold-linear neurons}, Neural Computation, 25 (2013),
  pp.~2858--2903.

\bibitem{Curto2016}
{\sc C.~Curto and K.~Morrison}, {\em Pattern completion in symmetric
  threshold-linear networks}, Neural Computation, 28 (2016), pp.~2825--2852.

\bibitem{dabrowski2012largest}
{\sc A.~Dabrowski}, {\em The largest transversal {L}yapunov exponent and master
  stability function from the perturbation vector and its derivative dot
  product ({TLEVDP})}, Nonlinear Dynamics, 69 (2012), pp.~1225--1235.

\bibitem{dahms2012cluster}
{\sc T.~Dahms, J.~Lehnert, and E.~Sch{\"o}ll}, {\em Cluster and group
  synchronization in delay-coupled networks}, Physical Review E, 86 (2012),
  016202.

\bibitem{deneve2016efficient}
{\sc S.~Den{\`e}ve and C.~K. Machens}, {\em Efficient codes and balanced
  networks}, Nature Neuroscience, 19 (2016), pp.~375--382.

\bibitem{bernardo2008piecewise}
{\sc M.~{di Bernardo}, C.~Budd, A.~R. Champneys, and P.~Kowalczyk}, {\em
  Piecewise-Smooth Dynamical Systems: Theory and Applications},
  Springer-Verlag, Heidelberg, Germany, 2008.

\bibitem{di2008bifurcations}
{\sc M.~{di Bernardo}, C.~J. Budd, A.~R. Champneys, P.~Kowalczyk, A.~B.
  Nordmark, G.~O. Tost, and P.~T. Piiroinen}, {\em Bifurcations in nonsmooth
  dynamical systems}, SIAM Review, 50 (2008), pp.~629--701.

\bibitem{di1999local}
{\sc M.~{di Bernardo}, M.~I. Feigin, S.~J. Hogan, and M.~E. Homer}, {\em Local
  analysis of {C}-bifurcations in n-dimensional piecewise-smooth dynamical
  systems}, Chaos, Solitons and Fractals, 11 (1999), pp.~1881--1908.

\bibitem{dumortier2006qualitative}
{\sc F.~Dumortier, J.~Llibre, and J.~C. Art{\'e}s}, {\em Qualitative Theory of
  Planar Differential Systems}, Springer-Verlag, Heidelberg, Germany, 2006.

\bibitem{ermentrout2019recent}
{\sc B.~Ermentrout, Y.~Park, and D.~Wilson}, {\em Recent advances in coupled
  oscillator theory}, Philosophical Transactions of the Royal Society A, 377
  (2019), 20190092.

\bibitem{Ermentrout1992}
{\sc G.~B. Ermentrout}, {\em Stable periodic solutions to discrete and
  continuum arrays of weakly coupled nonlinear oscillators}, SIAM Journal on
  Applied Mathematics, 52 (1992), pp.~1665--1687.

\bibitem{ermentrout1991multiple}
{\sc G.~B. Ermentrout and N.~Kopell}, {\em Multiple pulse interactions and
  averaging in systems of coupled neural oscillators}, Journal of Mathematical
  Biology, 29 (1991), pp.~195--217.

\bibitem{ermentrout2010mathematical}
{\sc G.~B. Ermentrout and D.~H. Terman}, {\em Mathematical Foundations of
  Neuroscience}, Springer-Verlag, Heidelberg, Germany, 2010.

\bibitem{Essen2013}
{\sc D.~C.~V. Essen, S.~M. Smith, D.~M. Barch, T.~E. Behrens, E.~Yacoub, and
  K.~Ugurbil}, {\em {The WU--Minn human connectome project: {A}n overview}},
  NeuroImage, 80 (2013), pp.~62--79.

\bibitem{feigin1994forced}
{\sc M.~I. Feigin}, {\em Forced oscillations in systems with discontinuous
  nonlinearities}, Nauka, (in Russian) (1994).

\bibitem{filippov2013differential}
{\sc A.~F. Filippov}, {\em Differential Equations with Discontinuous Righthand
  Sides: Control Systems}, Springer-Verlag, Heidelberg, Germany, 1988.

\bibitem{Forrester2019}
{\sc M.~Forrester, S.~Coombes, J.~J. Crofts, S.~N. Sotiropoulos, and R.~D.
  O'Dea}, {\em The role of node dynamics in shaping emergent functional
  connectivity patterns in the brain, submitted}, Network Neuroscience, 4
  (2020), pp.~467--483.

\bibitem{Labaree:1962}
{\sc B.~Franklin}, {\em The Collected Papers of Benjamin Franklin}, New Haven:
  Yale University Press, 1962.

\bibitem{fredriksson2000normal}
{\sc M.~H. Fredriksson and A.~B. Nordmark}, {\em On normal form calculations in
  impact oscillators}, Proceedings of the Royal Society of London. Series A,
  456 (2000), pp.~315--329.

\bibitem{freire1998bifurcation}
{\sc E.~Freire, E.~Ponce, F.~Rodrigo, and F.~Torres}, {\em Bifurcation sets of
  continuous piecewise linear systems with two zones}, International Journal of
  Bifurcation and Chaos, 8 (1998), pp.~2073--2097.

\bibitem{freire2012canonical}
{\sc E.~Freire, E.~Ponce, and F.~Torres}, {\em Canonical discontinuous planar
  piecewise linear systems}, SIAM Journal on Applied Dynamical Systems, 11
  (2012), pp.~181--211.

\bibitem{Glasser2013}
{\sc M.~F. Glasser, S.~N. Sotiropoulos, J.~A. Wilson, T.~S. Coalson, B.~Fischl,
  J.~L. Andersson, X.~Junqian, S.~Jbabdi, M.~Webster, J.~R. Polimeni, D.~C.~V.
  Essen, and M.~Jenkinson}, {\em The minimal preprocessing pipelines for the
  human connectome project}, NeuroImage, 80 (2013), pp.~105--124.

\bibitem{Glendinning2015}
{\sc P.~Glendinning}, {\em {II.18} {Hybrid Systems}}, in The Princeton
  Companion to Applied Mathematics, N.~J. Higham, ed., Princeton University
  Press, Princeton, NJ, USA, 2015, pp.~103--104.

\bibitem{golubitsky2023}
{\sc M.~Golubitsky and I.~Stewart}, {\em Dynamics and Bifurcation in Networks:
  Theory and Applications of Coupled Differential Equations}, Society for
  Industrial and Applied Mathematics, Philadelphia, PA, USA, 2023.

\bibitem{Golubitsky1988}
{\sc M.~Golubitsky, I.~Stewart, and D.~G. Schaeffer}, {\em Singularities and
  Groups in Bifurcation Theory, Volume II}, Springer-Verlag, Heidelberg,
  Germany, 1988.

\bibitem{Golubitsky1986}
{\sc M.~Golubitsky and I.~N. Stewart}, {\em Hopf bifurcation with dihedral
  group symmetry: coupled nonlinear oscillators}, in Multiparameter Bifurcation
  Theory, M.~Golubitsky and J.~Guckenheimer, eds., vol.~56 of Contemporary
  Mathematics, American Mathematical Society, Providence, RI, USA, 1986,
  pp.~131--173.

\bibitem{Golubitsky2016}
{\sc M.~Golubitsky and I.~N. Stewart}, {\em {Rigid patterns of synchrony for
  equilibria and periodic cycles in network dynamics}}, Chaos: An
  Interdisciplinary Journal of Nonlinear Science, 26 (2016), 094803.

\bibitem{Gonzales2019}
{\sc C.~J.~S. Gonzales}, {\em Analyzing the sensitivity of nonlinear
  oscillators toparametric perturbations using isostable and isochron
  coordinates}, master's thesis, Mechanical Engineering, University of
  California, Santa Barbara, 2019.

\bibitem{Grilli2017}
{\sc J.~Grilli, G.~Barab{\'a}s, M.~J. Michalska-Smith, and S.~Allesina}, {\em
  Higher-order interactions stabilize dynamics in competitive network models},
  Nature, 548 (2017), pp.~210--213.

\bibitem{Guckenheimer1975}
{\sc J.~Guckenheimer}, {\em Isochrons and phaseless sets}, Journal of
  Mathematical Biology, 1 (1975), pp.~259--273.

\bibitem{guckenheimer2013nonlinear}
{\sc J.~Guckenheimer and P.~Holmes}, {\em Nonlinear Oscillations, Dynamical
  Systems, and Bifurcations of Vector Fields}, Springer-Verlag, Heidelberg,
  Germany, 1983.

\bibitem{Guillamon2009}
{\sc A.~Guillamon and G.~Huguet}, {\em A computational and geometric approach
  to phase resetting curves and surfaces}, SIAM Journal on Applied Dynamical
  Systems, 8 (2009), pp.~1005--1042.

\bibitem{Han1995}
{\sc S.~K. Han, C.~Kurrer, and Y.~Kuramoto}, {\em Dephasing and bursting in
  coupled neural oscillators}, Physical Review Letters, 75 (1995),
  pp.~3190--3193.

\bibitem{Hansel1993}
{\sc D.~Hansel, G.~Mato, and C.~Meunier}, {\em Clustering and slow switching in
  globally coupled phase oscillators}, Physical Review E, 48 (1993),
  pp.~3470--3477.

\bibitem{Harris2015}
{\sc J.~Harris and B.~Ermentrout}, {\em {Bifurcations in the Wilson--Cowan
  equations with nonsmooth firing rate}}, SIAM Journal on Applied Dynamical
  Systems, 14 (2015), pp.~43--72.

\bibitem{Havlin2012}
{\sc S.~Havlin, D.~Y. Kenett, E.~Ben-Jacob, A.~Bunde, R.~Cohen, H.~Hermann,
  J.~W. Kantelhardt, J.~Kert\'esz, S.~Kirkpatrick, J.~Kurths, J.~Portugali, and
  S.~Solomon}, {\em Challenges in network science: {A}pplications to
  infrastructures, climate, social systems and economics}, The European
  Physical Journal Special Topics, 214 (2012), pp.~273--293.

\bibitem{Hlinka2012}
{\sc J.~Hlinka and S.~Coombes}, {\em Using computational models to relate
  structural and functional brain connectivity}, European Journal of
  Neuroscience, 36 (2012), pp.~2137--2145.

\bibitem{Hoppensteadt97}
{\sc F.~C. Hoppensteadt and E.~M. Izhikevich}, {\em Weakly Connected Neural
  Networks}, Springer-Verlag, Heidelberg, Germany, 1997.

\bibitem{Izhikevich2000}
{\sc E.~M. Izhikevich}, {\em Phase equations for relaxation oscillators}, SIAM
  Journal on Applied Mathematics, 60 (2000), pp.~1789--1804.

\bibitem{FN-scholar}
{\sc E.~M. Izhikevich and R.~FitzHugh}, {\em Fitzhugh--{N}agumo model},
  Scholarpedia, 1 (2006), 1349.

\bibitem{jeffrey2011sliding}
{\sc M.~R. Jeffrey}, {\em Sliding bifurcations and non-determinism in systems
  with switching}, IFAC Proceedings Volumes, 44 (2011), pp.~13275--13280.

\bibitem{Jeffrey2016}
{\sc M.~R. Jeffrey}, {\em Hidden bifurcations and attractors in nonsmooth
  dynamical systems}, International Journal of Bifurcation and Chaos, 26
  (2016), 1650068.

\bibitem{jeffrey2018hidden}
{\sc M.~R. Jeffrey}, {\em Hidden Dynamics}, Springer, 2018.

\bibitem{jeffrey2011geometry}
{\sc M.~R. Jeffrey and S.~J. Hogan}, {\em The geometry of generic sliding
  bifurcations}, SIAM Review, 53 (2011), pp.~505--525.

\bibitem{jordan2007nonlinear}
{\sc D.~Jordan and P.~Smith}, {\em Nonlinear Ordinary Differential Equations:
  An Introduction for Scientists and Engineers}, Oxford University Press,
  Oxford, UK, 4th~ed., 2007.

\bibitem{klausmeier2008floquet}
{\sc C.~A. Klausmeier}, {\em Floquet theory: {A} useful tool for understanding
  nonequilibrium dynamics}, Theoretical Ecology, 1 (2008), pp.~153--161.

\bibitem{Kong2023}
{\sc N.~J. Kong, J.~J. Payne, J.~Zhu, and A.~M. Johnson}, {\em Saltation
  matrices: The essential tool for linearizing hybrid dynamical systems}, arXiv
  preprint arXiv:2306.06862,  (2023).

\bibitem{Kopell2002}
{\sc N.~Kopell and G.~B. Ermentrout}, {\em Mechanisms of phase-locking and
  frequency control in pairs of coupled neural oscillators}, in Handbook of
  Dynamical Systems, B.~Fiedler, ed., Elsevier, Amsterdam, The Netherlands,
  2002, pp.~3--54.

\bibitem{Kuramoto1984}
{\sc Y.~Kuramoto}, {\em Chemical Oscillations, Waves and Turbulence},
  Springer-Verlag, Heidelberg, Germany, 1984.

\bibitem{kuznetsov2013elements}
{\sc Y.~A. Kuznetsov}, {\em Elements of Applied Bifurcation Theory},
  Springer-Verlag, Heidelberg, Germany, third~ed., 2004.

\bibitem{kvalheim2019existence}
{\sc M.~D. Kvalheim and S.~Revzen}, {\em Existence and uniqueness of global
  {K}oopman eigenfunctions for stable fixed points and periodic orbits},
  Physica D, 425 (2021), 132959.

\bibitem{Kyrychko2014}
{\sc Y.~N. Kyrychko, K.~B. Blyuss, and E.~Sch{\"{o}}ll}, {\em {Synchronization
  of networks of oscillators with distributed delay coupling}}, Chaos: An
  Interdisciplinary Journal of Nonlinear Science, 24 (2014), 043117.

\bibitem{Ladenbauer2013}
{\sc J.~Ladenbauer, J.~Lehnert, H.~Rankoohi, T.~Dahms, E.~Sch\"oll, and
  K.~Obermayer}, {\em Adaptation controls synchrony and cluster states of
  coupled threshold-model neurons}, Physical Review E, 88 (2013), 042713.

\bibitem{lai2019master}
{\sc Y.~M. Lai, J.~Veasy, S.~Coombes, and R.~Thul}, {\em A master stability
  function approach to cardiac alternans}, Applied Network Science, 4 (2019),
  pp.~1--16.

\bibitem{Lambiotte2019}
{\sc R.~Lambiotte, M.~Rosvall, and I.~Scholtes}, {\em From networks to optimal
  higher-order models of complex systems}, Nature Physics, 15 (2019),
  pp.~313--320.

\bibitem{Landstrom2017}
{\sc A.~P. Landstrom, D.~Dobrev, and X.~H. Wehrens}, {\em Calcium signaling and
  cardiac arrhythmias}, Circulation Research, 120 (2017), pp.~1969--1993.

\bibitem{Lehmann2018}
{\sc S.~Lehmann and Y.-Y. Ahn}, eds., {\em Complex Spreading Phenomena in
  Social Systems: Influence and Contagion in Real-World Social Networks},
  Springer International Publishing, Cham, Switzerland, 2018.

\bibitem{leine2013dynamics}
{\sc R.~I. Leine and H.~Nijmeijer}, {\em Dynamics and Bifurcations of
  Non-Smooth Mechanical Systems}, vol.~18, Springer-Verlag, Heidelberg,
  Germany, 2013.

\bibitem{llibre2019limit}
{\sc J.~Llibre and X.~Zhang}, {\em Limit cycles created by piecewise linear
  centers}, Chaos: An Interdisciplinary Journal of Nonlinear Science, 29
  (2019), 053116.

\bibitem{lodato2007synchronization}
{\sc I.~Lodato, S.~Boccaletti, and V.~Latora}, {\em Synchronization properties
  of network motifs}, EPL (Europhysics Letters), 78 (2007), 28001.

\bibitem{Makarenkov2012}
{\sc O.~Makarenkov and J.~S.~W. Lamb}, {\em Dynamics and bifurcations of
  nonsmooth systems: {A} survey}, Physica D, 241 (2012), pp.~1826--1844.

\bibitem{Mauroy2012}
{\sc A.~Mauroy and I.~Mezi\'c}, {\em {On the use of {F}ourier averages to
  compute the global isochrons of (quasi)periodic dynamics}}, Chaos: An
  Interdisciplinary Journal of Nonlinear Science, 22 (2012), 033112.

\bibitem{Mauroy2018}
{\sc A.~Mauroy and I.~Mezi{\'{c}}}, {\em {Global computation of phase-amplitude
  reduction for limit-cycle dynamics}}, Chaos: An Interdisciplinary Journal of
  Nonlinear Science, 28 (2018), 073108.

\bibitem{Mauroy2013}
{\sc A.~Mauroy, I.~Mezic, and J.~Moehlis}, {\em {Isostables, isochrons, and
  Koopman spectrum for the action-angle representation of stable fixed point
  dynamics}}, Physica D, 261 (2013), pp.~19--30.

\bibitem{McKean70}
{\sc H.~P. McKean}, {\em Nagumo's equation}, Advances in Mathematics, 4 (1970),
  pp.~209--223.

\bibitem{mendel2001}
{\sc M.~T. Mendl and S.~Held}, {\em Living in groups: an evolutionary
  perspective}, in Social Behaviour in Farm Animals, L.~Keeling and H.~Gonyou,
  eds., CABI Publishing, Wallingford, UK, 2001, pp.~7--36.

\bibitem{monga2021augmented}
{\sc B.~Monga and J.~Moehlis}, {\em Augmented phase reduction for periodic
  orbits near a homoclinic bifurcation and for relaxation oscillators},
  Nonlinear Theory and Its Applications, IEICE, 12 (2021), pp.~103--116.

\bibitem{monga2019phase}
{\sc B.~Monga, D.~Wilson, T.~Matchen, and J.~Moehlis}, {\em Phase reduction and
  phase-based optimal control for biological systems: {A} tutorial}, Biological
  Cybernetics, 113 (2019), pp.~11--46.

\bibitem{morris1981voltage}
{\sc C.~Morris and H.~Lecar}, {\em Voltage oscillations in the barnacle giant
  muscle fiber}, Biophysical Journal, 35 (1981), pp.~193--213.

\bibitem{Muller:1995}
{\sc P.~C. M\"uller}, {\em Calculation of {L}yapunov exponents for dynamical
  systems with discontinuities}, Chaos, Solitons and Fractals, 5 (1995),
  pp.~1671--1681.

\bibitem{Newman2018}
{\sc M.~E.~J. Newman}, {\em Networks}, Oxford University Press, Oxford,
  second~ed., 2018.

\bibitem{Nicks2023}
{\sc R.~Nicks, R.~Allen, and S.~Coombes}, {\em Insights into oscillator network
  dynamics using a phase-isostable framework}, in preparation,  (2023).

\bibitem{Nicks2018}
{\sc R.~Nicks, L.~Chambon, and S.~Coombes}, {\em Clusters in nonsmooth
  oscillator networks}, Physical Review E, 97 (2018), 032213.

\bibitem{Oldham2019}
{\sc S.~Oldham and A.~Fornito}, {\em The development of brain network hubs},
  Developmental Cognitive Neuroscience, 36 (2019), 100607.

\bibitem{Ott2008}
{\sc E.~Ott and T.~M. Antonsen}, {\em Low-dimensional behavior of large systems
  of globally coupled oscillators}, Chaos: An Interdisciplinary Journal of
  Nonlinear Science, 18 (2008), 037113.

\bibitem{Otto2018}
{\sc A.~Otto, G.~Radons, D.~Bachrathy, and G.~Orosz}, {\em Synchronization in
  networks with heterogeneous coupling delays}, Physical Review E, 97 (2018),
  012311.

\bibitem{panaggio2015}
{\sc M.~J. Panaggio and D.~M. Abrams}, {\em Chimera states: {C}oexistence of
  coherence and incoherence in networks of coupled oscillators}, Nonlinearity,
  28 (2015), R67.

\bibitem{park2018infinitesimal}
{\sc Y.~Park, K.~M. Shaw, H.~J. Chiel, and P.~J. Thomas}, {\em The
  infinitesimal phase response curves of oscillators in piecewise smooth
  dynamical systems}, European Journal of Applied Mathematics, 29 (2018),
  pp.~905--940.

\bibitem{Park2021}
{\sc Y.~Park and D.~D. Wilson}, {\em High-order accuracy computation of
  coupling functions for strongly coupled oscillators}, SIAM Journal on Applied
  Dynamical Systems, 20 (2021), pp.~1464--1484.

\bibitem{Parmelee2022a}
{\sc C.~Parmelee, J.~L. Alvarez, C.~Curto, and K.~Morrison}, {\em Sequential
  attractors in combinatorial threshold-linear networks}, SIAM Journal on
  Applied Dynamical Systems, 21 (2022), pp.~1597--1630.

\bibitem{Parmelee2022}
{\sc C.~Parmelee, S.~Moore, K.~Morrison, and C.~Curto}, {\em Core motifs
  predict dynamic attractors in combinatorial threshold-linear networks}, PLoS
  ONE, 17 (2022), e0264456.

\bibitem{RMP2015}
{\sc R.~Pastor-Satorras, C.~Castellano, P.~Van~Mieghem, and A.~Vespignani},
  {\em Epidemic processes in complex networks}, Rev. Mod. Phys., 87 (2015),
  pp.~925--979.

\bibitem{Pecora1998}
{\sc L.~M. Pecora and T.~L. Carroll}, {\em Master stability functions for
  synchronized coupled systems}, Physical Review Letters, 80 (1998),
  pp.~2109--2112.

\bibitem{Pecora2013}
{\sc L.~M. Pecora and T.~L. Carroll}, {\em Master stability function for
  globally synchronized systems}, in Encyclopedia of Computational
  Neuroscience, D.~Jaeger and R.~Jung, eds., Springer-Verlag, Heidelberg,
  Germany, 2013, pp.~1--13.

\bibitem{PecoraCarroll:1997}
{\sc L.~M. Pecora, T.~L. Carroll, G.~A. Johnson, D.~J. Mar, and J.~F. Heagy},
  {\em Fundamentals of synchronization in chaotic systems, concepts, and
  applications}, Chaos: An Interdisciplinary Journal of Nonlinear Science, 7
  (1997), pp.~520--543.

\bibitem{Pecora2014}
{\sc L.~M. Pecora, F.~Sorrentino, A.~M. Hagerstrom, T.~E. Murphy, and R.~Roy},
  {\em Cluster synchronization and isolated desynchronization in complex
  networks with symmetries}, Nature Communications, 5 (2014), 4079.

\bibitem{Cervera2018}
{\sc A.~P\'erez-Cervera, G.~Huguet, and T.~M. Seara}, {\em In Nonlinear
  Systems, Vol. 2, pages 63--81.}, Springer-Verlag, Heidelberg, Germany, 2018,
  ch.~Computation of invariant curves in the analysis of periodically forced
  neural oscillators.

\bibitem{perko2013differential}
{\sc L.~Perko}, {\em Differential Equations and Dynamical Systems},
  Springer-Verlag, Heidelberg, Germany, third~ed., 2001.

\bibitem{Pogromsky2002}
{\sc A.~Pogromsky, G.~Santoboni, and H.~Nijmeijer}, {\em Partial
  synchronization: from symmetry towards stability}, Physica D, 172 (2002),
  pp.~65--87.

\bibitem{Pogromsky2008}
{\sc A.~Y. Pogromsky}, {\em A partial synchronization theorem}, Chaos: An
  Interdisciplinary Journal of Nonlinear Science, 18 (2008), 037107.

\bibitem{ponce2014bifurcations}
{\sc E.~Ponce}, {\em Bifurcations in piecewise linear systems: Case studies},
  in VI Workshop on Dynamical Systems-MAT 70 An International Conference on
  Dynamical Systems celebrating the 70th birthday of Marco Antonio Teixeira,
  2014.

\bibitem{Porter2019}
{\sc M.~A. Porter}, {\em {Networks + Nonlinearity: {A} 2020 Vision}}, in
  Emerging Frontiers in Nonlinear Science, P.~G. Kevrekidis, J.~Cuevas-Maraver,
  and A.~Saxena, eds., vol.~32, Springer-Verlag, 2020, pp.~131--159.

\bibitem{Porter2014}
{\sc M.~A. Porter and J.~P. Gleeson}, {\em Dynamical systems on networks: A
  tutorial}, Frontiers in Applied Dynamical Systems: Reviews and Tutorials, 4
  (2016), pp.~1--79.

\bibitem{porter2016}
{\sc M.~A. Porter and J.~P. Gleeson}, {\em Dynamical Systems on Networks: A
  Tutorial}, vol.~4 of Frontiers in Applied dynamical Systems: Reviews and
  Tutorials, Springer International Publishing, Cham, Switzerland, 2016.

\bibitem{Qu2014}
{\sc Z.~Qu, G.~Hu, A.~Garfinkel, and J.~N. Weiss}, {\em Nonlinear and
  stochastic dynamics in the heart}, Physics Reports, 543 (2014), pp.~61--162.

\bibitem{Traud2011}
{\sc V.~Red, E.~D. Kelsic, P.~J. Mucha, and M.~A. Porter}, {\em Comparing
  community structure to characteristics in online collegiate social networks},
  SIAM Review, 53 (2011), pp.~526--543.

\bibitem{rook1991}
{\sc A.~J. Rook and P.~D. Penning}, {\em Synchronisation of eating, ruminating
  and idling activity by grazing sheep}, Applied Animal Behaviour Science, 32
  (1991), pp.~157--166.

\bibitem{roudi2007balanced}
{\sc Y.~Roudi and P.~E. Latham}, {\em A balanced memory network}, PLoS
  Computational Biology, 3 (2007), e141.

\bibitem{rubin2004high}
{\sc J.~E. Rubin and D.~Terman}, {\em High frequency stimulation of the
  subthalamic nucleus eliminates pathological thalamic rhythmicity in a
  computational model}, Journal of Computational Neuroscience, 16 (2004),
  pp.~211--235.

\bibitem{rubin2017balanced}
{\sc R.~Rubin, L.~Abbott, and H.~Sompolinsky}, {\em Balanced excitation and
  inhibition are required for high-capacity, noise-robust neuronal
  selectivity}, Proceedings of the National Academy of Sciences USA, 114
  (2017), pp.~E9366--E9375.

\bibitem{Sayli2021PWL}
{\sc M.~{\c{S}}ayli}, {\em {Piecewise Linear Dynamical Systems: From Nodes to
  Networks}}, PhD thesis, University of Nottingham, UK, 2021.

\bibitem{csayli2019synchrony}
{\sc M.~{\c{S}}ayli, Y.~M. Lai, R.~Thul, and S.~Coombes}, {\em {Synchrony in
  networks of {F}ranklin bells}}, IMA Journal of Applied Mathematics, 84
  (2019), pp.~1001--1021.

\bibitem{segel1976}
{\sc L.~A. Segel and S.~A. Levin}, {\em {Application of nonlinear stability
  theory to the study of the effects of diffusion on predator‐prey
  interactions}}, AIP Conference Proceedings, 27 (1976), pp.~123--152.

\bibitem{shaw1983periodically}
{\sc S.~W. Shaw and P.~J. Holmes}, {\em A periodically forced piecewise linear
  oscillator}, Journal of Sound and Vibration, 90 (1983), pp.~129--155.

\bibitem{shirasaka2017phaseAddLATER}
{\sc S.~Shirasaka, W.~Kurebayashi, and H.~Nakao}, {\em Phase reduction theory
  for hybrid nonlinear oscillators}, Physical Review E, 95 (2017), 012212.

\bibitem{ShirokyGendelman:2016}
{\sc I.~B. Shiroky and O.~V. Gendelman}, {\em Dicreate breathers in an array of
  self-excited oscillator: Exact solutions and stability}, Chaos: An
  Interdisciplinary Journal of Nonlinear Science, 26 (2016), p.~103112.

\bibitem{simpson2007andronov}
{\sc D.~Simpson and J.~Meiss}, {\em {Andronov--Hopf} bifurcations in planar,
  piecewise-smooth, continuous flows}, Physics Letters A, 371 (2007),
  pp.~213--220.

\bibitem{Simpson2011}
{\sc D.~J.~W. Simpson and R.~Kuske}, {\em {Mixed-mode oscillations in a
  stochastic, piecewise-linear system}}, Physica D, 240 (2011), pp.~1189--1198.

\bibitem{simpson2008unfolding}
{\sc D.~J.~W. Simpson and J.~D. Meiss}, {\em {Unfolding a codimension-two,
  discontinuous, Andronov--Hopf bifurcation}}, Chaos: An Interdisciplinary
  Journal of Nonlinear Science, 18 (2008), 033125.

\bibitem{Singh2013}
{\sc P.~Singh, S.~Sreenivasan, B.~K. Szymanski, and G.~Korniss}, {\em
  Threshold-limited spreading in social networks with multiple initiators},
  Scientific Reports, 3 (2013), 2330.

\bibitem{Skardal2019}
{\sc P.~S. Skardal and A.~Arenas}, {\em Higher-order interactions in complex
  networks of phase oscillators promote abrupt synchronization switching},
  Communications Physics, 3 (2020), 218.

\bibitem{Sorrentino2016}
{\sc F.~Sorrentino, L.~M. Pecora, A.~M. Hagerstrom, T.~E. Murphy, and R.~Roy},
  {\em {Complete characterization of the stability of cluster synchronization
  in complex dynamical networks}}, Science Advances, 2 (2016),
  e1501737-e1501737.

\bibitem{Sporns2016}
{\sc O.~Sporns and R.~F. Betzel}, {\em Modular brain networks}, Annual Review
  of Psychology, 67 (2016), pp.~613--640.

\bibitem{sun2009master}
{\sc J.~Sun, E.~M. Bollt, and T.~Nishikawa}, {\em Master stability functions
  for coupled nearly identical dynamical systems}, EPL (Europhysics Letters),
  85 (2009), 60011.

\bibitem{Sun2011}
{\sc J.~Sun, E.~M. Bollt, M.~A. Porter, and M.~S. Dawkins}, {\em A mathematical
  model for the dynamics and synchronization of cows}, Physica D, 240 (2011),
  pp.~1497--1509.

\bibitem{szand2021}
{\sc T.~Szanda{\l}a}, {\em Review and comparison of commonly used activation
  functions for deep neural networks}, in Bio-inspired Neurocomputing, A.~K.
  Bhoi, P.~K. Mallick, C.-M. Liu, and V.~E. Balas, eds., Springer Singapore,
  Singapore, 2021, pp.~203--224.

\bibitem{Tewarie2019}
{\sc P.~Tewarie, R.~Abeysuriya, A.~Byrne, G.~C. O'Neill, S.~N. Sotiropoulos,
  M.~J. Brookes, and S.~Coombes}, {\em {How do spatially distinct frequency
  specific MEG networks emerge from one underlying structural connectome? the
  role of the structural eigenmodes}}, NeuroImage, 186 (2019), pp.~211--220.

\bibitem{sagemath}
{\sc {The Sage Developers}}, {\em {S}ageMath, the {S}age {M}athematics
  {S}oftware {S}ystem ({V}ersion 9.5)}, 2022.
\newblock {\tt https://www.sagemath.org}.

\bibitem{Thul2010}
{\sc R.~Thul and S.~Coombes}, {\em Understanding cardiac alternans: {A}
  piecewise linear modeling framework}, Chaos: An Interdisciplinary Journal of
  Nonlinear Science, 20 (2010), 045102.

\bibitem{tonnelier2003mckean}
{\sc A.~Tonnelier}, {\em {The McKean's caricature of the {Fitzhugh--Nagumo}
  model {I}. {T}he space-clamped system}}, SIAM Journal on Applied Mathematics,
  63 (2003), pp.~459--484.

\bibitem{vanMieghem2012}
{\sc P.~van Mieghem}, {\em Graph Spectra for Complex Networks}, Cambridge
  University Press, Cambridge, UK, 2012.

\bibitem{van1996chaos}
{\sc C.~Van~Vreeswijk and H.~Sompolinsky}, {\em Chaos in neuronal networks with
  balanced excitatory and inhibitory activity}, Science, 274 (1996),
  pp.~1724--1726.

\bibitem{Veasy2019}
{\sc J.~Veasy, Y.~M. Lai, S.~Coombes, and R.~Thul}, {\em Complex patterns of
  subcellular cardiac alternans}, Journal of Theoretical Biology, 478 (2019),
  pp.~102--114.

\bibitem{wang2019shape}
{\sc Y.~Wang, J.~P. Gill, H.~J. Chiel, and P.~J. Thomas}, {\em Shape versus
  timing: linear responses of a limit cycle with hard boundaries under
  instantaneous and static perturbation}, SIAM Journal on Applied Dynamical
  Systems, 20 (2021), pp.~701--744.

\bibitem{Wang2022}
{\sc Y.~Wang, J.~P. Gill, H.~J. Chiel, and P.~J. Thomas}, {\em Variational and
  phase response analysis for limit cycles with hard boundaries, with
  applications to neuromechanical control problems}, Biological Cybernetics,
  116 (2022), pp.~687--710.

\bibitem{wilson2019isostable}
{\sc D.~Wilson}, {\em Isostable reduction of oscillators with piecewise smooth
  dynamics and complex {F}loquet multipliers}, Physical Review E, 99 (2019),
  022210.

\bibitem{wilson2018greater}
{\sc D.~Wilson and B.~Ermentrout}, {\em Greater accuracy and broadened
  applicability of phase reduction using isostable coordinates}, Journal of
  Mathematical Biology, 76 (2018), pp.~37--66.

\bibitem{wilson2019sirev}
{\sc D.~Wilson and B.~Ermentrout}, {\em Augmented phase reduction of (not so)
  weakly perturbed coupled oscillators}, SIAM Review, 61 (2019), pp.~277--315.

\bibitem{Wilson2016}
{\sc D.~Wilson and J.~Moehlis}, {\em {Isostable reduction of periodic orbits}},
  Physical Review E, 94 (2016), 052213.

\bibitem{wilson2019optimal}
{\sc D.~D. Wilson}, {\em An optimal framework for nonfeedback stability control
  of chaos}, SIAM Journal on Applied Dynamical Systems, 18 (2019),
  pp.~1982--1999.

\bibitem{Wilson72}
{\sc H.~R. Wilson and J.~D. Cowan}, {\em Excitatory and inhibitory interactions
  in localized populations of model neurons}, Biophysical Journal, 12 (1972),
  pp.~1--24.

\bibitem{Wilson73}
{\sc H.~R. Wilson and J.~D. Cowan}, {\em A mathematical theory of the
  functional dynamics of cortical and thalamic nervous tissue}, Kybernetik, 13
  (1973), pp.~55--80.

\bibitem{Winfree1967}
{\sc A.~T. Winfree}, {\em Biological rhythms and the behavior of populations of
  coupled oscillators}, Journal of theoretical biology, 16 (1967), pp.~15--42.

\bibitem{xu2013homoclinic}
{\sc B.~Xu, F.~Yang, Y.~Tang, and M.~Lin}, {\em Homoclinic bifurcations in
  planar piecewise-linear systems}, Discrete Dynamics in Nature and Society,
  2013 (2013), 732321.

\end{thebibliography}

\end{document}